\documentclass[11pt]{article}
\usepackage{amscd,amsmath,amsxtra,amssymb,theorem,latexsym,amsfonts,color,graphics}
\usepackage{bbold}
\usepackage[all]{xy}
\title{Microlocal Category}
\author{Dmitry Tamarkin}
\newtheorem{Theorem}{Theorem}[section]
\newtheorem{Lemma}[Theorem]{Lemma}

\newtheorem{Proposition}[Theorem]{Proposition}
\newtheorem{Corollary}[Theorem]{Corollary}

\newtheorem{Assumption}[Theorem]{Assumption}

\newtheorem{Claim}[Theorem]{Claim}

\usepackage{bbm}
\usepackage{rotating}

\reversemarginpar
\newcommand{\CMetR}{\mathbf{CMetR}}
\newcommand{\DMetR}{\mathbf{\bS}}
\newcommand{\SMetR}{\mathbf{SMetR}}

\newcommand{\ZQ}{\mathbb{ K}}
\newcommand{\bB}{\mathbf{B}}

\newcommand{\bF}{{\mathbf{F}}}
\newcommand{\pr}{\mathbf{pr}}

\newcommand{\bo}{{\mathbf{o}}}
\newcommand{\vu}{\nu}
\newcommand{\bbH}{{\mathbb{H}}}
\newcommand{\bDelta}{\mathbb{\Delta}}
\newcommand{\bbE}{\mathbb{E}}

\newcommand{\bZ}{\mathbb{Z}}
\newcommand{\bQ}{\mathcal{Q}}
\newcommand{\bH}{{\mathbf{H}}}
\newcommand{\bbI}{\mathbb{I}}
\renewcommand{\Re}{\mathbb{R}}
\newcommand{\cB}{\mathcal{B}}
\newcommand{\cC}{\mathcal{C}}
\newcommand{\cD}{\mathcal{D}}
\newcommand{\cE}{\mathcal{E}}
\newcommand{\cH}{\mathcal{H}}
\newcommand{\cU}{{\mathcal{U}}}
\newcommand{\cV}{{\mathcal{V}}}
\newcommand{\cW}{{\mathcal{W}}}
\newcommand{\cL}{{\mathcal{L}}}
\newcommand{\cM}{{\mathcal{M}}}
\newcommand{\cN}{{\mathcal{N}}}
\newcommand{\cK}{{\mathcal{K}}}
\newcommand{\cI}{{\mathcal{I}}}
\newcommand{\cJ}{{\mathcal{J}}}
\newcommand{\cG}{{\mathcal{G}}}
\newcommand{\cF}{{\mathcal{F}}}
\newcommand{\cO}{{\mathcal{O}}}
\newcommand{\cP}{{\mathcal{P}}}
\newcommand{\cQ}{{\mathcal{Q}}}
\newcommand{\cR}{{\mathcal{R}}}
\newcommand{\bR}{{\mathbf{R}}}
\newcommand{\cS}{{\mathcal{S}}}
\newcommand{\cT}{{\mathcal{T}}}
\newcommand{\cX}{{\mathcal{X}}}
\newcommand{\cY}{{\mathcal{Y}}}
\newcommand{\cZ}{{\mathcal{Z}}}

\newcommand{\MC}{{\mathbf{MC}}}
\newcommand{\ve}{{\varepsilon}}
\newcommand{\vs}{\sigma}
\newcommand{\op}{{\mathbf{op}}}

\newcommand{\red}{\mathbf{red}}

\newcommand{\bP}{\mathbf{P}}

\newcommand{\Hoch}{\mathbf{Hoch}}
\newcommand{\Hochcyc}{\mathbf{Hochcyc}}

\newcommand{\bbB}{\mathbb{B}}

\newcommand{\bbF}{{\mathbb{F}}}
\newcommand{\bbM}{{\mathbb{M}}}
\newcommand{\ad}{{\mathbf{Ad}}}

\DeclareMathOperator{\holim}{holim}
\DeclareMathOperator{\hocolim}{hocolim}
\DeclareMathOperator{\RHom}{RHom}
\DeclareMathOperator{\Ob}{\mathbf{Ob}}
\newcommand{\Rhom}{\RHom}

\newcommand{\Int}{\mathbf{int}}
\DeclareMathOperator{\psh}{{psh}}
\DeclareMathOperator{\sh}{{sh}}

\DeclareMathOperator{\Chech}{{Chech}}
\DeclareMathOperator{\swell}{{\mathbf{swell}}}
\DeclareMathOperator{\Open}{{Open}}
\DeclareMathOperator{\End}{End}
\newcommand{\Id}{\text{Id}}
\newcommand{\Hom}{\text{Hom}}
\newcommand{\Com}{\mathbf{Com}}

\renewcommand{\hom}{\Hom}
\newcommand{\ihom}{{\underline{\Hom}}}

\newcommand{\GZ}{\mathbf{GZ}}

\newcommand{\open}{\Open}

\newcommand{\Cov}{\mathbf{Cov}}
\newcommand\cyctrees{{\mathbf{cyctrees}}}
\DeclareMathOperator{\Cone}{Cone}

\newcommand{\cut}{{\mathbf{cut}}}

\DeclareMathOperator{\Loc}{Loc}
\newcommand{\const}{\text{const}}

\renewcommand{\SS}{\mathbf{\mu S}\;}
\newcommand{\Sp}{\text{Sp}}
\newcommand{\shinf}{\sh_q}

\newcommand{\IntB}{{\stackrel{\circ}{B}}}
\newcommand{\bS}{{\mathbf{S}}}
\newcommand{\pt}{{\mathbf{pt}}}

\newcommand{\cyc}{\mathbf{cyc}}

\renewcommand{\input}{{\mathbf{input}}}
\newcommand{\treecollec}{\mathbf{treecollec}}
\newcommand{\trees}{\mathbf{trees}}
\newcommand{\schur}{{\mathbb{S}}}

\newcommand{\assoc}{\mathbf{assoc}}
\newcommand{\unit}{{\mathbf{unit}}}

\newcommand{\gr}{\mathbf{gr}}

\newcommand{\val}{{\mathbf{val}}}

\newcommand{\Met}{{\mathbf{Met}}}
\newcommand{\Fr}{\mathbf{Fr}}

\newcommand{\MetR}{{\mathbf{MetR}}}

\newcommand{\ZZ}{\mathbf{ZZ}}

\newcommand{\Gr}{{\mathbf{Gr}}}

\newcommand{\bt}{\mathbf{t}}

\newcommand{\ovJ}{{\overline{\cJ}}}
\newcommand{\ovK}{{\overline{\cK}}}

\newcommand{\uncyc}{{\mathbf{noncyc}}}
\newcommand{\noncyc}{\uncyc}

\newcommand{\cA}{{\mathcal{A}}}
\newcommand{\bbA}{\mathbb{A}}

\newcommand{\ovQ}{{\overline{\cQ}}}

\newcommand{\bbG}{\mathbb{G}}
\newcommand{\bbS}{\mathbb{S}}
\newcommand{\ovSp}{\overline{\mathbf{Sp}}}
\newcommand{\ovU}{\overline{\mathbf{U}}}
\newcommand{\rad}{\mathbf{rad}}

\newcommand{\bbP}{{\mathbb{P}}}

\newcommand{\ho}{\text{ho}}

\newcommand{\bE}{{\mathbf{E}}}

\newcommand{\con}{\mathbf{con}}

\newcommand{\triv}{{\mathbf{triv}}}

\newcommand{\col}{\mathbf{col}}
\newcommand{\inp}{{\mathbf{inp}}}

\newcommand{\ZQmod}{{\ZQ\text{-mod}}}

\newcommand{\free}{\mathbf{free}}

\newcommand{\sets}{\mathbf{Sets}}
\newcommand{\GZtrunc}{{\mathbf{GZtrunc}}}
\newcommand{\GZtrunk}{\GZtrunc}

\newcommand{\bbQ}{{\mathbb{Q}}}

\newcommand{\TR}{\mathbf{TR}}
\newcommand{\tr}{{\mathbf{tr}}}

\newcommand{\coll}{\mathbf{coll}}
\newcommand{\Coll}{{\mathbf{Col}}}

\newcommand{\SU}{\mathbf{SU}}

\newcommand\Funct{\text{Funct}}

\newcommand{\loc}{{\mathbf{loc}}}

\newcommand{\frakH}{\mathfrak{H}}
\newcommand{\frakS}{\mathfrak{S}}
\newcommand{\tilU}{{\overline{U}}}

\newcommand{\bT}{{\mathbf{T}}}

\newcommand{\bGamma}{{\mathbf{\Xi}}}
\newcommand{\unA}{{\underline{A}}}

\newcommand{\ovFr}{{\overline{\Fr}}}

\newcommand{\ba}{{\mathbf{a}}}
\newcommand{\bm}{{\mathbf{m}}}
\newcommand{\Tr}{{\mathbf{Tr}}}

\newcommand{\bY}{{\mathbf{Y}}}
\newcommand{\bU}{\mathbf{U}}
\newcommand{\bbZ}{{\mathbf{Z}}}

\newcommand{\treesm}{\mathbf{treesm}}
\newcommand{\cTM}{{\mathcal{\cT}_{\mathcal{ \cM}}}}
\newcommand{\cCR}{\mathcal{\cC}{\cC_R}}
\newcommand{\cDR}{\mathcal{\cD}{\cD_R}}
\newcommand{\bbFR}{\bbF_R}

\newcommand{\ComAb}{{\mathbf{ComAb}}}
\newcommand{\Ker}{{\mathbf{Ker}}}
\newcommand{\bq}{\mathbf{q}}
\newcommand{\Quant}{\mathbf{Quant}}
\newcommand{\Quantum}{\Quant}
\newcommand{\Classic}{\mathbf{Classic}}
\newcommand{\alg}{\mathbf{a}}
\newcommand{\br}{\mathbf{res}}

\newcommand{\bg}{\mathbf{g}}
\newcommand{\Doplus}{({\mathbf{D}\bigoplus})}
\newcommand{\Dprod}{{\mathbf{D}\prod}}
\newcommand{\Doplusprod}{{\mathbf{D}\bigoplus\prod}}
\newcommand{\Dprodoplus}{{\mathbf{D}\prod\bigoplus}}
\newcommand{\bfD}{\mathbf{D}}
\newcommand{\ZQfreemod}{{\ZQ\mathbf{-freemod}}}
\newcommand{\ccA}{{\mathfrak{A}}}
\newcommand{\trunc}{{\mathbf{trunc}}}
\newcommand{\bbGamma}{{\mathbf{\Gamma}}}
\newcommand{\uGamma}{{\underline{\bbGamma}}}
\newcommand{\Fun}{{\mathbf{Fun}}}
\newcommand{\bbD}{\mathbb{D}}
\newcommand{\bbDhom}{{\bbD\hom}}
\newcommand{\bbAV}{{\bbA^V}}
\newcommand{\bbBV}{{\bbB^V}}
\newcommand{\bbBr}{{{  \mathbb{\cR_\bbB}  }}}
\DeclareMathOperator{\Der}{Der}
\newcommand{\bbf}{{\mathbf{f}}}
\newcommand{\into}{{\hookrightarrow}}

\begin{document}
\maketitle
\tableofcontents
\section{Introduction} In this paper we associate a dg-category  over Novikov ring with coefficients in $\bbQ$ to a compact symplectic manifold $M$
with its symplectic form having integral periods.  We generalize the approach in \cite{Tam2008}, where we 
associate a certain category ---denoted by $D_{>0}(X\times \Re)$--- to the cotangent bundle $T^*X$. No pseudoholomorphic curves are used in the construction.  The constructed category possesses properties similar
to the celebrated Fukaya category \cite{FOOO}.  For example, every object
has a support which is a closed subset of $M$;  the full sub-category of objects supported on a smooth
Lagrangian manifold admits ---modulo the maximal ideal of Novikov's ring--- the same desctiption as in the Fukaya category, and likewise 
for the hom space   between objects supported on transveral smooth Lagrangian manifolds. We also have an  invariance under Hamiltonian flow property.  These properties will be proven in  a subsequent paper.

A different approach to associating a Fukaya-like category to a symplectic manifold based on deformation quantization  is developed in \cite{Tsygan2015}.

Given a symplectic manifold $M$ as specified above,  we consider a family $F$ of Its Darboux balls,  $I:F\times B_R\to M$ --- given an $f\in F$, we have a symplectic embedding $I_f:B_R\to M$, where $B_R\subset T^*\Re^N$ is the symplectic ball of radius $R$.
We then consider a category $D_{>0}(F\times \Re^N\times \Re)$ and build a version of a  monad acting on it.
More precisely, we consider a monoidal category  $D_{>0}((F\times \Re^N)^2\times \Re)$  which acts
on $D_{>0}(F\times \Re^N\times \Re)$ by convolutions and define an $A_\infty$ algebra $A$ with curvature in 
$D_{>0}((F\times \Re^N)^2\times \Re)$.  The desired microlocal category is then defined as that of $A$-modules
in  $D_{>0}(F\times \Re^N\times \Re)$.   The algebra $A$ is first constructed in the '$\ve$-quasi-classical limit'.  That is
in a modification of $D_{>0}((F\times \Re^N)^2\times \Re)$ having the same objects but the hom definition changes
to 
\begin{equation}\label{classicve}
\hom_\ve(F,G):=\Cone(\hom(F,T_{-\ve}G)\to  \hom(F,G)),
\end{equation}
where $\ve>0$ is a small number.  The next part is the quantization, that is, lifting $A$ to an $A_\infty$-algebra
with curvature in $D_{>0}((F\times \Re^N)^2\times \Re)$.  The  quasi-classical part of the problem  can be carried over
$\bZ$ but the quantization requires passage to $\bbQ$: our approach is  based on obstruction theory which only vanish up-to torsion. 

Let us now dwell on the quantization problem.  We start with discussing:
\subsection{Ground rings/categories for quantization}

 First of all, let us draw a parallel with the more traditional
quantization setting, where we have a local complete ring $\bbQ[[q]]$ and the problem is to lift
a certain object over $\bbQ$ to $\bbQ[[q]]$. In our case it is more convenient to work in the category
of $\ve.\bZ$-graded objects in  the category  of complexes of $\bbQ$-vector spaces. Let us  denote
this category by $\cG$ (this notation will only be used in this introduction).  For such a
 $V\in \cG$, we denote by $\gr^{n\ve}V$ its corresponding graded component.

Let $\Lambda:=\bbQ[q]$ be a $\ve.\bZ$-graded ring in $\cG$, where the variable $q$ is of grading $-\ve$.
This ring is a graded analogue of the celebrated Novikov's ring.

Let us denote an appropriate  dg-version  of $D_{>0}(X\times \Re)$ by $\sh(X\times \Re)$ so that the latter
category is enriched over the category of complexes of $\bbQ$-vector spaces.
The category $\sh(X\times\Re)$ has an enrichment over $\cG$.  Let $F,G\in D_{>0}(X\times \Re)$.
Set 
$\gr^{-n\ve}\ihom(F,G):=\hom(F,T_{n\ve}G)$, where $T_{n\ve}$ denotes the shift along $\Re$ by $n\ve$.
As was explained in \cite{Tam2008},  we have natural maps
$T_{n\ve}G\to T_{m\ve G}$, $n\leq m$, whence induced maps
$$
\tau_{n,m}:\gr^{-n\ve}\ihom(F,G)\to \gr^{-m\ve}\ihom(F,G),\quad n\leq m.
$$
The action of $\bbQ[q]$ is as follows: the generator $q$ acts on $\gr^{-n\ve}\ihom(F,G)$  by $\tau_{n,n+1}$.
Observe that the properties of $D_{>0}(X\times \Re)$ do also imply the acyclicity  of the following homotopy projective limit:
$$
\gr^0\ihom(F,G)\stackrel q\leftarrow \gr^1\ihom(F,G)\stackrel q\leftarrow \gr^2\ihom(F,G)\leftarrow \cdots.
$$

This can be interpereted as the  $(q)$-adic completeness of a $\bbQ[q]$-module $\ihom(F,G)$.

Let us now consider 'the derived reduction mod $q$' which can be computed as follows
$$
\gr^{-n\ve}\ihom_\ve(F,G):=\Cone (q:\gr^{-(n-1)\ve}\ihom(F,G)\to \gr^{-n\ve}(F,G))
$$
which coincides with the definition in  (\ref{classicve}).

The category of $\bbQ[q]$-modules causes problems, for example, when we want to tensor multply over $\bbQ[q]$ ---
this requires derived tensor product.  One can replace this category with a full sub-category of 'semi-free' objects.
More precisely, we define a category $\Quant\langle\ve\rangle$ as follows.

  We first define a dg- category $Q$ whose every object is a 
$\bZ.\ve$-graded complex of $\bbQ$-vector spaces and we set
$$
\hom_Q(X,Y):=\prod\limits_{n\geq m} \hom(\gr^{m\ve}X;\gr^{n\ve}Y).
$$

An object of $\Quant\langle \ve\rangle$ is a pair $(X,D_X)$, where $X\in Q$, and $D_X\in \hom^1(X,X)$
satisfies the Maurer-Cartan equation $dD_X+D_X^2=0$.
Set
$$
\hom((X,D_X),(Y,D_Y)):=(\hom_Q(X,Y),D_{XY}),
$$
where we change the differential on the complex $\hom_Q(X,Y)$ as follows:  
$$D_{XY}f:=df+D_Yf-(-1)^{\deg f}D_X.
$$

Let $\cX:=(X,D_X)\in \Quant\langle\ve\rangle$.  Let $D_{nm}\in \hom^1(\gr^{n\ve}X;\gr^{m\ve}X)$,  $n\leq m$, be the components.
Every object in $\Quant\langle \ve\rangle$ is isomorphic to that with $D_{nn}=0$ for all $n$ which we will be assumed
from now on.

 Define a $\bbQ[q]$-module  $\cX[[q]]\in \cG$ as follows:
$$
\gr^{n\ve}\cX[[q]]:=(\prod\limits_{p\geq n} \gr^{p\ve} X, D_\cX^{(n)}),
$$
where
$$
D_{\cX}^{(n)}:=\sum\limits_{p,q|n\leq p\leq q} \; D_{pq}.
$$

The functor $\cX\mapsto \cX[[q]]$ is fully faithful; from now on we replace the category of $\bbQ[q]$-modules
in $\cG$ with  $\Quantum\langle \ve\rangle$.  
As the $\bbQ[q]$-module $\cX[[q]]$ is semi-free, we can define its classical reduction as the reduction mod $(q)$
which is nothing else but $X$ viewed as a $\bZ.\ve$-graded complex of $\bbQ$-vector spaces; in other words, we 
forget $D_X$ (since we assume that  all the diagonal components  $D_{nn}$ are equal to 0).

We thus denote by $\Classic\langle \ve\rangle$ the category of $\bZ.\ve$-graded complexes of $\bbQ$-vector
spaces (with $\hom(X,Y):=\prod_{n}(\gr^{n\ve}X;\gr^{n\ve}Y)$) and define a functor
$\red:\Quantum\langle \ve\rangle\to \Classic\langle \ve\rangle$ as explained  above: $\red(X,D):=X$.

Denote by $\sh_q(X)$ the category $\sh(X\times \Re)$ enriched over $\Quantum\langle \ve\rangle$
and by $\sh_0(X)$ the triangulated hull of $\red \sh_q(X)$.   Still denote by $\red:\sh_q(X)\to \sh_0(X)$ the natural
functor which has the meaning of the  $\ve$-quasi-classical reduction.

\subsection{$\bZ\times\bZ$-equivariance}  We have a $\bZ\times \bZ$-action on the category $\sh_q((\Fr\times \Re^N)^2)$ where $(m,n)$ acts by $T_n[2m]$, where $T_n$ is the shift along $\Re$ by $n$ units.
We define $A_0$ as the $\bZ\times \bZ$-equivariant algebra.  For simplicity, we neglect this aspect in this introduction.

\subsection{Passage to operads}     We have  monoidal categories  $\cM_q:=\sh_q((\Fr\times \Re^N)^2)$,
$\cM_0:=\sh_0((\Fr\times \Re^N)^2)$ and a monoidal functor $\red:\cM_q\to \cM_0$.
We also have an associative algebra $A_0$ in $\cM_0$.   Assume for simplicity that there exists an object $A\in \cM_q$
where $\red A\cong A_0$.     Denote by $\cO$ the full  asymmetric operad of $A$ which is defined in the category
$\Quantum\langle\ve\rangle$.    The associative algebra structure on $A_0$ reflects in a map of  asymmetric operads
$\assoc\to \red \cO$ (defined over $\Classic\langle\ve\rangle$), where $\assoc$ is the asymmetric operad
of associative  algebras,  $\assoc(n)=\bbQ$ for all $n$.
The quantization problem is now to lift this map  onto the level of $\Quantum\langle \ve\rangle$.
This can be conveniently formulated as the problem of finding a Maurer-Cartan element in $\cO$ whose classical
reduction is equivalent to the Maurer-Cartan element determined by $A_0$.
More precisely, we are to find elements
$M_n^k\in (\gr^{k\ve}\cO(n))^1$,   where:

---$
M_n^{<0}=0
$
for all $n$;

--- $M_2^0$ equals the binary product of $A_0$;  $M_n^0=0$ for all $n\neq 2$;

--- $$dM_n^k+\sum\limits_{m,l} M_m^l\{M_{n-m+1}^{k-l}\}=0,
$$

where $f\{g\}$ is the brace operation and we observe that the sum is essentially finite.
Such a Maurer-Cartan element is equivalent to the structure of an $A_\infty$ -algebra with curvature on $A$.

This type of deformation problems if well known to be controlled by the Hochschild complex of $A_0$ which does not vanish.  In order to achieve an unobstructedness we enrich the structure.

\subsubsection{Monoidal categories with a trace/Circular operads}  The notion of a trace on a monoidal category
is the categorification of the notion of a trace on an associative algebra.    Let $\cM$ be a monoidal category over
enriched over a SMC $\cC$.  Informally speaking, a contravariant trace on $\cM$ is a contravariant  functor $\Tr:\cM\to \cC$ such that
the polyfunctors $$(X_0,X_1,\ldots,X_n)\mapsto \Tr(X_0\otimes X_1\otimes \cdots\otimes X_n),\quad X_i\in \cM$$
are  invariant under cyclic permutations of $(X_0,X_1,\ldots,X_n)$.  See  Sec \ref{traces} for a more precise definition.
Given an object $A$ of a monoidal category with a trace, one associates to it the following collection of spaces:
$$
\cO_A^\uncyc(n):=\hom(A^{\otimes n};A); \quad \cO_A^\cyc(n):=\Tr(A_0\otimes A_1\otimes \cdots\otimes A_n).
$$

The objects $\cO_A(n)^\cyc$ carry a cyclic group of $(n+1)$-st  order action, we also have insertion maps
$$
\cO^\cyc_A(n)\otimes \cO^\uncyc_A(k_0)\otimes \cO^\uncyc_A(k_1)\otimes \cdots\otimes \cO^\uncyc_A(k_n)\to 
\cO^\cyc_A(k_0+k_1+\cdots+k_n).
$$
One also has an asymmetric operad structure on $\cO^\uncyc_A$.  These data satisfy certain associativity axioms.
We call such a structure {\em  a circular operad} abreviated as CO.  It is a restricted version of
a modular operad defined in \cite{GK}.

Given a monoidal category with a trace $\cM$ one defines  a trace on an algebra $A$ in $\cM$ as 
a map $t:\unit_\cM\to \Tr(A)$ satisfying:  let $m_n:A^{\otimes n}\to A$ be the $n$-fold product on $A$.
Then each element $m_n^*t\in \Tr(A^{\otimes n})$ must be cyclically invariant.

Equivalently, define a circular operad $\assoc$ whose every space is $\bbQ$, and all the composition  and cyclic group action maps preserve $1\in \bbQ$.  The structure of an algebra with a trace on an object $A$ is now
equivalent to a map of circular operads $\assoc\to \cO_A$.  

Back to our setting, we endow the category $\cM_q$, hence $\cM_0$, with a trace, we also endow the algebra $A_0$
with a trace resulting in a map of circular operads $\assoc\to \red \cO_A$, where $\cO_A$ is the full circular operad
of $A$.  One generalizes the definition of a MC-element to this setting which allows us to define the deformation
problem in the circular operad setting. This problem is still obstructed --- so we have to further enhance the structure
which is done by means of the following tool. For future purposes we will now describe the obstruction compex.
\subsubsection{Deformation complex of a map $a_0:\assoc\to \red \cO_A$} \label{defcom}   Given a map 
$a_0:\assoc\to \red \cO_A$, denote $\Hoch(a_0)^n:=\cO_A^\uncyc(n)$ and $\Hochcyc(a_0)^n:=\cO_A^\cyc(n)$.
We have a standard co-simplicial structure on $\Hoch(a_0)^\bullet$ and a co-cyclic structure on $\Hochcyc(a_0)^\bullet$.   Denote by $\Hoch(a_0)$ the total cochain complex of $\Hoch(a_0)^\bullet$.
Likewise denote by $\Hochcyc(a_0)$ the total  cyclic cochain complex of $\Hochcyc(a_0)^\bullet$ which computes
$$
R\hom_\Lambda(\underline{\bbQ};\Hochcyc(a_0)^\bullet),
$$
where $\underline{\bbQ}$ is the constant cyclic object.
The map $\red a_0:\assoc^\cyc\to \red \cO_A^\cyc$ defines a cocycle $\Omega$ in $\Hochcyc(a_0)$.  The Lie derivative $X\mapsto L_X\Omega$  defines a map
$\omega:\Hoch(a_0)\to \Hochcyc(a_0)[-1]$.   One can show that
the obstruction complex to lifting $a_0$ is as follows:
$$
\cD(a_0)=\Cone \omega.
$$

\subsection{$c_1$-Localization}    Let $\Lambda$ be the cyclic category.   Given a (co)-cyclic object $X$, we have
a natural map $c_1:X\to X[2]$, the first Chern class.  One generalizes to the case of an arbitrary circular operad as follows.   Let $\cE$ be an asymmetric operad.  Let us consider a category $C_\cE$ whose every object
is a CO $\cO$ with an identification $\cO^\uncyc=\cE$.  One can construct a category
$\bY(\cE)^\cyc$ whose every object is of the form $(n)$, $n\geq 0$, such that the structure of an  object  $\cO\in C_\cE$ is equivalent to that of a functor $F_\cO:\bY(\cE)^\cyc\to \cC$, where  $F((n))=\cO^\cyc(n)$.
For example,  $\bY(\assoc)^\cyc$  is the $\bbQ$-span of the cyclic category. 
It turns out that the $c_1$ map extends to the category of functors $\bY(\cE)^\cyc\to \cC$.   Let us now define
the $c_1$-localization in the standard way:  given $X:\bY(\cE)^\cyc\to \cC$, we define an inductive sequence
$$
X\stackrel{c_1}\to X[2]\stackrel{c_1}\to X[4]\stackrel{c_1}\to \cdots
$$
and define the hocolim of this seqence as an ind-object in the category of functors $\bY(\cE)^\cyc\to \cC$.
Denote this ind-object by $X_\loc$. A similar procedure is used for  the definition of the periodic cyclic homology.

  Let us get back to our setting, where we have
a circular operad $\cO$ in $\Quant\langle \ve\rangle$ and a map of circular operads
 \begin{equation}\label{aa0}
a_0:\assoc \to \red \cO.
\end{equation}
Consider
the object $\cO_\loc^\cyc$, which is an ind-functor $\bY^\cyc(\cO)\to \Classic\langle \ve\rangle$,
whence an induced ind $\bY(\assoc)^\cyc$-structure on $\red \cO_\loc^\cyc$, which is the same as the structure
of an ind-cocyclic object.  It turns out that this structure is rigid --- it admits a canonical lifting to the quantum level so
that we get an ind-cyclic object structure on $\cO_\loc^\cyc$.   Denote by $\cO_l$ the resulting circular operad
$(\assoc^\uncyc,\cO^\cyc_l)$.
Let  $\cU:=(\assoc^\uncyc\oplus \cO^\uncyc,\cO^\cyc)$, where $\assoc^\uncyc$ acts on $\cO^\cyc$ by zero.
We then construct a homotopy map of CO $\cU\to \cO_l$ which on the level of the underlying asymmetric operads
reduces to the projection $\assoc^\uncyc\oplus\cO^\uncyc\to \assoc^\uncyc$.  

\subsubsection{Deformation complexes and VCO }  

Let us define a map $i:\red \assoc \to \red \cU$,
where $\red \assoc^\uncyc\to \red(\assoc^\uncyc\oplus \cO^\uncyc)$  is the diagonal map $\Id\oplus a_0^\uncyc$
and the cyclic component coincides with $a_0^\cyc$, where $a_0$ is as in (\ref{aa0}).
We thus get a diagram
$$
 \xymatrix{\assoc\ar[r]\ar[dr]^\Id&\red \cU\ar[r]\ar[d]&\red\cO_l\ar[d]\\
                                    &\assoc\ar[r] & (\assoc^\uncyc,0)}.
$$
Whence an induced diagram of the deformation complexes
\begin{equation}\label{locd}
\xymatrix{
\cD(\assoc\to  \red \cU)\ar[d]^{\pi_1}\ar[r]&\cD(\assoc\to \red \cO_l)\ar[d]^{\pi_2}\\,
\cD(\assoc\to \assoc)\ar[r] & \cD(\assoc\to (\assoc^\uncyc,0)}
\end{equation}

Denote by $$\cD'(\assoc\to \red \cU):=\Cone \pi_{1}[-1];\quad  \cD'(\assoc\to\red \cO_l):=\Cone \pi_2[-1]$$.
We have an induced map 
\begin{equation}\label{wtrih}
\cD'(\assoc\to \red \cU)\to \cD'(\assoc\to \red \cO_l), 
\end{equation} which we will now describe.

  Let $a\in \Classic\langle \ve\rangle$ be as follows
$$
\gr^{n\ve} \ba=\bigoplus\limits_{k\in \bZ} \bbQ[2k],
$$
if $n\ve \in \bZ$, and $\gr^{n\ve} \ba=0$ otherwise. 
Let $C_M$ be the singular chain complex on $M$ with coefficients in $\bbQ$ and let $\Omega_M:=C_M\otimes a$  We have 
$$
\Hoch(\assoc\to\red \cO)\sim\Omega_M,
$$
$$
\Hochcyc(\assoc\to \red \cO)\sim \Omega_M[[u]],
$$
where $u$ is a variable of degree +2 --- the $c_1$-action is by the multiplication by $u$.  
The freeness of the $c_1$-action on $\Omega_M[[u]]$ comes from the triviality of the circle action
on the realization of the co-cyclic object $\Hochcyc^\bullet(\assoc\to \red \cO)$.

We have $$\cD'(\assoc\to \red \cU)\sim u^{-1}. \Omega_M[[u]][1];  \quad \cD'(\assoc\to \red \cO_l)\sim
\Omega_M[u^{-1},u]][1];$$
the map (\ref{wtrih}) is then homotopy equivalent to the obvious embedding
$$
u^{-1}. \Omega_M[[u]][1]\to \Omega_M[u^{-1},u]][1].
$$

It is crucial in our approach that this map admits a splitting.   That is, there exists an object $V_0\in \Classic\langle\ve\rangle$
and a map $V_0\to \cD(\assoc\to \red \cO_l)$ such that the induced map
$$
\cD(\assoc\to  \red \cU)\to\Cone(V_0\to \cD(\assoc\to \red \cO_l)).
$$

Because of rigidity,  the map $V_0\to  \cD(\assoc\to \red \cO_l)$ lifts as follows.  First we define $V\in \Quantum\langle\ve\rangle$ as $(V_0,0)$.  Next, we construct a homotopy map
of circular operads
$$
(\assoc^\uncyc;V\otimes \assoc^\cyc)\to \cO_l.
$$
 
We then have the following structure:

--- circular operads $\cU,\cO_l$ and an object $V$, all the data in the category $\Quant\langle\ve\rangle$;
 
--- maps of circular operads $\cU\to \cO\leftarrow (\assoc^\uncyc;V\otimes \assoc^\cyc)$.

We call such a structure {\em VCO} (V from the object $V$).     

Denote the resulting VCO by $W$.   Let also $\assoc$ be the VCO with its both CO being $\assoc$ and $V=0$,
the structure maps
$$
\assoc\to \assoc\leftarrow (\assoc^\uncyc,0)
$$
are the obvious map.  It follows that we have a map of VCO $\assoc \to \red W$.    The problem of lifting this map
turns out to be unobstructed, which concludes the solution of the quantization problem.

\subsection{Structure of the paper}

  First, we fix the notation  and review basic operations with dg-categories such as formal direct sums, products, twists by a Maurer-Cartan element.  
Next, we define dg-versions of $D_{>0}(X\times \Re)$, refered to as $\sh_q(X)$, and its quasi-classical reductions,
$\sh_\ve(X)$, $\ve>0$.
We then proceed to  the quasi-classical and the quantization parts of the construction  both having detailed 
 introductions (see Sec \ref{intro1} and Sec \ref{intro2})  to which we refer the reader.

{\em Acknowledgements} The author would like to thank M. Abouzaid,  A. d'Agnolo, A. Beilinson,  V. Drinfeld, S. Guillermou,  X. Jin,  D. Kazhdan,
  K. Kostello, Y.-G.  Oh, L. Polterovich,  N. Rosenblyum,   
P. Schapira, D. Treuman, B. Tsygan,   N. Vichery for fruitful discussions.
\section{DG Categories}
\subsubsection{Fixing a ground category}

 Let $\ZQ$ be either $\bZ$ or $\mathbb{Q}$.  Let $\ZQmod$ be the symmetric monoidal category of complexes 
of $\ZQ$-modules.
\subsection{ Category $\bfD\cC$}   Let $\cC$ be a category enriched over $\ZQmod$. 
Let $X\in \cC$.
{\em A MC-element on $X$} is an element $D\in \hom^1(X,X)$ satisfying $dD+D^2=0$.
Let $\bfD\cC$ be a category  enriched over $\ZQ$-mod whose every element is a pair $(X,D)$, where $X\in \cC$ and
$D$ is a MC-element on $X$.   Let 
$$
\hom_{\bfD\cC}((X_1,D_1),(X_2,D_2)):=(\hom_\cC(X_1,X_2),d'),
$$
where the modified differential $d'$ is defined by the formula
$$
d'f=df+D_2f-(-1)^{|f|}fD_1.
$$

A category $\cC$ enriched over $\ZQmod$ is called $D$-closed if the embedding $\cC\to \bfD\cC$,
$X\mapsto (X,0)$,
is an equivalence of categories.

 \subsection{Ground category} \label{groundcat} {\em A ground cateory } is a SMC $\ccA$ which is 

--- $D$-closed;

--- closed under direct sums and products;

--- the tensor product is compatible with direct sums and differentials,
the latter means that the natural map
$$(X,D_X)\otimes (Y,D_Y)\to (X\otimes Y,D_X\otimes 1+1\otimes D_Y)
$$
is an isomorphism for all $(X,D_X),(Y,D_Y)$;

--- has inner hom.

--- Let $\ZQfreemod\subset \ZQmod$ be the  full sub-category of complexes of free $\ZQ$-modules.
We then have a fully faithful embedding $\ZQfreemod\to \ccA$ which is a  strict tensor functor preserving differentials and direct sums.

We use the term {\em dg-category} to mean a category enriched over $\ccA$.
Throughout the paper one can assume everywhere that $\ccA=\ZQmod$.   

In the case one needs to use $\bZ$ as a ground ring one can still use the category of complexes of abelian groups
which satisfies all the properties of a ground category. However, one encounters problems
coming from inexactness of the tensor product and hom.
Alternatively, one can define a ground category $\swell(\bZ)$
which is free of these drawbacks. As this won't be used in the paper, we will sketch the construction very briefly.   We will use the notion of a co-filter on a set $S$, which, by definition, 
is a collection of subsets of $S$ satisfying: 

---if  $A\in \cF$; $B\subset A$, then $B\in \cF$;

--- if $A,B\in \cF$, then $A\cup B\in \cF$.

Given a filter $\cF$ on $S$ and $\cG$ on $T$, let us define  filters $\cF\times \cG, \hom(\cF,\cG)$ on $S\times T$, where
$\cF\times \cG$ is the smallest co-filter contaning all the sets $A\times B$, $A\in \cF$, $B\in \cG$.

The co-filter $\hom(\cF;\cG)$ consists of 
all subsets $U\subset S\times T$ satisfying: 

--- for every $t\in T$ and every $R\in \cF$, the set $(U\cap  R\times t$ is finite.

--- Let $R\in \cF$.  Define a  sub-set $H_R\subset T$ consisting of all $t\in T$  such that the set  $(U\cap  (R\times t)$ is non-empty. Then $H_R\in \cG$ for all $R\in \cF$.

Given a set $S$, a filter $\cF$ on $S$ and abelian groups $A_s$, $s\in S$, define the restricted product
$$
\prod_{s\in S}{}^\cF A_s\subset \prod_{s\in S} A_s,
$$
to consist of all families $(a_s)_{s\in S}$,  $a_s\in A_s$, satisfying: $$\{s\in S|a_s\neq 0\}\in \cF. 
$$
One defines the restricted product of complexes
of abelian groups  by setting
$$
(\prod_{s\in S} {}^\cF A_s)^k:= \prod_{s\in S}{}^\cF A_s^k.
$$

Let us now define a category $\swell_0(\bZ)$ is the following collection of data: $X:=(S_X,\cF_X,\{X_s\}_{s\in S})$, where
$\cF_X$ is a filter on $S_X$ and each $X_s$ is a finite complex of  free finitely generated abelian groups.
Set
$$
\hom(X,Y):=Z^0\prod_{(s,t)\in S\times T}{}^{\hom(\cF_X;\cF_Y)}\ \hom(X_s;Y_t),
$$
where $Z^0$ denotes the group of 0-cycles. 
Define a tensor product on $\swell_0(\bZ)$, where
$$
X\otimes Y:=(S\times T;\cF_X\times \cF_Y,\{X_s\otimes Y_t\}_{(s,t)\in S\times T}).
$$

We now set $\swell(\bZ):=D\swell_0(\bZ)$, one can check that $\swell(\bZ)$ satisfies all the properties of the ground category.

\subsection{Weak $t$-structure on a ground category $\ccA$} 
 We will define the following  weakening of $t$-structure on $\ccA$ which 
is given in terms of prescription of 
a full sub-category $D_{\leq  0}\cC$ with the following properties:

1)  $D_{\leq 0}\ccA$ contains the unit,  is $D$-closed,  closed under the tensor product, and is closed under direct sums.  $X\in D_{\leq 0} \ccA$ implies $X[n]\in D_{\leq 0}\ccA$, $n\geq 0$. 

Denote by $D_{>0}\ccA$ the full sub-category of $\ccA$ consisting of all objects $T\in D_{>0}\ccA$
satisfying $H^0\hom_\ccA(X,T)=0$ for all $X\in D_{\leq 0}X$.

2)    Let us define the catetory $\trunc(\ccA)$ consisting of all $X\in \ccA$ satisfying:
There exists a free $\ZQ$-module $A\in \cA$ and a map $A\to X$ in $\ccA$ such that
$\Cone(A\to X)\in D_{>0}\ccA$. 
   Then 
the functor $T_X:D_{\leq 0}\ccA^\op\to \sets$,  $T_X(U)=Z^0\hom_{\ccA}(U,X)$ is representable (we view here
$D_\leq \ccA^\op$ as a category over $\sets$, where the hom set is defined by $Z^0\hom_\ccA(-,-)$).
Denote the represented object by $\tau_{\leq 0}X$.
Furthermore, the induced map $A\to \tau_{\leq 0}X$ is then a homotopy equivalence and admits a unique retraction.

Instead of writing $X\in \trunc(\ccA)$ we will also say that {\em $X$ admits a truncation}.

Denote $A:=H^0(X)$, we can define $H^0(X)$ as the univeral free abelian group endowed with a map
$\tau_{\leq 0}X\to H^0(X)$ so that $H^0$ is a functor from the category $\trunc(\ccA)$ (enriched over $\sets$)
to that of free abelian groups.  We have a natural transformation $\tau_{\leq 0}\to H^0$.

It follows that $\trunc(\ccA)$ is closed under differential, products,  and positive cohomological shifts.

\begin{Lemma} \label{trunc1} Let $I$ be a category over $\sets$ and let $F:I\to \trunc\ccA$ be a functor and suppose
that $$
L^0:=\projlim\limits_{i\in I} H^0(F(i))
$$
is a free abelian group.   Then $L:=\holim_{i\in I} F(i)\in \trunc\ccA$  and we have an isomorphism $H^0(L)\cong L$.
\end{Lemma}

Let us fix throughout the paper  a ground category with a weak $t$-structure $\ccA$.
For example, the category of complexes of $\bbQ$-vector spaces.

\subsection{Tensor product}  Let $\cC,\cD$ be dg-categories.  Let
$\cC\otimes \cD$ be a dg-category whose every object is a pair $(X,Y)\in \Ob \cC\times \Ob \cD$ and
$$
\hom((X_1,Y_1);(X_2,Y_2)):=\hom_\cC(X_1,X_2)\otimes \hom_\cD(Y_1,Y_2).
$$

\subsubsection{Categories $\bigoplus \cC,\prod \cC$}  Let $\bigoplus \cC$ be the category of formal direct sums in $\cC$.   An object of $\bigoplus\cC$ is a collection of the following data:

--- a set $S$;

--- a family of objects $\{X_s\}_{s\in S}$

--- a function $n:S\to \bZ$

We refer to such an object as 
$$
\bigoplus\limits_{s\in S} X[n_s]
$$

Set $$
\hom(\bigoplus\limits_{s\in S} X[n(s)];\bigoplus\limits_{t\in T} Y[m(t)]):=
\prod\limits_{s\in S}\bigoplus\limits_{t\in T} \hom_\cC(X_s,Y_t)[m_t-n_s].
$$

Similar to above,  define the category $\prod \cC$, whose every object is defined by the same data as in $\bigoplus \cC$; an object of $\prod \cC$ is denoted as follows:
$$
\prod\limits_{s\in S} X[n(s)].
$$

Set $$
\hom(\bigoplus\limits_{s\in S} X[n(s)];\bigoplus\limits_{t\in T} Y[m(t)]):=
\prod\limits_{t\in T}\bigoplus\limits_{s\in S} \hom_\cC(X_s,Y_t)[m(t)-n(s)].
$$

We will often use the category $\Doplus \cC$ which is closed under direct sums and 
D-closed, or, briefly, $\Doplus$-closed.    Likewise the catefory $\Dprod \cC$ is  $\Dprod$-closed.

Let us list some properties of these operations.

a) Let $\cD$ be a  $\Doplus$-closed dg-category.  Let
$F:\cC\to \cD$ be a dg-functor.   The functor $F$ extends canonically to 
a functor
$$F:\Doplus \cC\to \cD.
$$

b)  Every $\Doplus$-closed category $\cC$ is tensored by the category of complexes
free $\ZQ$-modules.  

c)  We have a natural funtors $$
\Doplus \cC\otimes  \Doplus \cD\to \Doplus \cC\otimes \cD;
$$

We have a natural funtors $$
\Dprod \cC\otimes  \Dprod \cD\to \Dprod \cC\otimes \cD;
$$

d) Suppose $\cC$ is a SMC, then we have  induced SMC structures on $\Doplus \cC$, $D\prod \cC$.

e)  We have functors $$
\ihom:(\Doplus \cC)^\op\otimes \Doplus \cD\to \Dprodoplus (\cC^\op\otimes \cD);
$$
$$
\ihom:(\Doplus \cC)^\op\otimes \Dprod \cD\to \Dprod (\cC^\op\otimes \cD).
$$

\subsection{Operation $\otimes^L$}   Let $I$, $\cC$, $\cD$  be dg-categories.

Let $F:I\to \cC$; $G:I^\op\to \cD$ be dg-functors.  Define a functor
$$F\otimes^L_I G\in \Doplus (\cC\otimes\cA\otimes \cD)
$$ as follows.

Set
$$
F\otimes^L_I G:=(\bigoplus\limits_{n=0}^\infty (F\otimes^L_I G)_n,D)
$$
where
\begin{multline*}
(F\otimes^L_I G)_n=\\
\bigoplus\limits_{(X_0,X_1,\ldots,X_n\in I} F(X_0)\otimes \hom_I(X_0,X_1)\otimes \hom_I(X_1,X_2)\otimes \cdots \otimes \hom_I(X_{n-1},X_n)\otimes G(X_n)[n]\in \Doplus \cC\otimes \cD.
\end{multline*}

The MC-element  $D$ is the sum of its  components
$$
D_n:(F\otimes^L_I G)_n\to (F\otimes^L_I G)_{n-1}
$$
each of which is the standard bar differential.
\subsection{$R\hom_I(F,G)$}   Let now $F:I\to \cC$, $G:I\to cD$.   Define
$$
R\hom(F,G)\in \Dprod(\cC^\op\otimes \cA^\op\otimes \cD),
$$
where
$$
R\hom(F,G)=R\hom^0(F,G)\to F\hom^1(F,G)\to\cdots \to R\hom^n(F,G)\to\cdots,
$$
and
$$
R^n\hom(F,G)=\prod\limits_{X_0,X_1,\ldots,X_n}(F(X_0),\hom_I(X_0,X_1),\ldots,\hom_I(X_{n-1},X_n),F(X_n)).
$$

\subsection{hocolim,holim}  Let $\pt$ be the category with one object whose endomorphism
group is $\ZQ$.
Let $I$ be a category over $\sets$.   
  Let $J:=\ZQ[I]$.  Let
$\underline{\ZQ}:J^\op\to \pt$ be the constant functor.  

Let $F:J\to \cC$ be a dg-functor.
Set
$$
\hocolim_I F:=F\otimes^L_I \underline\ZQ\in \Doplus( \cC\otimes \pt)\cong
\Doplus \cC;
$$
$$
\holim_I F=R\hom_I(\underline{\ZQ},F)\in \Dprod(\cC).
$$

Suppose $\cC$ is $\Doplus$-closed and let $\tau:\Doplus \cC\to \cC$ be the canonical functor.  We will write by abuse of notation 
$\hocolim_I F\in \cC$ instead of $\tau(\hocolim_I F)$.  Likewise, whenever $\cC$ is $\Dprod$ closed,
we define $\holim_I F\in \cC$.

\section{Sheaves:  a dg-model}
\subsection{Pre-sheaves } Let $X$ be a locally compact topological space. Let $\open_X$
be the category of open subsets of $X$, where we have a unique arrow 
$U\to V$ iff $U\subset V$.    Denote by the same symbol the $\ZQ$-span of $\open_X$.
Let $\psh(X):=\Doplus  \open_X^\op$. 

Let $U\in \open_X$.   Denote by $h_U:\open_X^\op\to \Doplus \pt$ the Ioneda functor
$h_U(V):=\ZQ[\hom(V,U)]$.

Still denote by $h_U$ the extension:
$$
h_U:\Doplus \open_X^\op\to \Doplus \pt.
$$

Let $F\in \psh(X)$.   Define a functor $H_F:\open_X\to \Doplus \pt$, where
$H_F(U):=h_U(F)$.

\subsection{Sheaves}  We will define a full sub-category $\sh(X)\subset \psh(X)$
by the following conditions.  

a)  'stability':  let $\Id:\open_X^\op\to \open_X^\op$ be the identity functor.
Let $$
\rho:\Doplus \pt\otimes \Doplus \open_X^\op\to \Doplus \open_X^\op
$$ be the natural functor.

$$
G_F:=\rho(H_F\otimes^L_{\open_X} \Id)\in \Doplus \open_X^\op=\psh(X). 
$$
We have a natural map $G_F\to F$ in $\psh(X)$.

Call $F$ {\em stable} if this map is a homotopy equivalence.

b)  Coverings.

  Let $A$ be a set.  Let  $\{U_a\}_{a\in A}$ be a family of open sets in $X$.   Let $B$ be the set of all finite intersections $U_{a_1}\cap U_{a_2}\cap\cdots\cap U_{a_n}$.   For $b\in B$, denote by $U_b\subset X$ the corresponding open subset.   Let $U$ be the union of all $U_b$,  $b\in B$.
View $B$ as a sub-poset of $\open_X$.   
We have a natural map
$$
\hocolim_B H_F|_B \to H_F(U).
$$

We require this map to be a homotopy equivalence for every family $\{U_a\}_{a\in A}$.

\subsection{The sheaf $\ZQ_K$ }
\subsubsection{Coverings}
 Let $K\subset X$ be a compact subset. 
{\em A finite covering $\cU$ of $K$} is a collection $\cU_0,\cU_1,\cU_2,\ldots,\cU_N\in \open_X$ for some $N\in \bZ_{\geq 0}$.

Denote by $\Cov(K)$ the set of all finite coverings of $K$.
\subsubsection{Chech complex}
Denote by $\Chech(\cU)\in \Doplus \open_X^\op$ the following object:
$
\Chech(\cU)=(\bigoplus\limits_{n=0}^N \Chech_I^n,D)
$

$$
\Chech(\cU)^n:=(\bigoplus\limits_{0\leq i_0<i_1<\ldots i_n\leq N} \cU_{i_0}\cap \cU_{i_1}\cap\cdots\cap\cU_{i_n})[-n].
$$
The differential $D$ is the sum of its components
$D^n:\Chech(\cU)^n\to \Chech(\cU)^{n+1}$ defined in the well-known way.
  
We have a natural map 
\begin{equation}\label{chech}
X\to \Chech(\cU).
\end{equation}

\subsubsection{Cap-product}  Let $\cap:\open_X\times \open_X\to \open_X$
be the intersection bi-functor which naturally extends to a  bi-functor
$$
\cap:\psh(X)\otimes \psh(X)\to \psh(X).
$$

We have 
\begin{equation}\label{capa}
X\cap F\cong F \text{  for all } F\in \psh(X).  
\end{equation}

Let $S(K)$ be the poset of all finite subsets of $\Cov(K)$.   Set $$
\Chech(S):=\bigcap\limits_{\cU\in S} \Chech(\cU).
$$
Let $S_1\subset S_2$.  The maps (\ref{chech}), (\ref{capa}) induce a map
$\Chech(S_1)\to \Chech(S_2)$.  So that 
$$
\Chech:S\to \psh(X).
$$
Set
$$
\ZQ_K:=\hocolim_{S(K)} \Chech.
$$
\begin{Proposition} We have $\ZQ_K\in \sh(X)$.
\end{Proposition}

Let $K_1\subset K_2$.  Then $\Cov_{K_2}\subset \Cov_{K_1}$, whence 
an induced map 
\begin{equation}\label{resk}
\ZQ_{K_2}\to \ZQ_{K_1}.
\end{equation}
\subsection{ Sheaf $\ZQ_U$} 
\subsubsection{Compactification. Restriction functor}
 Let $\overline{X}$ be the one-point compactification of $X$.
Let $$r:\open_{\overline{X}}^\op\to \psh(X)
$$ 
be given by
$r(U)=U$ if $U\subset X$;  $r(U)=0$ otherwise.  The functor $r$ extends to a functor
$$
r:\psh(\overline{X})\to \psh(X).
$$
\subsubsection{Definintion of $\ZQ_U$} Let $U\in \open_X$.    Set
$$
\ZQ_U:=r\Cone(\ZQ_{\overline{X}}\to \ZQ_{\overline{X}\backslash U})[-1].
$$

The maps $\ref{resk}$ give the rule
$U\mapsto \ZQ_U$ the structure of a functor $\open_X\to \sh(X)$.

\subsection{Derived category of sheaves on $X$ versus $\sh(X)$}

Let $\Com(X)$  be the category of complexes of sheaves of $\ZQ$-modules on $X$
bounded from above, enriched over the category of  complexes of $\ZQ$-modules.

Define a functor $\zeta:\open_X\to \Com(X)$ by setting  $\zeta(U):=\ZQ_{\overline U}\in \Com(X)$.
We have an induced funcor
$$
\ho\zeta|_{\ho \sh(X)}:\ho \sh(X)\to D(X).
$$

Let us constuct the inverse functor.   Let us construct a  (non-additive) functorial free resolution functor
from the category $\ComAb^{-}$ of complexes of abelian  groups  bounded from above to  $\Doplus\pt$.

For an abelian group  $A$ set $F(A):=\bZ[A]/\bZ.0$.  We have a natural transformation
$F(A)\to A$.  Next, $F$ transfers every 0 arrow to a 0 arrow and preserves the 0 object.
Therefore, $F$ extends to a (non-additive) endofunctor  on  $\ComAb^{-}$.

   Let $X\in \ComAb^{-}$.
  Define a complex
$$
\cdots \to R^{-n}(X)\stackrel{D_n}\to R^{-n+1}(X)\stackrel{D_{n-1}}\to \cdots \stackrel{D_2}\to R^{-1}(X)\stackrel{D_1}\to R^0(X)\stackrel{D_0}\to 0,
$$
by induction.   Set $R^0(X):=F(X)$,  $D_0=0$.   Set $$
R^{-(n+1)}(X)=F(\Ker D_n).
$$
Set $$
D_{n+1}:R^{-(n+1)}(X)\to \Ker D_n\to R^{-n}(X).
$$

Let also $i_0: R^0(X)=F(X)\to X$ be the canonical map.
We have thus constructed a  bicomplex $R^\bullet(X)$; denote by $R(X)\in \Doplus\pt$ the total of $R^\bullet(X)$.   The map $i_0$ induces a natural transformation 
$R(X)\to X$ which is a quasi-isomorphism so that $R$ is a functorial free resolution.

Let now $F^\bullet$ be a complex of sheaves  on $X$ bounded from above.
We then get a functor $R(F):\open_X^\op\to \Doplus\pt$, where 
$$
R(F)(U):=R(F(U)).
$$
Define an object  $\lambda(F)\in \sh(X)$ by setting
$$
\lambda(F):=R(F)(U)\otimes^L_{U\in \open_X^\op}\ZQ_U.
$$

Let us extend $\lambda$ to the category of un-bounded complexes of sheaves on $X$.
By setting
$$
\lambda(F):=\hocolim_n \lambda(\tau_{\leq n}F).
$$

It follows that $\lambda$ converts quasi-isomorphisms  into homotopy equivalence, whence
an induced functor
$$
\ho\lambda:D(X,\ZQ)\to  \ho\sh(X).
$$

Consider the composition
$$
\zeta\lambda(F)=\hocolim_n \zeta\lambda(\tau_{\leq n}F)=
\hocolim_n R(\tau_{\leq n}F)(U)\otimes^{L}_{U\in \open_X^\op} \zeta(\ZQ_U).
$$
We have a natural transformation $\zeta(\ZQ_U)\to  \ZQ_U$ which is a quasi-isomorphism.
We therefore get a zig-zag of natural transformations which are term-wise quasi-isomorphisms between $\zeta\lambda$ and $\Id$.

Let us now consider 
$$
\lambda\zeta(F)\stackrel\sim \leftarrow F(U)\otimes^L \lambda(\zeta\ZQ_{\overline U}).
$$

We have a natural transformation $$
\ZQ_{\overline U}\to \lambda\zeta\ZQ_{\overline U}
$$
which is a homotopy equivalence.
This proves the statement.

\subsection{Base of topology}  Let $\cB\subset \open_X$ be a base of topology, that is:
$\cB$ is closed under intersections and
every open subset of $X$ can be represented as a union of sets from $\cB$.  
View $\cB$ as a poset, hence as a dg-category.
Let $\sh(\cB)\subset \Doplus(\cB^\op)$  be the full sub-category satisfying the stability conditiion and the covering condition with respect to all  coverings of elements of $\cB$
whose all terms  are in $\cB$.

Let us define a functor
$\beta_\cB:\sh(X)\to \sh(\cB)$, where  $$\beta_\cB(F)=H_F|_{\cB}\otimes^L \Id_{\cB^\op}.
$$
It follows that $\beta_\cB$ is a weak equivalence of dg-categories.

\subsection{Products} Set $\sh(X|Y):=\sh(\open_X\times \open_Y)$. 
We have a functor
$$
\boxtimes:\sh(X)\otimes \sh(Y)\to \Doplus(\open_X^\op\times \open_Y^\op).
$$
It follows that 
$$
\boxtimes:\sh(X)\otimes \sh(Y)\to \sh(X|Y).
$$

\subsubsection{Convolution}  Define a functor $$g_X:\open_X^\op\times \open_X^\op\to \Doplus\pt$$
via
$$
g_X(U,V):=\ZQ 
$$
if $U\cap V\neq \emptyset$;
$$
g_X(U,V):=0
$$
otherwise.

\begin{Lemma}  Let $X$ be compact.
We have a zig-zag homotopy equivalence
$$
g_X(U,V)\otimes^L_{(U,V)\in \open_X^\op\times \open_X^\op} \ZQ_{U\times V}\to \ZQ_{\Delta_{X\times X}},
$$
where $X$ is the diagonal.
\end{Lemma}

We have a bifunctor $$\circ:\sh(X|Y)\otimes  \sh(Y|Z)\to \sh(X|Z)
$$ 
defined as follows:
$$
\sh(X|Y)\otimes  \sh(Y|Z)\to \sh(X|Y|Y|Z)\stackrel{g_Y}\to \sh(X|Z).
$$
We have an induced map on homotopy categories
$$
\ho *:D(X\times Y)\times D(Y\times Z)\to D(X\times Z).
$$
\begin{Proposition} The bifunctor $\ho \circ $ is isomorphic to the following one:
$$
(F,G)\mapsto F*_Y G,
$$
where $*_Y$ is the operator of composition of kernels.
\end{Proposition}

\section{Quantum and Classical sheaves}
\subsection{More on categories}
\subsubsection{The  category of complexes}  Let $\cC$ be a category enriched over $\ccA$. 
Let $n\in \bZ_{\geq 0}$.   Let us define a category $\Com(\cC)\langle1/2^n\rangle'$
enriched over $\ccA$
whose every object is a $(1/2^n).\bZ$-graded object in $\cC$.  For such an $X$
denote by $\gr^{k/2^n} X$ the corresponding graded component.

Set $$
\Hom(X,Y):=\prod\limits_{m\leq n}\hom_{\cC}(\gr^m X;\gr^n Y)\in \ccA.
$$

Set 
$$
\Com(\cC)\langle1/2^n\rangle:=D\Com(\cC)\langle1/2^n\rangle'.
$$

One can represent every object of $\Com(\cC)\langle1/2^n\rangle$ as 
$(X,D)$, where $X\in \Com(\cC)\langle1/2^n\rangle_0$ and 
all the  diagonal components $D_{kk}:\gr^kX\to \gr^k X$  vanish.

Suppose the category is closed under direct products (resp. direct sums) then so is
$\Com(\cC)\langle 1/2^n\rangle$,  where $$
\gr^k(\prod\limits_a X_a)=\prod\limits_a \gr^k X_a;
$$
$$
\gr^k(\bigoplus\limits_a X_a)=\bigoplus\limits_a \gr^kX_a.
$$

Suppose $\cC$ is a SMC closed under direct sums and its tensor product is compatible with direct sums.  This structure is inherited by $\Com(\cC)\langle 1/2^n\rangle$,
where
$$
\gr^k (X\otimes Y):=\bigoplus\limits_{p\in \bZ} \gr^{k-p/2^n}X\otimes \gr^{p/2^n}Y.
$$

Suppose $\cC$ is, in addition to the  hypotheses from the  previous paragraph, closed under direct products and has an internal $\hom$.  Then so is $\Com(\cC)\langle 1/2^n\rangle$.

We have
$$
\gr^k\ihom(X,Y)=\prod\limits_{l\in \bZ}\ihom_\cC(\gr^{l/2^n}X;\gr^{(l/2^n)+k}Y).
$$ 

We have an obvious embedding functor $$
i_{nm}:\Com(\cC)\langle 1/2^n\rangle\to \Com(\cC)\langle 1/2^m\rangle,
$$
whenever $0\leq n\leq m$. The embedding $i_{nm}$ has a right adjoint, to be denoted by $c_{mn}$.  We have
$$
\gr^{p/2^n}(X,D)=(\bigoplus\limits_{2^{m-n}p\leq q<2^{m-n}(p+1)} \gr^{q/2^m}X,D'),
$$
where the differential $D'$ is induced by $D$.
In the case $\cC$ is an SMC,  $i_{mn}$ has a strict tensor structure, and $c_{mn}$ has a lax
tensor structure, that is we have a natural transformation
$$
c_{mn}(X)\otimes c_{mn}(Y)\to c_{mn}(X\otimes Y).
$$

Let also $\Gr(\cC)\langle 1/2^n\rangle$ be the category whose every object
is a $1/2^n\bZ$-graded object $X$ in $\cC$.  
Set
$$
\hom(X,Y):=\prod\limits_{n} \hom(\gr^n X;\gr^n Y).
$$

The category $\Gr(\cC)\langle 1/2^n\rangle$ has similar properties to that of $\Com(\cC)\langle 1/2^n\rangle$.

We have a functor 
$$\Gr:\Com(\cC)\langle 1/2^n\rangle\to \Gr(\cC)\langle 1/2^n\rangle.
$$

Let $(X,D)\in \Com(\cC)\langle1/2^n\rangle$, where $X\in \Com(\cC)\langle1/2^n\rangle_0$ and all the diagonal components of $D$ vanish.
Then $\Gr(X,D):= X$.

\subsubsection{The categories $\Classic(\cC)\langle 1/2^n\rangle$,   $\Quantum(\cC)\langle 1/2^n\rangle$.}

Let $\cC$ be a $D$-closed category enriched over $\ccA$.  Define a category $\Quantum(\cC)\langle 1/2^n\rangle'$ enriched over $\Com(\ccA)\langle 1/2^n\rangle$
 whose every object  $X$ is a  $\Re$-graded object in $\cC$.
Denote by $\gr^c X$ the corresponding graded component. 

Set $\hom(X,Y)=(\cH,0)$, where
$$
\gr^{k/2^n}\cH:=\prod\limits_{c\in \Re} \bigoplus\limits_{0\leq \delta<1/2^n}\hom(\gr^c X;\gr^{c+k/2^n+\delta} Y).
$$

Let $$
\Quantum(\cC)\langle 1/2^n\rangle:=D\Quantum(\cC)\langle 1/2^n\rangle'.
$$

Let us list some properties.

We have a functor
\begin{equation}\label{tquant}
T_\Quantum:\Quantum(\cC)\langle 1/2^n\rangle\otimes \Quantum(\cD)\langle 1/2^n\rangle\to \Quantum(\Doplus(\cC\otimes \cD)),
\end{equation}
where
$$
\gr^c T_\Quantum(X,Y)=\bigoplus\limits_{a\in \Re} (\gr^a X, \gr^{c-a}Y).
$$

Let $X,Y$ be $\Re$-graded objects in $\cC$.  We then have a natural isomorhphism
$$
\hom_{\Quantum(\cC)\langle 1/2^n\rangle}(X,Y)\cong c_{mn}\hom_{\Quantum(\cC)\langle 1/2^m\rangle}(X,Y),
$$
whence functors
$$
\iota_{nm}:i_{nm}\Quantum(\cC)\langle 1/2^n\rangle \to \Quantum(\cC)\langle 1/2^m\rangle,\quad n\leq m.
$$

The functor $\iota_{nm}$ is identical on objects.  Let us define
the action of $\iota_{nm}$ on hom's:
$$
\iota_{nm}:i_{nm}\hom(X,Y) =i_{nm}c_{mn}\hom(\iota_{nm}X;\iota_{nm}Y)\to \hom(\iota_{nm}X;\iota_{nm}Y),
$$
where the last arrow comes from the conjugacy.

Let $\Gamma:\Com(\ccA)\langle 1/2^n\rangle \to \ccA$ be given by
$\Gamma(X)=\hom(\unit_{\Com(\ccA)};X)$.   
It now follows that $\iota_{nm}$ induces an equivalence of categories enriched over $\ccA$:

$$
\iota_{nm}:\Gamma(\Quantum(\cC)\langle 1/2^n\rangle)\to \Gamma(\Quantum(\cC)\langle 1/2^m\rangle).
$$
Denote by  $\Quantum(\cC)$ any of $\Gamma(\Quantum(\cC)\langle 1/2^n\rangle)$.

Let us now define the category $\Classic(\cC)\langle 1/2^n\rangle:=D\Gr\Quantum(\cC)\langle 1/2^n\rangle$ enriched over
$\Gr(\ccA)\langle 1/2^n\rangle$.    

We have functors
$$
\Gr:\Gr\Quantum(\cC)\langle 1/2^n\rangle \to \Classic(\cC)\langle 1/2^n\rangle;
$$
$$
i_{nm}:i_{nm}\Classic(\cC)\langle 1/2^n\rangle \to  \Classic(\cC)\langle 1/2^m\rangle.
$$
In the case $\cC$ is an SMC, these functors have a tensor structure.

We have functors
\begin{equation}\label{tquant1}
T_{\Quant}:\Quant(\cC)\otimes \Quant(\cD)\to \Quant(\Doplus(\cC\otimes \cD);
\end{equation}
\begin{equation}\label{tclassic}
T_{\Classic}:\Classic(\cC)\langle 1/2^n\rangle\otimes \Classic(\cD)\langle 1/2^n\rangle\to \Classic(\Doplus( \cC\otimes\cD))\langle 1/2^n\rangle
\end{equation}
defined similarly to (\ref{tquant}).
\subsubsection{Functors $R_{\leq a}$,  $R_{>a}$}\label{rleq}
Let $\Quantum(\cC)_{\leq a}\subset \Quantum(\cC)$ consist of all objects
of the form
$
(X,D)
$
with $\gr^b X=0$ for all $b>a$.  

Suppose $\cC=\Doplus \cV$ for some category $\cV$ enriched over $\ccA$.
We then have an equivalence
$$
\Quantum(\cC)_{\leq a}\cong \Doplus(\cV\otimes (-\infty,a]),
$$
where  $(-\infty,a]$ is viewed as a poset, whence a category structure.

 The embedding $\Quantum(\cC)_{\leq a}\subset \Quantum(\cC)$ has a right adjoint, to be denoted by $R_{\leq a}$,
where $\gr^b(R_{\leq a}T)=0,b>a$;  $\gr^b(R_{\leq a}T)=\gr^b T$,  $b\leq a$. 
Denote $$R_{>0}X:=\Cone R_{\leq a}X\to X.
$$

\subsection{The category $\sh_q(X)$} Let us define a full subcategory
$$\sh_q(X)\subset  \Quantum(\Doplus \open_X^\op)$$ as follows.  

Let us define a functor
$$
\eta:\Quantum(\Doplus\open_X^\op)\to \Doplus(\open_{X\times \Re}^\op)
$$

Which is defined as follows:

--- let us first define $$\eta_{\leq a}:\Quantum(\Doplus \open_X^\op)_{\leq a}\cong
\Doplus(\open_X^\op\times (-\infty,a])\to  \Doplus(\open_{X\times \Re}^\op),
$$
where  the last arrow is unduced by the map
$$
(U,b)\mapsto U\times (b,\infty).
$$
Set 
$$
\eta(T):=\hocolim_{a\to \infty} \eta_{\leq a}R_{\leq a} T.
$$

Let $\sh_q(X)\subset \Quantum(\Doplus\open_X^\op)$ consist of all objects
$T$ with $\eta(T)\in \sh(X\times \Re)$.

Let $\sh(X\times \Re)_{>0}\subset \sh(X\times \Re)$ be the full sub-category
consisting of all objects $F$ with $F(U\times (-\infty,a))\sim 0$ for all $U\in \open_X$ and
all $a\in \Re$.
\begin{Proposition}  The functor $\eta$ establishes a weak equivalence 
$$\sh_q(X)\to \sh(X\times \Re)_{>0}.$$
\end{Proposition}

Let also $$\sh_{q,1/2^n}(X)\subset \Quantum(\open_X^\op)\langle 1/2^n\rangle
$$ be a full sub-category consisting of all the objects from $\sh_q(X)$.
\subsection{Full sub-category $\sh_{1/2^n}(X)\subset \Classic(\Doplus\open_X^\op)\langle 1/2^n\rangle$}

Let $\Gamma:\Classic(\ccA)\langle 1/2^n\rangle\to \ccA$ be given by
$\Gamma(X):=\hom(\unit;X)$.

We have  functors
 $$e_a:\Gamma(\Classic(\Doplus\open_X^\op)\langle 1/2^n\rangle)\to \Quantum(\Doplus \open_X^\op),\ a\in \Re,
$$
where 
$\gr^b e_a(X)=\gr^b X$, where $a\leq b<a+1/2^n$ and
$\gr^b e^a(X)=0$ otherwise.

Let $$\sh_{1/2^n}(X)\subset\Classic(\Doplus\open_X^\op)\langle 1/2^n\rangle
$$ be the full sub-category consisting of all objects $T$ with $e_a(T)\in \sh_q(X)$
for all $a\in \Re$.

The functor $$\red:\Quant(\open_X^\op)\langle 1/2^n\rangle\to \Classic(\open_X^\op)\langle 1/2^n\rangle
$$
descends to a functor
$$
\red:\sh_{q,1/2^n}(X)\to \sh_{1/2^n}(X)).
$$
We also have functors
$$
\red_{1/2^n,1/2^m}:i_{1/2^n1/2^m}\sh_{1/2^n}(X)\to \sh_{1/2^m}(X),\quad n\leq m.
$$

\subsection{Product}  
The functor $T_\Quant$ as in (\ref{tquant1}) gives rise to a functor
$$
\Quant(\Doplus\open_X^\op)\otimes \Quant(\Doplus\open_Y^\op)\to \Quant(\Doplus (\open_X^\op\times \open_Y^\op))
$$
which descends onto the level of sheaves so that we obtain the following functors
$$
\boxtimes:\sh_q(X)\otimes \sh_q(Y)\to \sh_q(X|Y);
$$
$$
\boxtimes:\sh_{q,1/2^n}(X)\otimes \sh_{q,1/2^n}(Y)\to \sh_{q,1/2^n}(X|Y).
$$
Likewise, the functors $T_\Classic$ as in (\ref{tclassic}) 
produce  functors:
$$
\boxtimes:\sh_{1/2^n}(X)\otimes\sh_{1/2^n}(Y)\to \sh_{1/2^n}(X|Y).
$$
\subsection{Convolution of kernels}
Denote 
$$\circ:\sh_q(X|Y)\otimes \sh_q(Y|Z)\to \sh_q(X|Y|Y|Z)\stackrel {g_Y}\to \sh_q(X|Z)
$$
The induced bi-functor on the homotopy categories is isomorphic to 
$$
(F,G)\mapsto Ra_!(F*_Y G),
$$
where $$
a:X\times \Re\times  Z\times \Re\to X\times Y\times \Re,
$$
$$
a(x,t_1,z,t_2)=(x,z,t_1+t_2).
$$

\subsection{Functors $f_!,f^{-1}$} Let $f:X\to Y$ be a map of locally compact topological spaces.
Let $\Gamma_f\subset X\times Y$ be the graph of $f$.   Set
$K_f\in \sh_q(X|Y)$,  $K_f=\zeta \ZQ_{\Gamma_f\times [0,\infty)}$.  Let $K_f^t\in \sh_q(Y|X)$
be the image of $K_f$ under the equivalence $\sh_q(X|Y)\sim \sh_q(Y|X)$.

Set $f_!:\sh_q(X)\to \sh_q(Y)$,  $f^{-1}:\sh_q(Y)\to\sh_q(X)$,  be given by
$$f_!(F)=F*_X K_f;\quad f^{-1}G:=F*_Y K^t_f.
$$

The functors $f_!,f^{-1}$ lift the correspondent functors between $D_{>0}(X\times \Re)$,  $D_{>0}(Y\times \Re)$.

\subsection{Microsupport} Let $F\in \sh_q(X)$.  We then have a sub-set 
$$\text{SS}\; \eta(F)\subset T^*(X\times \Re)=T^*X\times T^*\Re.
$$

Let us refer to a point in $T^*\Re$ as$ (t,k)$,  where $t\in \Re$ and $k\in T^*_k\Re$.   Let
$T^*_{>0}(X\times \Re)$ consist of all points whose $k$-coordinate is positive. 
The group  $R_{>0}^\times $ of positive dilations acts  freely on $T^*_{>0}(X\times \Re)$ fiberwise.   The quotient of
 this action is diffeomorphic to $T*X\times \Re$.   
The projection $\pi:T^*_{>0}(X\times \Re)\to T^*X\times \Re$ is  as follows
$$
\pi(x,\omega,t,k)=(x,\omega/k,t),
$$
where $x\in X$, $\omega\in T^*_x X$.

Set  $$\SS(F):=\pi(\text{SS}\; \eta(F)\cap T^*_{>0}(X\times \Re)).
$$
It follows that $$
\text{SS}\;\eta F\cap T^*_{>0}(X\times \Re)=\pi^{-1}\SS(F).
$$

The homogeneous symplectic structure  on $T^*_{>0}(X\times \Re)$ gives rise to a contact structure on
$T^*X\times \Re$, the corresponding contact form is  $\theta=dt+\alpha$, where $\alpha$ is the Liouville form on 
$T^*X$. We will always assume this contact structure on $T^*X\times \Re$.

\section{Classical theory: Introduction} \label{intro1}
\subsection{A category associated to an open subset of $T^*X\times\Re$}
 In the previous section, we have associated categories $\sh_q(X)$,
$\sh_\ve(X)$  to a smooth manifold $X$.  Given an  object $F\in \sh_q(X)$, one associates a closed
subset $\SS(F)\subset  T^*X\times \Re$.   In \cite{T2008}, it was shown that Hamiltonian symplectomorphisms
act by natural transformations on $D_{>0}(X\times \Re)$, an analogue of $\sh_q(X)$.  These observations suggest to
associate the category
$\sh_q(X)$   to a contact manifold $T^*X\times\Re$. 
Next, given an open subset $U\subset T^*X\times \Re$, we can build  categories $\sh_q(X)[U]$,  
(resp.  $\sh_\ve(X)[U]$),
which are quotients by the  full sub-category $C_q(U)$  (resp.  $C_\ve(U)$), of objects micro-supported away from $U$.

\subsection{A category associated to a symplectic ball}
The case of $
U=B_R\times \Re\subset T^*\Re^N\times\Re,
$ 
where $B_R\subset T^*\Re^N$ is the standard open ball of radius $R$, was treated in \cite{Chiu}. In particular, it was shown that the category
$\sh_q(\Re^N)[B_R\times \Re]$ is equivalent to the left orthogonal complement to $C_q(U)$.  Same is true for
$\sh_\ve(\Re^N)[B_R\times \Re]$.    These results are reviewed  in Sec \ref{symball}.

\subsection{Symplectic embeddings of a ball into $T^*\Re^N\times \Re$}
Let now $i:B_R\to T^*\Re^N$ be a symplectic embedding and consider the category 
$\sh_q(\Re^N)[i(B_R)\times \Re]$ This  category  turns out to be weakly equivalent to
$\sh_q(\Re^N)[B_R\times \Re]$, however, there is  no canonical way to choose a weak equivalence.
Such a choice requires  the following additional data

--- the lifting of $i$ to a $\Re$-equivariant contact embedding $i_c:B_R\times \Re\to T^*\Re$;

--- let $di_0\in \Sp(2N)$ be  the differential of $i$ at 0, we then have to specify a lifting $D\in \ovSp(2N)$ 
of $di_0$, where $\ovSp(2N)$ is the universal cover. 

Call these data {\em the grading  of $i$}.  Given such a graded $i$,  one constructs an object
\begin{equation}\label{individ}
\cP_i\in \sh_q(\Re^N\times \Re^N)
\end{equation}
 such that the convolution with $\cP_i$ induces an equivalence from
$\sh_q(\Re^N)[B_R\times \Re]$ 
and $\sh_q(\Re^N)[i(B_R)\times \Re]$.

\subsection{Families of symplectic embeddings of a ball into $T^*\Re^N$}

One  generalizes  to  families of grading embeddings as follows.  Let
$i:F\times B_R\to T^*\Re^N$ be a family of symplectic embeddings of $B_R$ into $T^*\Re^N$,
let $d_F:F\to \Sp(2N)$ be the map, where $d_F(f)$, $f\in F$, is the differential of $i|_{f\times B_R}$ at $(f,0)$.

Define a grading of $i$ as:

--- the lifting of $i$ to a family  $i_c:F\times B_R\times \Re\to T^*\Re^N\times \Re$   of $\Re$-equivariant   contact embeddings;

--- the lifting of $d_F$ to a  map $D_F:F\to \ovSp(2N)$.  

Given such a graded family, one constructs an object $\cP_i\in \sh_q(F\times \Re^N\times \Re^N)$, where 
the restriction $\cP_i|_{f\times \Re^N\times \Re^N}$ is homotopy  equivalent to $\cP_{i_{f\times B_R}}$ as in (\ref{individ}). 
\subsubsection{Example: parallel transitions}
Consider the simplest example,
 where $F=T^*\Re^N$.   Given an $f\in F$, set   $i_f(v)=f+v$,  where $v\in B_R$ and we use
the vector space structure on $T^*\Re^N$.

We will now pass to defining a grading on this family.
Let us refer to a point of $T^*\Re^N$ as $(q,p)$, $q\in \Re^N$;  $p\in (\Re^N)^*$.  
First,  lift our family to a $\Re$-equivariant family of contact embeddings
$
i_c:T^*\Re^N\times B_R\times \Re\to T^*\Re^N\times \Re.
$
Set
$$
i_c((q_f,p_f),(q_b,p_b),t)=(q_f+q_b,p_f+p_b,t-p_fq_b).
$$

Next, the  differential $d_0i:F\to \Sp(2N)$ is identically the identity, thus admitting an obvious lifting to $\ovSp(2N)$.

The quantization can be shown to be as follows:
$$
\cP_i:=\ZQ_{Z}*_{\Re^N} \cP_R
$$
where $Z\subset T^*\Re^N\times \Re^N\times\Re^N\times  \Re$ consists of all points
$(q_f,p_f,q,q_b',t)$ satisfying $q=q'_{b}+q_f$;  $t+p_fq'_b\geq 0$.

\subsubsection{The adjoint functor}  
 The convlution with $\cP_i$ defines a functor 
$C:\sh_q(F\times \Re^N)[T^*F\times B_R\times \Re]\to \sh_q(\Re^N)$.
Explicitly: let $A:F\times \Re^N\times \Re\to \Re^N\times \Re$ be given by
\begin{equation}\label{sdvigA}
A((q_f,p_f),q_b,t)=(q_f+q_b,t+p_fq_b),
\end{equation}
then we have
$
C(F)=A_!F.
$
Using the standard theorems from 6 functors formalism, one can show that  $C$ has a right adjoint, to be denoted by $D$.  One has $D(\cF)\sim A^!\cF\circ_{\Re^N}\cP_R$.  Next,  $A$ is a smooth trivial  fibration with its fiber being diffeomorphic
to $F$ so that $A^!\sim A^{-1}[2N]$ and $D\cF\sim A^{-1}\cF*_{\Re^N} \cP_R[2N]$ is a convolution with a kernel, 
to be denoted by $\cQ_i$, where
$$
\cQ_i:=\cP_R*_{\Re^N} \ZQ_{Z'}[2N].
$$ 
where $Z'\subset \Re^N\times T^*\Re^N\times \Re^N\times \Re$,
$$
Z':=\{(q'_b,q_f,p_f,q,t)|q=q'_b+q_f,t- p_fq_b\geq 0\}
$$

The functor $D$ has no right adjoint one, however,   there is a map
\begin{equation}\label{dcid}
DC\to \Id[2N], 
\end{equation}
defined as follows:

$$
DC\cF= \cP_R*_{\Re^N} A^{-1}A_!\cF[2N]\to \ZQ_{\Delta_{\Re^N}}*_{\Re^N} \cF[2N]\sim \cF[2N].
$$

The induced map
$$
\hom(\cF,C\cG[-N])\to \hom(D\cF;DC\cG[-N])\to \hom(D\cF;\cG)
$$
is a homotopy equivalence whenever $\cG$ is right orthogonal to all objects microsupported away from
$T^*F\times B_R\times \Re$ and the support of $\cG$ is proper along the projection $F\times \Re^N\times \Re\to \Re^N\times \Re.
$

As $D$ is a right adjoint to $C$,  we can form a corresponding monad: $\cA:=DC$ acting on $\sh_q(F\times \Re^N)$.
As the functor $C$ is essentially surjective,  we get an equivalence between the category of $\cA$-modules
in $\sh_q(F\times \Re^N)$ and
$\sh_q(\Re^N)$.

This observation suggests an idea of producing a category associated to an arbitrary compact symplectic manifold
$M^{2N}$:  choose a family $F$ of symplectic embeddings $B_R\hookrightarrow M$;  define an appropriate monad
$\cA_M$ acting on $\sh_\ve(F\times \Re^N)$, and, finally, consider a category of $\cA_M$-modules in $\sh_\ve(F\times \Re^N)$.
The problem thus  reduces to finding a monad $\cA_M$, which, in turn, reduces to finding an algebra $\bbA_M$
in the monoidal category $\sh_\ve(F\times \Re^N\times F\times \Re^N)$, so that the  $\cA_M$-action on
$\sh_\ve(F\times \Re^N)$ 
is by  the convolution with $\bbA_M$.

Unfortunately, the way, the monad $\cA$ was constructed above is not generalizeable to an arbitrary $M$ ---
the construction uses the existent category  $\sh_\ve(\Re^N)$.   However, there is an alternative --- more 'microlocal'--- way to produce
$\cA$. This alternative way works   for  more general symplectic manifolds.   We will first work out
a simplified version in the framework of the homotopy category $\ho\sh_\ve(F\times \Re^N\times F\times \Re^N)$.
   Next, we will generaize the constructuion to the case of a compact symplectic
$M$ with its symplectic form having integer periods.  We will finally carry this construction over
to the dg-setting.  We will start with formulating the key Lemma.

\subsubsection{Averaging Lemma}  Let $A:F\times \Re^N\times \Re\to \Re^N\times\Re$ be as  in (\ref{sdvigA}).
For $f\in F$, let $A_f:\Re^N\times \Re\to \Re^N\times \Re$ be the restriction of $A$ onto $f\times \Re^N\times \Re$.
The functor $A_{f!}:\ho\sh_q(\Re^N\times\Re)\to \ho\sh_q(\Re^N\times \Re)$ is then a quantization of the shift by $f$
in $T^*\Re^n$.

 Let $\gamma\in \ho\sh_q(\Re^N\times \Re^N)$ be an object satisfying
$\SS(\gamma)\subset B_R\times B_R\times \Re\subset T^*\Re^N\times T^*\Re^N\times \Re$.
Denote by $T_\gamma$ the following endofunctor on $\ho\sh_q(T^*\Re^N\times \Re^N)$,
$T_\gamma\cF:=\cF*_\Re \gamma$.

For $f\in F$, denote  $$
T_\gamma^f:=A_{f!}T_\gamma A_{-f!}:\ho\sh_q(\Re^N)\to \ho\sh_q(\Re^N).
$$
This endofunctor is representable by a kernel $\gamma^f:=\ho\sh_2(\Re^N\times \Re^N)$.
By letting $f$ vary, we get an object
$$
\gamma^F\in \ho\sh_q(F\times \Re^N\times \Re^N).
$$
Define the {\em averaging} of $\gamma$ as
 $$p_{F!}\gamma^F\in \ho\sh_q(\Re^N\times \Re^N),
$$
where 
$p_F:F\times \Re^N\times \Re^N\to \Re^N\times \Re^N$ is the projection along $F$.

Denote by $F_\gamma: \ho\sh_q( \Re^N)\to\ho \sh_q(\Re^N)$ the endofunctor determined by the kernel $p_{F!}\gamma^F$.   We have
$$
F_\gamma= C T_\gamma D.
$$

As a result of averaging, $F_\gamma$ has the following simple structure.  Let $\delta:\Re^N\to \Re^N\times \Re^N$
be the diagonal embedding  and let  $p:\Re^N\to \pt$ be the projection.      Set
$\bbf_\gamma:=p_!\delta^{-1} \gamma[N]\in \ho\sh_q(\pt)$.  By the conjugacy, we have a natural map
$\gamma\to \bbf_\gamma\boxtimes \ZQ_{\Delta_{\Re^N}}$.     Let $\cG\in \ho\sh(T^*\Re^N\times \Re^N)$.
We have an induced map $T_\gamma \cG\to \bbf_\gamma*\cG$.  Let $S_{\bbf_\gamma}$ be the operation of
convolution with $\bbf_\gamma$.  We have maps
\begin{equation}\label{ctd}
C T_\gamma D\to C S_{\bbf_\gamma} D\sim S_{\bbf_\gamma}CD\to S_{\bbf_\gamma},
\end{equation}
where the rightmost arrow is induced by the map $CD\to \Id$ coming from the conjugacy.

\begin{Claim} (Averaging Lemma) The through map $C T_\gamma D\to S_{\bbf_\gamma}$ is an isomorphism of functors.
\end{Claim}
The above construction descends onto the level of $\ho\sh_\ve$.

\subsubsection{The monoidal structure on the natural transformation (\ref{ctd})}\label{tenssovm}
Denote $\cA:=\ho\sh_q(\Re^N\times \Re^N)[B_R\times B_R]$
and  $\cB:=\ho\sh_q(\Re^N\times \Re^N)$.    We have convolution monoidal structures on both $\cA$ and $\cB$.
Denote $I(\gamma):=C T_\gamma D[-2N]$;  $J(\gamma):=S_{\bbf_\gamma}[-2N]$.  Vie replacing functors
with representing kernels, we interpret  both $I$ and $J$ as functors $\cA\to \cB$.
As it turns out,  the functors $I,J$ can be endowed with  lax tensor structure.   The natural transformation
$I\to J$  in (\ref{ctd}) upgrades to an isomorphism of lax tensor functors.

The natural transformation $I(\gamma_1)*I(\gamma_2)\to I(\gamma_1*\gamma_2)$
is as follows:
$$
CT_{\gamma_1}DCT_{\gamma_2} D[-4N]\to  CT_{\gamma_1}T_{\gamma_2}D[-2N]\sim CT_{\gamma_1*\gamma_2}D[-2N],
$$
where the first arrow is induced by the natural transformation $DC\to \Id[2N]$ as in (\ref{dcid}).

Let us define the natural transformation $J(\gamma_1)*J(\gamma_2)\to J(\gamma_1*\gamma_2)$.

We first define a natural transformation $\bbf_{\gamma_1}*\bbf_{\gamma_2}[-4N]\to \bbf_{\gamma_1*\gamma_2}[-2N]$.
Let $p:(\Re^N)^4\to\pt$ be the projection.  Denote the coordinates on $(\Re^N)^4$ by $(q_1,q_2,q_3,q_4)$. We then have:
$$
\bbf_{\gamma_1}*\bbf_{\gamma_2}[-4N]\sim p_! ((\gamma_1\boxtimes \gamma_2)\otimes \ZQ_{q_1=q_2,q_3=q_4})[-2N]\stackrel*\to
p_!((\gamma_1\boxtimes \gamma_2)\otimes \ZQ_{q_1=q_4,q_2=q_3}[-N])\sim \bbf_{\gamma_1*\gamma_2}[-2N],
$$
where the arrow $(*)$ is induced by the canonical map
\begin{equation}\label{1234}
\ZQ_{q_1=q_2,q_3=q_4}\to \ZQ_{q_1=q_2=q_3=q_4}\to \ZQ_{q_1=q_4,q_2=q_3}[N].
\end{equation}
 
We now define the tensor structure on $J$:
$$
J(\gamma_1)*J(\gamma_2)\sim S_{\bbf_{\gamma_1}*\bbf_{\gamma_2}[-4N]}\to S_{\bbf_{\gamma_1*\gamma_2}[-2N]}
\sim J(\gamma_1*\gamma_2).
$$
\subsubsection{Redefining a monad structure on $\cA$}
 It turns out that given $\ve<\pi R^2$, one can choose $\gamma$ so that $\bbf_\gamma$ is isomorphic 
to the unit of the monoidal category $\ho\sh_\ve(\pt)$ so that $S_{\bbf_\gamma}\sim \Id$.  Unfortunately there is no
analogous statement on the level of $\ho\sh_q$, that's whay we have  to switch to $\ho\sh_\ve$.

Thus, we have a map 
\begin{equation}\label{gammaunit}
\gamma\to \unit
\end{equation}
 in $\ho\sh_\ve(\Re^N\times \Re^N)$ which induces an isomorphism
$$
C T_\gamma D\to C \Id D\to \Id.
$$

Let us multiply by $D$ from the left and by $C$ from the right. This gives an isomorphism:
\begin{equation}\label{atga}
\cA T_\gamma \cA\to \cA \Id  \cA\to \cA,
\end{equation}
where $\cA=DC$, and  the last map is induced by the product on $\cA$.
On the other hand,  we have a natural transformation
$$
\cA\cA\to \cA=DC\to \Id[2N].
$$
Denote $A:=\cA$, $B:=\cA[-2N]$, so that we have a natural transformation $BA\to \Id$.   Let $\zeta:=\gamma[2N]$.  We can now rewrite (\ref{atga}) as an isomorphism
\begin{equation}\label{axiba}
AT_\zeta B\to \cA.
\end{equation}

 Observe that an associative algebra structure on $\zeta$ gives an induced associative algebra structure on $AT_\zeta B$:
$$
AT_\zeta BAT_\zeta B\to AT_\xi T_\zeta B\to AT_{\zeta*\zeta} B =AT_\zeta B\to A
$$

The goal is to find an associative structure on $\zeta$ so that  (\ref{axiba}) be an isomorphism of associative algebras.    One gets the following  sufficient condition for this.
We first observe that the associative product $\mu:\zeta*\zeta\to \zeta$ gives, by the conjugacy,  the following map
$$
\nu:\zeta\boxtimes \zeta\to p_{14}^{-1}\ZQ_{\Delta}\otimes  p_{23}^{-1}\zeta[N],
$$
where $p_{ij}:(\Re^N)^4\to (\Re^N)^2$ are the projections.

We then have the following diagram
\begin{equation}\label{d1234}
\xymatrix{
\zeta\boxtimes \zeta\ar[r]^\nu\ar[d]&  p_{14}^{-1}\ZQ_{\Delta}\otimes  p_{23}^{-1}\zeta[N]\ar[d]\\
\ZQ_{q_1=q_2,q_3=q_4}[4N]& \ZQ_{q_1=q_4,q_2=q_3}[3N]\ar[l]_{*}}
\end{equation}
where  the vertical arrows are induced by the map $\zeta\to \ZQ_{\Delta}[2N]$ induced by the map (\ref{gammaunit}).
The bottom horizontal arrow $(*)$ is induced by the map (\ref{1234}).

 \begin{Claim}  Suppose the diagram in (\ref{d1234}) commutes.  Then  (\ref{axiba}) is an isomorphism 
of algebras.
\end{Claim}

{\em Sketch of the proof} Follows from Sec \ref{tenssovm}.

\bigskip

An algebra $\zeta$ satisfying the specified conditions  is constructed (see Sec. \ref{xigamma}).
We, therefore, have an alternative construction of the monad $\bbA$.  It only  requires kernels
$A,B\in\ho\sh_\ve(F\times \Re^N\times F\times \Re^N)$ endowed with a natural transformation
$$B*_{F\times \Re^N} A\to \ZQ_{\Delta_{F\times \Re^N\times F\times \Re^N}}.
$$
One then gets a monad by the formula $\bbA:=A*T_\zeta*B$.   This approach applies to more general symplectic manifolds.
\subsection{Generalization to a compact symplectic $M$ with its symplectic form $\omega$ having integer periods} 
\label{catg}
   Choose a pseudo-Kaehler metric $g$ on $M$. 
Let $P_g$ be the bundle of $g$-unitary frames on $M$.   Given a point $m\in M$ and a $g$-unitary frame
$f\in P_g|_m$, there is a standard way to choose Darboux coordinates near $m$.  Therefore, for $R$ small enough,
one has a family of symplectic embeddings of a ball $B_R$ into $M$ parameterized by $P_g$:
$$
I:P_g\times B_R\to M.
$$
Let now $p_L:L\to M$ be the pre-quantization  bundle of $\omega$, that is a principal circle bundle  with connection whose
curvature is $\omega$ --- the existence of such a bundle follows from the integrality of periods of $\omega$.

We have the  connection 	1-form  $\theta$ on $L$   such that $p_L^*\omega=d\theta$ so that $(L,\theta)$ becomes
an $S^1$-equivariant contact manifold.  
  One can upgrade the family $I$ into a family of contact $S^1$-equivariant embeddings
$J_R:P_g\times_M L \times (B_R\times S^1)\to L$.  Denote $\Fr_g:=P_g\times_M L$.  

Let $r>0$ be a small enough number.   Denote by $J_r$ the restriction of $J_R$ onto $\Fr_g\times B_r\times S^1$.
 Consider  a subset $U\subset \Fr_g\times \Fr_g$ consisting of all pairs $(f_1,f_2)$ such that
$J_r(f_1\times B_r)\subset J_R(f_2\times B_R)$.   We then get a family of $S^1$-equivariant contact embeddings
$$
E:U\times B_r\times S^1\to B_R\times S^1.
$$
Let $u\in U$.
Let $d_E(u)\in \Sp(2N)$  be the differential  at $0\in B_R$ of the composition
$$
u\times B_R\times 0\into U\times B_r\times S^1\stackrel E\to B_R\times S^1\stackrel{\pr_{B_R}}\to B_R.
$$
Let $s_E(u)\in S^1$ be the image of the following map:
$$
u\times 0\times 0\into U\times B_r\times S^1\to B_R\times S^1\stackrel{\pr_{S^1}}\to  S^1.
$$
Let $V\subset U$ consist of all points $v$ with $s_E(v)=e_{S^1}$ and $d_E(v)$ being a Hermitian matrix in $\Sp(2N)$.
Observe that every element $g\in\Sp(2N)$  can be uniquely written as $g=hu$, where $h$ is hermitian  symplectic and $u$ is unitary.
 
One can canonically lift $E|_{V\times B_r\times S^1}$ to a $\Re$-equivariant family of graded symplectic embeddings
$$
\ve:V\times B_r\times \Re\to B_R\times \Re.
$$
We now have a quantization $\cP_\ve\in\ho \sh_q(V\times  \Re^N\times \Re^N)$.  Let $$
\cP:=j_!\cP_\ve\in \ho\sh_q(\Fr_g\times  \Re^N\times \Fr_g\times \Re^N),
$$
where 
$$j:V\times \Re^N\times \Re^N\subset \Fr_g\times \Re^N\times \Fr_g\times \Re^N
$$
is the embedding.
Let $\vs$ be the permutation  of the factors of $(\Fr_g\times \Re^N)^2$. Let $\cQ:=\vs \cP[2N]$.
In particular we have a natural transformation $\cQ*_{\Fr_g\times \Re^N} \cP\to \unit$.

We then have the following  replacements for the ingredients of the construction in the previous subsection.

--- replace $T^*\Re^N$ with $M$;

--- replace $\ho\sh_\ve(F\times \Re^N\times F\times \Re^N)$ with 
$\ho\sh_\ve(\Fr_g\times \Re^N\times\Fr_g\times \Re^N)$, where $\ve<\pi r^2$;

--- set $A$ to be the endofunctor determined by $\cP$ and $B$ to the endofunctor determined by $\cQ$.

Let us now define a generalization for $\zeta$.  We have an action of $\overline{U(N)\times S^1}$ on $\Fr_g$, where
$\overline{U(N)\times S^1}$ is the universal cover of $U(N)\times S^1$.  Therefore the category
$
\ho\sh_\ve(\overline{U(N)\times S^1}\times \Re^N\times \Re^N)
$
acts on $\ho\sh_\ve(\Fr_g\times \Re^N\times \Fr_g\times \Re^N)$ both from  the right and  the left. 
One then defines an $\overline{U(N)\times S^1}$-equivariant version of $\zeta$, to be denoted by $\xi$, where 
$\xi$ is an algebra in $\ho\sh_\ve(\overline{U(N)\times S^1}\times \Re^N\times \Re^N)$.  We now define the required monad 
$\bbA_M\in\ho\sh_\ve(\Fr_g\times \Re^N\times \Fr_g\times \Re^N)$  in the same way as above, namely,
$$
\bbA_M:=A*\xi*B.
$$
\subsubsection{Reformulation in terms of a sequence of tensor functors}  
It is convenient to break this construction into a sequence of maps of monoidal categories.

1)  Consider a monoidal   category  $\ho\sh_\ve(\Fr_g\times_M  \Fr_g\times  \Re^N\times \Re^N)$
which is identified with a full sub-category of $\ho\sh_\ve(\Fr_g\times \Fr_g\times \Re^N\times \Re^N)$
consisting of objects supported on $\Fr_g\times_M \Fr_g\times \Re^N\times \Re^N$.   This full sub-category
is closed under the tensor product, whence an embedding  tensor functor:
$$
I:\ho\sh_\ve(\Fr_g\times_M \Fr_g\times \Re^N\times \Re^N)\to \sh_\ve(\Fr_g\times \Fr_g\times \Re^N\times \Re^N).
$$

Let also $ JF:=A*I(F)*B$,  we have a tensor structure on $J:\ho\sh_\ve(\Fr_g\times_M \Fr_g\times \Re^N\times \Re^N)\to \sh_\ve(\Fr_g\times \Fr_g\times \Re^N\times \Re^N).$

Let now $\Phi$ be a pull-back of the following diagram
$$
\Fr_g\times_M \Fr_g\to U(N)\times S^1\leftarrow \overline{U(N)\times S^1}.
$$
We have a groupoid structure on $\Phi\rightrightarrows \Fr_g$ as well as a  map of groupoids
\begin{equation}\label{piq}
\pi:(\Phi\rightrightarrows \Fr_g)\to (\overline{U(N)\times S^1}\rightrightarrows \pt).
\end{equation}
The convolution product  gives a monoidal structure to the category $\sh_\ve(\Phi)$.

Let us define a tensor functor $$
K:\sh_\ve(\Phi)\to \sh_\ve(\Fr_g\times_M \Fr_g \times \Re^N\times\Re^N)
$$
as follows.

Let $q:\Phi\to \Fr_g\times_M \Fr_g$ be the covering map.    We have maps
$$\xymatrix{
\Phi && \Fr_g\times_M \Fr_g\times \Re^N\times \Re^N\\
&\Phi\times \Re^N\times \Re^N\ar[ul]^{p_1}\ar[ur]^{p_2}\ar[d]^{q}&\\
&\overline{U(N)\times S^1}\times \Re^N\times \Re^N&}
$$
where $p_1$ is the projection onto $\Phi$, $p_2=q\times \Id_{\Re^N\times \Re^N}$, and $q:\Phi\to \overline{U(N)\times S^1}$ is induced
by (\ref{piq}).
Set $K(F):=p_{2!}(p_{1}^{-1}F\otimes \pi^{-1}\xi)$.

\subsubsection{The action of $\bZ\times \bZ$}   For $a\in \bZ\times\bZ$, let $T^{\bZ\times \bZ}_a:\overline{U(N)\times S^1}\to 
\overline{U(N)\times S^1}$ be the action, where we consider $\bZ\times\bZ\subset \overline{U(N)\times S^1}$
as the kernel of the  homomorphism $\overline{U(N)\times S^1}\to U(N)\times S^1$.  
Let $a=(a_1,a_2)$, $a_1,a_2\in \bZ$. We have $T^{\bZ\times \bZ}_a \xi\cong T^{\Re}_{a_2}\xi[2a_1]$,
where, for every topological space $X$, $T^\Re_c$ denotes the endofunctor on $\sh_\ve(X)$ induced by
the parallel transition $X\times \Re\to X\times \Re$, $(x,t)\mapsto (x,t+c)$. 
One naturally extends the $\bZ\times \bZ$-action onto $\Phi$--- the action is by the deck transformations
of the covering map $\Phi\to \Fr_g\times_M \Fr_g$.

It now follows that we have an isomorphism of functors 
$$
KT^{\bZ\times \bZ}_a\stackrel\sim\to T^\Re_{a_2}K[2a_1].
$$
One can take this twisted equivariance into  account as follows.

Let $\Classic^{\bZ\times \bZ}\langle \ve\rangle$ be the SMC of $\bZ\times\bZ$-graded objects in $\Classic\langle \ve\rangle$.
For $X\in \bR_\ve^{\bZ\times \bZ}$ and $c\in \bZ\times \bZ$, let $\gr^{\bZ\times \bZ}_cX$ be the
$c$-th graded component of $X$.

Let us define the enrichment of $\sh_\ve(\Phi)$ over $\Classic^{\bZ\times \bZ}\langle \ve\rangle$, where
we set
$$
\gr^{\bZ\times \bZ}_a\ihom(F;G):=\ihom(F,T_{a_2}T^{\bZ\times\bZ}_aG[2a_1]).
$$

Denote the resulting category enriched over $\Classic^{\bZ\times \bZ}\langle \ve\rangle$ by $\sh_\ve^{\bZ\times \bZ}(\Phi)$.

We have a tensor functor $||:\Classic^{\bZ\times \bZ}\langle \ve\rangle\to \Classic\langle \ve\rangle$;  
 $$|X|:=\bigoplus\limits_{a\in \bZ\times \bZ}\gr^{\bZ\times \bZ}_a X.
$$
Whence a category $|\sh_\ve^{\bZ\times \bZ}(\Phi)|$. We have a monoidal structure, enriched over $\ho \Classic\langle \ve\rangle$,
on $\ho|\sh_\ve^{\bZ\times \bZ}(\Phi)|$.

The tensor functor $K$ now factors as the following sequence of tensor functors:
$$
\ho\sh_\ve(\Phi)\to\ho |\sh_\ve^{\bZ\times \bZ}(\Phi)|\to \ho\sh_\ve(\Fr_g\times_M \Fr_g\times \Re^N\times \Re^N).
$$

The rightmost arrow admits the following reformulation.  Let $\ZZ:=\bZ\times \bZ\cup\{\infty\}$ viewed as a partially ordered  commutative monoid;  where the elements of $\bZ\times \bZ$ are pairwise incomparable and $\infty$ is the greatest element. The addition on $\bZ\times \bZ$ is the group law and  we set $\infty+\ZZ=\infty$.
 This way, $\ZZ$ becomes a SMC, where the category structure is determined by the partial order and the tensor
product by the addition in the monoid. 

Let $\Classic^\ZZ\langle \ve\rangle$ be the SMC, enriched over $\Classic\langle\ve\rangle$, whose every object  $X$ is a collection of objects $\gr^\ZZ_c X\in \Classic\langle\ve\rangle$,
$c\in \ZZ$.   Set 
$$
\hom(X,Y):=\prod\limits_{a,b\in \ZZ| a\leq b} \hom(\gr^\ZZ_a X;\gr^\ZZ_b Y).
$$
The composition is well defined because for every $a\in \ZZ$ there are only finitely many elements $b$ such that $a\leq b$.

Next, we define the tensor product in $\Classic^\ZZ\langle \ve\rangle$ by setting $$
\gr^\ZZ_c(X\otimes Y):=\bigoplus\limits_{a,b|a+b=c} \gr^\ZZ_a X\otimes \gr^\ZZ_b Y.
$$

We have  tensor functors $i:\Classic^{\bZ\times \bZ}\langle \ve\rangle\to \Classic^\ZZ\langle \ve\rangle$; $i_\infty:\Classic\langle\ve\rangle\to \Classic^\infty\langle \infty\rangle$,
where 
$$\gr^{\ZZ}_a i(X)=\gr^{\bZ\times \bZ}_a X, a\in \bZ\times \bZ,\quad \gr^\ZZ_\infty i(X)=0;
$$
$$
\gr^\ZZ_a i_\infty(X)=0, a\in \bZ\times \bZ,\quad \gr^\ZZ_\infty(X)=X.
$$

We can now summarize our achievements in terms of a sequence of tensor functors of monoidal categories
enriched over $\ho \Classic^\ZZ\langle\ve\rangle$:
\begin{equation}\label{achiev}
\ho i\sh_\ve^{\bZ\times \bZ}(\Phi)\to \ho i_\infty \sh_\ve(\Fr_g\times_M \Fr_g\times \Re^N\times \Re^N)\to 
\ho i_\infty \sh_\ve(\Fr_g\times \Re^N\times \Fr_g\times \Re^N).
\end{equation}
\subsubsection{Building an associative algebra revisited}

Instead of building an associative algebra in the rightmost monoidal category we can do it in the leftmost one.
We have such a structure on the constant sheaf $A_\Phi:=\ZQ_{\Phi\times [0,\infty)}$.  In fact we have
a reacher structure on $A$: namely we have elements
$$\tau_a:T_{a_2}^\Re\ZQ[2a_1]\to  \gr^\ZZ_a \hom(A_\Phi,A_\Phi)$$
we have a commuatative diagram
\begin{equation}\label{achievalg}
\xymatrix
{T_{a_2+b_2}^\Re\ZQ[2a_1+2b_1]\ar[r]\ar[d]& \gr^\ZZ_{a}\hom(A_\Phi,A_\Phi)\otimes \gr^\ZZ_{b}\hom(A_\Phi,A_\Phi)\ar[d]\\
\gr^\ZZ_{a+b}(A_\Phi;A_\Phi)\ar[dr]&\gr^\ZZ_{a+b}\hom(A_\Phi\otimes A_\Phi;A_\Phi\otimes A_\Phi)\ar[d]\\
& \gr^\ZZ_{a+b}\hom(A_\Phi\otimes A_\Phi;A_\Phi)}
\end{equation}

\subsection{Lifting to a dg-level} 
\subsubsection{ As asymmetric operad as a lax version of  monoidal category}  As is usual when passing to dg structure, it is more convenient to
use lax versions of various algebraic structures.  
 Our lax version of a  monoidal category will be that 
of a colored  asymmetric operad $\cO$ over a ground SMC $\cC$. By such a structure we mean:

--- a collection of objects $C_\cO$;

--- a collection of 'operadic spaces'  $\cO(X_1,X_2,\ldots,X_n|X)\in \cC$, $X_i,X\in \cC$, $n\geq 0$ endowed with
the standard composition maps satisfying the associativity axioms.

Given a  small monoidal category $\cM$  enriched over $\cC$, with its set of objects $\Ob \cM$, we define  an asimmetric operad $\cO_\cM$
called the endomorphism operad of $\cM$, where $\cO_\cM(X_1,\ldots,X_n;X):=\hom_\cM(X_1\otimes X_2\otimes \cdots X_n;X)$.   In the case $n=0$, we set $\cO(|X):=\hom_\cM(\unit_\cM;X)$.

One has an inverse procedure of converting an asymmetric colored operad $\cO$ into a monoidal unital category $\cM^\cO$.  An object of $\cM^\cO$ is an ordered  collection $(X_1,X_2,\ldots,X_n)$ of objects of $C_\cO$.
Empty collection  is also also allowed --- it is the unit of $\cM^\cO$.   We set
$$
(X_1,X_2,\ldots,X_n)\otimes (Y_1,Y_2,\ldots,Y_m):=(X_1,\ldots,X_n,Y_1,\ldots,Y_m).
$$
Let us now define 
$$
\hom((X_1,X_2,\ldots,X_n);(Y_1,Y_2\ldots,Y_m))
$$
 
Let $S\subset \{1,2,\ldots,n\}$ be a segment, that is $S=[a,b]:=\{x|a\leq x\leq b\}$ for some $a,b\in S$, $a\leq b$, or
$S=\emptyset$.   Let $i\in m$.  Denote $([a,b]|i):=\hom(X_a,X_{a+1},\ldots,X_b|Y_i)$, set $(\emptyset|i):=\cO(|Y_i)$.

Set$$
\hom((X_1,X_2,\ldots,X_n);(Y_1,Y_2\ldots,Y_m))=\bigoplus\limits_{f:\{1,2,\ldots,n\}\to \{1,2,\ldots,m\}}
\bigotimes\limits_{i=1}^m (f^{-1}i|i),
$$
if $m>0$.  

In the case $m=0$, we set $\hom((X_1,X_2,\ldots,X_n);\unit)=0$ if $n>0$ and $\hom(\unit|\unit)=\unit_\cC$ if $n=0$.
Thus, $\cM^\cO$ is the associated  asymmetric PROP  of $\cO$.

\subsubsection{Further relaxing} We will need to further relax the definition of a colored asymmetric operad which can be done in the following standard way. First, the structure of an asymmetric operad  is, in turn, can be formulated as 
that of an algebra over a certain operad (not asymmetric anymore!), denote it temporarily by $T$.   Next, one can choose
a resolution $T'\to T$ thereby producing the notion of a $T'$-algebra as a lax version of a $T$-algebra, i.e an asymmetric operad.

Let us now discuss the operad $T$.   Of course, one can define it as a colored symmetric operad, but in fact one can do better by taking into account the specificity of the situation:   the arguments of
every meaningful  operation in $T$  correspond to vertices  of a certain planar tree. This leads us to the definition of a tree operad $\cT$  as:

--- a prescription of an object $\cT(t)$ for every isomorphism class  of   planar trees $t$;

--- Let $t$ be a planar tree with its set of   vertices $V_t$.  For each $v\in V_t$ let $E_v$ be the number of its inputs.
Let also $E_t$ be the number of inputs of $t$.   Given planar trees $t_v$ for each $v\in V_t$ such that
$E_v=E_{t_v}$, one can insert  the $t_v$'s into $t$ thus getting a new planar tree, to be denoted by
$t\{t_v\}_{v\in V_t}$.     In each such an ocasion one should have  a composition map
$$
\cT(t)\otimes \bigotimes\limits_{v\in V} \cT(t_v)\to \cT(t\{t_v\}_{v\in V_t}).
$$
These composition maps must satisfy associativity and unitality axioms.
Given a collection  of objects $O(n)\in \cC$, $n\geq 0$, one gets its endomorphism tree operad $E_O$, where
$$
E_O(t):=\hom_\cC(\bigotimes\limits_{v\in V_t}O(E_v);O(E_t)).
$$
Given a tree operad $\cT$, a $\cT$-algebra structure on $O$ is a map of tree operads $\cT\to E_O$.

One develops a colored version of a tree operad, in which case we should first fix the set of colors  $C$
and replace planar trees with $C$-colored planar trees (that is planar trees whose all edges, inputs, and the output, are $C$-colored).   Given a tree operad $\cT$ and a set of colors $C$, we define a $C$-colored tree opead
$\cT_C$ where $\cT_C(t)=\cT(\overline{t})$, where $\overline{t}$ is a planar tree with its coloring erased.
One has a notion of an algebra over a colored tree operad.  We call $\cT_C$-algebras simply $\cT$-algebras.

A $C$-colored asymmetric operad is then the same as an algebra over the trivial tree-operad $\triv$ whose all operadic spaces are set to the unit of the ground category, and  its every operadic composition sends the tensor product
of   a number of copies of the unit identically to the unit.

Define {\em a contractible tree operad} as  a tree operad $\cT$  endowed with a termwise homotopy equivalence 
$p_{\cT!}:\cT\to \triv$.  One then define a homotopy asymmetric operad as an algebra over a contractible $\cT$.
If the ground category is $\Doplus$-closed, then there is a straightening-out functor
from  the category of $\cT$-algebras to that of $\triv$-algebras so that a homotopy colored asymmetric operad
can be converted to a colored asymmetric operad proper.

\subsubsection{Homotopy asymmetric operad structures on categories of sheaves}\label{monxx} 
Our goal is to lift monoidal categories as in (\ref{achiev}) onto  dg level.  This can be done  according
the following scheme.

Let $X$ be a locally compact topological space and we have a monoidal structure on $\ho\sh_\ve(X)$, to be denoted by
$*$,
enriched over $\Classic^{\ZZ}\langle \ve\rangle$.

 GIven  objects $\gr^\ZZ_a K_n\in \sh_\ve(X^n| X)$, $n\geq 0$, $a\in \ZZ$,
we  have a $\sh_\ve(X)$-colored collection, to be denoted by $\cB^K$:
\begin{equation}\label{cbk}
\gr^\ZZ_b\cB^K(F_1,F_2,\ldots,F_n;F):=\hom_{\sh_\ve(X^n)}(F_1\boxtimes F_2\boxtimes \cdots \boxtimes F_n;\gr^\ZZ_b K_n\circ_X F).
\end{equation}
We assume that  for $n\geq 0$ we have natural transformations in  $\ho \sh_\ve(X)$:
\begin{equation}\label{agreetens}
\hom_{\ho\sh_\ve(X^n)}(F_1\boxtimes F_2\boxtimes \cdots \boxtimes F_n;\gr^\ZZ_b K_n\circ_X F)\to
\gr^\ZZ_b\hom_{\ho\sh_\ve(X)}(F_1* F_2* \cdots *F_n;F).
\end{equation}
which are isomorphisms for all $n\geq 1$.  As is seen from the examples,  for $n=0$, the situation is more delicate,
the natural transformation is only an isomorphism  for a certain full sub-category. 

Let us consider for example the category $\sh_\ve(F\times \Re^N\times F\times \Re^N)$.  Denote
$Y:=F\times \Re^N$ so that $X=Y\times Y$.  Let us refer to a point of $X$ as $(y,z)$, $y,z\in Y$.
Let $L_n\subset X^n\times X$ consist of all points $((y_1,z_1),(y_2,z_2),\ldots, (y_n,z_n),(y,z))$
satisfying  $z_1=y_2$, $z_2=y_3$, \ldots, $z_{n-1}=y_n$, $y_1=y$, $z_n=z$ if $n>0$. If $n=0$,
set $L_0\subset Y$ to be the diagonal.  Set $\gr^\ZZ_\infty K_n:=\ZQ_{L_n}[(n-1)N]$;
$\gr^\ZZ_a K_n=0$,  $a\in \bZ\times \bZ$ for all $n$.

Let $i_n:Y^{n+1}\to (Y\times Y)^n$ be given by
$$
i_n((y_0,y_1,\ldots,y_{n-1},y_n))=((y_0,y_1),(y_1,y_2),\ldots,(y_n,y_{n+1})).
$$
Let $p:Y^{n+1}\to Y\times Y$ be given by 
$$
i_n((y_0,y_1,\ldots,y_{n-1},y_n))=(y_0,y_n).
$$
Let also $q_n=pi_n$.

One now has  the following maps in the homotopy category $(n>0)$:
\begin{multline*}
\hom(F_1,F_2,\ldots,F_n,\gr^\ZZ_\infty K_n\circ_{X} F)\cong \hom(F_1\boxtimes F_2\boxtimes \cdots \boxtimes F_n;
i_{n!}i_n^{-1}p^{-1}F[(n-1)N])\\
\cong
\hom(i_n^{-1}(F_1\boxtimes \cdots\boxtimes F_n);q^!F)\cong \hom(q_!i_{n}^{-1}(F_1\boxtimes \cdots \boxtimes F_n);F)\\
\cong \hom(F_1*F_2*\cdots*F_n;F).
\end{multline*}

Consider the case $n=0$.  Let $\Delta:Y\to Y\times Y$ be the diagonal map.  We now have the following sequence of maps
$$
(\ZQ_{\Delta(Y)}\circ F)[-N]\cong \Gamma_c(Y;\Delta^{-1}F[-N])\to \Gamma(Y;\Delta^!F)\cong \hom(\Delta_!\ZQ_Y;F)
\cong \hom(\unit,F)
$$
which is an isomorphism provided that $F$ has a compact support and non-singular along the diagonal $\Delta(Y)$.

Let us now define a tree operad $\bbE$, defined over the ground category $\ccA$, which acts on a collection of objects $\sh_\ve(X)$ as follows.
Given a planar tree $t$, and elements $a_v\in \ZZ$ for each $v\in V$, one associates the object $\cK^{a_v}_{E(v)}$ to every vertex $v\in V$.  Taking convolution wih respect
to every edge of $t$, we get an object, to be denoted by $\cK_t(\{a_v\}_{v\in V})\in \sh(X^{E_t}\times X)$.  Set
\begin{equation}\label{e1}
\bbE(t)(\{a_v\}_{v\in V};a):=\hom(\cK_t(\{a_v\}_{v\in V});K_{E_t}^a);
\end{equation}
\begin{equation}\label{e2}
\bbE(t):=\prod\limits_{a\geq \sum_v a_v} \bbE(t)(\{a_v\}_{v\in V};a).
\end{equation}
  The action of $\bbE$ on $\sh_\ve(X)$ is straightforward. 
Suppose for simplicity that  $\ccA$ is the category of complexes of $\bbQ$-vector spaces.  One then can apply
the functor $\tau_{\leq 0}$ so as to get a tree operad $\tau_{\leq 0}\bbE_K$ mapped into $\bbE_K$, hence still
acting on $\cB$.   Suppose that
\begin{equation}\label{contr5}
 \tau_{\leq 0}\bbE_K\text{ has only the cohomology in degree 0.}
\end{equation}
  Therefore, we have
a homotopy equivalence of tree operads $\tau_{\leq 0}\bbE_K\to H^0(\bbE_K)$.  On the other hand,
the monoidal structure on $\ho(\sh_\ve(X))$ results in a map $\triv\to H^0(\bbE_K)$.   Let 
$C$ be a pull-back of the diagram
$
\tau_{\leq 0}\bbE_K\to H^0(\bbE_K)\leftarrow \triv.
$
As the left arrow is termwise surjective,  the structure map $C\to \triv$ is a term-wise homotopy equivalence.
On the other hand, $C$ maps to $\tau_{\leq 0}\bbE_K$, hence acts on $\cB^K$, whence
a homotopy colored asymmetric operad structure on $\sh_\ve(X)$.   Call $K$ {\em a quasi-contractible collection of kernels},
if  the condition (\ref{contr5}) is the case.

Let us denote
$\cC^F_\ve:=\sh_\ve(\Fr_g\times \Re^N\times \Fr_g\times \Re^N)$,
$\cC^\Psi_\ve:=\sh_\ve(\Fr_g\times_M \Fr_g\times \Re^N\times \Re^N)$, $\cC^\Phi_\ve:=\sh_\ve^{\bZ\times\bZ}(\Phi)$.
These categories are exactly those used in Sec \ref{catg}  The  above outlined methods 
 allow one to endow $\cC^F_\ve,\cC^\Phi_\ve,\cC^\Psi_\ve$ with a homotopy colored asymmetric operad structure,
as the corresponding collections of kernels happen to be quasi-contractible.

\subsubsection{Homotopy  maps of  homotopy colored asymmetric operads}
We follow the ideas from \cite{Fresse} to use various versions of operadic modules and bimodules.
First of all, we switch to the language of Schur functors, appropriately modified.  Fix a set of colors $C$ and a ground
SMC $\cC$.  Let $\col$ be the category, enriched over $\sets$, whose every object  $\cO$ is a collection
of objects $\cO(X_1,X_2,\ldots,X_n|X)\in \cC$ for all $n\geq 0$ and all $X_i,X\in C$. 

Let $t$ be a $C$-colored tree and $v$ a vertex of $T$.  Let $X_1^v,X_2^v,\ldots, X_{E(v)}^v$ be the colors of the inputs 
of $v$, listed according to the order of the inputs,  and let $X^v$ be the color of the output of $v$. 
Let $\cO\in \coll$.  Set $\cO_v:=\cO(X_1^v,X_2^v,\ldots,X_n^v|X^v)$.  Set $\cO(t):=\bigotimes\limits_{v\in V_t} \cO_v$.

Let also $\treecollec$ be the category, enriched over $\sets$, whose every object  $\cT$ is a collection of objects $\cT(t)\in\cC$ for every isomorphism class of $C$-colored planar trees $t$.

Let us define a functor $\schur:\treecollec\times \coll\to \coll$,  $$
\schur_\cT \cO:=\bigoplus\limits_t \cT(t)\otimes  \cO(t),
$$
where the direct sum is taken over the set of all isomorphisms classes of $C$-colored planar trees.

One has a monoidal structure on $\treecollec$, denoted by $\circ$, such that we have a natural isomorphism
$$
\schur_{\cT_1}\schur_{\cT_2}\cO\cong \schur_{\cT_1\circ \cT_2} \cO.
$$

One has 
$$
\cT_1\circ  \cT_2(t)=\bigoplus \cT_1(t_1)\otimes \bigotimes\limits_{v\in V_{t_1}} \cT_2(t_v),
$$
where the direct sum is taken over the set of all isomorphism classes of collections of $C$-colored planar trees $t_1,t_v$,  where $t=t_1\{t_v\}_{v\in V_{t_1}}$.

The structure of a tree operad on an object $\cT\in \treecollec$ is equivalent to that of a unital monoid with respect to $\circ$.
The structure of a $\cT$ algebra on an object $\cO\in \coll$ is equivalent to that of a  left $\cT$-module with respect
to the action $\schur$.

Let $\cT_1,\cT_2$ be  $C$-colored tree operads, that is  monoids in $\treecollec$.
Let $M\in \treecollec$ be a  $\cT_1-\cT_2$-bi-module and let $\cO\in \coll$ be a left $\cT_2$-module.
We have the following diagram
\begin{equation}\label{equ}
\schur_M\schur_{\cT_2}\cO_2\rightrightarrows\schur_M \cO_2
\end{equation}
where
the arrows on the left are as follows.   The upper arrow comes from the $\cT_2$-action on $M$:
$$
\schur_M\schur_{\cT_2}\cO_2\cong \schur_{M\circ \cT_2}\cO_2\to \schur_{M}\cO_2;
$$
the bottom arrow comes from the $\cT_2$-action on $\cO_2$:
$$
\schur_M(\schur_{\cT_2}\cO_2)\to \schur_M\cO_2.
$$ 
Suppose the diagram (\ref{equ}) admits a co-equalizer, to be denoted $\schur_M\circ_{\cT_2} \cO_2$.
This co-equalizer  then  carries a  left $\cT_1$-action.

One can show that if $M$ is good enough, namely if $M$ is a semi-free $\cT_2$-module, then 
$\schur_M\circ_{\cT_2}\cO_2$ does exist for any $\cO_2$. 

Let $\cT_1,\cT_2$ be contractible tree operads.
Define the notion of a contractible $\cT_1-\cT_2$-bimodule as:

--- a semi-free $\cT_1-\cT_2$-bimodule  $M$ and a map of triples ' a pair of tree operads and their bimodule':
$$
(\cT_1,\cM,\cT_2)\to (\triv,\triv,\triv),
$$
where the maps $\cT_1,\cT_2\to \triv$ are the ones coming from a contractible tree operad structure;
we vies $\triv$ as a bimodule over itself in the standard way, the map $\cM\to \triv$ is a term-wise homotopy equivalence. 

Let $\cO_k$ be a  $\cT_k$-algebra, $k=1,2$.  Define a homotopy map $\cO_2\to \cO_1$ as a prescription of
a contractible $\cT_1-\cT_2$-bimodule  $\cM$ and 
a map of $\cT_1$-algebras: 
$$
f:\cM\circ_{\cT_2} \cO_2\to \cO_1.
$$

As above, one has a straightening out procedure that converts  such a generalized  map into 
a pair $\cO_1',\cO_2'$ of strict  colored asymmetric operads and their map $f':\cO_1'\to \cO_2'$.

\subsubsection{Producing contracible bimodules}  Let $X^{(1)},X^{(2)}$ be locally-compact topological spaces and suppose
we have objects $\gr^\ZZ_a K^{(k)}_n\in \sh_\ve((X^{(k)})^{n}\times X^{(k)})$, $k=1,2$, $n\geq 0$,  $a\in \ZZ$,
similar to Sec. \ref{monxx},
in particular, we have monoidal structures on $\ho\sh_\ve(X^{(k)})$, $k=1,2$ which agree with 
the kernels  $\gr^\ZZ_a K^{(k)}$ in a way specified in (\ref{agreetens}).

We now have $X^{(k)}$-colored collections $\cB^{(k)}$ defined similar to (\ref{cbk}).
Next, we can define  objects $K^{(k)}_t\langle a_v\rangle_{v\in V_t}\in \sh_\ve((X^{(k)})^{E_t}\times X^{(k)})$, the objects
$$\bbE^{(k)}(t)( \{a_v\}_{v\in V_t};a):=\hom(K^{(k)}_t\langle a_v\rangle_{v\in V_t};\gr^\ZZ_a K^{(k)}_{E_t}),$$ and the tree operads $\bbE^{(k)}$ as in (\ref{e1}), (\ref{e2}).  Suppose we also have a tensor functor
$g:\ho\sh_\ve(X^{(1)})\to \ho\sh_\ve(X^{(2)})$ which is induced by a kernel 
$H\in \sh_\ve(X^{(2)}\times X^{(1)})$ so that $g=\ho h$, where \newline $h:\sh_\ve(X^{(1)})\to \sh_\ve(X^{(2)})$ is the convolution with $H$.

 Set $$M(t):=\hom_{(X^{(2)})^{E_t}\times X^{1}}\left(H^{\boxtimes E_t}\circ_{(X^{(1)})^{E_t}} K^{(1)}_{t};K^{(2)}_{E_t}\circ_{X^{(2)} }H\right).
$$
One gets a $\bbE^{(2)}-\bbE^{(1)}$-bimodule structure on $M$.    Let $N:=\tau_{\leq 0}M$.

 We have
a map 
\begin{equation}\label{neo}
\schur_{N}\circ_{\bbE^{(1)} }\cO^{(1)}\to h^{-1}\cO^{(2)}
\end{equation}
 of $\bbE^{(2)}$-modules.  If $\tau_{\leq 0}\bbE^{(k)}$, $k=1,2$, and $N$ only have cohomology in degree 0, we call a collection of kernels $(K^{(1)},K^{(2)},H)$ {\em quasi-contractible.}   

Similar to Sec \ref{monxx},  we now have a diagram of triples  ``a pair of operads and their bimodule"
$$
(\tau_{\leq 0}\bbE^{(1)}, N, \tau_{\leq 0}\bbE^{(2)})\to (H^0\bbE^{(1)}, H^0N, H^0\bbE^{(2)})\leftarrow (\triv,\triv,\triv)
$$
whose pull-back, to be denoted $(C_1,N_{12},C_2)$, is contractible. Let us choose a resolution $N'_{12}\to N_{12}$
of a $C_1-C_2$ bimodule $N_{12}$.  The map in (\ref{neo}) gives rise
to a map
$$
\schur_{N'_{12}}\circ_{C_2} \cO^{(1)}\to \cO^{(2)}.
$$
We get this way a  homotopy map of homotopy colored asymmetric operads $\cO^{(1)}\to \cO^{(2)}$.

 The  outlined construction applies to the monoidal categories 
in (\ref{achiev}) and their maps.   We wil get a sequence of  homotopy colored asymmetric operads  and their  homotopy maps over $\Classic^\ZZ\langle \ve\rangle$:
$$
i\cC^{\Phi}_\ve\to i_\infty \cC^{\Psi}_\ve\to i_\infty\cC^F_\ve
$$
The straightening out procedure produces a correponding sequence of strict colored asymmetric operads and their maps
\begin{equation}\label{noncycf}
\cO^{\Phi}_\ve\to  i_\infty \cO^\Psi\to i_\infty  \cO^F.
\end{equation}

\subsection{Associative algebra structure}  Let us now lift the structure from Sec \ref{achievalg} to the dg-level.
Let $\ba\in \Classic^\ZZ\langle \ve\rangle$ be given by 
$$
\gr^\ZZ_a\ba=T_{a_2}\unit[2a_1], a\in \bZ\times \bZ,\quad \gr^\ZZ_\infty \ba=0.
$$
We have a commutative algebra structure on $\ba$.
Let $\assoc$ be an asymmetric operad  with one color in $\Classic^\ZZ\langle \ve\rangle$, where $\assoc(n)=\unit$.
Let, finally, $\cO^{\bbA}_\ve:=\ba\otimes \assoc$, which is an asymmetric operad with one color in $\Classic^\ZZ\langle \ve\rangle$.

Using the above method, one can get a homotopy map $\cO^\bbA\to \cC^\Phi$ lifting the maps
from Sec \ref{achievalg}. The straightening out gives rise to  a diagram of colored asymmetric operads and their maps
$$
\cO^\bbA_\ve\stackrel\sim\leftarrow \overline{\cO}^\bbA_\ve\to \cO^\Phi_\ve.
$$
Passing to the associated monoidal categories allows one to produce maps
Let us pass to the associated PROPS, we then get a monoidal cateory
$\bbM^\Phi_\ve$, an object $X$ in it, and a map from $\overline{\cO}^\bbA_\ve$ to the full operad of $X$.
Let us pass to $\Doplus \bbM^\Phi_\ve$.  The straightening out  procedure  now produces
an object $X'\in \Doplus \bbM^\Phi_\ve$ homotopy equivalent to $X$ and a  map from $\cO^\bbA_\ve$
to the full operad of $X'$. All in all, we get a sequence of monoidal categories and their maps
\begin{equation}\label{monomap}
\bbM^{\bbA}\to \bbM^\Phi\to i_\infty \bbM^\Psi\to i_\infty \bbM^F,
\end{equation}
where we denote $\bbM^F:=\bbM^{\cO^F}$, $\bbM^\Psi:=\bbM^{\cO^\Psi}$, etc.

 We have constructed a $\bZ\times \bZ$-equivariant algebra in the monoidal category
$\bbM^F_\ve$.  In order to be able to quantize it we need to enrich a structure on it (via introduction of a trace)
and to prove some properties of this structure.  The rest of the Introduction will be devoted to  reviewing the related results.
\subsection{Traces} \label{traces}

 In this section we will enrich the structure on $\bbM^F$ and on the algebra by introducing {\em a trace} on
both of them.   

Let us define {\em a  (contravariant) trace} on a monoidal category $\cM$ over a ground category $\cC$ as:

--- a functor $\Tr:\cM^\op\to \cC$;

---  natural isomorphisms $\vs_{X,Y}: \Tr(X\otimes Y)\stackrel\sim \to \Tr(Y\otimes X)$.

These isomorphisms should satisfy:

--- $\vs_{X,Y}\vs_{Y,X}=\Id$;

--- the following diagram should commute:
$$\xymatrix{ \Tr(X\otimes \unit) \ar[rr]^{\sigma_{X,\unit}} && \Tr(\unit\otimes X)\\
& \Tr(X)\ar[ul]\ar[ur]&}
$$
where the diagonal arrows are induced by the structure maps in $\cM$:  $X\otimes \unit\to X$;
$\unit\otimes X\to X$.

---the following diagram should commute
$$
\xymatrix{\Tr(X\otimes Y\otimes Z)\ar[rr]^{\sigma_{X,Y\otimes Z}}&& \Tr(Y\otimes Z\otimes X)\ar[ddl]^{\sigma_{Y,Z\otimes X}}\\
&&\\
&\TR(Z\otimes X\otimes Y)\ar[uul]^{\sigma_{Z,X\otimes Y}}&}
$$

Let $(\cM,\Tr_\cM)$,  $(\cN,\Tr_\cN)$ be monoidal categories with a trace.  Define {\em a functor $F:(\cM,\Tr_\cM)\to (\cN,\Tr_\cN)$ of monoidal categories with trace} as:

--- a lax tensor functor $F:\cM\to \cN$ endowed with a natural transformation of functors
$s:\Tr_\cM\to \Tr_\cN\circ F$  fitting into the following commutative diagram

$$
\xymatrix{ \Tr_\cM(X\otimes Y)\ar[rr]^{\vs_{X,Y}}\ar[d]_s && \Tr_\cM(Y\otimes X)\ar[d]^s\\
                 \Tr_\cN(F(X\otimes Y))  \ar[d]_{(*) }                       &&\Tr_\cN(F(Y\otimes X))\ar[d]^{(*)}\\
            \Tr_\cN(F(X)\otimes F(Y)) \ar[rr]_{\vs_{F(X),F(Y)}}&& \Tr_\cN(F(Y)\otimes F(X))}
$$

where the arrows $(*)$ are induced by the natural maps $F(X)\otimes F(Y)\to F(X\otimes Y)$ and
$F(Y)\otimes F(X)\to F(Y\otimes X)$.

Let now $A$ be an associative algebra in $\cM$.  {\em A trace on $A$} is an element $\tr_A:\unit\to \Tr(A)$
such that the following diagram commutes
$$\xymatrix{
\unit \ar[r]^{\tr_A}\ar[d]_{\tr_A}&\Tr(A)\ar[rr]^{\Tr(m_A)}&& \Tr(A\otimes A)\ar[dl]^{\vs_{A,A}}\\
\Tr(A)\ar[rr]^{\Tr(m_A)}&& \Tr(A\otimes A)&}
$$

Define the trace on the  category $\ho\sh_\ve(Y\times Y)$, where  $Y=\Fr_g\times \Re^N$,
by setting $\Tr(F):=\hom(F;\unit)$.   One endows every category in (\ref{achiev}) with a trace in a similar way, 
all the tensor functors in (\ref{achiev}) can be defined as  functors of monoidal categories with a trace.

\subsection{Circular operads}  In order to transfer the  traces onto the dg level, we resort to operads.
As was discussed above,  every monoidal category gives rise to a colored asymmetric operad.
In a similar fashion,  every monoidal category  with trace  enriched over an SMC  $\cC$ produces {\em a  colored circular operad in $\cC$} which, by definition, is

--- a colored  asymmetric operad $\cO^\uncyc$  in $\cC$ with its  set of  colors $C$ 

--- objects $\cO^\cyc(c_0,c_1,\ldots,c_n)$ for all $n\geq 0$ and all $c_k\in C$.

--- let $f:\{0,1,\ldots,n\}\to \{0,1,\ldots,m\}$ be a cyclically monotone map.   For each $k\in \{0,1,\ldots,m\}$,
$f^{-1}k$ is a sub-interval of $\{0,1,\ldots,n\}$, including  an empty set.    Let $d_1,d_2,\ldots,d_n\in C$.
For every sub-interval $I\subset \{0,1,\ldots,n\}$,  $I=[a,b]$,  let $d_I$ be a sub-sequence 
$d_{i_a},d_{i_{a+1}},\ldots,d_{b}$.     We then shoud have an insertion map
$$
\cO^\cyc(c_0,c_1,\ldots,c_m)\otimes \bigotimes\limits_{k=0}^m \cO^\uncyc(d_{f^{-1}k}|c_k)\to \cO(d_0,d_1,\ldots,d_n).
$$

The insertion maps should obey the obvious associativity and unitality axioms.

Every  small monoidal category  $\cM$ with a trace enriched over $\cC$ defines such a circular operad  $\cO$ whose set of colors
is the set of objects of $\cM$,  where, as above,
$$
\cO^\uncyc(X_1,X_2,\ldots,X_n|X):=\hom_\cM(X_1\otimes X_2\otimes \cdots \otimes X_n;X). 
$$
Next, we set
$$
\cO^\cyc(X_0,X_1,\ldots,X_n):=\Tr(X_0\otimes X_1\otimes \cdots\otimes X_n).
$$

The above constructions carry over   onto the setting of circular operads straightforwardly so that
we get a sequence of circular operads and their maps similar to (\ref{noncycf}). The straightening out procedure
leads to a sequence of monoidal categories with a trace and their functors similar to (\ref{monomap}).

\subsection{Hochschild Complexes}  Let $A$ be an associative algebra a monoidal category with trace $\cM$
enriched over a SMC $\cC$.   One associates to $A$ two  standard objects. One of them is a co-simplicial object
$\Hoch^\bullet(A)$, where $\Hoch^n(A):=\hom(A^{\otimes n};A)$.  The totalization of this object is the standard Hochschild cochain complex of $A$.    The trace on $\cM$ allows one to define a co-cyclic object
$\Hochcyc^\bullet(A)$,  where $\Hochcyc^n(A):=\Tr(A^{\otimes n+1})$.   Let us denote
by $\Hoch(A)$ and $\Hochcyc(A)$ the totalizations.

Recall that we have constructed an associative algebra with trace in $\bbM^\Phi$, to be denoted by $A_\Phi$ (see (\ref{monomap}).
Denote the images of $A_\Phi$  by $A_\Psi\in \bbM^\Psi$ and $A_F\in \bbM^F$. The tensor functors
in (\ref{monomap}) define maps of Hochschild complexes
$$
\Hoch(A_\Phi)\to \Hoch(A_\Psi)\to \Hoch(A_F);\quad  \Hochcyc(A_\Phi)\to \Hochcyc(A_\Psi)\to \Hochcyc(A_F).
$$
One shows that each of these maps is a homotopy equivalence.

\subsection{$\bY(\bbB)$-modules}  Let us now take the $\bZ\times\bZ$-equivariant structure on our algebras into
account.  We will use the notion of a $\cO$-{\em bimodule}, where $\cO$ is an arbitrary circular operad over the set
of colors $C$.  It is convenient to define   a set $Y(C)^\noncyc$ whose every element is either a collection of the form
$(c_1,c_2,\ldots,c_n|c)$,  $n\geq 0$, $c_i,c\in C$, as well as a set $Y(C)^\cyc$, whose every object
is $(c_0,c_1,\ldots,c_n)$,   $n\geq 0$, $c_i\in C$.   Let finally set $Y(C):=Y(C)^\uncyc\sqcup Y(C)^\cyc$.  $Y(C)$ can be interpreted as the set of 'arities'. 

Suppose we have a family of objects $\cM(a)\in \cC$ for all $a\in Y(C)$   Define a $\cO$-bimodule structure
on $\cM$ as the structure of a circular operad on the direct sum $\cO\oplus \cM$ such that:
---
the natural inclusion/projecion maps $\cO\to \cO\oplus \cM\to \cO$ are  circular operad maps;

--- every composition map restricted onto $\cM(a_1)\otimes \cM(a_2)$, $a_1,a_2\in Y(C)$, vanishes.

One can define a category  $\bY(\cO)$ over $\cC$ whose set of objects is $Y(C)$ with the property that an $\cO$- bimodule structure 
on $\cM$ is equivalent to that of a functor $\cM:\bY(\cO)\to \cC$.  

Given a map of circular operads $\cO\to \cO_1$, $\cO_1$ has a natural structure of a $\cO$-bimodule on $\cO_1$, that is $\cO_1$ is a functor $\cO_1:\bY(\cO)\to \cC$.

In our case, we have  $\bZ\times \bZ$-equivariant algebras $A_\Phi\in \cM^\Phi$, $A_\Psi\in \cM^\Psi$,
$A_F\in \cM^F$.    Let $\cO^{A^\Phi}$ be a full circular operad of $A_\Phi$, and similar for $\cO^{A^\Psi}$ and
$\cO^{A^F}$. We now have maps  of circular operads
$$
\bbB\to\cO^{A^\Phi}\to \cO^{A^\Psi}\to \cO^F,
$$
whence an induced sequence of $\bY(\bbB)$-modules $\cO^{A^\Phi}\to \cO^{A^\Psi}\to \cO^F$ in $\Classic^\ZZ\langle \ve\rangle$.    The monoidal functor $$||:\Classic^\ZZ\langle \ve\rangle\to \Classic\langle \ve\rangle$$ gives rise to an induced sequence of $\bY(|\bbB|)$-modules in  $\Classic\langle \ve\rangle$:
\begin{equation}\label{phipsif}
|\cO^{A^\Phi}|\to |\cO^{A^\Psi}|\to |\cO^F|.
\end{equation}

We will now review some properties of these modules.   
\subsubsection{Condensation}   We have a map of circular operads $\assoc\to \bbB$   which gives a functor
$\bY(\assoc)\to \bY(\bbB)$.  Given a functor $\cM:\bY(\bbB)\to \cC$, we therefore have an induced structure
$\cM:\bY(\assoc)\to \cC$.  The latter structure give  rise
  to cosimplicial structures on $\cM^\uncyc$ and $\cM^\cyc$.     Call $\cM:\bY(\bbB)\to  \cC$  {\em quasi-constant}
if every arrow $f\in \Delta$ acts on both $\cM^\uncyc$ and $\cM^\cyc$ by quasi-isomorphisms.

Call $\cN:\bY(\bbB)\to \cC$ {\em right orhogonal to quasi-constant objects}  if  $R\hom(\cM;\cN)\sim 0$
for every quasi-constant object $\cM$.   

One can define an endofunctor $\con$ on the category of functors $\bY(\bbB)\to \cC$ (enriched over $\sets$)
as well as a natural transformation
$\con \to \Id$  such that $\con(\cM)$ is quasi-constant for every $\cM:\bY(\bbB)\to \cC$  and the 
cone of the induced map $\con(\cM)\to \cM$ is right orthogonal to quasi-constant objects.
In a homotopical sense,  $\con(\cM)$ is therefore the  universal quasi-constant object  mapping into $\cM$.

\begin{Claim}  Apply $\con$ to (\ref{phipsif}):
$$
\con |\cO^{A^\Phi}|\to \con |\cO^{A^\Psi}|\to\con |\cO^F|.
$$
All the arrows in this sequence are termwise equivalences.
\end{Claim}

This claims allows one to prove the propeties of $\con |\cO^F|$ via passing to $\con |\cO^{A^\Phi}|$.
\subsubsection{The $c_1$-action} {\em  In this section our basic category is that of complexes of $\bbQ$-vector spaces.}

 It is well known that every cyclic object $X$ has a natural endomorphism
$c_1:X\to X[2]$,  'the first Chern class'.       One has a similar property for every functor $\bY(\bbB)\to \cC$.
In particular, we have a $c_1$-action
on the object  $$
H:=R\hom(\con |\cO^{A^F}|;\con |\cO^{A^F}|)\in  \Classic\langle \ve\rangle.
$$
we show that such an action on every non-zero graded component $\gr^{n\ve}H$, $n\neq 0$ is homotopy nilpotent.
That is there exists a number $N(n)\geq 0$ such that the map $c_1^{N(n)}:\gr^{n\ve}H\to \gr^{n\ve}H[2N(n)]$
is homotopy equivalent to 0.

Finally, one studies  an object $$T:=R\hom_{\bY(|\bbB|)}(|\bbB|;\con|\cO^{A^F}|)\stackrel\sim\to 
R\hom_{\bY(|\bbB|)}(|\bbB|;|\cO^{A^F}|)$$ which also carry a $c_1$-action.
We prove that there is a homotopy equivalence of objects with $c_1$-action
$$
G[u]\to T,
$$
where $u$ is a variable of degree $+2$, and $c_1$ acts on the LHS by the multiplication by $u$.

\subsection{Dependence on a pseudo-Kaehler metric}  The above outlined approach depends on the choice of
a pseudo-Kaehler metric $g$.  In order to make the construction invariant under symplectomorphisms of $M$,
one coniders the set of all pseudo-Kaehler metrics on $M$ and carries over all the steps in this settings.
One cannot use the category $\Classic\langle \ve\rangle$ anymore, as one cannot choose a $\ve$ uniformly
for all  the metrics.   One however can build  ground categories $\bR_q,\bR_0$ as a replacement.

\section{Action of $\Sp(2N)$}\label{group} Let $G$ be the  universal cover of $\Sp(2N)$.   Let $V=\Re^{2N}$ be the standard symplectic
vector space with the coordinates $(q,p)$ and let $E=\Re^N$ so that $V=T^*E$.    The group $\Sp(2N)$, hence $G$,
acts on $V$.
\subsection{Graph of the  $G$-action on $T^*E$}.     Let  $a:T^*E\to T^*E$ be the antipode map
$(q,p)\mapsto (q,-p)$.  
   Let $\Gamma\subset G\times V\times V$
 consist of all points of the form  $\{(g,v,gv^a)|\ g\in \Sp(2N);v\in V\}$.  It follows that there exists 
a unique
  Legendrian sub-manifold $\cL\subset T^*(G\times E\times E)\times \Re$ which:

--- diffeomorphically projects onto $\Gamma$ under the projection $$
T^*(G\times  E\times E)\times \Re\to G\times T^*(E\times E).
$$

--- contains all the points of the form $(e,v,v^a,0)$, where $e$ is the unit of $G$ and $v\in V$.

Let $\cC$ be the full sub-category of 
$\sh_q(G\times E\times E)$ consisting of all objects $F$ with $\SS(F)\subset \cL$
such that there exists a homotopy equivalence in $\sh_q(E\times E)$:
$$F|_{e\times E\times E}\stackrel\sim\to \ZQ_{\Delta_E\times [0,\infty)},
$$
where $\Delta_E\subset E\times E$ is the diagonal.

Let $\cA_\unit\subset \cA$ be the full sub-category of all objects isomorphic to $\unit_\cA$.
We have a functor
$$
\Psi:\cC\to \cA_\unit,
$$
where $$\Psi(F):=\hom(\ZQ_{\Delta_E\times [0,\infty)};F|_{e\times E\times E}).
$$

\begin{Theorem}    The functor $\Psi$ is an equivalence of categories.
\end{Theorem}

Sketch of the proof.   {\em Part 1: Let us construct at least one object $\bbS$ of $\cC$} 

1)  For an open subset $U\subset G$, let $\cL_U\subset T^*(U\times E\times E)\times \Re)$ be the restriction of $\cL$.
Let $\cC_U$ be
the full sub-category of $\sh_q(U\times E\times E)$ consisting of all objects $F$
such that $\SS(F)\subset \cL_U$ and  there exists a homotopy equivalence $F|_{e\times E\times E} \sim \ZQ_{\Delta_E\times [0,\infty)}$.

 2) Let $\cU$ be a small enough geodesically convex neighborhood of unit in $\Sp(2N)$ satisfying:   for each $g\in \cU$ we have:
$(q,p')$ is a non-degenerate system of coordinates, where $(q',p')=g(q,p)$.     $\cU$ lifts uniquely to $G$, to be denoted by the same letter.

 3) Let $\bF,\bF'\in \sh_q(E\times E)$  be the Fourier kernels:  $$
\bF:=\ZQ_{\{(q_1,q_2,t)|t+\langle q_1,q_2\rangle\geq 0\}};\quad \bF':=\ZQ_{\{(q_1,q_2,t)|t-\langle q_1,q_2\rangle\geq 0\}}[2N].
$$

 Let $R:T^*E\times \Re\to T^*E\times \Re$ be  the 'Fourier' contactomorphism given by
$$R(q,p,t)=(-p,q,t+\langle p,q\rangle).
$$

  Let $$
R_1:T^*E\times T^*E\times \Re\to T^*E \times  T^*E\times \Re,
$$
be defined by $R_1(u_1,u_2,t)=(u_1,R^{-1}(u_2,t))$.
  Let $\cC'_\cU\subset \sh_q(\cU\times E\times E)$ consist of all objects $F$ such that

---  there exists a homotopy equivalence
$F|_{e\times E\times E}\sim \bF'$.

--- $\SS F\subset R_1(\cL_\cU)$.

It follows that the functor $G\mapsto G*_E \bF$ induces a homotopy equivalence of categories
$\cC'_\cU\to \cC_\cU$.

4) The Legendrian manifold $R_1\cL_\cU\subset T^*(G\times E\times E)\times \Re$ projects uniquely onto the base
$G\times E\times E$, therefore, $R_1\cL_\cU$ is of the form $\cL_f$ for some smooth function $f$ on $G\times E\times E$. 

Let $\bbB\subset \sh_q(\cU\times E\times E)$ be the full sub-category of objects $F$ satisfying:

--- $\SS(F)\subset (T^*_{\cU\times E\times E} \cU\times E\times E)\times 0$;

--- there exists a homotopy equivalence $F|_{e\times E\times E}\sim \ZQ_{E\times E\times [0,\infty)}$.

It follows that $\bbB$ is the category consisting of all objects homotopy equivalent to $\ZQ_{\cU\times E\times E\times [0,\infty)}$.

  The convolution with $\ZQ_{\{e,e,t|,t+f(e)\geq 0\}}$ gives a homotopy equivalence of 
categories $\bbB\to \cC'_\cU$.

 Fix an object $S_\cU\in \cC_\cU$ along with a homotopy equivalence $$S_\cU|_{0\times E\times E}\sim \ZQ_{\Delta_E\times [0,\infty)}.
$$

5) For $h\in \cU$, set $S_h:=S_\cU|_{h\times E\times E}$.   Every $g=G$ can be written as 
$g=g_1g_2\cdots g_n$, where $g_i,g_i^{-1}\in \cU$.

  Set $S_{g_1,\ldots ,g_n}=S_{g_1}*_E S_{g_2}*_E\cdots *_E S_{g_n}$.     

For each $g$, choose an object $S_{g\cU}$ which is homotopy equivalent to one of $S_{g_1,\ldots, g_n}*_E S_\cU$ for
$g_1\cdots g_n=g$.    Observe that the objects $S_{g_1,\ldots,g_n}$ and $S_{g'_1,\ldots,g'_m}$, where
$g_1\cdots g_n=g'_1\cdots g'_m=g$ are homotopy equivalent.    It suffices to show
that $$S_{g_1,\cdots, g_m,(g'_{m})^{-1},\cdots,(g'_1)^{-1}}\sim \ZQ_{[\Delta_E,0]}$$ that is
$S_{g_1g_2\cdots g_n}=\ZQ_{\Delta_E,0}$ whenever $g_1g_2\cdots g_n=e$.  
As $U$ is geodesically closed, there is a unique shortest geodesic line joining $g_1\cdots g_k$ and 
$g_1\cdots g_{k+1}$.  We will thus get a  broken  geodesic line starting and terminating at $e$.
As $G$ is simply connected, this line can be contracted to a point.   By possibly adding intermediate points, 
one can reduce the problem to the case when there exist smooth paths $h_k:[0,1]\to U$ such that
$h_1(t)\cdots h_n(t)=e$, $h_k(1)=e, h_k(0)=g_k$ for all $k$.   Let $S_k\in \sh([0,1]\times E\times E)$,
$S_k:=h_k^{-1} S_\cU$.   Consider 
$$
\Sigma:=S_1*_E S_2*_E\cdots *_E S_n\in \sh_q(I^n\times E\times E)|_{\Delta_I\times E\times E},
$$
where $\Delta_I\subset I^n$ is the diagonal.

It follows that $$\Sigma_{1\times E\times E}\sim\ZQ_{[\Delta_E,0]};\quad \Sigma_{0\times E\times E}\sim S_{g_1,g_2,\cdots,g_n}.
$$
Next,  the singular support estimate shows that $\Sigma$ is locally constant along  $\Delta_I$, which implies the statement.

6) Choose a covering $G=\bigcup_n g_n\cU$.  Let $I$ be the poset consisting of all non-empty intersections
$g_{i_1}\cU\cap \cdots g_{i_k}\cU$.  Each element of $I$ is geodesically convex.  It follows that
all the restictions $S_{g_{i_l}\cU}|_{g_{i_1}\cU\cap \cdots \cap g_{i_k}\cU}$ are homotopy equivalent. 
Indeed,  choose a point $h\in g_{i_1}\cU\cap \cdots \cap g_{i_k}\cU$;  4) implies that there is a homotopy equivalence of restrictions $S_{g_{i_l}\cU}|_{h\times E\times E}$  with $S_h$.  The statement now follows from 4).

  For every $A\in I$, $A=g_{i_1}\cU\cap g_{i_2}\cU\cap \cdots\cap g_{i_k}\cU$,  choose an object
$S_{A}\in \cC_A$
to be homotopy equivalent  to each of the restrictions $S_{g_{i_l}\cU}|_{A\times E\times E}$.

7) For each $V\in I$ let $j:V\to G$ be the embedding. Let $T_V:=j_{\cV!} S_\cV$.

8) Whenever $A\subset B$, $A,B\in I$,  we have a homotopy equivalence $\ZQ\sim \hom(T_A,T_B)$.  
Let  $r_{AB}:T_A\to T_B$ be the image of $1\in \ZQ$.

9) We have $r_{BC}r_{AB}$ is homotopy equivalent  to $E_{ABC}r_{AC}$ for some  $E_{ABC}\in \ZQ^\times $.

10)  $E_{ABC}$ is a 2-cocycle on $I$. Since $H^2(G,\ZQ^\times)=0$,  $E_{ABC}$ is exact.  Therefore, wlog we can assume that $E_{ABC}=1$.

11)  Denote $\cJ(A,B):=\tau_{\leq 0}\hom(T_A,T_B)$, see Sec \ref{groundcat} for the definition of $\tau_{\leq 0}$.

  We have a  functor $\cJ\to I$ which is a homotopy equivalence
of categories so that we have the constant functor
$
Z:\cJ^\op\to I^\op\to \cA,
$
$Z(A)=\ZQ$ for all $A$.

Finally, we set $\bbS:=S_G:=\cZ\otimes^L_{\cJ^\op} S$.

{\em Part 2. Uniqueness}   The convolution with $\bbS$ gives a pair of quasi-inverse maps between $\cC_G$ and the full sub-category
of objects $S\in \sh_q(G\times E\times E)$ with $\SS  S\subset T^*_G G\times T^*_E(E\times E)\times \{0\}$, where there exists
an isomorphism
$$S|_{e\times E\times E}\sim \ZQ_{\Delta_E\times [0,\infty)}.
$$

  Passing  to $\tau_{\leq 0}$ yields the statement.

\subsubsection{The object $\bbS$}\label{bbs} Fix an object $\bbS\in \cC$ endowed with a homotopy equivalence
$$\bbS|_{e\times E\times E}\to \ZQ_{\Delta_E\times [0,\infty)}.
$$
\section{Objects supported on a symplectic ball}\label{symball}
\subsection{Projector onto the ball}\label{xigamma} Let $i_0:\Re/2\pi\bZ\to \Sp(2N)$ be a one-parametric subgroup  consisting of all transformations
$$
q'=q\cos(2a)+p\sin(2a);
$$
$$
p'=-q\sin(2a)+p\cos(2a).
$$
Let 
$i:\Re\hookrightarrow G$ be the lifting.  Denote $\bbI:=i(\Re)$.    Let $\cT\in \shinf(\bbI\times E\times E)$ be the restriction of $\bbS$.
The object $\cT$ is microsupported within the set
\begin{equation}\label{sig}
\Sigma=\Sigma_0\cup \{(a,-(q^2+p^2),q,-p,q',p',-S(q,p,a))|(q,p)\in V;a\in \Re,\sin(2a)\neq 0\}\subset T^*\bbI\times T^*E\times T^*E\times \Re
\end{equation}
where
$$
\Sigma_0=\{(\pi n,-(q^2+p^2),q,-p,q,p,0)|(q,p)\in V,n\in \bZ\}\cup \{(\pi(\frac 12 +n);-(q^2+p^2),q,-p,-q,-p,0)|(q,p)\in V,n\in \bZ\};
$$
$$
S(q,p,a)=\frac{\cos (2a)(q^2+(q')^2)+2qq'}{2\sin( 2a)}.
$$

Let $\cB=\Re$ with the coordinate $b$.  Let $p_{\cB}:\cB\times E\times E\to E\times E$ be the projection.
Set $$\cP_R:=p_{\cB!}\cT*_\bbI \ZQ_{\{(a,b,t)\in \bbI\times \cB|b<R^2,t+ab\geq 0\}}[1]\in \sh_q(E\times E).$$

We have $$\cP_R\sim p_{\bbI!} \cT*_{\bbI} \ZQ_{\{(a,a,t)|a\leq 0,t+aR^2\geq 0\}},
$$
where $p_\bbI:\bbI\times E\times E\to E\times E$ is the projection.
\subsubsection{The map $\alpha:T_{-\pi R^2}\cP_R[2N]\to \cP_R$}
We have a homotopy equivalence
$$
T^a_{-\pi}\cT[-2N]\sim  \cT,
$$
where $T^a_{-\pi}$ is the translation along $\bbI$ by $-\pi$ units.

Thus, we have a map
 \begin{multline*}
\cP_R\sim p_{\bbI!} (T^a_{-\pi}\cT)*_{\bbI}\ZQ_{\{(a,a,t)|a\leq 0,t+aR^2\geq 0\}}[-2N])
\sim p_{\bbI!}(\cT*_{\bbI}\ZQ_{\{(a_1,a_2,t)| a_2\leq 0,a_1=a_2+\pi,t+a_2 R^2\geq 0\}})[-2N]\\
\sim p_{\bbI!}(\cT*_\bbI \ZQ_{\{(a_1,a_2,t)|a_1\leq \pi;a_1=a_2, t\geq \pi R^2-a_1R^2\}})[-2N]\\
\sim  T_{\pi R^2}p_{\bbI!}(\cT*_\bbI  \ZQ_{\{(a_1,a_2,t)|a_1\leq \pi;a_1=a_2,t+a_1R^2\geq 0\}})[-2N]\\
\to T_{\pi R^2}p_{\bbI!}(\cT*_\bbI  \ZQ_{\{(a_1,a_2,t)|a_1\leq 0;a_1=a_2,t+a_1R^2\geq 0\}})[-2N]\sim T_{\pi  R^2}\cP_R[-2N].
\end{multline*}
This map can be rewritten as
$$
\alpha:T_{-\pi R^2}\cP_R[2N]\to \cP_R.
$$

\subsubsection{$\hom(T_c\cP_R;\cP_R)$}   Let $(\nu-1)\pi R^2<c\leq \nu\pi R^2$, where $\nu\in \bZ$. Let $G_c:=\hom(T_c\cP_R;\cP_R)$. Then
$$
G_c\sim\bZ[-2N\nu] \text{  if } \nu\leq 0,  \quad   G_c=0 \text{ if } \nu>0.
$$

 The natural map  $G_{\nu \pi R^2}\to G_c$ is a homotopy equivalence.   The generator of $G_{\nu \pi R^2}$, $\nu<0$ is given by
$\alpha^{* \nu}$.

The map $\cP_R\to \ZQ_{[\Delta_E,0]}$ induces a homotopy equivalence
$$
\hom(T_c\cP_R;\cP_R)\to \hom(T_c\cP_R;\ZQ_{\Delta_E\times [0,\infty)}).
$$
\subsubsection{$\cP_R$ is a projector} We have a natural map 
\begin{equation}\label{proek}
\pr:\cP_R\to \ZQ_{\Delta_E\times [0,\infty)}.
\end{equation}   Let $\cC_R\subset \sh_q(E)$ be the full subcategory of objects
supported away from $\IntB_R\times \Re\subset T^*E\times \Re$.   Let $\sh_q[\IntB_R]   \subset \sh_q(E)$ be the left orthogonal complement
to $\cC_R$.  We have $\cP_R*_E F\in \sh_q[\IntB_R]$; $\Cone \cP_R*_E F\to F\in \cC_R$ so that $\cP_R$ gives a semi-orthogonal decomposition.

\subsubsection{Generalization}\label{orthog} Denote by $\sh_{1/2^n}[T^*X\times \IntB_R\times \Re]\subset \sh_{1/2^n}(X\times E)$ be the left orthogonal complement to the  full category
of objects supported away from $T^*X\times \IntB_R\times \Re$.  The convolution with $\cP_R$ gives a semi-orthogonal decomposition.

\subsubsection{The object $\gamma=\Cone \alpha$}\label{gamma}  Let $\gamma:=\Cone \alpha$.  
We have 
$$
\gamma\sim \cT*_\bbI \ZQ_{\{(a_1,a_2)|a_1=a_2; -\pi R^2<a_1\leq 0,t+aR^2\geq 0\}}
$$

We have a homotopy equivalence

\begin{multline*}
E_c:=\hom(T_c\gamma,\cP_R)\stackrel\sim\to \hom(T_c\gamma;\ZQ_{\Delta_E\times [0,\infty)})
\sim \hom(\Cone(T_{c-\pi R^2}\gamma[2N]\to T_c\gamma);\ZQ_{\Delta_E\times [0,\infty)})\\
\sim
\Cone (\hom(T_c\gamma;\ZQ_{\Delta_E\times [0,\infty)})\to \hom(T_{c-\pi R^2}\gamma;\ZQ_{\Delta_E\times [0,\infty)})[-2N])[-1]\\
\sim  
 \Cone( G_c\to  G_{c-\pi R^2}[-2N])[-1],
\end{multline*}
where the map is incuced by the multiplication by $\alpha$.

Therefore,

---$E_c=\ZQ[-2N-1]$,  $0<c\leq \pi R^2$;

---$E_c=0$ otherwise.

\subsubsection{Singular support of $\gamma$} We have
$$
\SS \cT*_\bbI \ZQ_{\{(a_1,a_2)|a_1=a_2, -\pi R^2<a_1\leq 0,t+aR^2\geq 0\}}\subset 
\{(a,R^2+k,q,-p,q',p',t-aR^2)\in \Sigma|-\pi<a<0\}\cup  S,
$$
where $\Sigma$ is as in (\ref{sig}) and 
$$
S=\{(-\pi,R^2+k,q,-p,q,p,-\pi R^2)|k\leq -p^2-q^2\}\cup \{(0,R^2+k,q,-p,q,p,0)|k\leq -p^2-q^2\}.
$$
Therefore,  we have
\begin{multline*}
\SS \gamma\subset \{(q,-p,q',p',-aR^2-S(a,q,q'))|p^2+q^2= R^2;-\pi<a<0\}\cup
\{(q,-p,q,p,-\pi R^2)| q^2+p^2\leq R^2\}\\
\cup \{(q,-p,q,p,0)|q^2+p^2\leq R^2\}.
\end{multline*}
It follows that $0\leq -aR^2-S(a,q,q')\leq \pi R^2$  if $-\pi<a<0$.
\subsubsection{Singular support of $\cP$}
Similarly, one can find
\begin{multline*}
\SS \cP\subset \{(q,-p,q',p',-aR^2-S(a,q,q')|p^2+q^2=R^2;a<0\}\cup \{(q,-p,q,p,0)|q^2+p^2\leq R^2\}.
\end{multline*}
\subsubsection{Singular support of $\Cone \cP\to \ZQ_{\Delta_E\times [0,\infty)}$}
We have 
$$
\Cone (\cP\to \ZQ_{\Delta_E\times [0,\infty)})\approx p_{\bbI!}\cT*_\bbI \ZQ_{\{(a_1,a_2)|a_1=a_2,a_1\leq 0,t+aR^2\geq 0\}}
$$
so that 
\begin{multline*}
\SS \cT*_\bbI \ZQ_{\{(a_1,a_2)|a_1=a_2,a_1\leq 0,t+aR^2\geq 0\}}\subset   \{(a,R^2+k,q,-p,q',p',t-aR^2)\in 
\Sigma|a<0\}\cup  S',
\end{multline*}
where $\Sigma$ is as in (\ref{sig}) and 
$$
S'=\{(0,R^2+k,q,-p,q,p,0)|k\geq -p^2-q^2\}.
$$

Therefore,
$$
\SS \Cone(\cP\to \ZQ_{\Delta_E\times [0,\infty)})\subset \{(q,-p,q',p',-aR^2-S(a,q,q')|p^2+q^2=R^2;a<0\}\cup \{(q,-p,q,p,0)|q^2+p^2\geq R^2\}.
$$
\subsubsection{Corollaries} 
\begin{Corollary} \label{konusp} 
 We have $$ R_{\leq c}\Cone(\cP\to \ZQ_{\Delta_E\times [0,\infty)})\approx 0;$$ 
$$R_{\leq c}\Cone(\cP\boxtimes \cP\to \ZQ_{\Delta_E\times \Delta_E\times [0,\infty)})\approx 0.$$
for all $c\leq 0$.
\end{Corollary}
\begin{Corollary}
 Let $F\in \sh(E\times E)$.  Then the natural maps
$$
\hom(\ZQ_{\Delta_E\times \Delta_E};F)\stackrel\sim\to \hom(\cP\boxtimes \cP;F\boxtimes \ZQ_{[0,\infty)});
$$
$$
\hom(\ZQ_{\Delta_E\times \Delta_E};F)\stackrel\sim\to \hom(T_{2\pi R^2}\gamma\boxtimes \gamma[-4N];F\boxtimes\ZQ_{[0,\infty)})
$$
are homotopy equivalences.
\end{Corollary}

\subsubsection{Convolution of $\gamma$ with itself}\label{obrez}
 We have a homotopy equivalence
$$
\gamma*_E \gamma\sim \gamma \oplus T_{-\pi R^2}\gamma[2N].
$$

Denote by $\mu:\gamma*_E \gamma\to \gamma$ the projection.

We now have the following  homotopy equivalence
$$
\hom(T_c\gamma,\ZQ_{\Delta_E\times [0,\infty)})\stackrel{\mu}\to \hom(T_c\gamma*_E\gamma; \ZQ_{\Delta_E\times [0,\infty)}),
$$
for all $c$ except those in $(\pi R^2,2\pi R^2]$.

In particular,  for $0<c\leq \pi R^2$,  we have:
$$
\hom(T_c\gamma*_E \gamma; \ZQ_{\Delta_E\times [0,\infty)})\sim\ZQ[-2N-1];
$$
For $c\leq 0$, the above expresssion is homotopy equivalent to 0.

 Let $\Lambda\in \sh_q(\pt)$;  $\Lambda=\Cone (\ZQ_{[-\pi R^2,\infty)}\to \ZQ_{[0,\infty)})$.

We  have a chain of homotopy equivalences
$$
\hom(\gamma;\ZQ_{\Delta_E}\boxtimes \Lambda)\stackrel{\mu}\sim \hom (\gamma*_E \gamma; \ZQ_{\Delta_E}\boxtimes \Lambda)\sim \ZQ[-2N].
$$
In particular, we have a homotopy equivalence
$$
\hom(\gamma,\Lambda\boxtimes \ZQ_{\Delta_E}[2N])\sim \ZQ.
$$
Let \begin{equation}\label{nu}
\nu:\gamma\to \Lambda\boxtimes \ZQ_{\Delta_E}[2N]
\end{equation}
be the generator.

One also has a map $\ve:\Lambda\boxtimes\ZQ_{\Delta_E}\to \gamma$ which has a homotopy unit property with respect to $\mu$,
the through map
$$
\gamma\sim \ZQ_{\Delta_E}*_E\gamma\to \Lambda\boxtimes \ZQ_{\Delta_E}*_E \gamma\to \gamma*_E \gamma\to \gamma
$$
is homotopy equivalent to the  Identity.

The induced map
\begin{equation}\label{ekviva}
\hom(\gamma,\Lambda\boxtimes \ZQ_{\Delta_E})\stackrel\ve\to \hom(\gamma,\gamma)
\end{equation}
is a homotopy equivalence.
The map $\nu$ on the LHS corresponds to $\Id$ on the RHS.
\subsubsection{Lemma on $\nu\boxtimes \nu$} \label{digamma}  Consider the following maps
\begin{equation}\label{oto1}
\gamma\boxtimes \gamma\stackrel{\nu\boxtimes \nu}\to \Lambda\boxtimes \ZQ_{\Delta_E}\boxtimes \Lambda\boxtimes \ZQ_{\Delta_E}[4N]
\to \Lambda\boxtimes \ZQ_{\Delta_E\times \Delta_E}[4N];
\end{equation}

\begin{equation}\label{oto2}
\gamma\boxtimes \gamma\stackrel{\overline{\mu}} \to
p_{14}^{-1}\gamma \boxtimes p_{23}^{-1}\ZQ_{\Delta_E}\stackrel\nu\to \Lambda\boxtimes p_{14}^{-1}\ZQ_{\Delta_E}\boxtimes p_{23}^{-1}\ZQ_{\Delta_E}[3N]\to \ZQ_{\Delta_E\times \Delta_E}[4N].
\end{equation}

Here the maps $\overline{\mu}$ is obtained from $\mu$ by conjugation.  The last arrow
is the generator of  $$\hom(p_{23}^{-1}\ZQ_{\Delta_E}\otimes p_{14}^{-1}\ZQ_{\Delta_E}; \ZQ_{\Delta_E\times \Delta_E}[N]).
$$
\begin{Lemma} The maps (\ref{oto1}) and (\ref{oto2}) are homotopy equivalent.
\end{Lemma}
{\em  Sketch of the proof}  One  reformulates the statement as follows:

By the conjugacy, the  map $\vu$ correponds to a homotopy equivalence
$$
\xi:\Lambda\to \gamma*_{E^2} \ZQ_{\Delta_{\Re^n}\times [0,\infty)}[n] 
$$

The problem reduces to showing that the map
\begin{equation}\label{otb1}
\Lambda\to \Lambda\approx (\gamma\boxtimes \gamma) *_{E^4} \ZQ_{\Delta\times \Delta}[2n]\to
(\gamma\boxtimes \gamma)*_{E^4}
\ZQ_{(v_1,v_2,v_3,v_4)\in E^4|v_1=v_4;v_2=v_3}[n]\approx (\gamma*_E \gamma)*_{E^2}\ZQ_{\Delta}[n]\to 
\gamma*_{E^2} \ZQ_\Delta[n]
\end{equation}

is homotopy equivalent to
\begin{equation}\label{otb2}
\Lambda\otimes \Lambda \to \Lambda\to \gamma*_{E^2} \ZQ_\Delta[n].
\end{equation}

We have a homotopy equivalence,
$$
\gamma*_{E_2} \ZQ_{\Delta}[n]\cong \hom(\ZQ_\Delta;\gamma).
$$

The map $\xi$ rewrites as $\xi':\Lambda\to \hom(\ZQ_{\Delta};\gamma)]$ which produces a
map $e:\Lambda\otimes \ZQ_{\Delta}\to \gamma$. 

The map (\ref{otb1}) rewrites as 
$$
\Lambda\otimes \Lambda \to \hom(\ZQ_\Delta;\gamma)\otimes \hom(\ZQ_\Delta;\gamma)\to 
\hom(\ZQ_\Delta;\gamma*_E \gamma)\to \hom(\ZQ_\Delta;\gamma).
$$

The map (\ref{otb2}) rewrites as
$$
\Lambda\otimes \Lambda\to \Lambda\to \hom(\ZQ_{\Delta};\gamma).
$$
Homotopy equivalence of the two maps follows from the following maps being homotopy equivalent:
$$
\Lambda\otimes\ZQ_{\Delta}*_E \Lambda\otimes \ZQ_{\Delta}\stackrel{e*e}\to  \gamma*_E\gamma\to \gamma
$$
and
$$
\Lambda\otimes \Lambda\to \Lambda \to \gamma.
$$
The latter statement follows from Sec. \ref{obrez}.

\subsubsection{$\gamma$ as an object of $\sh_{\pi R^2}(E\times E)$} \label{gam} 
We assume that $\pi R^2=1/2^n$ for some $n\in\bZ_{\geq 0}$.
It follows that $\gamma$ is supported within the set $E\times E\times [-\pi R^2;0]$.   Therefore, $\gamma$ determines
an object of $\sh_{\pi R^2}(E\times E)$, to be denoted by $\Gamma$.

Using the bar-resolution for $\Gamma*_E \Gamma$,   we see that it is glued of $\gamma*_E\Lambda^{*_E^ n}*_E \gamma$.   
We therefore have the following homotopy equivalences (all the  hom's are in $\sh_{\pi R^2}(E\times E)$:
$$
\hom(\Gamma;\ZQ_{\Delta_E})\stackrel{\xi}\sim \hom (\Gamma*_E \Gamma; \ZQ_{\Delta_E})\sim \ZQ[-2N].
$$
  
\subsection{Study of the cateory $\sh_q(F\times E\times E)[T^*F\times \Int B_R\times \Int B_R\times \Re]$}
\subsubsection{The category $\bbB_I$}
Let $I\subset \Re$ be an open subset.  Denote by $\bbB_I$ the full sub-category of 
$$\sh_q(F\times E\times E)[T^*F\times \Int B_R\times \Int B_R\times \Re]$$ consisting of all objects $X$, where
$$
\SS(X)\cap T^*F\times \Int B_R\times \Int B_R\times I=\emptyset.
$$

\subsubsection{Study of  $\bbB_{(a,\infty)}$ }      

Let $F\in \bbB_{(a,\infty)}$.

We have a natural map
$$
 F* (P_R\boxtimes P_R)\to (R_{\leq a} F)* (P_R\boxtimes P_R),
$$
where $R_{\leq a}$ is as in Sec \ref{rleq}.
\begin{Lemma} The above map is a homotopy equivalence.
\end{Lemma}
{\em Sketch of the proof}   Equivalently,  we are to show $$ (R_{>a} F)*(P_R\boxtimes P_R)\sim 0.
$$
We have $$\hocolim_{c\downarrow a} R_{>c}F\stackrel\sim \to  R_{>a} F,$$
therefore, it suffices to show that

\begin{equation}\label{prpr}
(R_{>c}F)* (P_R\boxtimes P_R)\sim 0.
\end{equation}

Let us study $\SS R_{>c} F.$  Denote $\Phi:=F\times E\times E$.
  The functor $R_{>c}$ descends onto $D_{>0}(F\times E\times E\times \Re)$,
where it is isomorphic to the functor 
$$
F\mapsto \Cone\left(i_{c!}F|_{\Phi\times [c,\infty)}\to F|_{\Phi\times c}\boxtimes \ZQ_{[c,\infty)}\right)[-1]
$$
where $$
i_{c!}:\Phi\times [c,\infty)\to  \Phi\times\Re
$$
is the embedding.    $F$ is therefore homotopy equivalent (as an object of $\sh(\Phi\times \Re)$)
to:
$$
\Cone(j_{c!}F|_{\Phi\times (c,\infty)}\to  F|_{\Phi\times c}\boxtimes \ZQ_{(c,\infty)})[-1],
$$
where 
$$
j_c:\Phi\times (c,\infty)\to \Phi\times \Re
$$
is the embedding.

It follows that $$\text{SS}( F|_{\Phi\times c}\boxtimes \ZQ_{(c,\infty)})\cap T^*_{>0}(\Phi\times \Re)=\emptyset.
$$

The object on the left hand side is isomorphic in $D(\Phi\times \Re)$ to
$F\otimes \ZQ_{t>c}$.  Let us use the SS estimate from KS.

We have $\text{SS}(F)$ contains no points of the form $(f,\eta,v_1,\zeta_1,v_2,\zeta_2,t,k)$,
where $k>0$ and $|(v_i,\zeta_i/k)|<R$, $i=1,2$.       Next, $\SS(\ZQ_{t>c})=\{(f,0,v_1,0,v_2,0,t,k)|k\leq 0,t\geq c, t>c\Rightarrow k=0\}$.  Therefore,  $\Phi\boxtimes \ZQ_{t>c}\in D(\Phi\times \Re\times \Re)$ is
non-singular along the diagonal $\Phi\times \Delta_\Re$. Hence, 
$
\text{SS}(F\otimes \ZQ_{t>c})$ is obtained via fiberwise adding of $\SS(F)$ and $\SS(\ZQ_{t>c})$. The resulting
sum  contains no points  of the form $(f,\eta,v_1,\zeta_1,v_2,\zeta_2,t,k)$, 
where $k>0$ and $|(v_i,\zeta_i/k)|<R$, $i=1,2$,
which implies (\ref{prpr}).

\subsubsection{Study of $\bbB_{(-\infty,a)} $}
Let $\tau_{\geq a},\tau_{<a}:\sh(\Phi\times \Re)\to \sh(\Phi\times \Re)$
be  given by $$
\tau_{\geq a}F=i_{a!}F|_{\Phi\times [a,\infty)};\quad  \tau_{< a}F=j_{a!}F|_{\Phi\times (-\infty,a)}.
$$

\begin{Lemma} Let $F\in \bbB_{(-\infty,a)}$.  Then $\tau_{<a} F\sim 0$.
\end{Lemma}
{\em Sketch of the proof}  It suffices to show  that $R_{\leq c} F\sim 0$ for all $c<a$. 
Similar to the previous Lemma, we deduce that $R_{\leq c}F $ is  non-singular on  the set
$$
T^*F\times \Int B_R\times \Int B_R\times \Re.
$$

Next, we have homotopy equivalences
$$
 R_{\leq c} (F *(P_R\boxtimes P_R))\stackrel\sim\to R_{\leq c}( R_{\leq c} F*(P_R\boxtimes P_R))\sim 0.
$$
This proves the statement.

\subsubsection{Study of $\bbB_{\Re\backslash a}$}
Let $b_R\subset E$ be the open ball of radius $R$ centered at 0.
We have functors 
$$
\alpha:\sh(F\times b_R\times b_R)\to \bbB_{\Re\backslash a},$$
where
$$
\alpha(S)=(S\boxtimes \ZQ_{[a,\infty)})*(P_R\boxtimes P_R);
$$
$$
\beta:\bbB_{\Re\backslash a}\to \sh(F\times b_R\times b_R),
$$
where
$$
\beta(T)=T|_{t=a}.
$$

\begin{Proposition}\label{OOO} The functors $\alpha, \beta$ establish homotopy inverse homotopy  equivalences 
of categories.
\end{Proposition}
{\em Sketch of the proof} 
Let  $S\in \bbB_{\Re\backslash a}$.  
 According to the  two previous subsections we have homotopy equivalences:
$$
S\sim R_{\leq a} S*(P_R\boxtimes P_R)\sim  (\tau_{\geq a}R_{\leq a} S)*(P_R\boxtimes P_R)\sim
((S\bullet \ZQ_{t\geq a})\boxtimes \ZQ_{[a,\infty)})*(P_R\boxtimes P_R),
$$
which implies the statement.

\subsubsection{$\SS(\alpha(F))$} \begin{Proposition}
 Let $C$ be a closed conic subset of $T^*F\times T^*b_R\times T^*b_R$.
 $F\in \bbB_{\Re\backslash a}$ and suppose $\SS(F)\cap T^*F\times B_R\times B_R\times a \subset C.$
Then $\SS (\alpha(F)\boxtimes \ZQ_{[a,\infty)})\subset C\times a$.
\end{Proposition}

\subsubsection{The category $\bbB_{\Re\backslash a,\Delta}$}  Let $\alpha:B_R\to B_R$ be the antipode map,
$\alpha(q,p)=(q,-p)$.
Let 
$$
\Delta^\alpha=\{(\alpha(v),v)|v\in \Int B_R\}\subset \Int B_R\times \Int B_R.
$$
Let $\bbB_{\Re\backslash a,\Delta}\subset \bbB_{\Re,\backslash a}$ be the full sub-category of objects
$X$ where 
$$
\SS(X)\cap T^*F\times \Int B_R\times \Int B_R\times \Re\subset
T^*_F F\times T^*_{\Delta^\alpha}(\Int B_R\times \Int B_R)\times a. 
$$

Let $A_F\subset \sh(F\times b_R\times b_R)$ be the full sub-category of objects $T$
where \begin{equation}\label{af}
\SS(T\boxtimes \ZQ_{[a,\infty)})\subset T^*_{F\times \Delta_{b_R}}(F\times b_R\times b_R\times a).
\end{equation}
According to the previous subsection,  we have a homotopy equivalence
$$
\beta:A_F\to \bbB_{\Re\backslash a,\Delta}.
$$
Furthermore, let  $\Loc(F)\subset \sh(F)$ be the full sub-category of objects $T$ where
$$\SS(T\boxtimes \ZQ_{[a,\infty)})\subset T^*_F F\times a.
$$
Let $
\gamma:\Loc(F)\to A_F
$  
be given by
$$
\gamma(S)=S\boxtimes \ZQ_{\Delta_{b_R}}.
$$
\begin{Lemma} $\gamma$ is a homotopy equivalence of categories.
\end{Lemma}

Therefore, 
\begin{Proposition}\label{OO}
the functor $\zeta:=\beta\gamma:\Loc(F)\to \bbB_{\Re\backslash a,\Delta}$ is a homotopy equivalence of
categories.
\end{Proposition}
 \section{Families of symplectic embeddings of a ball into $T^*E$}
Let $$\alpha:=\sum\limits_i p_idq_i$$
be the Liouville form on $B_R$.    Let $\theta=dt+\alpha$ be the contact form on $B_R\times \Re$.
Let $F$ be a smooth family and let 
$$
I:F\times B_R\to T^*E\times \Re
$$
be a smooth map such that, for all $f\in F$, the restriction $$
I_f:=\cI|_{f\times B_R}:f\times B_R\to T^*E\times \Re
$$
satisfies 
$$
I_f^*\theta=\alpha.
$$

Let $\Gamma_I\subset F\times \Int B_R\times T^*E\times \Re$ consist of all points of the form
$$
(f,(q,-p),\cI(f,q,p)).
$$

There is a unique  Legendrian manifold $\cL\subset T^*F\times \Int B_R\times T^*E\times \Re$ which 
projects uniquely onto $\Gamma$ under the projection
$$
T^*F\times \Int B_R\times T^*E\times \Re\to F\times \Int B_R\times T^*E\times \Re.
$$

Let $f\in F$.  Let $$J_f:f\times B_R\to T^*E\times \Re\to T^*E$$ be the through map, which is a symplectomorphic
embedding.  

Let $DI(f)\in \Sp(2N)$ be the  differential of $J_f$ in $0\in B_R$.  This way we get a map $DI:F\to \Sp(2N)$.
Suppose we have a lifting $\overline{DI}:F\to \ovSp(2N)$  of $DI$.      
Call such a collection of data $\cI:=(I,\overline{DI})$  {\em a graded family of symplectic embeddings of a ball
into $T^*E$}.

\subsection{The category $\cC_\cI$}  Let $\sh_q(F\times \Re^n\times \Re^n)[T^*F\times\Int  B_R\times T^*\Re^n\times \Re]
$ be the full sub-category of $\sh_q(F\times \Re^n\times\Re^n)$ consisting of all objects $F$ which are left
orthogonal to all objects non-singular on $T^*F\times \Int B_R\times T^*\Re^n\times \Re$, same as in
Sec \ref{orthog}.

Below we will study the full sub-category
$$\cC_\cI\subset \sh_q(F\times \Re^n\times \Re^n)[T^*F\times\Int  B_R\times T^*\Re^n\times \Re]
$$
consisting of all objects $T$ satisfying $\SS(T)\cap T^*F\times \Int B_R\times T^*\Re^n\times \Re\subset \cL$.
\subsection{Main Theorem}
Let $A_F$ be a category as in (\ref{af}).
\begin{Theorem}
  We have a homotopy equivalence between  the categories $\cC_\cI$ and
$A_F$.
\end{Theorem}
The proof of this theorem occupies the rest of the subsection.

 1) Extend $I$ to  $F\times [-1,1]\times B_R$ as follows.    For $t\in [-1,1]\backslash 0$, set
$$J(f,t,x)= \frac{I(f,tx)-I(f,0)}{t}+I(f,0).
$$
This map extends uniquely to a smooth map $J:F\times [-1,1]\times B_R\to T^*\Re^n$. 
The grading of $I$ extends uniquely to a grading $\cJ$ of $J$.

Let  $K=J|_{F\times 0}$.  It follows that $K$ is a family of linear symplectomorphisms of $T^*\Re^n$ restricted to
$B_R$.     The grading $\cJ$ determines uniquely a map 
\begin{equation}\label{mumu}
\mu:F\to \ovSp(2N)\times \Re.
\end{equation}

2)  For every $(f,t)\in F\times B_R$ we have a Hamiltonian vecor field  on $B_R$, namely
$\displaystyle\frac{dJ(f,t)}{dt}$.   Let $H_{(f,t)}$ be a smooth function on $B_R$ correpsonding
to this vector field and satisfying $H_{(f,t)}(0)=0$.   It follows that $H: F\times I\times B_R\to \Re$
is a smooth function.   It extends to a smooth function on $F\times I\times T^*E$ whose support projects
properly onto $F\times I$.   

3) Let $\chi:\Re\to [-1,1]$ be a non-decreasing  smooth function such that $\chi(t)=-1$ for all $t\leq -1$,
$\chi(t)=1$ for all $t\geq 1$, and $\chi(0)=0$.  Let 
$K(f,t)=J(f,\chi(t))$ and $h(f,t,v)=H(f,\chi(t),v)\chi'(t)$ so that $h(f,t,-)$  is the Hamiltonian function
of the vector field $\displaystyle\frac {dK(f,t)}{dt}$.   It follows that there exists a unique
family of symplectomorphisms $M:F\times \Re\times E\to E$ such that

a) $M|_{F\times 0}$ is the family of linear symplectomorphisms coinciding with $J|_{F\times 0}=K_{F\times 0}$;

b) $\displaystyle\frac{dM(f,t)}{dt}$ is the Hamiltonian vector field of $h(f,t,-)$.

It also follows that $M|_{F\times \Re\times B_R}=K$

4)  The family $M$ defines a Legendrian sub-manifold $\cL_M\subset T^*(F\times E\times E\times \Re)$
such that $\cL_M\cap T^*F\times B_R\times T^*E\times \Re=\cL_K$.  

5)  According to the theorem of Guillermou-Kaschiwara-Schapira,  there exists 
a quantization of $\cL_M$:  an object $Q\in \sh_q(F\times E\times E)$ such that $\SS Q\subset \cL_M$
and $Q|_{t=0}=\mu^{-1}\bS$, where $\mu$ is as in (\ref{mumu}).

6) Similarly, one defines a quantization $Q'$ of the family $M^{-1}$ of inverse symplectomorphisms.

7)   Let $\Delta:F\times E\times E\to F\times I\times F\times I\times E\times E$ be the following embedding
$$
\Delta(f,v_1,v_2)=(f,1,f,1,v_1,v_2).
$$
We  have    endofunctors
$$
S\mapsto S*_{F\times E} \Delta_{!} Q;\quad S\mapsto S*_{F\times E} \Delta_! Q'
$$
of $\sh_q(F\times E\times E)$ which descends to homotopy inverse homotopy equivalences 
between $
\cC_O$ and $ \cC_I$, where $O:F\times B_R\to B_R\stackrel\iota\to  E$ is the constant family, where
$\iota$ is the standard embedding.

By definition, $\cC_O=\bbB_{\Re\backslash 0,\Delta}$.  By Proposition \ref{OO} we have a homotopy equivalence
$\zeta:\Loc(F)\to \cC_O$.  We thus have constructed a zig-zag homotopy equivalence
between $\Loc(F)$ and $\cC_I$.     Denote by  $\cP_I\in \cC_I$ the object corresponding
to $\ZQ_F\in \Loc(F)$.

\subsubsection{Inverse functor} We have $\cP_I\in \sh_q(F\times b_R\times E)$.

Let $I':F\times B_R\to T^*E$ be given by $I'(f,v)=\alpha I(f,\alpha(v))$, where $\alpha: T^*E\to T^*E$,
 $\alpha(q,p)=\alpha(q,-p)$.
Let $\cQ_I:=\sigma_!\cP_{I'}\in \sh_q(F\times E\times b_R)$, where $\sigma:b_R\times E\to E\times b_R$ is
the permutation.

Let $\Delta_F:F\to F\times F$ be the diagonal embedding.
\begin{Proposition} We have $$\cQ_I*_{F\times E }\Delta_{F!}\cP_I\approx \ZQ_F\boxtimes \cP_R\in
\sh_q(F\times b_R\times b_R).
$$
\end{Proposition}

\subsubsection{} \label{semeistvo} Let $\pi:T^*F\to F$ be the projection.
Let $G_I\subset T^*F\times T^*E$ be an open subset defined as follows 
$$
G_I=\{(\phi,v)|v\in I(\pi(f)\times \Int B_R)\}.
$$

Let us also define functors $$
\bbP:\sh_q(F\times b_R)[T^*F\times \Int B_R]\to \sh_q(F\times E)[G_I];\quad \bbQ:\sh_q(F\times E)[G_I]\to
\sh_q(F\times b_R)[T^*F\times \Int B_R],
$$
where
$$
\bbP_I(S)=S*_{F\times b_R}\Delta_{F!} \cP_I;\quad \bbQ_I(T)=T*_{F\times E} \Delta_{F!} \cQ_I.
$$
\begin{Proposition} The functors $\bbP_I,\bbQ_I$ establish homotopy mutually inverse homotopy equivalences
between the categories $\sh_q(F\times b_R)[T^*F\times \Int B_R]$ and $\sh_q(F\times E)[G_I]$.
\end{Proposition}

\subsection{Pair of consequitive families} Let $u:F\times B_r\to B_R$, $v:F\times B_R\to E$
be graded families of symplectic embeddings.   Let
$w:F\times B_r\to E$ be defined by $w(f,b)=v(f,u(f,b)).$ 
The gradings define liftings $g_u:F\times B_r\to \ovSp(2N)$; $g_v:F\times B_R\to \ovSp(2N)$ of the corresponding
differential maps.

Let $g_w:F\times B_r\to \ovSp(2N)$ be given by  $g_w(f,b)= g_v(f,u(f,b))g_u(f,b)$. 
It follows that $g_w$ lifts the differential map $F\times B_r\to E$ determined by $w$.  Therefore, $g_w$ is
a grading of $w$.

\begin{Proposition} We have a homotopy equivalence $\bbP_v\circ \bbP_u\stackrel\sim\to \bbP_w$.
\end{Proposition}
{\em Sketch of the proof}  As above, let us extend the family $v$ to
a family 
$$v_t:F\times [-1,1]\times B_R\to E,$$
where 
$$v_t(f,t,b)=\frac{v(f,tb)-v(f,0)}{t}+v(f,0).
$$

Let $w_t:F\times [-1,1]\times B_r\to E$, where $w_t(f,t,b)=v_t(f,t,u(f,b))$. 
The gradings from $v$ and $w$ extend to $v_t,w_t$.
We will show that there exists a homotopy equivalence
\begin{equation}\label{restt}
\bbP_{v_t}\circ \bbP_u\stackrel\sim\to  \bbP_{w_t}.
\end{equation}
Restriction to $t=1$ will then show the Proposition.

To show the existence of (\ref{restt}), it suffices to establish the homotopy equivalence 
of the restriction to $t=0$.      Observe that $v_0$ comes from  a family of linear symplectomorphisms
$F\to \Sp(2N)$ whose grading defines a lifting $V_0:F\to \ovSp(2N)$.
Let $V\in \sh_q(F\times \times E\times E)$ be the corresponding object.
We have a homotopy equivalence
$$
\bbP_{v_0}\circ \bbP_u \sim  V\circ \bbP_u
$$
so the problem reduces to establishing a homotopy equivalence 
$V\circ \bbP_u \stackrel\sim\to \bbP_{v_0u}$.

In a similar way (via considering the family $u_t$), one reduces the problem to the case when the family $u$ is linear.
The grading then defines an object $U\in \sh_\infty(F\times E\times E)$.    Similarly, the linear family $v_0u$, along with
its grading, defines an object $W\in \sh_\infty(F\times E\times E)$.

 Next, we have
 homotopy equivalences
$U\circ \bbP_{B_r}\stackrel\sim\to \bbP_u$;  $W\circ \bbP_{B_r}\stackrel\sim\to \bbP_{v_0u}$ so that the problem
reduces to establishing a homotopy equivalence
$$
V\circ U\stackrel\sim\to W,
$$
which follows from Sec \ref{group}.

 \subsection{Mobile families} \label{mobile}
\subsubsection{Definition}
Let $U\subset T^*E$ be an open subset let $j:U\to T^*E$ be the corresponding open embedding.  Let $I:U\times B_R\to T^*E$ be a family of symplectic embeddings, where we assume
$I|_{U\times 0}=j$.   

The family $I$ defines a Lagrangian sub-manifold
$$
L_I\subset T^*U\times \Int B_R\times T^*E.
$$
Set $F=E\oplus E^*$.

We have a natural identification $T^*U=U\times F$.    For each $\xi\in U$ let $L_\xi:= T^*_\xi U \times \Int B_R\times T^*E
\cap L_I\subset F\times \Int B_R\times T^*E$.  Let $P_\xi\subset  F\times\Int  B_R$ be the image of $L_\xi$ under the projection along $T^*E$
Call $I$ {\em mobile} if for every $\xi$,  $P_\xi$ is a graph of an embedding $\Int B_R\to F$.

\subsubsection{Main proposition}

We have objects $\cP_I,\cQ_I\in \shinf(U\times E\times E)$.        Let $p_1,p_2:U\times E\times E\times E\times E\to U\times E\times E$
be the projections
$$
p_1(u,e_1,f_1,e_2,f_2)=(u,e_1,f_1);\quad p_2(u,e_1,f_1,e_2,f_2)=(u,e_2,f_2).
$$

Consider
$$
R_I:=p_1^{-1}\bbP_I\circ p_2^{-1} \bbQ_I.
$$

Let $i:E^3\to E^4;p:E^3\to E^2$ be given by $i(a,b,c)=(a,b,b,c)$; $p(a,b,c)=(a,c)$.   According to the previous subsection, we have
 a map
$$
p_!i^{-1} R_I\to \ZQ_{U\times \Delta_E\times [0,\infty)}
$$
where $\Delta_{E}\subset E\times E$ is the diagonal.

By the conjugacy, we have a map
$$
R_I\to \ZQ_{U\times \Delta_{14}\times \Delta_{23}\times [0,\infty)}[N],
$$
where $N=\dim E$ which, in turn, gives rise to a map
$$
\alpha:\pi_{U!}R_I\to  \ZQ_{\Delta_{14}\times \Delta_{23}\times [0,\infty)}[N],
$$
where
$
\pi_U:U\times E^4\to E^4$ is the  projection along $U$.

Let $V\subset U$ be an open subset satisfying:   for every $u\in U$,  if $I(u\times B_R)\cap V\neq \emptyset$,   then $I(u\times B_R)\subset U$.

Let $p_i:T^*E^4\to T^*E$ be the projections $i=1,2,3,4$.  Let $p_{ij}:=p_i\times p_j:T^*E^4\to T^*E^2$.

\begin{Proposition}       Let  $A,B\in \sh_q(E\times E)$ and assume that $\SS A\subset B_R\times B_R\times \Re$;  $\SS B\subset V$.
Then  $H:=(\Cone \alpha)*_{E^4} (p_{23}^{-1}A\circ p_{14}^{-1}B)\sim 0$.
\end{Proposition}

Sketch of the proof.    Let us define a family of symplectic embeddings 
$$
J:U\times (-1,1)\times B_R\to T^*E
$$
by means of dilations, same as above.     One then defines  an object $\pi_{U!}R_J\in \sh_q((-1,1)\times E^4)$, a map
$$
\alpha_J:\pi_{U!}R_J\to  \ZQ_{(-1,1)\times\Delta_{14}\times \Delta_{23}\times [0,\infty)}[N],
$$
and an object
$$
H_J:=(\Cone \alpha_J)*_{E^4} (p_{23}^{-1}A\otimes  p_{14}^{-1}B)\in \sh_q((-1,1)).
$$

Singular support estimate (see below) shows that
$$
\SS H_J\subset  T^*_I I\times \Re.
$$

Therefore, it suffices to show that $H_J|_0\sim 0$, in other words, the problem reduces to the case when $I$ is a family of linear symplectic embeddings.
The latter case can be reduced to the case when every embedding is a parallel transfer which is straightforward.

{\em Estimate of $\SS H_J$}.   It suffices to show that
 $$
\SS(\pi_{U!} R_J*_{E^4}(A\boxtimes B))\subset  T^*_{(-1,1)}(-1,1).
$$
  Let us identify
$$T^*(U\times \Re\times E^4)\times\Re= (U\times \Re)\times (F\oplus \Re)\times F^4\times \Re.
$$
We have
\begin{multline*}
\SS(R_J)\subset \{(\tau,\eta_J(\tau,v_1)-\eta_J(\tau,v_2), v_1^a,J(\tau,v_1),v_2^a,J(\tau,v_2))|\tau\in U\times\Re, v_i\in F,
|v_i|<R\}\times \Re\\
\cup \{(\tau,\zeta,v_1,w_1,v_2,w_2)| |v_1|,|v_2|\leq R; \max( |v_1|,|v_2|)=R.\}\times \Re.
\end{multline*}

Consider now $\SS(R_J*_{E\times E} A)$.  As $\SS(A)\subset \{(v_1,v_2)||v_1|,|v_2|<R\}$, it follows that
$$
\SS(R_J*_{E\times E} A)\subset \{(\tau,\eta_J(\tau,v_1)-\eta_J(\tau,v_2),J(\tau,v_1),J(\tau,v_2))||v_1|,|v_2|<R\}\times \Re.
$$
Let us estimate
$$
\SS((R_J*_{E\times E} A)*_{E\times E} B).
$$

It follows that  there exists a compact subset $K\subset U$ such that 

$$
\SS((R_J*_{E\times E} A)*_{E\times E} B)\subset \{(\tau,\eta_J(\tau,v_1)-\eta_J(\tau,v_2))|\tau\in K\times (-1,1),|v_1|,|v_2|<R\}\times \Re.
$$

Namely,  one can choose $K=\overline{\{u\in U| I(u,B_R)\cap V\neq \emptyset\}}$.

Let now $\tau=(u,x)\in U\times (-1,1)$.   We have $ \eta_J(\tau,v)\in F\oplus \Re$.
Let $f(\tau,v)$ be the  $F$-component and $x(\tau,v)$ be the $\Re$-component.
Let us now estimate 
$$
\SS(\pi_{U!} (R_J*_{E^4} (A\boxtimes B))).
$$

As $\pi_U$ is proper on the support of $R_J*_{E^4}(A\boxtimes B)$,  the singular support in question
is determined by the condition $f(\tau,v_1)-f(\tau,v_2)=0$.  As the family $I$ is mobile, this condition
implies $v_1=v_2$, which implies $\eta_J(\tau,v_1)-\eta_J(\tau,v_2)=0$ and
$$
\SS(\pi_{U!} (R_J*_{E^4} (A\boxtimes B)))\subset T^*_{(-1,1)} (-1,1)\times \Re.
$$

\section{Tree operads and multi-categories }

\subsection{Planar/cyclic trees} Let us introduce a notation for a  tree $\bt$.   Denote by $\inp(\bt)$  the set of inputs
of $\bt$, $V_\bt$ the set of inner vertices of $\bt$,   for $v\in V_\bt$, denote by $E_v$ the set of inputs of $v$.
Let $p_\bt$ be the principal vertex of $\bt$. 
\subsubsection{Planar trees}
Define {\em a planar tree} as a tree with a  total order on every set $E_v$; we then have an induced total order on $\inp(\bt)$.  

We have a unique identification of ordered sets $E_v=\{1,2,\ldots,n_v\}$, where $n_v=\# E_v$;
$\inp_\bt=\{1,2,\ldots,n_\bt\},$ where $n_\bt=\# \inp_\bt$.

\subsubsection{Cyclic trees}
Define {\em a cyclic tree} as a tree with a total order on every set $E_v$, $v\neq p_\bt$, and a cylic order on $p_\bt$.
We then have an induced cycic order on $\inp_\bt$, in particular, we assume $\inp_\pt\neq \emptyset$.

{\em A rigid cyclic tree} is a cyclic tree along with identifications $E_{p_\bt}=\{1,2,\ldots,n_{p_{\bt}}\}$;
$\inp_{p_\bt}=\{1,2,\ldots,n_\bt\}$ which agree with the cyclic order on both sets.

\subsubsection{Inserting trees into a tree} Let $\bt$ be a planar tree.  Let $\bt_v$, $v\in V_\bt$ be planar trees
where $n_{\bt_v}=n_v$.  One then can insert the trees $\bt_v$ into $\bt$.  Denote the resulting tree
by $\bt\{\bt_v\}_{v\in V_\bt}$.   

Similarly, let $\bt$ be a rigid cyclic tree.  Let $\bt_v$, $v\in V_\bt\backslash p_\bt$ be planar trees
with $n_{\bt_v}=n_v$; let $\bt_{p_\bt}$ be a rigid cyclic tree with $n_{\bt_{p_\bt}}=n_{p_\bt}$.
One then can define a similar insertion, to be denoted by $\bt\{\bt_v\}_{v\in V_\bt}$.
\subsubsection{Isomorphism classes of trees}
Let $\trees$ be the set of isomorphism classes of planar trees
and $\trees^\cyc$ be the set of isomorphism classes of rigid cyclic trees. 

Let also $\trees_n\subset \trees$ be the subset consisting of all isomorphism classes of trees
with $n_\bt=n$ and likewies for $\cyctrees_n$.   The above defined insertions are defined on the level of isomorphism classes.

\subsubsection{Famililes parameterized by isomorphism classes of trees} Let $\ccA$ be a $\bigoplus$-closed 
dg SMC.   Let $\cT(\ccA)$ be a category, enriched over sets, whose every object is 
a family of objects $X_\bt\in \ccA$, $\bt\in \trees\sqcup \cyctrees$.
Let $X,Y\in \cT(\ccA)$.  Let us define a new family $X\circ Y\in \cT(\ccA)$ as follows:
\begin{equation}\label{kruzhochek}
X\circ Y(\bT)=\bigoplus\limits_{\bT=\bt\{\bt_v\}_{v\in V_\bt}} X(\bt)\otimes \bigotimes\limits_{v\in V_\bt} Y(\bt_v).
\end{equation}
This way, $\cT(\ccA)$ becomes a monoidal category.  The unit object  $\unit\in \cT(\ccA)$ is defined
by setting $\unit(\bt)=\unit_\ccA$ for all isomorphism classes of planar trees with  one vertex (corollas)
and all isomorphism classes $\bt$ of rigid cyclic trees with one vertex and matching numberings of $E_p$ and $\inp_\bt$. Otherwise, $\unit(\bt)=0$.

\subsubsection{Planar trees with marked right branch}   Let $\treesm\subset  \trees$ be a subset
consisting of all planar trees $\bt$ with $\inp \bt\neq \emptyset$.  Let $r_\bt\in \inp \bt$ be the rightmost input.
   For $\bt\in \treesm$, let $V^R_\bt\subset V_\bt$
consist of all vertices for which $r_\bt$ is an output.

Let $\bt\in \treesm$; let $\bt_v\in \trees$, $v\in V_\bt\backslash V^R_\bt$ and  $\bt_w\in \treesm$, $w\in V^R_\bt$.
Suppose $\# \inp \bt_v =\# E_v$ for all $v\in V_\bt$.   We then have  a well defined insertion
$\bt\{\bt_v\}_{v\in V_\bt}$.   

Let $\cTM(\ccA)$  be a category, enriched over sets, whose every object is 
a family of objects $X_\bt\in \ccA$, $\bt\in \trees\sqcup \cyctrees\sqcup \treesm$.
Let $X,Y\in \cT(\ccA)$.  We define a new family $X\circ Y\in \cT(\ccA)$ by the same formula  (\ref{kruzhochek}).

This way, $\cTM(\ccA)$ becomes a monoidal category.  The unit object  $\unit\in \cT(\ccA)$ is defined
by setting $\unit(\bt)=\unit_\ccA$ for all isomorphism classes of planar trees or planar trees with marked
right branch with  one vertex (corollas)
and all isomorphism classes $\bt$ of rigid cyclic trees with one vertex and matching numberings of $E_p$ and $\inp_\bt$. Otherwise, $\unit(\bt)=0$.
\renewcommand{\swell}{{\Doplus}}
\subsection{Collections of functors}  Let $\cC,\cD$ be categories  enriched over $\ccA$ and  tensored by $\ccA$.
Suppose we are also given a functor $h:\cC\otimes \cD\to \ccA$.

 Let us define a category over $\sets$, $\cF(\cC,\cD)$, as follows $$
\cF(\cC,\cD):=\prod\limits_{n=0}^\infty\swell(\cC^n\otimes \cD) \times \prod\limits_{n=1}^\infty \cC^{n}
$$
so that an object $F\in \cF(\cC,\cD)$ is a collection of objects 
   $F^{[n]}\in\swell(\cC^{\otimes n}\otimes \cD)$, $n\geq 0$, and 
$F^{(n)}\in \swell(\cC^{\otimes n})$, $n\geq 1$.

Let $\bt$ be a planar tree.   Define an object$$F(\bt)\in \swell(\cC^{\otimes n_\bt}\otimes \cD).$$

A) We have an equivalence of categories
$$
 \bigotimes\limits_{v\in V_\bt}( \cC^{\otimes n_v}\otimes \cD)\cong \left(\bigotimes\limits_{v\in V_\bt\backslash p_\bt}
\cC\otimes \cD\right)\otimes (\cC^{\otimes n_\bt}\otimes \cD),
$$
coming  from the bijection
$$
\bigsqcup_{v\in V_\bt} E_v\cong  V_\bt\sqcup \inp_\bt\backslash p_\bt
$$
which associates to an edge its target.  

As a result we have a through map (via applying the functor $h$):
$$
\circ_\bt:
 \bigotimes\limits_{v\in V_\bt}\swell( \cC^{\otimes n_v}\otimes \cD)\to \swell\left(\left(\bigotimes\limits_{v\in V_\bt\backslash p_\bt}
\cC\otimes \cD\right)\otimes (\cC^{\otimes n_\bt}\otimes \cD)\right)\to \swell(\cC^{\otimes n_\bt}\otimes \cD).
$$

C) Set $F(\bt):=\circ_\bt\left(\bigotimes\limits_{v\in V_\bt} F^{[n_v]}\right)$.

Let now $\bt$ be a rigid cyclic tree.  Define a functor $F(\bt)\in \swell(\cC^{n_\bt})$ in a similar way.
Let
$$\circ_\bt: \cC^{\otimes n_{p_\bt}}\otimes \bigotimes\limits_{v\in V_\bt\backslash p_\bt} (\cC^{\otimes n_v}\otimes \cD)\to \cC^{\otimes n_\bt}$$
be defined similar to above and set 
$$
F(\bt):=\circ_\bt\left(F^{(n_{p_\bt})}\otimes \bigotimes\limits_{v\in V_\bt\backslash p_\bt} F^{[n_v]}\right).
$$

\subsubsection{Extended collection  of functors}  Let $\cC,\cD,\cCR,\cDR$ be  dg categories tensored by $\ccA$.
Suppose we are also given functors $h:\cC\otimes \cD\to \ccA$, $h_R:\cCR\otimes \cDR\to \bbB$.

 Let us define a category over $\sets$, $\cF(\cC,\cD.\cCR,\cDR)$, as follows \begin{multline*}
\cF(\cC,\cD,\cCR,\cDR):=\prod\limits_{n=0}^\infty\swell(\cC^n\otimes \cD) \times \prod\limits_{n=1}^\infty \cC^{n}
\times \prod\limits_{n=0}^\infty\swell(\cC^n\otimes \cCR\otimes \cDR).
\end{multline*}
so that an object $F\in \cF(\cC,\cD)$ is a collection of objects 
   $F^{[n]}\in\swell(\cC^{\otimes n}\otimes \cD)$, $n\geq 0$, 
$F^{(n)}\in \swell(\cC^{\otimes n})$, $n\geq 1$, and $F^{[n,1]}\in \swell(\cC^{\otimes n}\otimes \cCR\otimes \cDR)$.
\subsubsection{Definition of $F(\bt)$}  Let $F$ be a collection of functors.
Let $\bt$ be a planar tree.   Define an object
$$
F(\bt)\in \swell(\cC^{\otimes n_\bt}\otimes \cD).
$$

A) We have an equivalence of categories
$$
 \bigotimes\limits_{v\in V_\bt}( \cC^{\otimes n_v}\otimes \cD)\cong \left(\bigotimes\limits_{v\in V_\bt\backslash p_\bt}
\cC\otimes \cD\right)\otimes (\cC^{\otimes n_\bt}\otimes \cD),
$$
coming  from the bijection
$$
\bigsqcup_{v\in V_\bt} E_v\cong  V_\bt\sqcup \inp_\bt\backslash p_\bt
$$
which associates to an edge its target.  

As a result we have a through map (via applying the functor $h$):
$$
\circ_\bt:
 \bigotimes\limits_{v\in V_\bt}\swell( \cC^{\otimes n_v}\otimes \cD)\to \swell\left(\left(\bigotimes\limits_{v\in V_\bt\backslash p_\bt}
\cC\otimes \cD\right)\otimes (\cC^{\otimes n_\bt}\otimes \cD)\right)\to \swell(\cC^{\otimes n_\bt}\otimes \cD).
$$

C) Set $F(\bt):=\circ_\bt\left(\bigotimes\limits_{v\in V_\bt} F^{[n_v]}\right)$.

Let now $\bt$ be a rigid cyclic tree.  Define a functor $F(\bt)\in \swell(\cC^{n_\bt})$ in a similar way.
Let
$$\circ_\bt: \cC^{\otimes n_{p_\bt}}\otimes \bigotimes\limits_{v\in V_\bt\backslash p_\bt} (\cC^{\otimes n_v}\otimes \cD)\to \cC^{\otimes n_\bt}$$
be defined similar to above and set 
$$
F(\bt):=\circ_\bt\left(F^{(n_{p_\bt})}\otimes \bigotimes\limits_{v\in V_\bt\backslash p_\bt} F^{[n_v]}\right).
$$

Let  now $F$ be  an extended collection of functors.  We then define $\bt(F)$ in a similar way.

\subsection{Schur functors} 
 
Suppose $\cC,\cD$ are  enriched and tensored  tensored over $\ccA$.  Let $X\in \cT(\ccA)$  and $F\in \cF(\cC,\cD)$.
Define an object
$
\schur_X(F)\in \cF(\cC,\cD)
$
as follows
$$
\schur_X(F)^{[n]}:=\bigoplus_{\bt\in \trees_n} \bt(F);\quad \schur_X(F)^{(n)}=\bigoplus_{\bt\in \cyctrees_n} \bt(F).
$$
We have  natural isomorphisms
$$
\schur_X\schur_Y F\cong \schur_{X\circ Y} F;\quad \schur_\unit F\cong F.
$$
In fact, we have a $\cT(\ccA)$-action on $\cF(\cC,\cD)$.

One defines a $\cTM(\ccA)$-action on $\cF(\cC,\cCR,\cD,\cDR)$ in a similar way.

\subsection{Tree operads} {\em A tree operad in $\cT(\ccA)$} is the same as  a unital monoid  in $\cT(\ccA)$.
Respectively, an extended tree operad in $\cTM(\ccA)$ is the same as a unital monoid in $\cTM(\ccA)$.
\subsubsection{A tree operad $\triv$}    Let $\triv\in \cT(\ccA)$ be given by $\triv(\bt)=\unit_\ccA$ for all $\bt$.
Define $\triv\in \cTM(\ccA)$ in a similar way.
\subsubsection{Endomorphism tree operad}  Let $\cC,\cD$ be enriched and tensored over $\ccA$. 
Let $F,G\in \cF(\cC,\cD)$.
Consider a functor $H_{F,G}:\cT(\swell \ccA)\to \sets$,
$$
H_{F,G}(X)=\hom(\schur_X F;G)
$$
The functor $H_{F,G}$ is representable.  Denote the representing object by $\cH_{F,G}$. We have ($\bt$ is planar):
$$
\cH_{F,G}(\bt)=\ihom_{\swell (\cC^{\otimes n_\bt}\otimes \cC^\op)}(F(\bt);G^{[n_\bt]});
$$
if $\bt$ is a rigid cyclic tree, we have:
$$
\cH_{F,G}(\bt)=\ihom_{\swell \cC^{\otimes n_\bt}}(F(\bt);G^{(n_\bt)}).
$$

Set $\End_F:=\cH_{F,F}$. We have a natural tree operad structure on $\End_F$.    Furthermore, we
have an $\End_F-\End_G$-bi-module structure on $\cH_{F,G}$ (where we interpret tree operads $\End_F,\End_G$
as monoids in  $\cT(\swell \ccA)$.

Let $F\in \cF(\cC,\cCR,\cD,\cDR)$.  We define an extended tree operad  $\End_F$ in a similar way.
\subsection{Pull backs from $\cF(\cC',\cD')$ to $\cF(\cC,\cD)$} Let $\ccA$ have internal hom.
 Let $\cC,\cC',\cD,\cD'$ be categories  enriched over $\ccA$; let $h:\cC\otimes \cD\to \ccA$;
$h':\cC'\otimes \cD'\to \ccA$ be functors.

Let $G\in \cF(\cC',\cD')$.  Let $L\in \swell(\cD\otimes \cC')$.

Consider the following functor $H:\cF(\cC,\cD)^\op\to \sets$ as follows.   

1)  We have  functors $$
e_L:\cC^{\otimes n}\otimes \cD\otimes (\cD\otimes \cC')^{\otimes n}\to 
(\cC')^{\otimes n}\otimes \cD,$$
via using the hom-functor $h^{\boxtimes n}:\cC^{\otimes n}\otimes (\cD)^{\otimes n}\to \ccA$, as well as
$$
f_L:(\cC')^{\otimes n}\otimes \cD'\otimes \cD\otimes \cC'\to (\cC')^{\otimes n}\otimes \cD.
$$
via the hom functor $h':\cC'\otimes\cD'\to \ccA$.

Similarly, one defines a  cyclic version:
set $$
e^\cyc_L: \cC^{\otimes n}\otimes (\cD\otimes \cC')^{\otimes n}\to (\cC')^{\otimes n}.
$$

2) set 
$$
H^{[n]}(F,G):=\hom(e_L(F^{[n]}\otimes \cL^{\otimes n});f_L(G^{[n]}\otimes \cL));
$$
$$
H^{(n)}(F,G):=\hom(e_L^\cyc(F^{(n)}\otimes  \cL^{\otimes n});  G^{(n)}).
$$

Set $$
H(F,G)=\prod\limits _{n\geq 0} H^{[n]}(F,G)\times \prod_{n>0} H^{(n)}(F,G).
$$
It follows that the functor $F\mapsto H(F,G)$ is representable.  Denote the representing object by $L^{-1}G$.
Let $X\in \cT(\ccA)$.  We have a natural map $\schur_X L^{-1}G\to L^{-1}\schur_X G.$

Therefore,   $\cH_{F,L^{-1}G}$ is a $\End_F-\End_G$-bimodule.

One  generalizes this construction to the extended case straightforwardly.

\subsubsection{Quasi-contracible tree operads}    Let now $\ccA=\pt$ so that $\swell \ccA=\GZ$.
Call a tree operad $\cO\in \cT(\GZ)$ {\em  pseudo-contractible} if

1)  $\cO(\bt)\in \GZ$ admits a truncation for every  $\cO(\bt)$.  We therefore have an induced 
tree operad structure on $\tau_{\leq 0} \cO$ and a map of  tree operads $\tau_{\leq 0}\cO\to \cO$.

2) Every object $\tau_{\leq 0} \cO$ admits a trunctation $\tau_{\geq 0}$, to be denoted $H^0\cO(\bt)$ which
is a finitely generated free $\ZQ$-module; we have an induced map of tree operads
$\tau_{\leq  0}\cO\to H^0 \cO$.  We require this map to be a term-wise homotopy equivalence.

{\em A quasi-contractible tree operad} is a pseudo-contractible operad $\cO$ endowed with a map of tree
operads $\triv\to H^0(\cO)$. 

In this case there exists a splitting of the map $\tau_{\leq 0} \cO(\bt)\to H^0\cO(\bt)$, hence a 
  pull-back of the diagram $$
\triv\to H^0(\cO)\leftarrow \tau_{\leq 0} \cO,
$$
to be denoted  by  $\triv_\cO$ so that we have a diagram
$$
\triv\stackrel\sim\leftarrow \triv_\cO\to \cO.
$$
Let $\cO_1,\cO_2$ be quasi-contractible operads  and $\cM$ a $\cO_1-\cO_2$-bi-module.
Call $\cM$ {\em pseudo-contractible} if     there exist truncations
$\tau_{\leq 0} \cM(\bt)$ and $\tau_{\geq 0}\tau_{\leq 0}\cM(\bt)=:H^0\cM(\bt)$, where each $H^0\cM(\bt)$ is a finitely generated
free $\ZQ$-module.   

{\em A quasi-contractible $\cO_1-\cO_2$-bi-module $\cM$} is a pseudo-contractible $\cO_1-\cO_2$-bi-module
$\cM$ endowed with a map 
$$
(\triv,\triv,\triv)\to (H^0\cO_1,H^0\cM,H^0\cO_2)
$$
of triples:  a pair of tree-operads and their bi-module.

Similar to above,  we have a pull-back of the diagram
$$
(\triv,\triv,\triv)\leftarrow (\tau_{\leq 0} \cO_1,\tau_{\leq 0}\cM,\tau_{\leq 0}\cO_2)\to (H^0\cO_1,H^0\cM,H^0\cO_2),
$$
to be denoted by $(\triv_{\cO_1},\triv_{\cM};\triv_{\cO_2})$ so that we have a diagram
$$
(\triv,\triv,\triv)\stackrel\sim\leftarrow (\triv_{\cO_1},\triv_{\cM},\triv_{\cO_2})\to (\cO_1,\cM,\cO_2).
$$

\section{Straightening out} \label{straight} Fix a ground category $\ccA$.
\subsection{Pseudo-contractible sequences}\label{pcs}
Fix categories $\cC_i,\cD_i$, $i=1,n$;  functors $h_i:\cC_i\otimes \cD_i\to \ccA$,  objects $F_i\in \cF(\cC_i,\cD_i)$ and
$L_{i+1,i}:\swell(\cD_{i+1}\otimes \cC_i)$.  Call such a collection of data {\em a sequence}.
Call such a sequence {\em  pseudo-contractible} if
 every triple $$
(\End_{F_{i}},\cH_{F_{i+1},L_{i+1,i}^{-1}F_i},\End_{F_{i+1}}), 
$$
$0\leq i\leq n-1,$ is pseudo-contractible.

 Denote by  $(\cO_{i},\cO_{i,i+1},\cO_{i+1})$ the pull-back of the diagram
$$
\xymatrix
{\tau_{\leq 0}(\End_{F_{i+1}},\cH_{F_{i+1}F_i},\End_{F_i})\ar[r] &H^0(\End_{F_{i+1}},\cH_{F_{i+1}F_i},\End_{F_i})\\
                                                                                      &(\triv,\triv,\triv)\ar[u]}
$$

Denote by $$\pi_{i+1,i}:(\cO_{i},\cO_{i,i+1},\cO_{i+1})\stackrel\sim\to (\triv,\triv,\triv)$$
the projection.

\subsection{Straightening}  Let us define new sequences of functors $G_i\in \cF(\cC_i,\cD_i)$ 
such that we have  a $\triv$ -action on each $G_i$ as well as maps  of $\triv$-modules $G_{i+1}\to L_{i+1,i}^{-1}G_i$.
Namely,  fix canonical quasi-free resolutions  $\cR_i\to \cO_i$ of $\cO_i$ as of a bimodule over itself.

Set
$$
G_0:=\triv\circ_{\cO_0}\cR_0 \circ_{\cO_0} F_0;
$$
$$
G_1:=\triv\circ_{\cO_0}\cR_0\circ_{\cO_0} \cO_{01}\circ_{\cO_1} \cR_i \circ_{\cO_1} F_1;
$$
$$
G_2:=\triv\circ_{\cO_0}\cR_0\circ_{\cO_0} \cO_{01}\circ_{\cO_1} \cR_i \circ_{\cO_1} \cO_{12}\circ_{\cO_2}\cR_2\circ_{\cO_2} F_2;
$$
etc.
 
The maps of $\triv$-modules $G_{i+1}\to L_{i+1,i}^{-1}G_i$  are induced by the maps
$$
\cO_{i,i+1}\circ_{\cO_{i+1}} \cR_{i+1}\circ_{\cO_{i+1}} F_{i+1}\to  \cO_{i,i+1}\circ_{\cO_{i+1}} F_{i+1}\to F_i.
$$

\subsubsection{Extended case} Let $F$ be an extended collection of functors acted upon by a quasi-contractible
extended tree operad.  One then constructs an extended collection $F'$ with a $triv$-action in a similar way.
\subsection{Definition of a monoidal category with trace} 
\subsubsection{The category $\tau(\cM)$}  Let $\cM$ be a monoidal category over $\ccA$.  Define a category
$\tau(\cM)$ enriched over $\ccA$ as follows.  An object of $\tau(\cM)$ is an object 
of $\cM^{\otimes n}$ for some $n\geq 1$.    

Let $\cX:=(X_1, X_2, \ldots ,X_n),\  \cY:=(Y_1, Y_2,\ldots, Y_m)\in \tau(\cM)$, where $X_i,Y_j\in \cM$.
Let $f:\{1,2,\ldots,n\}\to \{1,2,\ldots,m\}$ be a cyclically non-decreasing map.
For $p\in  \{1,2,\ldots,m\}$  we then have a total order on $f^{-1}p$.
 Let
$$
X_p:=\bigotimes\limits_{i\in f^{-1}p} X_i\in \cM.
$$
Let 
$$
\hom_f(\cX,\cY):=\bigotimes\limits_{1\leq p\leq m} \hom_\cM(X_p;Y_p).
$$
Let $$\hom_{\tau(M)}(\cX,\cY):=\bigoplus\limits_f \hom_f(\cX,\cY),
$$
where $f$ runs through the set of all cyclically non-decreasing maps $\{1,2,\ldots,n\}\to \{1,2,\ldots,m\}$.

\subsubsection{Definition of a  monoidal category with trace} {\em A monoidal category with trace} 
$(\cM,\cM^\cyc)$ is a  pair consisting of a monoidal category $\cM$ and a category $\cM^\cyc$ endowed with  a functor $\TR:\tau(\cM)\to \cM^\cyc$.

\subsection{A monoidal category  with trace  $\bU(F)$ from a collection  $F\in \cF(\cC,\cD)$ with a $\triv$-action}

 As above, let $F\in \cF(\cC,\cD)$.  Suppose $F$ carries a $\triv$-action.
Let us construct a monoidal category  with trace $\bU(F)$. 

--- An object of $\bU(F)$ is an object of $\swell(\cC^{\boxtimes n})$ for some $n\geq  0$.

--- Let $X\in \swell(\cC^{\boxtimes n})$ and $Y\in \swell(\cC^{\boxtimes m})$, $m>0$.
Let $f:\{1,2,\ldots,n\}\to \{1,2,\ldots,m\}$ be a non-decreasing  map.    Let $\phi_k:=\# f^{-1}m$. 
Let $$
F(f):=\boxtimes_{k=1}^m F(\phi_k)\in \swell(\cC^{\boxtimes n}\boxtimes \cD^{\boxtimes m}).
$$  

Set $$
\hom_f(X,Y):=\hom(X;h^{\boxtimes m}(F(f)\boxtimes Y)).
$$

Set $$\hom(X,Y):=\bigoplus_{f} \hom_f(X,Y).$$

If $n>0$ and $m=0$, then we set $\hom(X,Y)=0$.  Finally, if $n=m=0$, we set
$$\hom(X_n,Y_m)=\unit_\ccA.
$$

The $\triv$-action on $F$ induces maps
$$
\hom_f(X,Y)\otimes \hom_g(X,Y)\to \hom_{gf}(X,Y),
$$
whence a $\ccA$-category structure on $\bU(F)$.

Define the monoidal structure on $\bU(F)$  by means of $\boxtimes$ and natural maps
$$
\hom_f(X_1,Y_1)\otimes \hom_g(X_2,Y_2)\to \hom_{f\times g}(X_1\boxtimes X_2;Y_1\boxtimes Y_2).
$$
The unit is a unique object in $\cC^{\boxtimes 0}$.

Define a trace on $\bU(F)$ via prescribing restrictions
$$
\TR_n: \cC^{\boxtimes n}\to \ccA,\quad \TR_n(X)=\hom(X;F^{(n)}).
$$

Suppose we have  a map of $\triv$-modules $\pi_{i+1,i}:G_{i+1}\to L_{i+1,i}^{-1}G_i$, as in the previous subsection.
We then have an induced tensor functor $\pi_{i+1,i}:\bU(G_{i+1})\to \bU_{G_i})$, where for $X\in \cC_{i+1}^{\boxtimes n}$
we set
$$\pi_{i+1,i}(X):=h^{\boxtimes n}.(\cL_{i+1,i}^{\boxtimes n}\boxtimes X) \in \cC_i^{\boxtimes n}.
$$

\subsubsection{The category $\bU_1(F)$} Denote by $\bU_1(F)\subset\bU(F)$ the full sub-category
of objects homotopy equivalent to an object from $\cC\subset \bU(F)$.

\subsubsection{The category $\bU^\cyc(\cO)$}   Let $\bU^\cyc(\cO)\subset \tau(\bU(F))$ be the full
sub-category consisting of all objects of the form $(X_1,X_2,\ldots,X_n)$, where $n>0$ and $X_i\in \cC$.

We have a functor $I:\bU(F)\to \bU^\cyc(F)$ which is identical on objects.

We have a  functor $$
\TR:\bU^\cyc(F)^\op\to \ccA.
$$

Let $\bU^\cyc_{\geq 0}(F)$ be a category which is obtained from $\bU^\cyc(F)$ by adding an extra object,
to be denoted by  $(0)$,  where $\hom((0),X)=0$,  $\hom(X,(0))=\TR(X)$, for all $X\in \bU^\cyc(F)$,
$\hom((0),(0))=\unit_\ccA$.

\subsubsection{Extended case: a monoidal category with a trace and its module}
Let now $F$ be an extended collection of functors with $\triv$-action.  We then get a monoidal category with traces $\bU(F)$ and
its module $\bU_M(F)$, that is, $\bU_M(F)$ is a category  endowed with a functor
$\bU(F)\otimes \bU_M(F)\to \bU_M(F)$ satisfying  the associativity axiom.

\subsection{Circular operads}\label{circularoperad}  {\em A circular operad} is a $\triv$-module on an object $\cO\in \cF(p,p)$,
where $p$ is a   disjoint  union  of  a  finite number of copies of  $\pt$, the category with one object , to be denoted by $O$ and $\hom(O,)=\unit_\ccA$.
Set  $\cO(n):=\cO^{[n]}(\pt,\ldots,\pt)$;  $\cO^\cyc(n):=\cO^{(n)}(\pt,\ldots,\pt)$.

A circular operad structure is equivalent to an asymmetric operad structure on a collection $\cO$ and
an $\bU^\cyc(\cO)$-module structure on $\cO^\cyc$.
In particular, we have a  $\bZ/n\bZ$-action on $\cO^\cyc(n)$ for every $n\geq 1$.
Call $\cO$ {\em cyclically semi-free} if such is each $\bZ/n\bZ$-module $\cO^\cyc(n)$. 

Let $\Lambda$ be the cyclic category.  Denote its objects by $(n)$, $n>0$ (observe the shift by 1 unit with respect
to the traditional numbering!).
     Let $\bbZ^\cyc(\cO)$ be the category with the same set of 
objects as in $\bU^\cyc(\cO)$ and we set 
$$\bbZ^\cyc(\cO)((m),(n))=\bU^\cyc((m),(n))\otimes \Lambda^\op((m),(n)). 
$$ 

We have a functor $\bU^\cyc(\cO)\to \bbZ^\cyc(\cO)$.
Let $U:\bU^\cyc(\cO)^\op\to \ccA$,  $L:\Lambda\to \ccA$ be functors.   Set $U\diamond L((n)):= U(n)\otimes L(n)$.
We have 
$$
U\diamond L: \bU^\cyc\to \bbZ^\cyc\to \ccA.
$$
Hence,  we have an induced circular operad structure on $(\cO,\cO^\cyc\diamond L)$.

Fix a free   resolution $\cR_\Lambda\to \underline{\ZQ}$   of the constant $\Lambda$-module $\underline{\ZQ}$.
It follows that $\cO_\cR:=(\cO,\cO^\cyc\diamond \cR_\Lambda)$ is a cyclically semi-free circular operad.  We have a natural
map of circular operads $\cO_\cR\to \cO$.

\subsection{Algebras over circular operads}   Let $\cO$ be a circular operad 
Let $(\cM,\cM^\cyc)$ be a monoidal category with a  trace.  {\em A $\cO$-algebra in $\cM$} is a strict tensor functor
of monoidal categories with a trace:
$$
(\bU(\cO),\bU^\cyc_{\geq 0}(\cO))\to (\cM,\cM^\cyc).
$$

Suppose $\cM$ is $\bigoplus $-closed and $\cO$ is  cyclically semi-free. Then the structure  of a $\cO$-algebra is equivalent to 
that of an algebra over a certain monad, to be denoted by  $\schur_\cO$, acting on
the category $K:=\cM\times \cM^\cyc$  over $\sets$.  
 By definition,
$$
\schur_\cO(X,U):=(X',U'),
$$
where 
$$
X':=\bigoplus\limits_{n\geq 0} \cO(n)\otimes X^{\otimes n};
$$
$$
U':=U\oplus\bigoplus\limits_{n\geq 1} \cO^\cyc(n)\otimes_{\bZ/n\bZ} \TR(\underbrace{X,\ldots,X}_{n \text{ times}}).
$$

There exists a monoidal structure on $\cF(\pt,\pt)$, to be denoted by $\circ$,  so that we have an isomorphism
$$
\schur_X\schur_Y\cong \schur_{X\circ Y},
$$
for all cyclically free $X$, $Y$.   A structure of a circular operad on $\cO$ is equivalent to that of a unital monoid
in $\cF(\pt,\pt)$ whence a monad structure on $\schur_\cO$. 

Let $\cO_2\to \cO_1$ be a map of circular operads.    We then have a $\cO_1-\cO_2$-bimodule structure on
$\cO_2\in \cF(\pt,\pt)$.    Choose its  semi-free resolution  $R$.  

Let $(X,Y)$ be a $\cO_2$-algebra in $\cM$.
One then has a $\cO_1$-algebra structure
on
$$
R\circ_{\cO_2} (X,Y),
$$
which is well-defined as long as $(\cM,\cM^\cyc)$ is $\bigoplus D$-closed.

\centerline{\bf PART 3. MICROLOCAL CATEGORY: CLASSICAL LEVEL}

\section{Geometric setting} 
\subsection{Principal bundles}
---Let $M$ be a compact symplectic  $2N$-dimensional manifold whose symplectic form $\omega$ has integral periods.

---Let $L$ be a circle bundle whose first Chern class equals $\omega$.  

---Let $P_0\to M$ be the principal bundle of symplectic frames on $M$ with its structure group $\Sp(2N)$. Let  $P:=P_0\times_M L$.  $P$
is a $\Sp(2N)\times S^1$-principal bundle  over $M$.

--- Let $\frakH:=P_0/U(N)$.   We have a smooth fibration $\frakH\to M$ with contractible fiber.
We have a princilal $U(N)\times S^1$-bundle $P\to \frakH$.

\subsection{Pseudo-Kaehler metrics}

---Let $\Met$ be the set of all pseudo-Kaehler metrics on $M$ with their symplectic form being $\omega$. 
One identifies $\Met$ with the set of sections of the bundle $\frakH\to M$.  For $g\in \Met$, denote
by $i_g:M\to \frakH$ the corresponding section and $\Fr_g:=i_g^{-1}P$.

Fix retractions $\pi_g: P\to \Fr_g$ so that we have the following retraction of principal
$U(N)\times S^1$-bundles
\begin{equation}\label{retract}
\xymatrix{ 
\Fr_g\ar[d]\ar[r]& P\ar[d]\ar[r]^{\pi_g}& \Fr_g\ar[d]\\
 M\ar[r]&\frakH\ar[r]& M}
\end{equation}

The retraction $\pi_g$ is constructed as follows.  Fix a $U(N)\times U(N)$-equivariant retraction
 \begin{equation}\label{sp2u}
U(N)\to \Sp(2N)\stackrel r\to U(N).
\end{equation}

We then set
$$
\pi_g:P=\Fr_g\times_{U(N)} \Sp(2N)\stackrel{r}\to \Fr_g\times_{U_N} U(N)=\Fr_g.
$$

Denote by $i_\Met: \Met\times M\to \frakH$ the union of all $i_g$, $g\in \Met$.
 Let $\Fr:=\bigsqcup_{g\in \Met} \Fr_g$.   We have $\Fr=i_\Met^{-1}P$. One has the following 
retraction of principal $U(N)$-bundles
\begin{equation}\label{retracmet}
\xymatrix{ 
\Fr\ar[d]\ar[r]& P\times \Met\ar[d]\ar[r]& \Fr\ar[d]\\
 M\times \Met\ar[r]&\frakH\times \Met\ar[r]& M\times \Met}
\end{equation}

\subsection{Groupoids}
\subsubsection{The groupoid $P\times_M P$}
Let us start with a groupoid
$$
P\times_M P\rightrightarrows P.
$$
\subsubsection{The groupoid $Q$, a covering of $P\times_M P$}\label{kaen}
We have a natural map $P\times_M P\to \Sp(2N)\times S^1$.  Let  $\cH$ be the universal
cover of $\Sp(2N)\times S^1$.   Set
$$
Q:=(P\times_M P)\times_{\Sp(2N)\times S^1} \cH.
$$

We have a groupoid structure on $Q\rightrightarrows P$  and a map of gropoids
$$
(Q\rightrightarrows P)\to (P\times_M P\rightrightarrows P).
$$
The $\bZ\times \bZ$-action on $\cH$ carries over to $Q$. Let us denote this action by $T$.
Let $\vs:Q\to Q$ be the  inversion map.

One also has iterations
$$
C_n:Q\times_P Q\times_P \cdots \times_P Q\to Q,
$$
where there are $n$ occurences of $Q$ on the LHS.

Let $I_n: Q\times_P Q\times_P \cdots \times_P Q\to Q^n$ be the embedding.

Let  $a\in \bZ\times \bZ$.  Let $K_n^a$ be the image of the map
$$
I_n\times a.C_n:Q\times_P Q\times_P \cdots \times_P Q\to Q^{n+1}.
$$

Let also $K_n^{\cyc,a}$ be obtained from $K_n^a$ by applying the permutation $\vs$ to the last factor of $Q$.

\subsubsection{The groupoid $\Phi$}\label{groupphi} The deformation retraction (\ref{retracmet}) induces a deformation
retraction of groupoids 
\begin{equation}\label{deform2}\xymatrix{
\Fr\times_M \Fr\ar@<1ex>[d]\ar@<-1ex>[d]\ar[r]^(0.3){j}& P\times_M P\times \Met\times \Met\ar@<1ex>[d]\ar@<-1ex>[d]\ar[r]& \Fr\times_M\times \Fr\ar@<1ex>[d]\ar@<-1ex>[d]\\
\Fr\ar[r]& P\times \Met  \ar[r] &\Fr}
\end{equation}

Let now $\Phi$ be the $j$-pull-back of the covering $$Q\times \Met\times \Met\to P\times_M P\times \Met\times \Met.
$$

We have an induced groupoid structure on $\Phi\rightrightarrows \Fr$;  the map $j$ then naturally extends 
to a map of groupoids
$$
(\Phi\rightrightarrows \Fr)\to (Q\times \Met\rightrightarrows P\times \Met).
$$

As $j$ is a deformation retraction, we have a 
 deformation retraction of groupoids:
\begin{equation}\label{deform2}\xymatrix{
\Phi\ar@<1ex>[d]\ar@<-1ex>[d]\ar[r]^(0.3){j}& Q\times \Met\times \Met\ar@<1ex>[d]\ar@<-1ex>[d]\ar[r]^(0.7){p_\Phi}& \Phi\ar@<1ex>[d]\ar@<-1ex>[d]\\
\Fr\ar[r]& P\times \Met  \ar[r]& \Fr}
\end{equation}

The arrows on the top commute with the $\bZ\times \bZ$-actions on all spaces.

Let $\Phi_{g_1g_2}$ be the pre-image of $(g_1,g_2)\in \Met\times\Met$ under the map
$$
\Phi\to \Met\times \Met.
$$

We then have groupoid composition maps
$$
C(g_1g_2\cdots g_n): \Phi_{g_1g_2}\times_{ \Fr_{g_2} }\Phi_{g_2g_3}\times_{\Fr_{g_3}} \cdots \times_{\Fr_{g_{n-1}}} \Phi_{g_{n-1}g_n}\to \Phi_{g_1g_n}.
$$

Let us define   subspaces
$$
K(g_1g_2\cdots g_n)^a\subset \Phi_{g_1g_2}\times \Phi_{g_2g_3}\times \cdots \times \Phi_{g_{n-1}g_n} \times \Phi_{g_1g_n}
$$
and
$$
K(g_1g_2\cdots g_n)^{\cyc,a}\subset \Phi_{g_1g_2}\times \Phi_{g_2g_3}\times \cdots \times \Phi_{g_{n-1}g_n} \times \Phi_{g_ng_1}
$$

similar to Sec. \ref{kaen}

The retraction (\ref{deform2}) induces  deformation retractions:
$$\xymatrix{
K(g_1g_2\cdots g_n)^a\ar[rr]\ar[d]&& K_n^a\ar[rr] \ar[d]&&K(g_1g_2\cdots g_n)^a\ar[d]\\
\Phi_{g_1g_2}\times \Phi_{g_2g_3}\times \cdots \times \Phi_{g_{n-1}g_n}\ar[rr]&&
Q^{n-1}\ar[rr]&&\Phi_{g_1g_2}\times \Phi_{g_2g_3}\times \cdots \times \Phi_{g_{n-1}g_n}}
$$
\subsubsection{The groupoid $\frakS$}\label{frakss} Let $S:=P\times_\frakH P$ we have a homotopy equivalence of groupoids
$$
(S\rightrightarrows P)\to (P\times_M P\rightrightarrows P),
$$

Let $$\frakS:=S\times_{P\times_M P} Q.$$ Equivalently, one can define $\frakS$ as follows.
Since $P\to \frakH$ is a principal $U(N)$-bundle, we have  a map
$$
S=P\times_\frakH P\to U(N).
$$
Let $\tilU(N)\to U(N)$ be the universal cover. Then $\frakS=S\times_{U(N)} \tilU(N).$

We  have an induced groupoid structure on $\frakS\rightrightarrows  P$ as well as a $\bZ\times \bZ$-equivariant groupoid maps

\begin{equation}\label{S2Phi}
(\frakS\rightrightarrows P)\to (Q\rightrightarrows P)\stackrel{p_{\Phi_{g_1g_2}}}\longrightarrow (\Phi_{g_1g_2}\rightrightarrows M).
\end{equation}

We define spaces $K^{a,\frakS}_n\subset \frakH^{n}$,  $K^{a,\frakS,\cyc}_n\subset \frakH^n$
similar to above.  
\subsubsection{The groupoid $\Sigma$} \label{sigmo}  Let $\cL:=P/\SU(N)$.  We then get an $S^1\times S^1$- bundle $\cL\to \frakH$, hence,
a map $\cL\times_\frakH \cL\to S^1$.  Let $$\Sigma:=(\cL\times_\frakH \cL)\times_{S^1\times S^1}\Re\times \Re$$
We have a groupoid structure on
$
\Sigma\rightrightarrows \cL
$

as well as a  $\bZ\times \bZ$-equivariant groupoid  map
$$
(\frakS\rightrightarrows P)\to (\Sigma\rightrightarrows \cL).
$$

One defines the spaces $K^{\Sigma}_n,K^{\Sigma,\cyc}_n\subset \Sigma^n$ in the same way as above.

\subsection{Symplectic balls} \label{radius}
Let $g\in \Met$.

Let $B_{g,R}M\subset TM$ be the sub-bundle of balls of $g$-radius $R$. Let $\pi_{g,R}:B_{g,R}M\to M$ be the projection.
 
One can define a function $R'':\Met\to \Re_{>0}$ satisfying:   

(1) for every $g\in \Met$, 
we have a map  $I^g:B_{R''(g)}^gM\to M$, where for each $m\in M$, the induced map
$I^g_m:(\pi_{g,R''(g)})^{-1}m\to M$ is a symplectic embedding such that $I^g_m(0)=m$. For $R\leq R''(g)$, denote $B_{R,g}(m):=I^g_m((\pi_{g,R})^{-1}m)$.

(2) for all $R\leq R''(g)$ and all $m\in M$,  
$B_{g,R}(m)\subset M$
is $g$-geodesically convex.

(3) Let $m_0\in M$ and fix an identification $\psi_{m_0}:B_{R''(g)}\cong B_{g,R''}(m_0)$.
 For each $n\in B_{R''(g)}(m_0)$  let $\phi_n: B_{R''(g)}\cong B_{g,R''}(n)\to M$
be the embedding such that the local coordinates near $n$  coming from $\psi_{m_0}$ and $\phi_n$
have coincident differentials at $n$.  This way, we have  a family of symplectic embeddings
$\phi:B_{R''(g)}\times B_{R''(g)}\to M$.   We require this family to be mobile.

For any function $\phi:\Met\to \Re_{>0}$ which is invariant under the action of symplectomorphisms, one can  define another such a function $\mu_\phi)$ satisfying:

--- let  $m_i,i=1,100$ be points in $M$ such that
$B_{g,\mu_\phi(g)}(m_i)\cap B_{g,\mu_\phi(g)}(m_{i+1})\neq \emptyset$, $1\leq i<100$.  Then for all $i$,  $m_i\in B_{\phi(g)}(m_1)$.

Let $R'=\mu_{R''}$,  $R=\mu_{R'}$, and  $r=\mu_R$.

Let $\MetR$ be the set of pairs $(g,R)$, where $g\in \Met$ and  we assume$R<r(g)$.
 We have a partial order on $\MetR$: $(g_1,R_1)\leq (g_2,R_2)$
iff $B_{g_1,R_1}(m)\subset B_{g_2,R_2}(m)$ for all $m\in M$.
 Write $(g_1,R_1)<<(g_2,R_2)$ if  for all $m_1,m_2\in M$, $B^{g_1}_{R_1}(m_1)\cap B^{g_2}_{R_2}(m_2)\neq \emptyset$ implies
that $B^{g_2}_{R_2}(m_2)\subset \Int B_{R(g_1)}^{g_1}(m_1)$.  Note that $<<$ is not a partial order.

It follows immediately that for all $i\in \MetR$ we have
$i<<i$.
 
Let $i_k=(g_k,R_k)\in \MetR$.   We set $\Phi_{i_1i_2}:=\Phi_{g_1g_2}$.   

\section{More on categories} 
\subsection{A poset  $\SMetR$}   Fix a subset $\MetR'\subset \MetR$. We assume that
$\MetR'$ is filtered and that for every $g\in \MetR'$ there
exists an $i\in \MetR$ such that $g<<i$.

 Define a poset $\SMetR$ as follows.
An element of $\SMetR$ is
 linearly ordered  sub-set  $S\subset \MetR'$.
Write $S_1\geq S_2$ if $S_1\subset S_2$ (sic!).

Let $\SMetR'\subset \SMetR$ consist of  all non-empty subsets $S$.

Let $\mu:\SMetR'\to \MetR'$ be the following monotone map 
$\mu(S)=\min(S)$. 

Suppose the minimal element of $S$ is of the form $(g,R)$.
We then write $R_{\min}(S)=R$.
   Let also $\val:\SMetR'\to \Re_{>0}$ be
given by $\val(S)=\pi R^2_{\min}$.   Let us also set $\val(\emptyset)=\infty$.

\subsection{Categories $\bR_q$, $\bR_0$}   Let $\cC$ be a  category enriched over the ground category $\cA$ and
$\Doplus$-closed.

An object of $\bR_q(\cC)'$, enriched over $\ccA$, is a family $\{X_S\}_{S\in \SMetR}$, where 
$$
X_S\in \Com(\cC)\langle \val(S)\rangle.
$$
Set
$$
\hom(X,Y):=\prod\limits_{T\subset S}\hom(i_{\val T\val S} X_T;Y_S)\in \ccA.
$$
The composition is well defined because, given $T,S\in \SMetR$, there are only finitely many 
$R\in \SMetR$, where $T\subset R\subset S$.

Set $\bR_q(\cC):=D\bR_q(\cC)'$.

An object of $\bR_0'(\cC)$ is a collection $\{X_S\}_{S\in \SMetR}$, where $X_S\in \Gr(\cC)\langle \val S\rangle$.
Set
$$
\hom(X,Y):=\prod\limits_{T\subset S}\hom(i_{\val T\val S} X_T;Y_S)\in \ccA.
$$
Set $\bR_0(\cC):=D\bR_0(\cC)'$.

Suppose $\cC$ is a ground category,  then so are $\bR_q(\cC)$, $\bR_0(\cC)$.  Indeed,  
we have
$$
(\bigoplus\limits_{\alpha}X^\alpha)_S=\bigoplus\limits_{\alpha} (X^\alpha)_S;
$$
$$
(\prod\limits_{\alpha}X^\alpha)_S=\prod\limits_{\alpha} (X^\alpha)_S.
$$

Let $X,Y\in \bR_q(\cC)$ (resp $\bR_0(\cC)$).  Define the tensor product by
$$
(X\otimes Y)_S:=\bigoplus\limits_{S_1,S_2|S_1\cup S_2=S} i_{\val(S_1)\val(S)}X_{S_1}\otimes i_{\val(S_2)\val(S)}X_{S_2}.
$$
Define the inner hom
$$
\ihom(X,Y)_S:=\prod\limits_{T|T\subset S}\ihom(i_{\val T\val S} X_T;Y_S).
$$

In both cases above, the differential is determined by those on $X$ and $Y$.

\subsection{Category $\bS_q(\cC)$ enriched over $\bR_q(\ccA)$} An object of $\bS_q(\cC)'$ is a collection
$\{X^S\}_{S\in \SMetR}$, where 
$$
X^S\in \Quant(\cC)\langle \val(S)\rangle.
$$
Define $\hom(X,Y)\in \bR_q(\ccA)'$:
$$
\hom(X,Y)_S:=\prod\limits_{T\subset S}\hom(X^S;i_{\val T\val S} Y^T)\in \Com(\ccA)\langle \val(S)\rangle.
$$
The composition is well defined because, given $T,S\in \SMetR$, there are only finitely many 
$R\in \SMetR$, where $T\subset R\subset S$.

Set $\bS_q(\cC):=D\bS_q(\cC)'$.

An object of $\bS_0'(\cC)$ is a collection $\{X^S\}_{S\in \SMetR}$, where $X^S\in \Classic(\cC)\langle \val S\rangle$.
Define $\hom(X,Y)\in \bR_0(cA)$.
$$
\hom(X,Y)_S:=\prod\limits_{T\subset S} \ihom(X^S;i_{\val T\val S} Y^T)\in \Gr(\ccA)\langle \val(S)\rangle.
$$
Set $\bS_0(\cC):=D\bS_0(\cC)'$.

Suppose $\cC$ is closed under direct products/sums,  then so are $\bS_q(\cC)$, $\bS_0(\cC)$.  Indeed,  
we have
$$
(\bigoplus\limits_{\alpha}X^\alpha)^S=\bigoplus\limits_{\alpha} (X^\alpha)^S;
$$
$$
(\prod\limits_{\alpha}X^\alpha)^S=\prod\limits_{\alpha} (X^\alpha)^S.
$$

Suppose $\cC,\cD$ are categories enriched over $\ccA$.  Define a functor
$$
\boxtimes \bS_q(\cC)\otimes \bS_q(\cD)\to \bS_q(\Doplus \cC\otimes \cD).
$$

$$
(X\boxtimes Y)^S:=\bigoplus\limits_{S_1,S_2|S_1\cup S_2=S}( i_{\val(S_1)\val(S)}(X^{S_1}\boxtimes i_{\val(S_2)\val(S)}Y^{S_2}),
$$

where $\boxtimes$ on the RHS denotes the functor
$$
\boxtimes:\Classic(\cC)\langle \val(S)\rangle\otimes \Classic(\cD)\langle \val(S)\rangle \to
\Classic(\Doplus \cC\otimes \cD)\langle \val(S)\rangle.
$$
The differential is determined by those on $X$ and $Y$.

One defines a functor $$\boxtimes:\bS_0(\cC)\otimes \bS_0(\cD)\to \bS_0(\Doplus \cC\otimes \cD).
$$
in a similar way.  

It now follows that given an SMC $\cC$ with direct sums, $D$-closed, such that the direct sums and differentials
are compatible with tensor product,  we have an induced  SMC structure on $\bS_q(\cC)$, $\bS_0(\cC)$.

We also have a functor

$$
\bS_q(\cC)\otimes \bR_q(\cD)\to \bS_q(\Doplus\cC\otimes \cD),
$$

where 
$$
(X\otimes  U)^S=(X^S\otimes\bigoplus\limits_{T|T\subset S} U_T,D_S),
$$
where the differential $D_S$ is induced by the differential on $U$.

One defines a functor $$\bS_0(\cC)\otimes \bR_0(\cD)\to \bS_0(\Doplus \cC\otimes \cD).
$$

In particular,  $\bS_q(\cC)$ is tensored over $\bR_q(\cD)$ if $\cC$ is tensored over $\cD$, and likewise
for $\bS_0(\cC)$.

Let $X$ be a locally compact topological space.  Let $\cD$ be a category.   Denote
$$\psh(X,\cD,\bS_q):=\bS_q(\cD\otimes\open_X),\  \psh(X,\cD,\bS_0):=\bS_0(\cD\otimes \open_X).
$$

\subsubsection{Lemma on truncation} 
The following Lemma follows from Lemma \ref{trunc1}
\begin{Lemma}\label{trunc2}
Let $X\in \bR_0(\ccA)$ and suppose  that

1)  for each $S\in \SMetR$,  we have $\gr^0 X_{\leq S}\in \trunc \ccA$;

2)  $L^0:=\projlim_{S\in \SMetR} H^0\gr^0 X_{\leq S}$ is a free abelian group.

Then $L:=\hom(\unit_{\bR_0(\ccA)};X)\in \trunc \ccA$ and $H^0(L)\cong L^0$.
\end{Lemma}
\subsection{The category $\ZZ(\cC)$} 
\subsubsection{A partially ordered monoid $\ZZ$}
 Let $\ZZ$ be a partially ordered commutative monoid, where $\ZZ=\{\infty\}\sqcup \bZ\times \bZ$,
where we endow $\bZ\times \bZ$ with the discrete poset structure and declare $\infty$ to be the greatest object.
The addition  on $\ZZ$ is the group addition, when restricted  on $\bZ\times \bZ$, and we set $\infty+a=\infty$
for all $a\in \bZ\times \bZ$.
 
\subsubsection{ A category $\ZZ(\cC)$}   Let $\cC$ be a SMC enriched over another SMC $\cA$ which is
$D,\oplus,\prod$-closed.

  Let $\ZZ(\cC)'$
be a category, enriched over $\ccA$, whose every object is a $\ZZ$-graded object  $X_\bullet$ in $\cC$ and
we set
$$
\hom(X,Y)=\prod_{s\in \bZ\times \bZ}(\hom(X_s,Y_s)\oplus \hom(X_s,Y_\infty))\oplus \hom(X_\infty,Y_\infty).
$$
That is,
$$
\hom(X,Y)=\prod\limits_{s\geq t} \hom(X_s,Y_t).
$$

Set
$$
(X\otimes Y)_s=\bigoplus_{t_1,t_2\in \ZZ|t_1+t_2=s} X_{t_1}\otimes Y_{t_2}.
$$

The category $\ZZ(\cC)'$  is enriched over $\ZZ(\cA)'$:

 --- for $a\in \bZ\times \bZ$, we have
$$
\ihom(X,Y)_a=\prod\limits_{b\in \bZ\times \bZ}(X_b,Y_{a+b});
$$
$$
\ihom(X,Y)_\infty=\prod_{b\in \ZZ}\ihom(X_b;Y_\infty).
$$

Set $\ZZ(\cC):=D\ZZ(\cC)'$.   The cateory $\ZZ(\cC)$ has a tensor structure and is enriched over $\ZZ(\cA)$.
\section{Constructing the categories}
Throughout this section,  the categories are assumed to be enriched over $\bR_0(\ZZ(\ccA))$ unless 
otherwise specified.
\subsection{The categories $\bbF,\bbFR$}  \label{bbfr1}
 Let $\Fr:=\bigsqcup_{i\in \MetR} \Fr_i$.   
Let
$$
\bbF:=\psh(\Fr\times E\times \Fr\times E,\bS_q(\ZZ(\ccA)));
$$

$$
\bbFR:=\psh(\Fr\times E,\bS_q(\ZZ(\ccA))).
$$

Let $h:\bbF\boxtimes \bbF\to \bS_q(\ZZ(\ccA))$ be given by $h(U,V)=\unit$ if $U\cap V\neq \emptyset$
and $h_\bbF(U,V)=0$ otherwise.  Define the functor 
$h_R:\bbFR\boxtimes \bbFR\to \bS_q(\ZZ(\ccA))$ in a similar way.

Denote $\cE:=\Fr\times E$.
Let $$F_{\bbFR}^{[n]}\in \psh(\cE\times \cE)^n\times (\cE\times \cE),\bS_q(\ZZ(\ccA))$$
be given by $$F_\bbF^{[n]}:=\ZQ_{\Delta^{[n]}}\otimes \unit ,$$
where $\Delta^{[n]}\subset (\cE\times \cE)^n\times (\cE\times \cE)$ consists of all points of the form
$$
((e_1,e_2),(e_2,e_3),\ldots,(e_{n-1},e_n),(e_1,e_n))
$$

Define 
$$
F_{\bbFR}^{[n,1]}=F_\bbF^{[n]}\in \psh((\cE\times \cE)^n\times \cE\times \cE),\bS_q(\ZZ(\ccA))).$$

Let $F_{\bbFR}^{(n)}\in \psh((\cE\times \cE)^{\times n+1},\bS_q(\ZZ(\ccA)))$ be given by 
$$F_\bbF^{(n)}:=\ZQ_{\Delta^{(n)}}\otimes \unit [-2N],$$
where $\Delta^{(n)}\subset (\cE\times \cE)^{n+1}$ consists of all points of the form
$$
((e_0,e_1),(e_1,e_2),\ldots,(e_{n-1},e_n),(e_n,e_0)).
$$

Via the functor $\red:\DMetR_q\to \DMetR_0$, the category $\bbF$ is also defined  over $\DMetR_0$.
\subsection{The category $\bbG$}

 Set $\Phi:=\bigsqcup\limits_{i,j\in \MetR} \Phi_{ij}$.
Let $K^{[n]}\subset \Phi^n\times \Phi$ be the union
$$
(K^{[n]})^{ab}:=\bigsqcup\limits_{i_0,i_1,\ldots,i_n\in \MetR} K^{ab}(i_0,i_1,\ldots,i_n);
$$
$$
(K^{(n)})^{ab}:=\bigsqcup\limits_{i_0,i_1,\ldots,i_n\in \MetR} (K^\cyc)^{ab}(i_0,i_1,\ldots,i_n);
$$
Let $\Delta_E^{[n]}\subset (E\times E)^n\times (E\times E)$ consist of all points
of the form $((e_0,e_1),(e_1,e_2),\ldots,(e_{n-1},e_n);(e_0,e_n))$ and
$\Delta_E^{(n)}\subset (E\times E)^{n+1}$ consist of all points of the form
$((e_0,e_1),(e_1,e_2),\ldots,(e_{n-1},e_n),(e_n,e_0))$.

Set $$\bbG:=\psh(\Phi\times E\times E, \bS_0(\ZZ(\ccA))).
$$
Set $h_\bbG:\bbG\otimes \bbG\to  \bS_0(\ZZ(\ccA))$ be defined by the same rule as $h_\bbF$.
Let  $$F_{\bbG}^{[n]}\in 
\psh((\Phi\times E\times E)^{n+1},\bS_0(\ZZ(\ccA)))
$$
be given by
$$
F_{\bbG}^{[n]}:=\bigoplus\limits_{ab}\ZQ_{(K^{ab})^{[n]}\times \Delta_E^{[n]}}\otimes (a,b)\otimes \unit_{\DMetR_0};
$$
let
$$F_{\bbG}^{(n)}\in 
\psh(\open_{(\Phi\times E\times E)}^{n+1}, \DMetR_0)
$$
be given by
$$
F_{\bbG}^{(n)}:=\bigoplus\limits_{ab}\ZQ_{(K^{ab})^{(n)}\times \Delta_E^{(n)}}\otimes (a,b)\otimes \unit_{\DMetR_0};
$$

\subsubsection{Quasi-contractibility}  
The quasi-contractibility of both $F_\bbF$, $F_\bbG$ can be easily verified.
\subsection{Connecting $\bbG$ and $\bbF$}\label{uii}
\subsubsection{Constructing an object $\bGamma\in \psh(\Phi\times E\times E\times  \Fr\times \Fr\times E\times E,\bS_0(\ZZ(\ccA))$}

Let $\pi:\Phi\to \Fr\times_M \Fr$ be the projection.  Let $$\sigma:
\Phi\times E\times E\times \Fr\times \Fr\times E\times E\to \Fr\times_M \Fr\times E\times E\times \Fr\times\Fr\times E\times E
$$ be the induced map.      Let $\pr:\Fr\to \MetR$  be the projection. 
Let 
$$
i:\Fr\times_M \Fr\times E\times E\times \Fr\times \Fr\times E\times E\to\Fr\times \Fr\times E\times E \times \Fr\times\Fr\times E\times E
$$
be the closed embedding.

Let 
$$
j:\Fr\times_\MetR \Fr\times E\times E\times \Fr\times_\MetR\Fr\times E\times E\to \Fr\times \Fr\times E\times E \times \Fr\times\Fr\times E\times E
$$
be  the closed embedding defined as follows:
$$
j(f_1,f_2,v_1,v_2,f'_1,f'_2,v'_1,v'_2)=(f_1,f'_1,v_1,v'_1,f_2,f'_2,v_2,v'_2).
$$

We therefore have the following diagram
$$\xymatrix{
\Phi\times E\times E\times \Fr\times \Fr\times E\times E\ar[d]^\vs&\\
\Fr\times_M \Fr\times E\times E\times \Fr\times \Fr\times E\times E\ar[r]^{i}&\Fr\times \Fr\times E\times E \times \Fr\times\Fr\times E\times E\\
&\Fr\times_\MetR \Fr\times E\times E\times \Fr\times_\MetR \Fr\times E\times E\ar[u]^{j}}
$$

Let $i=(g,r)\in \MetR$ let $R=\rad g/2$.
Let $V_{ii}\subset \Fr_i\times \Fr_i$ consist of all pairs $(f_1,f_2)$ where $B_{(g,r)}(f_1)\subset B_{(g,R)}(f_2)$.
We then have a family of symplectic embeddings  $V_{ii}\times B_r\to B_R$.  Let $D:V_{ii}\to \Sp(2N)$
be the differential map at $0\in B_r$.   Let $U_{ii}\subset V_{ii}$ be  $D^{-1}e_{\Sp(2N)}$.

We now get a graded family of symplectic embeddings
$$
U_{ii}\times B_r\to B_R
$$
Let $$
\cP_{i}\in \sh_q(U_{ii}\times \Re^N\times \Re^N)
$$ be the object determined by this family.

Let $\iota:U_{ii}\hookrightarrow \Fr_i\times \Fr_i;$ 
be the embedding.

   Let $$
\bbP_i:=\iota_!\cP_i\in \sh_q(\Fr_i\times \Fr_i\times E\times E).
$$
  The objects $\bbP_i$ give rise to  an object
Let $\bbP_\MetR\in \sh_q(\Fr\times_\MetR \Fr\times E\times E).$

Set $$\bGamma':=\vs^{-1}i^{-1}j_!(\bbP_\MetR\boxtimes \bbP_\MetR)
$$

Set $$\bGamma:=\bGamma'\otimes \unit_{\bS_0} \in
\psh(\Phi\times E\times E\times \Fr\times \Fr\times E\times E,\bS_0(\ZZ(\ccA))).
$$
\subsubsection{Quasi-Contractibility}
Let $\bt$ be a planar  tree  with $n$ inputs. One sees that we have a homotopy equivalence
 $F_\bbG(\bt)\sim F_\bbG^{[n]}$.   Similarly, for a cyclic tree with $n+1$ inputs we have
a homotopy equivalence
$
F_\bbG(\bt)\sim F_\bbG^{(n)}.
$

The problem now reduces to showing the quasi-contractibility of
\begin{equation}\label{urav}
\hom(F_\bbG^{[n]}\circ \cL^{\boxtimes n};F_\bbF^{[n]}\circ \cL);
\end{equation}
\begin{equation}\label{uravcyc}
\hom(F_\bbG^{(n)}\circ \cL^{\boxtimes (n+1)};F_\bbF^{(n)}).
\end{equation}

Let us first show the quasi-contractibility of the object in  (\ref{urav}).  
Consider the following diagram
$$   \xymatrix{          (\Fr\times \Fr\times E\times E)^{2n}&\\
(\Fr\times \Fr\times E\times E)^{n}\times (\Fr\times \Fr\times E\times E)^n\ar[u]^\vs&\\
K^{[n]}\times E^{n+1}\times (\Fr\times \Fr\times E\times E)^n\ar[r]^{d}\ar[u]^{\tilde{i}}\ar[d]^{\tilde{p}}& \Fr^{\times_M^{n+1}}\times E^{n+1}\times (\Fr\times \Fr\times E\times E)^n\ar[ul]^i\ar[d]^p\\
\Phi\times E\times E\times (\Fr\times \Fr\times E\times E)^n\ar[r]^a& \Fr\times_M\Fr\times E\times E\times(\Fr\times \Fr\times E\times E)^n\\
\Phi\times E\times E\times (\Fr\times E)^{n+1}\ar[r]^b \ar[u]^{\tilde{j}}\ar[d]^{\tilde\rho}&\Fr\times_M\Fr\times E\times E\times  
(\Fr\times E)^{n+1}\ar[u]^j\ar[d]^\rho\\
\Phi\times E\times E\times \Fr\times \Fr\times E\times E\ar[r]^c& \Fr\times \Fr\times E\times E\times \Fr\times \Fr\times E\times E}
$$

Here:

---\begin{multline*}
\vs\big(\left((f_1,f_2,e_1,e_2),(f_3,f_4,e_3,e_4),\ldots,(f_{2n-1},f_{2n},e_{2n-1},e_{2n})\right)\\
\times \left((f'_1,f'_2,e'_1,e'_2),(f'_3,f'_4,e'_3,e'_4),\ldots,(f'_{2n-1},f'_{2n},e'_{2n-1},e'_{2n})\right)\big)\\
=((f_1,f_1',e_1,e'_1),(f_2,f_2',e_2,e_2'),\ldots,(f_{2n},f_{2n'},e_{2n},e_{2n'})).
\end{multline*}

--- $i=i_1\times \Id$, where  $$
i_1:\Fr^{\times_M^{n+1}}\times E^{n+1}\to (\Fr\times \Fr\times E\times E)^{n}
$$
is given by
$$
i_1(f_0,f_1,f_2,\ldots,f_{n-1},f_n,e_0,e_1,\ldots,e_n)=(f_0,f_1,f_1,f_2,f_2,\ldots,f_{n-1},f_{n-1},f_n,e_0,e_1,e_1,\ldots,e_{n-1},e_{n-1},e_n).
$$

--- The map $d$ is induced by the covering map $K^{[n]}\to \Fr^{\times_M^{n+1}}$;

--- set $\tilde{i}:=i\circ d$;

--- The map $p$ is induced  by the projection onto the marginal factors $\Fr^{\times_M^{n+1}}\to \Fr\times_M \Fr$;
$E^{n+1}\to E\times E$;

---  The map $\tilde{p}$ is defined in a similar way to $p$;

--- The map $a$ is induced  by the covering map $\Phi\to \Fr\times_M \Fr$.

--- The maps $j,\tilde{j}$ are induced by the map $$
j_1:(\Fr\times E)^{n+1}\to (\Fr\times \Fr\times E\times E)^{n},
$$
where
 $$
j_1(f_0,e_0,f_1,e_1,\ldots,f_n,e_n)=(f_0,e_0,f_1,f_1,e_1,e_1,\ldots,f_{n-1},f_{n-1},e_{n-1},e_{n-1},f_n,e_n).
$$

--- The map  $b$ is induced by the covering map $b:\Phi\times \Fr\times_M \Fr$.

--- The map $\rho$ is induced by the closed embedding $\Fr\times_M \Fr\to \Fr\times \Fr$ and
by the projection onto the marginal factors $$
(\Fr\times E)^{n+1}\to \Fr\times \Fr\times E\times E.
$$
The map $\tilde{\rho}$ is induced by the above projection;

---the map $c$ is a through map
$$
\Phi\to \Fr\times_M \Fr\to \Fr\times \Fr,
$$
where the left arrow is the covering map and the left arrow is the closed embedding.

Let  now $$\cP_{2n}:=\bbP^{\boxtimes 2n}\in \sh_q((\Fr\times \Fr\times E\times E)^{2n}).
$$
Let $$\cP_2:=\bbP\boxtimes \bbP\in \sh_q((\Fr\times \Fr\times E\times E)^{2}).
$$
We then rewrite (\ref{urav}) as follows:
\begin{multline*}
\hom(\tilde{p}_!(\tilde{i})^{-1}(\sigma^{-1}\cP_{2n});(\tilde{j})^{-1}\tilde{\rho}^{-1}c^{-1}\cP_2)\sim
\hom(\tilde{p}_!(\tilde{i})^{-1}(\sigma^{-1}\cP_{2n});(\tilde{j})^{-1}b^{-1}\rho^{-1}\cP_2)\\
\sim
\hom(\tilde{p}_!(\tilde{i})^{-1}(\sigma^{-1}\cP_{2n});a^{-1}j_!\rho^{-1}\cP_2)\sim
\hom(a_!\tilde{p}_!(\tilde{i}^{-1}\sigma^{-1}\cP_{2n});j_!\rho^{-1}\cP_2)\\
\sim\hom(p_!d_!d^{-1}i^{-1}\sigma^{-1}\cP_{2n};a^{-1}j_!\rho^{-1}\cP_2)\sim
\hom(p_!(i^{-1}\sigma^{-1}\cP_{2n}\otimes \lambda);j_!\rho^{-1}\cP_2)\\
\sim\hom(j^{-1}p_!(i^{-1}\sigma^{-1}\cP_{2n}\otimes \lambda);\rho^{-1}\cP_2),
\end{multline*}
where
$$
\lambda=d_!\ZQ_{K^{[n]}\times E^{n+1}\times (\Fr\times \Fr\times E\times E)^n}.
$$

Let us now consider the following diagram:
$$ \xymatrix{            (\Fr\times E\times \Fr\times E)^{2n}&&\\
(\Fr\times E)^{2n}\times (\Fr\times E)^{2n}\ar[u]^\vs&&\\
\Fr^{\times_M^{n+1}}\times E^{n+1}\times (\Fr\times E )^{2n}\ar[d]^{p_1}\ar[u]^{i}&\Fr^{\times_M^{n+1}}\times E^{n+1}\times (\Fr\times E)^{n+1}\ar[r]^{l_1}\ar[l]^{k_1}\ar[ul]^m\ar[d]^{q_1}& (\Fr\times E)^{n+1}\times (\Fr\times E)^{n+1}\ar[ull]_n\ar[d]^r\\
\Fr^{\times_M^{n+1}}\times E^2\times (\Fr\times  E)^{2n}\ar[d]^{p_2}&\Fr^{\times_M^{n+1}}\times E^2\times (\Fr\times E)^{n+1}\ar[l]^{k_2}\ar[d]^{q_2}\ar[r]^{l_2}&
\Fr^{n+1}\times E^2\times (\Fr\times E)^{n+1}
\\
\Fr\times_M \Fr\times E^2\times (\Fr\times E)^{2n}
&\Fr\times_M \Fr\times E^2\times (\Fr\times E)^{n+1}\ar[l]^j\ar[d]^\rho&\\
&\Fr \times \Fr\times E\times  E\times \Fr\times \Fr\times E\times E&}
$$

Here:

--- $\vs$ and $i$ are as above;

--- $p_1,q_1,r$ are induced by the projection onto the marginal factors $E^{n+1}\to E^2$;

--- $p_2,q_2$ are induced by the projection onto the marginal factors $\Fr^{\times_M n+1}\to \Fr\times_M \Fr$.

---the maps $k_1,k_2,j$ is induced by the following map $$
k':(\Fr\times E)^{n+1}\to (\Fr\times E\times \Fr\times E)^{2n}:
$$
$$
k'(f_0,e_1,f_1,e_1,\ldots,f_n,e_n)=(f_0,e_0,f_1,e_1,f_1,e_1,f_2,e_2,\ldots,f_{n-1},e_{n-1},f_{n-1},e_{n-1},f_n,e_n).
$$
The map $j$ coincides with that on the previous diagram.

--- The maps $l_1,l_2$ are induced by the closed embedding $\Fr^{\times_M^{n+1}}\to \Fr^{\times {n+1}}$

--- set $n=k'\times k'$,  set $m=nl_1$.

We have $p=p_2p_1$.    As all squares in the diagram are Cartesian, we have:
\begin{multline*}
j^{-1}p_!(\lambda\otimes i^{-1}\sigma^{-1}\cP_{2n})\sim j^{-1}p_{2!}(\lambda\otimes p_{1!}i^{-1})\sim q_{2!}(\lambda\otimes k_2^{-1}p_{1!}i^{-1}\sigma^{-1}\cP_{2n})\sim q_{2!}(\lambda\otimes q_{1!}k_1^{-1}i^{-1}\sigma^{-1}\cP_{2n})
\\ \sim
q_{2!}(\lambda\otimes q_{1!}m^{-1}\sigma^{-1}\cP_{2n})=q_{2!}(\lambda\otimes q_{1!}l^{-1}n^{-1}\sigma^{-1}\cP_{2n})\sim q_{2!}(\lambda\otimes(l_2^{-1}r_!n^{-1}\sigma^{-1}\cP_{2n}))
\end{multline*}

Let $$n_0:\Fr\times E\times \Fr\times E\to (\Fr\times E\times \Fr\times E)^2$$ be the diagonal map
$$
n_0(f_1,e_1,f_2,e_2)=((f_1,e_1,f_2,e_2),(f_1,e_1,f_2,e_2)).
$$
Let 
$$r_0:\Fr\times E\times \Fr\times E\to \Fr\times \Fr\times E$$
be the projection along the first factor of $E$.

Let $$\tau: (\Fr\times E)^{n+1}\times (\Fr\times E)^{n+1}\to (\Fr\times E\times \Fr\times E)^{n+1};
$$
be given by \begin{multline*}
\tau((f_0,e_0,f_1,e_1,\ldots,f_n,e_n),(f'_0,e'_0,f'_1,e'_1,\ldots,f'_n,e'_n))\\
\to ((f_0,e_0,f'_0,e'_0),(f_1,e_1,f'_1,e'_1),\ldots,(f_n,e_n,f'_n,e'_n))
\end{multline*}

We have 
$$
\sigma n=(\Id_{\Fr\times E\times \Fr\times E}\times n_0^{\times n-1}\times \Id_{\Fr\times E\times \Fr\times E})\tau
$$

Let 
$$
\tau_2:\Fr\times E\times \Fr\times E\times (\Fr\times \Fr\times E)^{n-1}\times \Fr\times E\times\Fr\times E
\to \Fr^{n+1}\times E^2\times (\Fr\times E)^{n+1}
$$
be  given by
\begin{multline*}
\tau_2(f_0,e_0,f'_0,e'_0),(f_1,f'_1,e'_1),(f_2,f'_2,e'_2),\ldots,(f_{n-1},f'_{n-1},e'_{n-1}),(f_n,e_n,f'_n,e'_n))\\
=((f_0,f_1,\ldots,f_n),(e_0,e_n),(f'_0,e'_0,f'_1,e'_1,\ldots,f'_{n+1},e'_{n+1})).
\end{multline*}

We now have 
$$
r= \tau_2    (\Id_{\Fr\times E\times\Fr\times E}\times r_0^{n-1}\times \Id_{\Fr\times E\times \Fr\times E})   \tau^{-1}
$$

so that
$$
 r_!n^{-1}=\tau_{2!}( \Id_{\Fr\times E\times \Fr\times E}\boxtimes (r_{0!}n_0^{-1})^{\boxtimes n-1}\boxtimes \Id_{\Fr\times E\times \Fr\times E}).
$$
Therefore,
$$
r_!n^{-1}\cP_{2n}\sim \tau_{2!}(\bbP\boxtimes(\boxtimes_{i=1}^{n-1} r_{0!}n_0^{-1}(\bbP\boxtimes 
\bbP))\boxtimes \bbP).
$$
Let $P_\MetR\in \sh_q(\Fr\times E\times E)$ be defined as follows:
$$
P_\MetR|_{\Fr_i\times E\times E}:= \ZQ_{\Fr_i}\boxtimes P_{r(i)}.
$$
Let $\delta:\Fr\times E\to\Fr\times E\times E$ be the diagonal embedding.
Let $V_{ii},U_{ii}\subset \Fr_i\times \Fr_i$  be as in Sec. \ref{uii}.  Let
 $$
\cU:=\bigsqcup\limits_{i\in \MetR} U_{ii}\subset \bigsqcup\limits_{i\in \MetR} \Fr_i\times \Fr_i\subset
\Fr\times \Fr.
$$

Let $p_1:\Fr\times \Fr\times E\to E$;  $p_2:\Fr\times \Fr\times E\to \Fr\times \Fr$ be projections, let $\delta:\Fr\times E\to\Fr\times  E\times E$ be the diagonal embedding.
We have $$
r_{0!}n_0^{-1}(\bbP\boxtimes \bbP)\sim   p_1^{-1}\delta^{-1} P_\MetR\otimes p_2^{-1}\ZQ_{\cU}.
$$

Let $\kappa:(\Fr\times \Fr)^{n+1}\to \Fr^{n+1}\times \Fr^{n+1}$ be the natural map.

Let $A:=l_2^{-1}(\kappa(\cU^{ {n+1}})\times E^{2}\times E^{n+1})$.
It follows that $q_2|_A:A\to \Fr\times_M \Fr\times E^2\times (\Fr\times E)^{n+1}$ ;
$l_2|_A:A\to \Fr^{n+1}\times E^2\times (\Fr\times E)^{n+1}
$
are embeddings of
a locally closed subset.  

Denote $i_A:=q_2|_A$. 
 Let us decompose 
$$
\Fr\times_M\Fr\times E^2\times (\Fr\times E)^{n+1}=\Fr\times_M \Fr\times E^2\times \Fr^2\times E^2\times
(\Fr\times E)^{n-1}
$$

Let $$\iota: \Fr\times_M \Fr\times E^2\times \Fr^2\times E^2\to \Fr\times \Fr\times E^2\times \Fr\times \Fr\times E^2$$ be the embedding. Let $q:A\to \Fr^{\times_M n+1}$ be the projection.
We then have
$$
q_{2!}(\lambda\otimes(l_2^{-1}r_!n^{-1})\cP_{2n})\sim  i_{A!}((\iota^{-1}(\bbP\boxtimes \bbP)\boxtimes 
(\delta^{-1} P_\MetR))\otimes q^{-1}\lambda).
$$

Next, we have
$\rho^{-1}\bbP_2\sim (\bbP\boxtimes \bbP)\boxtimes \ZQ$.

Let $B:=\rho^{-1}(\cU\times E^2\times \cU\times E^2)$.  Observe that $i_A(A)\subset B$ is an open subset.
Therefore, we have
$$
\hom(i_{A!}(\iota^{-1}(\bbP\boxtimes \bbP)\boxtimes 
(\delta^{-1} P_\MetR)\otimes q^{-1}\lambda);\rho^{-1}\bbP_2)\sim
\hom((\iota^{-1}(\bbP\boxtimes \bbP)\boxtimes 
(\delta^{-1} P_\MetR))\otimes q^{-1}\lambda;i_A^{-1}\rho^{-1}(\bbP\boxtimes \bbP))
$$
The convolution with $i_A^{-1}\rho^{-1}(\bbQ\boxtimes\bbQ)$ gives rise to a homotopy equivalence 
of the above hom with
\begin{multline*}
\hom(i_A^{-1}\Delta^{-1}(P_\MetR\boxtimes P_\MetR)\boxtimes (\delta^{-1}P_\MetR\boxtimes \delta^{-1}P_\MetR)^{n-1};
i_A^{-1}(\Delta^{-1}(P_\MetR\boxtimes P_\MetR)\boxtimes (\ZQ_{\Fr\times E})^{\boxtimes (n-1)}))\\
\stackrel{(*)}\sim
\hom(i_A^{-1}\Delta^{-1}(P_\MetR\boxtimes P_\MetR)\boxtimes (\delta^{-1}P_\MetR\boxtimes \delta^{-1}P_\MetR)^{n-1};i_A^{-1}(\Delta^{-1}(\ZQ_{[\Fr\times\Delta_E\times \Fr\times \Delta_E,0]})
\boxtimes (\ZQ_{\Fr\times E})^{\boxtimes (n-1)}))
\end{multline*}
where $$\Delta:\Fr\times_M\times \Fr\times E\times E\times \Fr\times \Fr\times E\times E\to
\Fr\times \Fr\times E\times E\times \Fr\times \Fr\times E\times E,$$
$\Delta_E\subset E^2$ is the diagonal
and the map $(*)$ is induced by the universal map
$$
P_\Met\to \ZQ_{[\Fr\times \Delta_E,0]}.
$$
The universality of this map implies that $(*)$ is a homotopy equivalence.
The latter hom is homotopy equivalent to
$$
\hom(\cut_{\leq 0}i_A^{-1}\Delta^{-1}(P_\MetR\boxtimes P_\MetR)\boxtimes (\delta^{-1}P_\MetR\boxtimes \delta^{-1}P_\MetR)^{n-1};i_A^{-1}(\Delta^{-1}(\ZQ_{\Fr\times\Delta_E\times \Fr\times \Delta_E})
\boxtimes (\ZQ_{\Fr\times E})^{\boxtimes (n-1)})
$$
The resulting space is  pseudo-contractible, as both arguments of the hom are constant sheaves on
locally closed constructible sub-sets of $A$.

Consider  now the case (\ref{uravcyc}).  

$$
\xymatrix{(\Fr\times \Fr\times E\times E)^{2n+2}&\\
K^{(n)}\times E^{n+1}\times (\Fr\times \Fr\times E\times E)^{n+1}\ar[r]^{d}\ar[u]^{\tilde{i}}\ar[d]^{\tilde{p}}& \Fr^{\times_M^{n+1}}\times E^{n+1}\times (\Fr\times \Fr\times E\times E)^{n+1}\ar[ul]^i\ar[dl]^p\\
 (\Fr\times \Fr\times E\times E)^{n+1}&\\
(\Fr\times E)^{n+1}\ar[u]^j&
}
$$

Let $\delta:K^{(n)}\to \Fr^{\times_M^{n+1}}$ be the projection.  Let $\lambda^\cyc:=\delta_!\ZQ_{K^{(n)}}$.

The space in (\ref{uravcyc}) can be rewritten as follows:
\begin{multline}\label{cepko}
\hom((\tilde{p})_!(\tilde{i})^{-1}\cP_{2n+2};j_!\ZQ_{[(\Fr\times E)^{n+1},0]}[-2N])\sim
\hom(j^{-1}p_!d_!d^{-1}i^{-1}(\cP_{2n+2})[2N];\ZQ_{[(\Fr\times E)^{n+1},0]})\\
\sim \hom(j^{-1}p_!((\lambda^\cyc\boxtimes \ZQ_{E^{n+1}\times (\Fr\times \Fr\times E\times E)^{n+1}})\otimes i^{-1}\cP_{2n+2})[2N];\ZQ_{[(\Fr\times E)^{n+1},0]}),
\end{multline}

Let us now consider the following diagram

$$
\xymatrix{&(\Fr\times E\times \Fr\times E)^{2n+2}&\\
\Fr^{\times_M^{n+1}}\times E^{n+1}\times (\Fr\times E\times \Fr\times E)^{n+1}\ar[d]^{p_1}\ar[ur]^{i}&\Fr^{\times_M^{n+1}}\times E^{n+1}\times (\Fr\times E)^{n+1}\ar[r]^{l_1}\ar[l]^{k_1}\ar[u]^m\ar[d]^{q_1}& (\Fr\times E)^{n+1}\times (\Fr\times E)^{n+1}\ar[ul]^n\ar[d]^r\\
\Fr^{\times_M^{n+1}}\times (\Fr\times \Fr\times E\times E)^{n+1}\ar[d]^{p_2}&\Fr^{\times_M^{n+1}} \times (\Fr\times E)^{n+1}\ar[l]^{k_2}\ar[d]^{q_2}\ar[r]^{l_2}&
\Fr^{n+1}\times (\Fr\times E)^{n+1}
\\
(\Fr\times \Fr\times E\times E)^{n+1}
&(\Fr\times E)^{n+1}\ar[l]^j&}
$$
As every square in this diagram is Cartesian, and $p=p_2p_1$, we have
\begin{multline*}
j^{-1}p_!((\lambda^\cyc\boxtimes \ZQ_{E^{n+1}\times (\Fr\times \Fr\times E\times E)^{n+1}})\otimes i^{-1}\cP_{2n+2})\\
\sim j^{-1}p_{2!}p_{1!}((\lambda^\cyc\boxtimes \ZQ_{E^{n+1}\times (\Fr\times \Fr\times E\times E)^{n+1}})\otimes i^{-1}\cP_{2n+2})\\
\sim q_{2!}k_2^{-1}p_{1!}((\lambda^\cyc\boxtimes \ZQ_{E^{n+1}\times (\Fr\times \Fr\times E\times E)^{n+1}})\otimes i^{-1}\cP_{2n+2})\\
\sim  q_{2!}q_{1!}k_1^{-1}((\lambda^\cyc\boxtimes \ZQ_{E^{n+1}\times (\Fr\times \Fr\times E\times E)^{n+1}})\otimes i^{-1}\cP_{2n+2})\\
\sim q_{2!}q_{1!}((\lambda^\cyc\boxtimes \ZQ_{E^{n+1}\times (\Fr\times E)^{n+1}})\otimes m^{-1}\cP_{2n+2})\\
\sim q_{2!}q_{1!}((\lambda^\cyc\boxtimes \ZQ_{E^{n+1}\times (\Fr\times E)^{n+1}})\otimes l^{-1}n^{-1}\cP_{2n+2})\\
\sim   q_{2!}((\lambda^{\cyc}\boxtimes \ZQ_{(\Fr\times E)^{n+1}})\otimes l_2^{-1}  r_!n^{-1}\cP_{2n+2})
\end{multline*}

We can now continue the chain of homotopy  equivalences  (\ref{cepko})
\begin{multline*}
\hom(j^{-1}p_!((\lambda^\cyc\boxtimes \ZQ_{E^{n+1}\times (\Fr\times \Fr\times E\times E)^{n+1}})\otimes i^{-1}\cP_{2n+2})[2N];\ZQ_{[(\Fr\times E)^{n+1},0]})\\
\sim \hom(q_{2!}((\lambda^{\cyc}\boxtimes \ZQ_{(\Fr\times E)^{n+1}})\otimes l_2^{-1}  r_!n^{-1}\cP_{2n+2})[2N];
\ZQ_{[(\Fr\times E)^{n+1},0]})\\
\sim \hom((\lambda^{\cyc}\boxtimes \ZQ_{(\Fr\times E)^{n+1}})\otimes l_2^{-1}  r_!n^{-1}\cP_{2n+2};q_2^!\ZQ_{[(\Fr\times E)^{n+1},0]}[-2N])
\end{multline*}

Similar to the previous proof, we have
$$
r_!n^{-1}\cP_{2n+2}\sim\tau^{-1} (r_{0!}n_0^{-1}\cP_2)^{\boxtimes (n+1)}.
$$
where 
$$
\tau:\Fr^{n+1}\times (\Fr\times E)^{n+1}\to (\Fr\times \Fr\times E)^{n+1}
$$
is a permutation.

Let $\cW_1:=\tau^{-1} (\cU\times E)^{n+1}$; $\cW_2:=l_2^{-1}\cW_1$;  $\cW_3:=q_2(\cW_2)$.
Let $l_2':\cW_2\to \cW_1$; $q_2':\cW_2\to \cW_1$ be the induced maps.   It follows that
$q_2'$ is a  smooth fibration  whose fiber is diffeomorphic to $\Re^{2n}$ which implies
$(q_2')^![-2N]\sim (q_2')^{-1}$. We have
\begin{multline*}
\hom((\lambda^{\cyc}\boxtimes \ZQ_{(\Fr\times E)^{n+1}})\otimes l_2^{-1}  r_!n^{-1}\cP_{2n+2};q_2^!\ZQ_{[(\Fr\times E)^{n+1},0]}[-2N])\\
\sim
\hom((\lambda^{\cyc}\boxtimes \ZQ_{(\Fr\times E)^{n+1}})|_{W_2}\otimes (l'_2)^{-1} (\tau^{-1}(r_{0!}n_0^{-1}\cP_2)^{\boxtimes (n+1)}|_{\cW_1});(q_2')^!\ZQ_{[W_2,0]}[-2N])\\
\sim \hom((\lambda^{\cyc}\boxtimes \ZQ_{(\Fr\times E)^{n+1}})|_{W_2}\otimes (l'_2)^{-1} (\tau^{-1}(r_{0!}n_0^{-1}\cP_2)^{\boxtimes (n+1)}|_{\cW_1});(q_2')^{-1}\ZQ_{[W_2,0]})
\end{multline*}
The statement now follows.
\subsection{The category $\bbH$}
Let $$\bbH:=\psh(\Sigma,\bS_0(\ZZ(\ccA))).$$
Let 
$$
F_{\bbH}^{[n]}:=\bigoplus\limits_{ab}\ZQ_{[K^{\Sigma,[n],a,b},a]}\otimes (a,b)\otimes\unit_{\DMetR_q} [2b]
$$
$$
F_{\bbH}^{(n)}:=\bigoplus\limits_{ab}\ZQ_{[K^{\Sigma,(n),a,b},a]}\otimes (a,b)\otimes \unit_{\DMetR_q} [2b]
$$
The quasi-contractibility is easy to check

\subsection{Connecting $\bbH$ and $\bbG$}
\subsubsection{The objects $\cP_{ki}$}  Let $k,i\in \MetR$, $k\leq i$.   Let $k=(g_k;r_k)$; $i=(g_i,r_i)$.
We have a family of symplectic embeddings 
\begin{equation}\label{fam1}
\Fr_k\times_M \Fr_i\times B_{r_k}\to B_{r_i}.
\end{equation}
The differential at $0\in B_{r_k}$ gives rise to a map
$\Fr_k\times_M \Fr_i\to \Sp(2N)$.   Let $\Fr_{ki}\subset \Fr_k\times_M \Fr_i$ be the pre-image 
of the set of all symmetric matrices in $\Sp(2N)$.   The subset $\Fr_{ki}$ projects diffeomorphically onto
both $\Fr_k$ and $\Fr_i$.   This way we get  identifications
 $$
b_{ki}:\Fr_k\stackrel\sim\to \Fr_i;\quad c_{ki}:\Fr_k\stackrel\sim\to \Fr_{ki}
$$
The restriction of (\ref{fam1}) onto $\Fr_{ki}\cong \Fr_k$ gives rise to a graded family
of symplectic embeddings
$$
\phi_{ki}:\Fr_k\times B_{r_k}\to B_{r_i}
$$
Whence  objects $\cP_{ki}\in \sh_q(\Fr_k\times E\times E)$.

Let $s:\Fr_k\times E\times E\to \Fr_i\times E\times E$ be given by
$s(f,e_1,e_2)=(b_{ki}(f),e_2,e_1)$.

Let
$\cQ_{ik}:=s_!\cP_{ki}\in \sh_q(\Fr_i\times E\times E)$

Denote by $D'_{ki}:\Fr_k\to \Sp(2N)$ the differential of $\phi_{ki}$ at $0\in B_{r_k}$.  $D'_{ki}$ takes values in
the space of symmetric matrices, therefore, we have a lifting $D_{ki}:\Fr_k\to \ovSp(2N)$.

Let $k_1\leq k_2\leq i$.  Let
$b^{i}_{k_1k_2}:\Fr_{k_1}\to \Fr_{k_2}$ be given by $$b^i_{k_1k_2}= b_{k_2i}^{-1} b_{k_1i}.$$
We have  a unique graded family of symplectic embeddings 
$$
\phi_{k_1k_2}^i:\Fr_{k_1}\times B_{r_{k_1}}\to B_{r_{k_2}}
$$
satisfying  $$
\phi_{k_2i}(b^i_{k_1k_2}(f),\phi_{k_1k_2}^i(f,b))=\phi_{k_1i}(f,b).
$$
The family $\phi_{k_1k_2}^i$ defines  objects 
$$
\cP^i_{k_1k_2}\in \sh_q(\Fr_{k_1}\times E\times E);\quad \cQ^i_{k_2k_1}\in \sh_q(\Fr_{k_2}\times E\times E).
$$

Let $$
p_E:\Fr_{k_1}\times E\times E\times E\to \Fr_{k_1}\times E\times E
$$
be the projection,  where $p_E(f,e_1,e_2,e_3)=(f,e_1,e_3)$;
let
$$
p:\Fr_{k_1}\times E\times E\times E\to \Fr_{k_1}\times E\times E,
$$
$p(f,e_1,e_2,e_3)=(f,e_1,e_2)$;

let
$$
q:\Fr_{k_1}\times E\times E\times E\to \Fr_{k_1}\times E\times E,
$$
$q(f,e_1,e_2,e_3)=(b_{k_1k_2}^i(f);e_2,e_3).$

For $U,V\in \sh_q(\Fr_{k_1}\times E\times E)$, denote 
$$
U\circ^{\Fr_{k_1}}_{E} V:=p_{E!} (p^{-1}U\otimes q^{-1}V).
$$

Let also $\beta_{k_1k_2}^i:\Fr_{k_1}\times E\times E\to \Fr_{k_2}\times E\times E$ be the map induced by $b_{k_1k_2}^i$.
We now have $$
\cP_{k_1i}\sim      (\beta_{k_1k_2}^i)^{-1} \cP_{k_2i}\circ^{\Fr_{k_1}}_E \cP_{k_1k_2}^i
$$

\subsubsection{The objects $\Gamma_{ij}^k$}  Let $i,j,k\in \MetR$,  $i,j\geq k$.  
Define an object 
$$\Gamma_{ij}^k\in \sh_{\pi r_k^2}(\Phi_{ij}\times E\times E).$$

Let $\vs:\Phi_{ij}\to \ovSp(2N)$ be the structure map.   Let $$p_i:\Phi_{ij}\to \Fr_i, \quad p_j:\Phi_{ij}\to \Fr_j;
$$
 $$
\pi_i:\Phi_{ij}\times E\times E\to \Fr_i\times E\times E; \quad \pi_j:\Phi_{ij}\times E\times E\to \Fr_j\times E\times E
$$
be the projections.

  Define a map $\tau:\Phi_{ij}\to \ovSp(2N)$
by the condition
$$
\sigma(f)D_{ki}(b_{ki}^{-1}p_i(f))=D_{kj}(b_{kj}^{-1}p_j(f))\tau(f).
$$

Let $$t:\Phi_{ij}\times E\times E\to \ovSp(2N)\times E\times E$$
be the map induced by $\tau$.

Set
$$
\Gamma_{ij}^k:=\pi_i^{-1}\cP_{ki}\circ^{\Phi_{ij}}_E \Gamma_{r_k} \circ^{\Phi_{ij}}_E  (\beta_{kj})_! \cQ_{jk}\in \sh_{\pi r_k^2}(\Phi_{ij}\times E\times E);
$$
$$
\gamma_{ij}^k:=\pi_i^{-1}\cP_{ki}\circ^{\Phi_{ij}}_E \gamma_{r_k} \circ^{\Phi_{ij}}_E  (\beta_{kj})_! \cQ_{jk}\in \sh_q(\Phi_{ij}\times E\times E),
$$
where $\beta_{kj}:\Fr_k\times E\times E\to \Fr_i\times E\times E$ is induced by $b_{kj}$.

Let $k_1\leq k_2\leq i,j$ and consider

\begin{multline*} G'_{ij}(k_2,k_1):=
\hom(\Gamma_{ij}^{k_2};\Gamma_{ij}^{k_1})\sim
\hom(\pi_i^{-1}\cP_{k_2i}\circ^{\Phi_{ij}}_E \gamma_{r_{k_2}} \circ^{\Phi_{ij}}_E  (\beta_{k_2j})_! \cQ_{jk_2};\pi_i^{-1}\cP_{k_1i}\circ^{\Phi_{ij}}_E \gamma_{r_{k_1}} \circ^{\Phi_{ij}}_E  (\beta_{k_1j})_! \cQ_{jk_1})\\
\sim\hom((\cP_{k_2i}\boxtimes \cP_{k_2j})\gamma_{r_{k_2}};(\cP_{k_2i}\cP_{k_1k_2}^i\boxtimes \cP_{k_2j}\cP_{k_1k_2}^j)\gamma_{r_{k_1}})\\
\sim \hom(\gamma_{r_{k_2}}; (\cP_{k_1k_2}^i\boxtimes \cP_{k_1k_2}^j)\gamma_{r_{k_1}})\\
\sim \hom((\cQ_{k_2k_1}^i\boxtimes \cQ_{k_2k_1}^i)\gamma_{r_{k_2}};\gamma_{r_{k_1}})\\
\sim \hom((P_{r_{k_1}}\boxtimes P_{r_{k_1}})(\cQ_{k_2k_1}^i\boxtimes \cQ_{k_2k_1}^i)\gamma_{r_{k_2}};\gamma_{r_{k_1}})\\
\sim\hom(\gamma_{r_{k_1}};\gamma_{r_{k_2}}).
\end{multline*}

It follows that $G'_{ij}(k_2,k_1)$ admits a truncation.   Let 
$$
\cG_{ij}(k_2,k_1):=\tau_{\leq 0} G'_{ij}(k_2,k_1).
$$

We have a homotopy equivalence $\pi:\cG_{ij}(k_2,k_1)\to \ZQ$.
 
We can now define  categories $\cG_{ij},\SMetR_{ij}$ whose objects consist of all $S\subset \SMetR$, where $$k\in S\Rightarrow k<<i,j.
$$
and $$
\hom_{\cG_{ij}}(S_2,S_1)=\cG_{ij}(\mu(S_2),\mu(S_1)); \quad\hom_{\SMetR_{ij}}(S_2,S_1)=\ZQ
$$ if $S_2\subset  S_1$  and $$\hom_{\cG_{ij}}(S_2,S_1)=\hom_{\SMetR_{ij}}(S_2,S_1)=0
$$ 
otherwise.

We have a functor $\pi:\cG_{ij}\to \SMetR_{ij}$ which is a weak equivalence of categories.

Let $\tau:(\SMetR_{ij})^\op\to \DMetR_0(\ccA)$ be  given by:

--- $\tau(S)^T=0$ if $S\neq T$ and $\tau(S)^S=\ZQ$.

  Let $\tau_{\cG_{ij}}:=\pi^{-1}\tau$,
$$
\tau_{\cG_{ij}}:(\cG_{ij})^\op\to \DMetR_0(\ccA).
$$

Let $S\in \SMetR$.  We then have a functor
$i_S:\Com(\ccA)\langle \val(S)\rangle\to \bR_0(\ccA)$, where
$i_S(X)^T=0$ if $T\neq S$ and $i_S(X)^S=X$ otherwise.  It induces a functor
$$
\iota_S:\psh_{\val(S)}(X,\ZZ(\ccA))\to \psh(X,\bR_0(\ZZ(\ccA)))
$$
in the obvious way.
We now have a functor $$\uGamma_{ij}:\cG_{ij}\to \psh({\Phi_{ij}\times E\times E}, \bR_0(\ZZ(\ccA))),
$$ where
$$
\uGamma_{ij}(S)=\iota_S(\Gamma_{ij}^{\mu(S)})
$$

Set
$$
\Gamma_{ij}:=\uGamma_{ij}\otimes_{\cG_{ij}}^L \tau_{\cG_{ij}}\in \psh(\Phi_{ij}\times E\times E, \DMetR_0(\ZZ(\ccA))),
$$
where we have used the tensor product $\bR_0(\ZZ(\ccA))\otimes \bS_0(\ZZ(\ccA))\to  \bS_0(\ZZ(\ccA))$.
\subsection{The linking object $\bbGamma$} Let us now construct an object
$$\bbGamma\in \psh(\Sigma\times \Phi\times E\times E,\bS_0(\ZZ(\ccA)))$$
by means of prescription its components
$$
\bbGamma_{ij}\in\psh(\Sigma\times \Phi_{ij}\times E\times E, \bS_0(\ZZ(\ccA))).
$$
As was explained above, we have a map
$$
a_{ij}:\Phi_{ij}\to \Sigma.
$$
Let $A_{ij}\subset \Phi_{ij}\times \Sigma$ be the graph of $a_{ij}$.  Let $p_A:A_{ij}\times E\times E\to \Phi_{ij}\times E\times E$; $i_A:A_{ij}\times E\times E\to \Phi_{ij}\times E\times E$ be the corresponding maps.
Set
$$
\bbGamma_{ij}:=i_{A!}p_A^{-1}\Gamma_{ij}[-N].
$$
\subsubsection{Quasi-contractibility}
Let us first consider
$$
H^{[n]}=\hom(F_\bbH^{[n]}\circ (\bbGamma^{\boxtimes n});F_\bbG^{[n]}\circ \bbGamma)
$$

Consider the following diagram
$$
\xymatrix{
(\Sigma\times \Phi\times E\times E)^n&&\\
K^{\Sigma,[n],ab}\times  \Phi^n\times E^{n+1}\ar[u]^i\ar[r]^{\tilde{p}}& \Sigma\times (\Phi\times E\times E)^{n}&\\
K^{\Sigma,[n],ab}\times  K^{[n],ab}\times E^{n+1}\ar[u]^j\ar[r]^p&\Sigma\times K^{[n],ab}\times E^{n+1}\ar[u]^{\tilde{j}}\ar[r]^q& \Sigma\times \Phi\times E\times E}
$$
We have
\begin{multline*}
H^{[n]}=\hom(\tilde{p}_!i^{-1}(\bbGamma^{\boxtimes n}); \tilde{j}_!q^{-1}\bbGamma)
\sim \hom(\tilde{j}^{-1}\tilde{p}_!i^{-1}(\bbGamma^{\boxtimes n});q^{-1}\bbGamma)
\sim \hom(p_!j^{-1}i^{-1}(\bbGamma^{\boxtimes n});q^{-1}\bbGamma)
\end{multline*}

Let $$a:=\bigsqcup\limits_{ij} a_{ij}:\Phi\to \Sigma.$$  The map $a$ induces maps
$$
a_n:K^{[n],ab}\to K^{\Sigma,[n],ab}.
$$
We now have the following diagram:
$$
\xymatrix{ &(\Phi\times E\times E)^n\ar[ld]^{\alpha}&\\
(\Sigma\times \Phi\times E\times E)^n&\Phi^n\times E^{n+1}\ar[ld]^{\tilde{\alpha}}\ar[u]^{\tilde{i}}&\\
K^{\Sigma,[n],ab}\times  \Phi^n\times E^{n+1}\ar[u]^i&  K^{[n],ab}\times E^{n+1}\ar[ld]^b\ar[d]^c\ar[u]^{\tilde{j}}\ar[r]^{\tilde{q}}&\Phi\times E\times E\ar[d]^{\tilde{c}}\\
K^{\Sigma,[n],ab}\times  K^{[n],ab}\times E^{n+1}\ar[u]^j\ar[r]^p&\Sigma\times K^{[n],ab}\times E^{n+1}\ar[r]^q& \Sigma\times \Phi\times E\times E.}
$$

Here:

---\begin{multline*}
\alpha((\phi_1,e_1,e_2),(\phi_2,e_3,e_4),\ldots,(\phi_n,e_{2n-1},e_{2n}))\\
= ((a(\phi_1),\phi_1,e_1,e_2),(a(\phi_2),\phi_2,e_3,e_4),\ldots,(a(\phi_n),\phi_n,e_{2n-1},e_{2n}))
\end{multline*}

--- $$\tilde{\alpha}=a_n\times\Id_{\Phi^n\times E^{n+1}}.
$$

---Let $\pi_{k,k+1}:K^{\Sigma,[n],ab}\to \Phi$ be the projection onto the corresponding factor, $0\leq k\leq n$. Then
\begin{multline*}
i(\kappa,\phi_1,\phi_2,\ldots,\phi_n,e_0,e_1,e_2,\ldots,e_{n})\\
=((\pi_{01}\kappa,\phi_1,e_0,e_1),(\pi_{12}\kappa,\phi_2,e_1,e_2),\ldots,(\pi_{n-1,n}\kappa,\phi_n,e_{n-1,n})).
\end{multline*}

---\begin{multline*}
\tilde{i}((\phi_1,\phi_2,\ldots,\phi_n),(e_0,e_1,\ldots,e_n))\\
=((\phi_1,e_0,e_1),(\phi_2,e_1,e_2),\ldots,(\phi_n,e_{n-1},e_n)).
\end{multline*}

--- Let $$
\pi:=\pi_1\times\pi_2\times\cdots\times \pi_n:K^{[n],ab}\to  \Phi^n.
$$
Then
$$
j(\kappa,\kappa',e)=(\kappa,\pi(\kappa'),e);
$$
$$\tilde{j}:=\pi\times \Id_{E^{n+1}}.
$$

--- $b=a_n\times \Id_{K^{[n],ab}\times E^{n+1}}$;

--- Let $\pi_{0n}:K^{\Sigma,[n],ab}\to \Sigma$ be the projection onto the first and the last factors.
Then  $$p(\kappa_\vs,\kappa,e)=(\pi_{0n}(\kappa_\vs),\kappa,e).
$$

--- Set $c=pb$;

--- Set $\tilde{c}(\phi,e_1,e_2)=(a(\phi),\phi,e_1,e_2)$;

--- Set $q(\sigma,\kappa,e_0,e_1,\ldots,e_n)=(\sigma,\pi_{0n}\kappa,e_0,e_n)$;

--- Set $\tilde{q}(\kappa,e_0,e_1,\ldots,e_n)=(\pi_{0n}\kappa,e_0,e_n)$.
All the squares in this diagram are Cartesian.
We have $$
(\bbGamma)^{\boxtimes n}\sim(\alpha)_!(\Gamma)^{\boxtimes n}; \quad  \bbGamma \sim  (\tilde{c})_!\Gamma,$$
therefore,
$$
H^{[n]}\sim \hom(p_!j^{-1}i^{-1}(\alpha)_!(\Gamma)^{\boxtimes n};q^{-1}{\tilde c}_!\Gamma)\sim
\hom(p_!b_!(\tilde{i}\tilde{j})^{-1}(\Gamma)^{\boxtimes n};{ c}_!(\tilde{ q})^{-1}\Gamma)
$$

Observe that $c=pb$ and that $c$ is a closed embedding,  therefore, we can continue:
$$
H^{[n]}\sim \hom((\tilde{i}\tilde{j})^{-1}(\Gamma[N+1])^{\boxtimes n};(\tilde{q})^{-1}\Gamma[N+1]).
$$

We have a  decomposition $$
K^{[n],ab}=\bigsqcup\limits_{i_0i_1\ldots i_n\in \MetR} K^{ab}_{i_0i_1\ldots i_n}.
$$
Accordingly $H^{[n]}$ splits into a direct  product of its components
$$
H^{[n]}=\prod_{i_0i_1\ldots i_n} H_{i_0i_1\dots i_n}
$$
The problem reduces to considering each such a component.

Denote by $$\iota:K^{ab}_{i_0i_1\ldots i_n}\times E^{n+1}\to (\Phi\times E\times E)^n$$
 the  restriction of  $\tilde{i}\tilde{j}$.  Let us  denote by $q'$ the  restriction of $\tilde{q}$ onto
$K^{ab}_{i_0\ldots i_n}$. Let  us decompose $q'=q_\Phi q_E$, where 
$$q_E:K^{ab}_{i_0\ldots i_n}\times E^{n+1}\to K^{ab}_{i_0\ldots i_n}\times E^2;$$
$$
q_\Phi:K^{ab}_{i_0\ldots i_n}\times E^2\to \Phi_{i_0i_n}\times E^2.
$$

Let us consider $$
\boxtimes_{m=0}^{n-1} \Gamma_{i_mi_{m+1}}
\in \psh\left(\prod\limits_{m=0}^{n-1}\Phi_{i_mi_{m+1}}\times E\times E,\bS_0(\ZZ(\ccA))\right)
$$

We have
$$
\boxtimes_{m=0}^{n-1} \Gamma_{i_mi_{m+1}}
\sim
(\boxtimes_{m=0}^{n-1} \uGamma_{i_{m}i_{m+1}})\otimes^L (\boxtimes_{m=0}^{n-1} \tau_{\cG_{i_{m}i_{m+1}}}),
$$
where $\otimes^L$ is taken over the category
$$
\bigotimes\limits_{m=0}^{n-1} \cG_{i_mi_{m+1}}.
$$

We then have a homotopy equivalence
\begin{multline}\label{hioin}
H_{i_0\ldots i_n}\sim
\hom\left(
(
  \iota^{-1} (\boxtimes_{m=0}^{n-1}       \uGamma_{i_{m}i_{m+1}}[-N]
                   )
\otimes^L 
                  (\boxtimes_{m=0}^{n-1}    \tau_{\cG_{i_{m}i_{m+1}}}
                  )
              );
    \tilde{q}^{-1}(\uGamma_{i_0i_n}[-N]\otimes^L    \tau_{\cG_{i_0i_n}}
              )
        \right)\\
\sim
\Rhom\left(\iota^{-1}\boxtimes_{m=0}^{n-1}       \uGamma_{i_{m}i_{m+1}}[-N];
  \hom_{\DMetR_0} (\boxtimes_{m=0}^{n-1}    \tau_{\cG_{i_{m}i_{m+1}}}; \tilde{q}^{-1}(\uGamma_{i_0i_n}[-N]\otimes^L    \tau_{\cG_{i_0i_n}})\right).
\end{multline}

For $S\in \SMetR$,  define an object $f^S_0\in \bS_0$, $(f^S_0)^T=0$ if $T\neq S$, and $(f^S_0)^S=\ZQ$.

Observe that for every $S\in \SMetR$, the natural map
$$
\tilde{q}^{-1}\left(\uGamma_{i_0i_n}[-N]\otimes^L _{\cG_{i_0i_n}}   \hom_{\DMetR_0}(f^S_0; \tau_{\cG_{i_0i_n}})\right)
\to \hom_{\DMetR_0}\left( f^S_0; \tilde{q}^{-1}(\uGamma_{i_0i_n}[-N]\otimes^L_{\cG_{i_0i_n}}    \tau_{\cG_{i_0i_n}})\right)
$$
is a homotopy equivalence (because there are only finitely many $T\in \cG_{i_0i_n}$ satisfying $T\subset S$.
We can now continue (\ref{hioin}):
$$
H(i_0i_1\cdots i_n)\sim
\RHom\left(\iota^{-1}\boxtimes_{m=0}^{n-1}       \uGamma_{i_{m}i_{m+1}}[-N];\tilde{q}^{-1}(\uGamma_{i_0i_n}[-N]\otimes^L_{\cG_{i_0i_n}} \lambda)\right),
$$
where 
$$\lambda:\bigotimes\limits_{m=0}^{n-1} \cG(i_mi_{m+1})\otimes \cG(i_0i_n)^{\op}\to \ccA;
$$
$\lambda(S_1,S_2,\ldots,S_n;S)=\ZQ$
if $S_1\cup S_2\cup\cdots \cup S_n\subset S$ and
$\lambda(S_1,S_2,\ldots,S_n;S)=0$ otherwise.

We have
$$
 (\tilde{q}^{-1}\Gamma_{i_0i_n}[-N]\otimes^L_{\cG_{i_0i_n}} \lambda)(S_1,S_2,\ldots,S_n)\sim \tilde{q}^{-1}\Gamma_{i_0i_n}^{S_1\cup S_2\cup \cdots\cup S_n}[-N]
$$
if $S_1\cup S_2\cup\dots\cup S_n$ is linearly ordered
and
$$
 (\tilde{q}^{-1}\Gamma_{i_0i_n}[-N]\otimes^L_{\cG_{i_0i_n}} \lambda)(S_1,S_2,\ldots,S_n)\sim 0 
$$
otherwise.

The problem now reduces to showing  quasi-contractibility of
$$
\hom(\iota^{-1}(\boxtimes_m \Gamma_{i_{m-1i_m}}^{k_m}[-N]);
\tilde{q}^{-1}(\Gamma_{i_0i_n}^k[-N])).
$$
under assumptions that $k_m\in G(i_{m-1},i_m)$ and $k\leq k_m$ for all $m$.

Let $$
q_E:K^{ab}_{i_0i_1\ldots i_n}\times E^{n+1}\to K^{ab}_{i_0i_1\ldots i_n}\times E^{2}
$$
be the projection onto the first and the last factor of $E^{n+1}$.
The above space can be rewritten as 
\begin{multline*}\hom(q_{E!}\iota^{-1}(\boxtimes_m \Gamma_{i_{m-1}i_m}^{k_m}[-N])[(n-1)N];\tilde{q}^{-1}_\Phi\Gamma_{i_0i_n}^k[-N])\\
=\hom(q_{E!}\iota^{-1}(\boxtimes_m \Gamma_{i_{m-1}i_m}^{k_m});\tilde{q}^{-1}
\Gamma_{i_0i_n}^k)
\end{multline*}

Step 1.    Let $$\red:\sh_q(\Phi_{ij}\times E\times E)\to \sh_{\val k}(\Phi_{ij}\times E\times E)
$$ be the projection.
The object $\Gamma_{ij}^k\in \sh_{\val k}(\Phi_{ij}\times E\times E)$ is glued of objects  $(\gamma_{ij}^k)_l$of the form

$$
(\gamma_{ij}^{k})_l=\red T_{-l\pi r_k^2}\gamma_{ij}^k[2l].
$$  

The object $\gamma_{ij}^k$ is glued of $\bP_{ij}^k$ and $T_{-\pi r_k^2}\bP_{ij}^k[-2N-1]$, 
The  object $(\gamma_{ij}^k)_{l}$  is therefore glued of the following objects
$$
X^1_{ijk|l}:=T_{-(l+1)\val k}\bP_{ij}^k[2l+2N+1]
$$
and
$$
X^0_{ijk|l}:=T_{-l\val k}\bP_{ij}^k[2l].
$$

Step 2. Let $L=q_{E!}\iota^{-1}(\boxtimes_m  \bP_{i_{m-1}i_m}^{k_m})$ and  $M=q_{E!}\iota^{-1}(\boxtimes_m \bP_{i_{m-1}i_m}^k)$.

We then have a natural map $M\to L$ because $k\leq  k_m$ for all $m$.
Singular support consideration shows that for any $a\in \Re$,
$$ 
\hom(T_aL;\tilde{q}^{-1}\gamma_{i_0i_n}^k[N])\to R\hom(T_aM;\tilde{q}^{-1}\gamma_{i_0i_n}^k[N])
$$
is a homotopy equivalence.

Step 3 We have a homotopy equivalence $M\sim \tilde{q}^{-1}\bP_{i_0i_n}^k$.

Step 4   $\hom(\bP_{ij}^k;T_aT_{-\val(k)}\gamma_{ij}^k)=\ZQ$ if $0\leq a<\val k$ and 0 otherwise (use the representation of $\gamma_{ij}^k$
as the cone of  $\bP_{ij}^k[2N+1]\stackrel{+1}\to T_{\pi r_k^2}\bP_{ij}^k$.

Step 5 $\tilde{q}_{!}\tilde{q}^{-1}F[\nu]=F\otimes \cL$, where $\cL\in \sh(\Phi_{i_0i_n}\times E\times E)$
is a local system obtained by gluing   constant sheaves of the form $\ZQ_{\Phi_{i_0i_n}\times E\times E}[a]$, $a\geq 0$. Here $\nu$ is the relative dimension of $p$.

Step 6 Count.  
It suffices to show that all objects
$$
\hom(q_{E!}\iota^{-1}\boxtimes_m (X^{a_m}(i_{m-1}i_{m}k_m|l_m));\tilde{q}^{-1}T_{-\val(k)}\gamma_{i_0i_n}^k)\in \GZtrunk.
$$
We have:
\begin{multline*}
q_{E!}\iota^{-1}\boxtimes_m (X^{a_m}(i_{m-1}i_{m}k_m|l_m))\\
=q_{E!}\iota^{-1}(\boxtimes_m T_{-(l_m+a_m)\val(k_m)}\bP_{i_{m-1}i_m}^{k_m}[2l_m+(2N+1)a_m])\\
=T_{u'}q_{E!}\iota^{-1}(\boxtimes_m \bP_{i_{m-1}i_m}^{k_m})[v']
\end{multline*}
where
$$
u'=\sum -(l_m+a_m)\val(k_m)
$$
and
$$v'=2\sum_m l_m+(2N+1)\sum_m a_m.
$$

So that we have
\begin{multline*}
\hom(q_{E!}\iota^{-1}\boxtimes_m (X^{a_m}(i_{m-1}i_{m}k_m|l_m));\tilde{q}^{-1}T_{-\pi r_k^2}\gamma_{i_0i_n}^k)\\
=\hom(T_{u'}\tilde{q}^{-1}\bP_{i_0i_n}^k[v'];\tilde{q}^{-1}T_{-\pi r_k^2}\gamma_{i_0i_n}^k)\\
=\hom(\tilde{q}^{-1}\bP_{i_0i_n}^k;T_{-\pi r_k^2-u'}\tilde{q}^{-1}\gamma_{i_0i_n}^k[2N+1-v'])\\
=\hom(\bP_{i_0i_n}^k\otimes \cL;T_{-\pi r_k^2-u'}\gamma_{i_0i_n}^k[2N+1-v']).
\end{multline*}

As $\cL$ is glued of $\ZQ_{\Phi_{i_0i_n}\times E\times E}[a]$, $a\geq 0$,
it suffices to show
that
$$
\hom(\bP_{i_0i_n}^k;T_{-\pi r_k^2-u'}\gamma_{i_0i_n}^k[2N+1-v'])\in D_{\geq 0} \ccA.
$$

 This is true in each of the following cases

1)  $-u'\notin [0,\pi r_k^2)$;

2) $-u'\in [0,\pi r_k^2)$ and $-v'\leq 0$.

However,  $v'\geq 0$ which proves the statement.

\subsection{The category $\bbM$}   Let us first define a category $\bbM_0$ enriched over $\ccA$.

Set $\bbM_0:=\open_M$.   
Let  $h:\bbM_0\otimes \bbM_0\to \ccA$,
$h(U,V)=\ZQ$ if $U\cap V\neq \emptyset$;  $h(U,V)=0$ otherwise.

  Let $\Delta_M^{[n]}\subset M^n\times M$; $\Delta^{(n)}\subset M^{n+1}$ be diagonals.
Set
$$
F^{[n]}_{\bbM_0}:=\ZQ_{\Delta_M^{[n]}}
\in  \sh(M^n\times M);
$$
$$
F^{(n)}_{\bbM_0}=\ZQ_{\Delta_M^{(n)}}  \in  \sh(M^{n+1}).
$$

Let $\forall_0\in \Classic(\ZZ(\ccA))\langle 1\rangle$ be defined  by
$$
\gr^a \forall_0:=\bigoplus\limits_b  \langle a,b\rangle[2b].
$$
Set 
$$
\forall:=i_\emptyset(\forall_0)\in \bR_0(\ZZ(\ccA)).
$$
We have an algebra structure on $\forall$.   Denote by $\bbB$ the category, enriched over $\bR_0(\ZZ(\ccA))$, with one object, whose
endomorphism space is $\forall$.    Denote this  object  by $e$. As $\forall$ is a commutative alegbra,
we have a symmetric monoidal structure on $\bbB$, where $e\otimes e=e$.

Let now  $\bbM:=\bbM_0\otimes \bbB$. The functors $F^{[n]},F^{(n)}$ extend to $\bbM$.

\subsubsection{Connecting $\bbM$ and $\bbH$}  Let $p_\Sigma:\Sigma\to M$ be the projection.
Let $G_\Sigma\subset \Sigma\times M$ be the graph of $p_\Sigma$.
Let $$
\bDelta:=\ZQ_{G_\Sigma}\otimes \unit_{\bS_0(\ZZ(\ccA))}\in \Doplus(\bbM\otimes \bbH \otimes \bS_0(\ZZ(\ccA))).
$$

\subsection{The categories $\bbB_\cT$} \label{triang} Fix a triangulation $\cT$ of $M$.  Let $\open_\cT\subset \open_M$
consist of $M$ and all stars of $\cT$.
Let $\bbB_\cT:=\open_\cT$ viewed as a discrete set.  Let $h_\bbB:\bbB\times \bbB\to \ccA$ be given by $h_\bbB(U,V)=\ZQ$ if $U=V$ and $h_\bbB(U,V)=0$ otherwise.

Let $$
F^{[n]}:=\bigoplus\limits_{U_1\cap U_2\cap \cdots\cap U_n\subset U} U_1\otimes U_2\otimes \cdots \otimes U_n\otimes U\otimes \forall\otimes \unit_{\bS_0(\ccA)}$$

Let $$F^{(n)}_{\bbB_\cT}:=\bigoplus\limits_{U_0\cap U_1\cap \cdots\cap U_n\neq \emptyset} U_0\otimes U_1\otimes \cdots\otimes U_n\otimes \unit_{\bS_0(\ccA)}$$

Let $\bbB\subset \bbB_\cT$ be the full sub-category consisting of a single object $M$.

\subsubsection{Connecting $\bbB_\cT$ and $\bbM$} 
Let $I:\SMetR^\op\to \DMetR_0$, 
$I(S)=f^S_0$.   Let $\cI\in \swell \DMetR_0$,
$$\cI:=\hocolim_{S\in \SMetR} I(S).
$$

Set $\bbE\in \psh(M,\bS_0(\ZZ(\ccA)))$,
$$\bbE|_{U\times \bbM}:=\ZQ_U\langle 0,0\rangle\otimes \cI.
$$
\newcommand{\bM}{\mathbf{M}}
Observe that the restriction of $\bbB_\cT$ and $\bbE$ onto $M\in \open_\cT$ is independent of the choice of $\cT$.
\subsection{Straightening out}
\subsubsection{Changing the ground category}  We have the constant  tensor functor $\ZZ\to \pt$ which
induces a tensor functor  $\bR_0(\ZZ( \ccA))\to \bR_0(\ccA)$ by means of which  the structure on
the categories $\bbB,\bbM,\ldots$ carries over to $\bR_0(\ccA)$ .  
All the categories below are  thus enriched over
$\bR_0$. The category $\bbF$ is enriched over $\bR_q$.

\subsubsection{Straightening out}\label{strt} Applying the procedures from Sec \ref{straight} we produce
a sequence of monoidal categories  and their functors
\begin{equation}\label{strght}
\bU(\bbB)\to \bU(\bbB_\cU)\to\swell \bU(\bbM)\to \swell\bU(\bbH)\to \swell\bU(\bbG)\to \swell\bU(\bbF).
\end{equation}

The category $\bU(\bbF)$ is enriched over $\bR_q$.   We also have a category $\bU(\bbF_R)$, where 
$\swell \bU(\bbF_R)$  enriched over $\Dprodoplus\bU(\bbF)$.

   We have a zig-zag map
from $\bU({\bbM_0})\otimes \bU(\bbA)$  to $ \swell\bU(\bbM)$.
We have  sub-categories $\bU(*)_1$, where $*=\bbB,\bbM,\bbH,\bbG,\bbF$. These sub-categories are 
preserved by the above functors.

\section{Hochschild complexes} \label{Hochsch}

\subsection{Hochschild complexes of an algebra  in a monoidal category with trace}    Let $\cM$ be a monoidal category with trace enriched over  a ground category $\ccA$.
Let  $A$ be an algebra in $\cM$.   We then can build  the Hochschild cochain complex $\Hoch^\bullet(A,A)$ in $\ccA$, where
 $\Hoch^n(A)=\hom_\cM(A^{\otimes n};A)$ as well as Hochschild chain complex $\Hochcyc^{-\bullet}(A)$,
 where  $\Hochcyc^{-n}(A)=\TR_\cM(A^{\otimes {n+1}})$, in the usual way.

Let $N$ be an $A$-bimodule in $\cM$. One  then defines Hochschild complexes $\Hoch^\bullet(M,A)$ and
$\Hochcyc^\bullet(N)$, where $$
\Hoch^{n}(M,A)=\bigoplus\limits_{i+j=n,i,j\geq 0} \hom_\cM(A^{\otimes i}\otimes N\otimes A^{\otimes j};A);
$$
$$
\Hochcyc^{-n}(M) =M\otimes A^{\otimes n}.
$$

Let  $\unA$ be $A$ viewed as a bimodule over itself.
we have  a natural map
\begin{equation}\label{hoh}
\Hoch(A)\to \Hoch(\unA,A).
\end{equation}

We will now formulate a sufficient condition for this map to be a homotopy equivalence.
Let $\unit\in\cM$ be the unit.  A map $\unit\to A$ in $\cM$ is called a homotopy unit of $A$
 if the induced maps $A=A\otimes \unit \to A\otimes A\to A$ and $A=\unit\otimes A\to A\otimes A\to A$
are homotopy equivalent to $\Id_A$.   
\begin{Proposition} Suppose $A$ admits a homotopy unit. Then the map (\ref{hoh}) is a homotopy equivalence.
\end{Proposition}

\subsubsection{A map $I:\Hoch(A)\to \Hochcyc(A)$, where $A$ is an algebra with trace} \label{mapI}
Let $A$ be an algebra with trace $\Tr$ in a monoidal category with trace $\cM$,
where $\Tr:\unit\to \TR(A)$.

Let us  define  maps
$$
I^n: \Hoch^n(A)\to \Hochcyc^n(A).
$$

To this end we first define maps $d_0:\Hoch^n(A)\to \Hoch^{n+1}(A)$
as follows:
$$
d_0:\Hoch^n(A)=\hom(A^{\otimes n};A)\to \hom(A\otimes A^{\otimes n};A\otimes A)\stackrel m\to  \hom(A^{n+1};A).
$$
We now set
\begin{multline*}
I:\Hoch^n(A):=\hom(A^{\otimes n};A)\otimes \unit\stackrel{d_0\otimes \Tr}\longrightarrow \hom(A^{\otimes {n+1}};A)\otimes \TR(A)\\
\to \TR(A^{\otimes {n+1}})=\Hochcyc^n(A).
\end{multline*}

It follows that $I$ gives rise to a map of complexes.
\subsection{Hochschild complexes in $\bM_\bbF$ and $\bM_\bbH$}
 We have an algebra $M\in \bM_\bbA$ and its bimodules $M_U$ for
every  element $U\in \open_\cT$.    We have a functor  $U\mapsto M_U$ from the cateory $\open_\cT$ to the category
of $M$-bimodules in $\bM_\bbB$.

 Denote by $M^\bbH,M^\bbH_U,M^\bbF,M^\bbF_U$ their images in $\bM_\bbH$ and $\bM_\bbF$.

We have  natural maps
\begin{equation}\label{utrian}
\hocolim_{U\in \cT} M^\bbH_U\to M^\bbH;\quad  \hocolim_{U\in \cT} M^\bbF_U\to M^\bbF.
\end{equation}
\begin{Lemma} \label{lutrian} The map (\ref{utrian}) is a homotopy equivalence.
\end{Lemma}

We have induced maps
\begin{equation}\label{cepm1}
\Hoch(M^\bbH)\to \Hoch(M^\bbF);\quad \Hochcyc(M^\bbH)\to \Hoch(M^\bbF);
\end{equation}
\begin{equation}\label{cepmu1}
\Hoch(M^\bbH_U,M^\bbH)\to \Hoch(M^\bbF_U,M^\bbF);\quad \Hochcyc(M^\bbH_U,M^\bbH)\to \Hochcyc(M^\bbF_U,M^\bbF)
\end{equation}
\begin{Proposition} \label{hchsch}

 The maps in (\ref{cepm1}),(\ref{cepmu1})  are homotopy equivalences.
\end{Proposition}

The rest of the section is devoted to proving the Proposition.

\subsubsection{Passage to the category $\Classic\langle \val S\rangle$} \label{refinecov} By vitrue of Lemma \ref{limitbq}, 
it suffices to show that the induced maps

\begin{equation}\label{cepm}
\tau_{\subset S}\Hoch(M^\bbH)\to \tau_{\subset S}\Hoch(M^\bbF);\quad \tau_{\subset S}\Hochcyc(M^\bbH)\to \tau_{\subset S}\Hoch(M^\bbF);
\end{equation}
\begin{equation}\label{cepmu}
\tau_{\subset S}\Hoch(M^\bbH_U,M^\bbH)\to \tau_{\subset S}\Hoch(M^\bbF_U,M^\bbF);\quad \tau_{\subset S}\Hochcyc(M^\bbH_U,M^\bbH)\to \tau_{\subset S}\Hochcyc(M^\bbF_U,M^\bbF)
\end{equation}

are homotopy equivalences in $\swell Q_{\val S}$ for every
$S\in \SMetR$.
From now on we fix an element $S\in \SMetR$.

Observe that it suffices to prove  $(\ref{cepmu})$ for  any  refinement $\cU_S$ of the cover $\cU$, where $\cU_S$ may depend on $S$.   We therefore make the following assumption without loss of generality: 

\begin{Assumption}\label{asb} Let
$\vs$ be the minimal element of $S$, then every element $U$ from $\cU_S$ is contained in a ball
$B_\vs(m_0)$ for some point $m_0\in M$.
\end{Assumption}
\subsubsection{Changing the set $\MetR'$}.   Let
us replace $\MetR'$ with its subset $S$.  It is clear that 
the induced map
$$
\tau_{\subset S}\Hoch^{\MetR'}(M^\bbH_U,M^\bbH)\to
\tau_{\subset S}\Hoch^{S}(M^\bbH_U,M^\bbH)
$$
is an isomorphism and similar for all the remaining ingredients of (\ref{cepm}), (\ref{cepmu}).  We therefore
will assume $\MetR'=S$ from now on.
\subsubsection{Passage to homotopy categories}   Let $K$ be a category enriched over the category of complexes
of $\ZQ$-modules.  Denote by $\ho K$ the category, enriched over the category of   $\ZQ$-modules.
where we set $\ho K(X,Y):=H^0\hom(X,Y)$.    If $\cM$ is a monoidal category  with a trace enriched over  a ground  SMC $\ccA$,
then $\ho \cM$ is  a monoidal category with a trace enriched over a SMC $\ho \ccA$.

\subsubsection{When  a homotopy equivalence in $\ho \ccA$ implies that in $\ccA$?}
Let $X^n\in \ccA$ and let $d_{nm}\in \hom^0(X^n,X^m),\quad n<m$ be elements satisfying
$$dd_{nm}+\sum_{k|n<k<m} d_{km}d_{nk}=0.$$ Call such data {\em a complex in $\ccA$}

Let us define an object $|X^\bullet|\in \bbB$ as the representing object
of 
$$
(\bigoplus\limits_{n<0} X^n[-n]\oplus \prod\limits_{n\geq 0} X^n[-n],\sum_{nm} d_{nm}).
$$
 
{\em A map of complexes } $f:X^\bullet\to Y^\bullet$ is  a collection of maps
$f_{nm}:X_n\to X_m$, $m\geq n$ satisfying:  let  $f:=\sum_{nm} f_{nm}\in \hom(|X^\bullet|,|Y^\bullet|)$,
then $df=0$.

One builds a complex $\Cone f$ in $\ccA$ in a standard way.

Let us now pass to $\ho \ccA$.  Let $X^n,d_{nm}$ be as above.  We have
  $$d d_{n,n+1}=0;\quad d_{n,n+1}d_{n+1}d_{n+2}+dd_{n,n+2}=0,
$$
which implies that
we have the following complex in $\ho \ccA$, to be denoted by $\ho_{X^\bullet}:$
$$
0\to X^0\stackrel{d_{01}}\to X^1\stackrel{d_{12}}\to \cdots.
$$

It also follows that a map of complexes $f:X^\bullet\to Y^\bullet$ induces a map  
$$
\ho_f: \ho_{X^\bullet}\to \ho_{Y^\bullet}
$$
of complexes in $\ho_\bbB$.  One define a complex $\Cone \ho_f$ in $\ho\bbB$ in a standard way.

Call $\ho_{X^\bullet}$ {\em acyclic}   if there exist
elements $\bo_i \in \hom_{\ho \bbB}(X^i;X^{i-1})$ such that 
$d_{i-1,i}\bo_i+\bo_{i+1}d_{i,i+1}=\Id_{X^i}$ in $\ho \bbB$.

Call a map $\ho_f$ {\em a homotopy equivalence} if $\Cone \ho_f$ is acyclic.

\begin{Proposition}  1) Suppose $\ho_{X^\bullet}$ is acyclic. Then so is $X^\bullet$.

2) Let $f:X^\bullet\to Y^\bullet$ be a map of complexes in $\bbB$ and suppose $\ho_f$ is a homotopy equivalence.
Then so is $f$.
\end{Proposition}
{\em Sketch of the proof}  It is clear that 1) implies 2).  Let us prove 1).
 Choose representatives of $\bo_i$ in $\hom^{-1}_\bbB(X^i,X^{-1})$,
to be denoted by $h_i$.  We have $dh_i=0$;   
$$
d_{i-1,i}h_i+h_{i+1}d_{i,i+1}=f_{ii}-du_i$$ 
for some $u_i\in \hom^{-1}_\bbB(X^i,X^i)$.

The sum of all $h_i$ and $d_i$ defines an element $H\in \hom^{-1}(\cX,\cX)$
such that
$dH=\Id+K$, where $K:X^i\to \prod\limits_{j\geq i+2} X^j$ so that   $\Id+K$ is an  automorphism of $\cX$.
This finishes the proof.

\subsubsection{Hochschild complexes in the homotopy category} Let $\cM,\cN$ are monoidal categories enriched over $\ccA$,
let $A$ be an algebra in $\cM$ admitting a homotopy unit and let $K$ be an $A$-bimodule in $\cM$.
Let $f:\cM\to \cN$ be a tensor functor.    We then have  maps of Hocshild complexes
\begin{equation}\label{hocha}
\Hoch(K,A)\to \Hoch(f(K),f(A));\quad \Hochcyc(K)\to \Hochcyc(f(K)).
\end{equation}

On the other hand we have an algebra $\ho A$ and its bi-module $\ho K$ in $\ho \bbB$
 and Hochschild complexes in  $\ho \bbB$:
\begin{equation}\label{hochho}
\Hoch(\ho K,\ho A)\to \Hoch(f(\ho K),f(\ho A));\quad \Hochcyc(\ho K)\to \Hochcyc(f(\ho K))
\end{equation}
It is clear that
$$
\ho_{\Hoch(K,A)}\cong \Hoch(\ho K,\ho A);\quad \ho_{\Hochcyc(K)}\cong \Hochcyc(\ho K).
$$
Therefore, 
\begin{Claim} If maps in (\ref{hochho}) are homotopy equivalences, then so are maps from (\ref{hocha}).
\end{Claim}

In particular our problem reduces to showing that the following maps are homotopy equivalences
\begin{equation}\label{cepu}
\Hoch(\ho M^\bbH_U,\ho M^\bbH)\to \Hoch(\ho M^\bbF_U,\ho M^\bbF);\quad \Hochcyc(\ho M^\bbH_U,\ho M^\bbH)\to \Hochcyc(\ho M^\bbF_U,\ho M^\bbF)
\end{equation}

\subsection{Describing $\ho (\bM_\bbH)_1$} 
We have an equivalence of categories $i:\ho \bbH\to \ho(\bM_\bbH)_1$. 
For each $n\geq 0$, consider the functor $T_n:((\ho \bbH)^\op)^n\times (\ho \bbH)\to \ho\bR_0$, where
$$T_n(X_1,X_2,\ldots,X_n;X):=\hom_{\ho \bbM_\bbH}(i(X_1)\otimes i(X_2)\otimes \cdots \otimes i(X_n);i(X)).
$$
We have an isomorphism of functors
\begin{multline*}
T_n(X_1,X_2,\ldots,X_n;X)\cong\hom_{\ho\psh_{\val S}(\Sigma^n)}(X_1\boxtimes X_2\boxtimes \cdots\boxtimes X_n;\bigoplus\limits_{ab}T_{ab!}T_a\ZQ_{(K^\Sigma_n)^{ab}}\circ_{\Sigma} X[2b])
\end{multline*}

Recall that we have a groupoid structure on $\Sigma\rightrightarrows \cL$.

Let $\mu_n:\Sigma^{\times_\cL^{n}}\to \Sigma$ be the composition map of this groupoid.

We have $(K^\Sigma_n)^{ab}\cong\Sigma^{\times_\cL^n}$ for all $a,b\in \bZ$.    The projection
$p_{1n}^{ab}:(K^\Sigma_n)^{ab}\to \Sigma$, under this identification reads as
$$
\cK_n\stackrel {\mu_n}\rightarrow \Sigma\stackrel{T_{ab}}\to \Sigma,
$$

Observe that $\mu_n$ is a smooth fibration.   Therefore, $\mu_{n}^{-1}=\mu_{n}^![-(n-1)D]$,
where $$
(n-1)D=
(n-1)(\dim \Sp(2N)-\dim \SU(N)),
$$
so that $D=\dim \Sp(2N)-\dim \SU(N)$.
Let also $i_n:\cK_n\to \Sigma^n$ be the embedding. 
so that we can rewrite
\begin{multline*}
T_n(X_1,X_2,\ldots,X_n;X)=\hom(X_1\boxtimes \cdots \boxtimes X_n;\bigoplus\limits_{ab} i_{n!}\mu_{n}^{-1} T_{ab!}T_aX[2b])\\
=\hom(i_n^{-1}(X_1\boxtimes \cdots \boxtimes X_n);\bigoplus\limits_{ab}\mu_{n}^{!}[-(n-1)D] T_{ab!}T_aX)[2b] )\\
\cong \hom(\mu_{n}^{-1}i_n^{-1}(X_1[D]\boxtimes \cdots \boxtimes X_n[D]);\bigoplus\limits_{ab}T_{ab!}T_aX[D][2b])).
\end{multline*}

The groupoid structure on $\Sigma\rightrightarrows \cL$ produces a convolution monoidal structure
on $\ho \sh_q(\Sigma)$, where the tensor product, to be denoted by $*$ is given by
$$X_1*X_2*\cdots *X_n=\mu_{n!}i_n^{-1}(X_1\boxtimes \cdots \boxtimes X_n),$$
so that we finally have
\begin{multline}\label{svertgroup}
T_n(X_1,X_2,\ldots,X_n;X)\cong \bigoplus\limits_{ab} T_{a}\hom(X_1[D]*\cdots*X_n[D];\bigoplus\limits_{ab}T_{ab!}T_aX[D][2b])\\=
\hom_{\ho (\bM_{\bbH})_1}(X_1[D]*\cdots * X_n[D];X[D]).
\end{multline}

This result has the following corollaries:

1) Let $\bbH'$ be the category enriched over $\ho\bR_0$, 
whose objects are the same as in $\ho \psh(\Sigma,\DMetR_0)$   and we 
set
$$
\hom_{\bbH'}(Y,X)=\hom_{\ho \sh_{\val S}(\Sigma)}(Y,\bigoplus\limits_{ab}T_{ab!}T_aX[2b])
$$

We then have a monoidal structure on $\bbH'$ given by $*$.  
Next, we have a  lax tensor functor  
$$
\Sigma_\bbH:(\bbH',*)\to \ho \bM_\bbH,
$$
where $\Sigma_\bbH(X)=i(X)[-D]$.  By lax tensor functor structure on $S$ we mean a natural transformation
$\Sigma_\bbH(X)\otimes \Sigma_\bbH(Y)\to \Sigma_\bbH(X*Y)$ which has an associativity property but need not be a homotopy equivalence.
The induced map
$$
\hom_{\bbH'}(X_1*X_2*\cdots* X_n;X)\to\hom_{\bM_\bbH}(\Sigma_\bbH(X_1)\otimes \Sigma_\bbH(X_2)\otimes \cdots \otimes \Sigma_\bbH(X_n);\Sigma_\bbH(X)).
$$
is an isomorphism for all $X_i,X\in \bbH'$.

\subsection{The categories $\ho(\bM_\bbG)_1$, $\ho(\bM_\bbF)_1$} One gets similar results for $\bM_\bbG,\bM_\bbF$.   
Let $$\bbG':=\ho \psh_{\val S}(\Phi\times E\times E);\quad \bbF':=\ho\psh_{\val S}((\Fr\times E)^2).
$$

Set $$\hom_{\bbG'}(A,B):=\hom_{\ho\psh_{\val S}(\Phi\times E\times E)}(A;\bigoplus_{ab} T_{ab!}T_aB[2b]).
$$
Set
$$\hom_{\bbF'}(A,B)=\hom_{\psh_{\val S}((\Fr\times E)^2)}(A,B).
$$

Let us denote by $*$ the convolution on $\bbG'$ coming from the groupoid
$$
\Phi\times E\times E\rightrightarrows\Fr\times E.
$$
We thereby get a monoidal structure on $\bbG'$.

Let us also denote by $*$ the convolution on $\bbF'$ coming from the groupoid
$$
(\Fr\times E)^2\rightrightarrows \Fr\times E.
$$
This way we get a monoidal structure on $\bbF'$.   

Let  $D_\bbG=\dim \Fr\times E-\dim M$;
$D_\bbF:=\dim \Fr\times E$.       Then we have  lax monoidal functors
$\Sigma_\bbG:\bbG'\to \bM_\bbG$ and
$\Sigma_\bbF:\bbF'\to \bM_\bbF$, where 
$\Sigma_\bbG(X)=i(X)[-D_\bbG]$
and
$\Sigma_\bbF(X)=i(X)[-D_\bbF]$.

The induced maps
$$
\hom_{\bbG'}(X_1*X_2*\cdots* X_n;X)\cong\hom_{\bM_\bbG}(\Sigma_\bbG(X_1)\otimes \Sigma_\bbG(X_2)\otimes \cdots \otimes \Sigma_\bbG(X_n);\Sigma_\bbG(X));
$$
$$
\hom_{\bbF'}(X_1*X_2*\cdots* X_n;X)\cong\hom_{\bM_\bbF}(\Sigma_\bbF(X_1)\otimes \Sigma_\bbF(X_2)\otimes \cdots \otimes \Sigma_\bbF(X_n);\Sigma_\bbF(X)).
$$
are  isomorphism for all $X_i,X\in \bbG'$ (resp.  all $X_i,X\in \bbF'$).

\subsection{Reduction to $\bbH'$, $\bbG'$, and $\bbF'$}

As the essential image of $\Sigma_\bbH$ is $\ho  (\bM_\bbH)_1$ and likewise for $\Sigma_\bbG,\Sigma_\bbF$, we have   induced 
lax tensor functors
$$
\bE_{\bbH'\bbG'}:\bbH'\to \bbG'\quad \bE_{\bbG'\bbF'}:\bbG'\to \bbF'.
$$
Denote $\bE_{\bbH'\bbF'}:=\bE_{\bbG'\bbF'}\bE_{\bbH'\bbG'}$.

The algebra $M^{\bbH}\in \bM_\bbH$ as well as its bimodule $M_U^\bbH$ determines an algebra in $\bbH'$ and its
bimodule,   denoted by $A^{\bbH'}$,  $M_U^{\bbH'}$.
Likewise we have  an algebra and its bimodule $A^{\bbF'},M_U^{\bbF'}$ in $\bbF'$ which are isomorphic to
$\bE_{\bbH'\bbF'}(A^{\bbH'}),\bE_{\bbH'\bbF'}(M_U^{\bbH'})$.  The Hochschild complexes for $M_U^\bbH,M_U^\bbF$ are isomorphic
to those for $M_U^{\bbH'},M_U^{\bbF'}$, the map $\bE_{\bbH'\bbF'}$ induces a map of Hochschild complexes isomorphic to
that induced by the map $\bM_\bbH\to \bM_\bbF$.  Thus,  the problem is reduced to showing that
the map $\bE_{\bbH'\bbF'}$ induces an isomorphism of Hochschild complexes
$$
\Hoch(M_U^{\bbH'},A^{\bbH'})\to \Hoch(M_U^{\bbF'};A^{\bbF'});\quad \Hochcyc(M_U^{\bbH'})\to \Hochcyc(M_U^{\bbF'}).
$$

\subsection{Description of the functor $\bE_{\bbG'\bbF'}:\bbG'\to \bbF'$} 

Let $$\pi:\Phi\times E\times E\to \Fr\times_M \Fr\times E\times E\to (\Fr\times E)^2.
$$

We also  have  objects $\bbP_\MetR,\bbQ_\MetR\in \bbF'$.   
    We now have
$$
\bE_{\bbG'\bbF'}(T)=\bbQ_\MetR*\pi_! T* \bbP_\MetR.
$$
\subsection{The functor $\bE'_{ji}$ and a natural transformation $
P_{r_j}*_E \pi_! T*_E P_{r_i}[-2N]\to \bE'_{ji}(T).
$}
Let $\Delta:\Fr_i\times E\times E\to \Fr_i\times \Fr_i\times E\times E$ be the  embedding.  Let
us calculate
$
\Delta^!\bbP_\MetR
$

The object $\bbP_\MetR$ is supported on the subset $U_{ii}\times E\times E\subset \Fr_i\times \Fr_i$.  The embedding
$\Delta$ factors through $U_{ii}\times E\times E$ so that we have a closed embedding of smooth manifolds
$\delta:\Fr_i\times E\times E\to U_{ii}\times E\times E$.       The embedding $\delta$ is of codimension $2N$.

It follows that $\bbP_{\MetR}|_{U_{ii}\times E\times E}$ is non-characteristic on $T^*_{\Fr_i\times E\times E} U_{ii}\times E\times E$, therefore we have
 $$
\delta^!\bbP_{ii}\sim \Delta^{-1}\bbP_{ii}[-2N]\sim p_{\Fr_i}^{-1}P_{r_i}[-2N],
$$
where $p_{\Fr_i}:\Fr_i\times E\times E\to E\times E$ is the projection.

Let $$\pi: \Phi_{ij}\times E\times E\to \Fr_i\times_M \Fr_j\times E\times E\to \Fr_i\times \Fr_j\times E\times E$$
be the through map.
Let $$(\bE_{\bbG'\bbF'})_{ji}(T):=\bE_{\bbG'\bbF'}(T)|_{\Fr_j\times \Fr_i\times E\times E}$$ 
be the corresponding component.

Let $$
\bE'_{ji}(T):=P_{r_j}*_E \pi_! T*_E P_{r_i}[-2N].
$$
We have a natural transformation 
$$
\bE'_{ji}\to  (\bE_{\bbG'\bbF'})_{ji},
$$
where we have taken into account that $\bbQ_\MetR=\bbP_\MetR[2N]$.

\subsubsection{Splitting the groupoid map $\Phi_{ff}\times_M U\to \Sigma_{ff}\times_M U$} 
Let $f\in S$ be the maximal element.  In particular, we have
$f>>f$.

Let $\cL_f:=\Fr_f/\SU(N)$;  $$
\Sigma_{ff}=(\cL_f\times_M \cL_f)\times_{S^1\times S^1} \Re\times \Re,
$$
so that we have an embedding $i_f: \Sigma_{ff}\to \Sigma$.

We have maps of groupoids
$$
(\Phi_{ff}\rightrightarrows\Fr_f)\to (\Sigma_{ff}\rightrightarrows \cL_f)\to (\Sigma\rightrightarrows \cL).
$$

Let us split the map of groupoids
$$
(\Phi_{ff}\rightrightarrows\Fr_f)|_U\to (\Sigma_{ff}\rightrightarrows \cL_f)|_U.
$$
First, fix a $U(N)\times S^1$-equivariant  trivialization of  the bunlde $\Fr_f|_U=U\times (U(N)\times S^1)$.
Let  $$s:(U(N)\times S^1)\times (U(N)\times S^1)\to S^1\times S^1$$ 
be given by $s(g,h)=\det(g^{-1}h)$. 

Let  $$s_0:(S^1\times S^1)\times (S^1\times S^1)\to S^1\times S^1$$ 
be given by $s(g,h)=g^{-1}h$. 

 Let 
$$
\overline{(U(N)\times S^1)\times( U(N)\times S^1)}:=((U(N)\times S^1)\times (U(N)\times S^1))\times_{s,S^1\times S^1}\Re\times \Re;$$
$$
\overline{(S^1\times S^1)\times (S^1\times S^1)}:=((S^1\times S^1)\times (S^1\times S^1))\times_{s_0,S^1\times S^1}\Re\times \Re
$$
The above map now can be rewritten as
$$
\left(U\times \overline{(U(N)\times S^1)\times  (U(N)\times S^1)}\rightrightarrows U\times (U(N)\times S^1)\right)\to 
\left(U\times \overline{(S^1\times S^1)\times (S^1\times S^1)}\rightrightarrows U\times (S^1\times S^1)\right).
$$

Let $z:S^1=U(1)\to U(N)$ be the embedding induced by the embedding $\mathbb{C}^1\to \mathbb{C}^N$.
 It defines a splitting
$\overline{(S^1\times S^1)\times (S^1\times S^1)}\to \overline{(U(N)\times S^1)\times (U(N)\times S^1)}$, whence the desired splitting
$$
\zeta:\Sigma_{ff}|_U\to \Phi_{ff}|_U.
$$

Denote by $D_\zeta$ the codimension of a smooth embedding $\zeta$.

We have the following natural transformation of functors
$$
\bE_{\bbH'\bbF'}(T)_{ff}\cong \bE_{\bbG'\bbF',ff}\bE_{\bbH'\bbG'}(T)_{ff}\leftarrow 
\bE'_{ji}(\zeta_!\zeta^! \bE_{\bbH'\bbG'}(T))=:\cK(T).
$$

\subsubsection{Applying the  natural transformation to the Hochschild complexes}
Consider the following commutative diagram in $\ho \swell \CMetR$:
$$
\xymatrix{ \Hoch(M_U^{\bbH'};A^{\bbH'})\ar[r]^\alpha\ar[d]^\beta & \Hoch(M_U^{\bbF'};A^{\bbF'})\ar[d]^{\lambda}\\
                   \Hom_{\bbH'}(M_U^{\bbH'};A^{\bbH'})\ar[r]^{\mu} &\Hom_{\bbF'}((\bE_{\bH'\bF'}M_U^{\bbH'})_{ff};A^{\bbF'}_{ff})\ar[d]^{\nu}\\
&\Hom_{\bbF'}(\cK(M_U^{\bbH'});A^{\bbF'}_{ff})}
$$

Our goal is to show that $\alpha$ is an isomorphism.  We will deduce it from the following statements:

---$\nu\mu$ is an isomorphism;

--- $\nu\lambda$ is an isomorphism;

---$\beta$ is an isomorphism.

\subsection{ $\nu\mu$  is an  isomorphism}

  Let $$e:\Fr_f\times_M \Fr_f\times E\times E\to \Fr_f\times \Fr_f\times E\times E
$$
be the embedding.

We have a natural transformation $e_!e^{-1}A^{\bbF'}_{ff}[-2N]\to A^{\bbF'}_{ff}$ which induces
a homotopy equivalence 
$$
\Hom_{\bbF'}(\cK(M_U^{\bbH'});e_!e^{-1}A^{\bbF'}_{ff}[-2N])\stackrel\sim\to\Hom_{\bbF'}(\cK(M_U^{\bbH'});A^{\bbF'}_{ff})
$$
The natural transformation  $\cK\to \bE_{\bbH'\bbF'}$ factors as
$$
\cK\to e_!e^{-1}\bE_{\bbH'\bbF'}[-2N]\to \bE_{\bbH'\bbF'},
$$
so that we have a map
$$
\hom_{\bbH'}(M_U^{\bbH'};A^{\bbH'})\to \hom_{\bbF'}(\cK(M_U^{\bbH'});e_!e^{-1}A^{\bbF'}_{ff}[-2N])
$$
the problem reduces to showing this map to be a homotopy equivalence.

Let us rewrite the definition of the functor $e^{-1}\bE_{\bbH'\bbF'}$.

Let  $f'=(g_f,R_f/2)$.      Let $W\subset M\times M$ consist of all pairs
$(m_1,m_2)$, where $B_{f}(m_2)\subset B_{f'}(m_1)$.  Let $\pi_1,\pi_2:W\to M$ be the projections
onto the factors.
We have a projection
$$
U_{f'f}\to W.
$$
Let 
$$
\cU_0:=U_{f'f}\times_W U_{f'f}.
$$
The projection $U_{f'f}\to \Fr_{f'}$ induces an identification
$$
\cU_0\to \Fr_{f'}\times_M \Fr_{f'}\times_{M,\pi_1} W.
$$

We have  projections
$$\cU_0\to\Fr_{f'}\times_M \Fr_{f'};\quad  \cU_0\to \Fr_f\times_M \Fr_f.
$$

Let  also $\cU:=\Phi_{f'f'}\times_M W$.   The above projections lift canonically to projections
$$
p_{f'}:\cU\to \Phi_{f'f'};\quad p_{f}:\cU\to \Phi_{ff}.
$$

Let $u_1,u_2:\cU\to U_{f'f}$ be the projections.

Let $v_1,v_2:\cU\times E\times E\times E\times E\to U_{f'f}\times E\times E$ be given by
$v_1(x,e_1,e_2,e_3,e_4)=(u_1(x),e_1,e_2);\quad v_2(x,e_1,e_2,e_3,e_4)=(u_2(x),e_3,e_4)$;
Let $w:\cU\times E\times E\times E\times E\to \cU\times E\times E$ be the projection
along the second and the third factors of $E$:
$$
w(x,e_1,e_2,e_3,e_4)=(x,e_1,e_4).
$$

Let $\vs:\Phi_{ff}\to \overline{U(N)\times S^1}$ be the projection.

 Let $$\ve: \cU\times E\times E\times E\times E\to 
\overline{U(N)\times S^1}\times E\times E
$$ 
be the map induced by $\sigma p_{f}:\cU\to \overline{U(N)\times S^1}$
and the
 projection onto  the second and the third factors of $E$.

Let us define an object $\cZ_1\in \psh(\cU\times E\times E,\DMetR_0)$ as follows.
Set
$$
\cZ_1:=w_!(v_1^{-1}\bbP_{f'f}\otimes v_2^{-1}\bbQ_{ff'}\otimes \ve^{-1}\bGamma).
$$

Set 
$$
\cZ:= P_{r_f}*_E \cZ_1*_E P_{r_f}.
$$
Let $q:\cU\stackrel{p_f}\to \Phi_{ff}\to \Sigma_{ff}$
be the through map.
We have
$$
e^{-1}\bE_{\bbH'\bbF'}(T)\sim  \rho_!p_{f'!}(\cZ\otimes  q^{-1}T),
$$
where $\rho:\Phi_{f'f'}\to \Fr_{f'}\times_M \Fr_{f'}$ is the covering map.

We have a natural transformation
$$p_{f'!}(\cZ\otimes q^{-1}T)\to p_{f'!}\cZ\otimes p_{f'!} q^{-1}T[2N].
$$
which is a homotopy equivalence as long as $T$ is a  locally constant object.

We have an identification $\Phi_{f'f'}=\Phi_{f'}\times \overline{U(N)\times S^1}$ which induces
an identification
$$\cU=U_{f'f}\times \overline{U(N)\times S^1}.
$$
Let 
is rewrite the definition of $\cZ_1$ under this identification.

Let $p_{12},p_{34},p_{14}:U_{f'f}\times E\times E\times E\times E\to U_{f'f}\times E\times E$
be the projections along the last two factors of $E$ (resp. the first two factors of $E$,resp. along
the 2-nd and the 3-rd factors).

Let $p_{23}:U_{f'f}\times E\times E\times E\times E\to E\times E$ be the projection
onto the 2-nd and the 3-rd factors of $E$.

Set
$$
\cY:=p_{14!}(p_{12}^{-1}\bbP_{f'f}\otimes p_{23}^{-1}\bGamma\otimes p_{34}^{-1}\bbQ_{ff'})
$$
We have
$$
\cZ\cong \cY*_E \bbS.
$$


Sec. \ref{mobile} implies a homotopy equivalence
$p_{f'!}\cY\sim  t^{-1}P_{r_f}[2N]$, where $t:\Fr_{f'}\times E\times E\to E\times E$ is the projection.

Therefore, we have a natural transformation of functors
$$
e^{-1}\bE_{\bbH'\bbF'}(T)\to  P_{r_f}*_E\rho_!(p_{f!}q^{-1}T*_E \bbS)*_E [2N]=:\bF(T)
$$

which is a homotopy equivalence provided that $T$ is locally constant.

We have a natural transformation $\cK\to \bE_{\bbH'\bbF'}\to \bF$, and the problem reduces to
showing that the induced map
\begin{equation}\label{isomm}
\hom_{\bbH'}(M_U^{\bbH'}; A^{\bbH'})\to \hom_{\bbF'}(\cK(M_U^{\bbH'});\bF(A^{\bbH'}))
\end{equation}
is a homotopy equivalence.

Let $$\iota:\Re\times \Re=\overline{U(1)\times S^1}\to \overline{U(N)\times S^1}$$ be the embedding.

Let $r:\Fr_f\times \overline{U(N)\times S^1}\to \Sigma$ be the projection and
let $$
\zeta:\Sigma=\Fr_f\times \overline{U(1)\times S^1}\to \Fr_f\times \overline{U(N)\times S^1}
$$ 
be the embedding  induced by $\iota$.
We have 
$$
\cK(T)=\rho_!((\zeta_! T)\boxtimes \bGamma)*_E \bbS)
$$

We now have the following chaing of homotopy equivalences:
\begin{multline*}
\hom_{\bbF'}(\cK(M_U^{\bbH'});\bF(A^{\bbH'}))\sim 
\hom_{\bbF'}(((\zeta_! M_U^{\bbH'})\boxtimes \bGamma)*_E \bbS[D_\zeta];\rho^!P_{r_f}*_E\rho_!(p_{f'!}q^{-1}A^{\bbH'}*_E \bbS)*_E P_{r_f}[2N])\\
\sim \hom_{\bbF'}(((\zeta_! M_U^{\bbH'})\boxtimes \bGamma)*_E \bbS[D_\zeta];\rho^{-1}\rho_!(p_{f'!}q^{-1}A^{\bbH'}*_E \bbS)[2N])\\
\sim
\hom_{\bbF'}((\zeta_! M_U^{\bbH'})\boxtimes \bGamma)*_E \bbS[D_\zeta];\bigoplus\limits_{ab}(T_a p_{f'!}q^{-1}A^{\bbH'}*_E \bbS[2b][2N])\\
\sim \hom_{\bbF'}(M_U^{\bbH'}\boxtimes \bGamma;\zeta^{-1}\bigoplus\limits_{ab} T_a(p_{f'!}q^{-1}A^{\bbH'}[2b][2N])\boxtimes \ZQ_{\Delta_E}),
\end{multline*}
 We have used a homotopy equivalence 
$$
\rho^{-1}\rho_!\bbS\sim \bigoplus\limits_{ab} T_b\pi^{-1}\bbS[2a],
$$
where $$
\pi:\Phi_{f'f'}\times E\times E\to \overline{U(N)\times S^1}\times E\times E
$$
is the projection.
The map $\bGamma\to \ZQ_{\Delta_E}$ in  $\sh_{\val S}(E\times E)$ induces a homotopy equivalence
\begin{multline*}
\hom_{\Sigma}(M_U^{\bbH'};\bigoplus\limits_{ab} T_a\zeta^{-1}(p_{f'!}q^{-1}A^{\bbH'})[2b][2N])\sim 
\hom_{\bbF'}(M_U^{\bbH'}\boxtimes \bGamma;\zeta^{-1}\bigoplus\limits_{ab} T_a(p_{f'!}q^{-1}A^{\bbH'}[2b][2N])\boxtimes \ZQ_{\Delta_E})
\end{multline*}

Let the embedding $$
\eta:\Fr_{f'}\times \overline{U(N)\times S^1}\to U_{f'f}\times \overline{U(N)\times S^1}
$$
 be  induced by the diagonal embedding $\Fr_{f'}\to U_{f'f}$. 
We have an equivalence (provided that $T$ is locally constant):
$$
T\sim \zeta^{-1}\eta_!\eta^{-1}p_{f'!}q^{-1} T
$$
as well as a map
$$
\zeta^{-1}\eta_!\eta^{-1}p_{f'!}q^{-1} T\to \zeta^{-1}p_{f'!}q^{-1}T[2N]
$$
which results in the following natural transformation
\begin{equation}\label{natur}
T\to \zeta^{-1}p_{f'!}q^{-1}T[2N].
\end{equation}
This transformation induces a map
$$
\hom_{\Sigma}(M_U^{\bbH'};\bigoplus\limits_{ab} T_a A^{\bbH'}[2b])\to 
\hom_{\Sigma}(M_U^{\bbH'};\bigoplus\limits_{ab} T_a\zeta^{-1}(p_{f'!}q^{-1}A^{\bbH'})[2b][2N])
$$
It follows that this map is isomorphic to the map in (\ref{isomm}).
The problem therefore reduces to showing that the map
$$
\ZQ_{\Sigma,0}\to  \zeta^{-1}p_{f'!}q^{-1}\ZQ_{[\Sigma,0]}[2N]
$$
induced by the natural transformation in (\ref{natur}) is a homotopy equivalence which is immediate.

\subsection{$\nu\lambda$ is an isomorphism}

 The second statement reduces to the following one.  We have a map
$\cK(M_U^{\bbH'})_{ff}\to (M_U^{\bbF'})_{ff}$.    We are to show that the following composition 
is an isomorphism in $\bbF'$:
\begin{equation}\label{eqB}
A^{\bbF'}_{if}*\cK(M_U^{\bbH'})_{ff}*A^{\bbF'}_{fj}\to  A^{\bbF'}_{if}*M_U^{\bbF'}*A_{fj}^{\bbF'}\to
(M^{\bbF'}_U)_{ij}.
\end{equation}

A. Let $T\subset S$, $T\in \SMetR'$.  Let $$h_T:\bS_0(\ccA)\to \Classic(\ccA)\langle \val T\rangle
$$ be the projection.

 It suffices to show that
\begin{multline} \label{ukv}h_T(A^{\bbF'}_{if})*h_T(\cK(M_U^{\bbH'})_{ff})*h_T(A^{\bbF'}_{fj})\\
\to
h_T(A^{\bbF'}_{if}*\cK(M_U^{\bbH'})_{ff}*A^{\bbF'}_{fi})\cong  h_T(A^{\bbF'}_{jf}*M_U^{\bbF'}*A_{fj}^{\bbF'})\\
\to
h_T((M^{\bbF'}_U)_{ij}).
\end{multline}
is an isomorphism in $\ho \Classic(\ccA)\langle \val T\rangle.$

Let $\vs$ be the minimum of $T$.   We have a natural map $h_T\to h_{\{\vs\}}$ which induces an isomorphism
on both sides of the above through map. Therefore,  we reduce the problem to the case
when $T=\{\vs\}$ is a one element set.

B.  Let us reformulate (\ref{ukv}).    Let $p_M:\Fr_f\to M$, $p_M:\Fr_\vs\to M$ be the projections.
  Let  $$
\cJ\subset \Fr_i\times \Fr_\vs\times \Fr_\vs\times \Fr_f\times \Fr_\vs\times \Fr_\vs\times \Fr_f\times \Fr_\vs\times \Fr_\vs\times\Fr_j
$$
be the subset consisting of all points $(f_1,f_2,\ldots,f_{10})$, where $(f_2,f_1)\in U_{\vs i_{R_i}}$,
$(f_2,f_3)\in \Fr_\vs\times_M \Fr_\vs$,  $(f_3,f_4)\in  U_{ \vs f_{R_f}}$,    $(f_5,f_4)\in  U_{\vs f_{R_f}}$,
$(f_5,f_6)\in \Fr_\vs\times_M \Fr_\vs|_U$,  $(f_6,f_7)\in U_{\vs f_{R_f}}$,  $(f_8,f_7)\in U_{\vs f_{R_f}}$,
$(f_8,f_9)\in \Fr_\vs\times_M \Fr_\vs$, $(f_9,f_{10})\in U_{\vs j_{R_j}}$.

Let $p_{mn}$,  $1\leq m<n\leq 10$ be the projection from $$
\Fr_i\times \Fr_\vs\times \Fr_\vs\times \Fr_f\times \Fr_\vs\times \Fr_\vs\times \Fr_f\times \Fr_\vs\times \Fr_\vs\times\Fr_j
$$
onto the $m$-th and the $n$-th factors.  Still denote by $p_{mn}$ the restriction of these projections 
onto $\cJ$.

Let  $\pi:\cJ\to (U(N)\times S^1)^3$ be the through map
$$\cJ\stackrel{p_{23}\times  p_{56}\times p_{89}}\longrightarrow (\Fr_\vs\times_M \Fr_\vs)^3\to \overline{U(N)\times S^1}^3$$
be the composition.

Let 
$$
\ovJ:=\cJ\times_{\pi,(U(N)\times S^1)^3} \overline{U(N)\times S^1}^3.
$$

Let $q_{23},q_{56},q_{89}:\ovJ\to \overline{U(N)\times S^1}$ be the projections.

We then rewrite the LHS of (\ref{ukv}) as
\begin{multline*}
p_{1,10!} (p_{12}^{-1}(P_{r_i}\cP_{i_{R_i}\vs})*_E q_{23}^{-1}\gamma_{r_\vs}*_E q_{34}^{-1} \cQ_{\vs f}*_E
q_{45}^{-1}\cP_{f \vs}*_E q_{56}^{-1} \zeta_!\zeta^!\gamma_{r_\vs}*_E p_{67}^{-1} \cQ_{\vs f}*_E p_{78}^{-1}\cP_{f \vs}*_E
q_{89}^{-1} \gamma_{r_\vs}*p_{9,10}^{-1} \cQ_{\vs j_{R_j}}P_{r_j}).
\end{multline*}

Let 
$$
\cK\subset \Fr_i\times (\Fr_\vs\times \Fr_\vs)^3\times \Fr_j.
$$
consist of all points $(f_1,f_2,f_3,f_5,f_6,f_8,f_9,f_{10})$, where $(f_1,f_2)\in U_{i_{R_i}\vs}$;
$(f_2,f_3),(f_8,f_9)\in \Fr_\vs\times_M \Fr_\vs$,  $(f_5,f_6)\in \Fr_\vs\times_M \Fr_\vs|_{U}$,
$(f_9,f_{10})\in U_{\vs j_{R_j}}$.   Let   $q_{23},q_{56},q_{89}:\cK\to U(N)$ be the projections.
Let $\pi:q_{23}\times q_{56}\times q_{89}$,  $\pi:\cK\to (U(N)\times S^1)^3$.  Let
$$\ovK:=\cK\times_{(U(N)\times S^1)^3}(\overline{U(N)\times S^1})^3.$$

Let $p_\cK:\ovJ\to \ovK$ be the projection.    Let $p^\cK_{ij}$ be the projection from $\ovK$ to the corresponding
pair of factors.

The map $\cP_{\vs f}*\cQ_{f\vs}\to \ZQ_{\Delta_{\Fr_\vs\times E}}$ and the algebra structure on $\gamma_{r_\vs}$
induce the following maps
\begin{multline*}
p_{1,10!} (p_{12}^{-1}(P_{r_i}\cP_{i_{R_i}\vs})*_E q_{23}^{-1}\gamma_{r_\vs}*_E q_{34}^{-1} \cQ_{\vs f}*_E
q_{45}^{-1}\cP_{f \vs}*_E q_{56}^{-1} \zeta_!\zeta^!\gamma_{r_\vs}*_E p_{67}^{-1} \cQ_{\vs f}*_E p_{78}^{-1}\cP_{f \vs}*_E
q_{89}^{-1} \gamma_{r_\vs}*p_{9,10}^{-1} \cQ_{\vs j_{R_j}}P_{r_j})\\
\to p_{1,10!} (p_{12}^{-1}(P_{r_i}\cP_{i_{R_i}\vs})*_E q_{23}^{-1}\gamma_{r_\vs}*_E p_{35}^{-1}\ZQ_{\Delta_{\Fr_\vs\times E}}*_E q_{56}^{-1} \zeta_!\zeta^!\gamma_{r_\vs}*_E p_{68}^{-1}\ZQ_{\Delta_{\Fr_\vs\times E}}*_E
q_{89}^{-1} \gamma_{r_\vs}*p_{9,10}^{-1} \cQ_{\vs j_{R_j}}P_{r_j}\\
\to p_{1,10!} (p_{12}^{-1}(P_{r_i}\cP_{i_{R_i}\vs})*_E (q_{235689}^{-1}\gamma_{r_\vs}|_{\Delta^{(6)}_{\Fr_\vs}}) *_E
p_{9,10}^{-1} \cQ_{\vs j_{R_j}}P_{r_j})
\end{multline*}

The resulting map is isomorphic to the map (\ref{ukv}).

C.  Let $\vs=(g_\vs,r_\vs)$.   For an $r>0$ denote $\vs_r:=(g_\vs,r)$; $i_r:=(g_i,r)$, etc.

 As $\vs<<i,j$,  
there exists a number $\rho>r_\vs$ such that  $\vs_\rho<<i,j$.
Fix such a $\rho$.

In particular we have   the following families of symplectic embeddings
$$ I_{\vs_\rho i}:U_{\vs_\rho i}\times  B_\rho\to B_{R_i};\quad  I_{\vs_\rho j}:U_{\vs_\rho j}\times  B_\rho\to B_{R_j}.
$$

2) There exists a number $r>0$ satisfying:  
$$
\forall n_1,n_2 \in M: B_{\vs}(n_1)\cap B_{\vs_{r}}(n_2)\neq \emptyset \Rightarrow  B_{\vs_{r}}(n_2)\subset \Int B_{\vs_\rho}(n_1).
$$

C.  Define
$\cJ_r\subset \Fr_{\vs_r}\times \cJ\times \Fr_{\vs_r}$ to consist of all points
$(f_0,f_1,\ldots,f_{10},f_{11})$ satisfying $B_{\vs_r}(f_0)\subset B_{i}(f_1)$,  $B_{\vs_r}(f_0)\cap B_{\vs}(f_2)\neq \emptyset$,
$B_{\vs}(f_9)\cap B_{\vs_r}(f_{11})\neq \emptyset$, $B_{j}(f_{10})\supset B_{\vs_r}(f_{11})$.

Let us define $\cK_r\subset \Fr_{\vs_r}\times \cK\times \Fr_{\vs_r}$ in a similar way.

The inclusions $B_{\vs_r}(f_0)\subset B_{i}(f_1)$,  $B_{j}(f_{10})\supset B_{\vs_r}(f_{11})$
define, via their differential at 0, the maps
$$
q_{01},q_{10,11}:\cJ_r\to \cK_r\to \Sp(2N).
$$
Let $$q_{\cJ}:\cJ_r\to \cK_r\stackrel{q_{\cK}} \to  \Sp(2N)\times \Sp(2N)$$ 
be the product.

Let $$\ovJ_r:=(\cJ_r\times_{\cJ} \ovJ)\times_{q_\cJ,\Sp(2N)\times \Sp(2N)} \overline{\Sp(2N)\times \Sp(2N)};  
$$
$$
\ovK_r:=(\cK_r\times_{\cK} \ovK)\times_{q_\cK,\Sp(2N)\times \Sp(2N)} \overline{\Sp(2N)\times \Sp(2N)}.
$$
We have a projection $p_{012}:\ovJ'_r\to K_r$.  Let $\ovJ_r:=\ovJ'_r\times_{K_r} K'_r$;
$\ovK_r:=\ovK'_r\times_{K_r}^{K'_r}$.

We have the following inclusions:   $B_{\vs_r}(f_0)\subset   B_{\vs_{\rho}}(f_2)\subset  B_{i}(f_1)$.
Whence the following families of symplectic embeddings
$$
a:\ovJ_r\times B_{\vs_r}(f_0)\to B_{\vs_\rho}(f_2);
$$
$$
b:\ovJ_r\times B_{\vs_\rho}(f_2)\to B_{i}(f_1).
$$
Let also
$$
b_r:\ovJ_r\times B_\vs(f_2)\to B_i(f_1)
$$
be the restriction of $b$.   By the construction, $b_r$ and $ba$ are graded.  Therefore, so is $b$ and, hence $a$.
Let $c=ba$.  Let $\cA^i,\cB^i$ be the quantizations of $a,b$.   The quantization of $b_r$ is then
$\cB^i*_E P_{r_\vs}.
$

Let $\cC^i:= \cB^i*_E \cA^i$ be the quantization of $ba$.   Let $\cA^j,\cB^j,\cC^j$ be the similar objects for $j$.

It suffices to show that the following map is an isomorphism
\begin{multline}\label{ukvs}
p_{0,1,10,11!} \Big((\cC^i)^t*_E  \cB^i*_E P_{r_\vs} *_E q_{23}^{-1}\gamma_{r_\vs}*_E q_{34}^{-1} \cQ_{r_\vs f}*_E
q_{45}^{-1}\cP_{f r_\vs}\\
*_E q_{56}^{-1} \zeta_!\zeta^!\gamma_{r_\vs}*_E p_{67}^{-1} \cQ_{r_\vs f}*_E p_{78}^{-1}\cP_{f r_\vs}*_E
q_{89}^{-1} \gamma_{r_\vs}*_E  P_{r_\vs}*_E (\cB^j)^t*_E( \cC^j)^t\Big)\\
\to p_{0,1,10,11!} ((\cC^i)^t*_E  (\cB^i)*_E P_{r_\vs} *_E (q_{235689}^{-1}\gamma_{r_\vs}|_{\Delta^{(6)}_{\Fr_\vs}}) *_E
 P_{r_\vs}*_E (\cB^j)^t*_E (\cC^j)^t)
\end{multline}

We have $\cC^t*_E \cB =\cA^t*_E\cB^t*_E\cB*_E P_{r_\vs}\cong \cA^t*_E P_{r_\vs}$. Accordingly, we rewrite:

\begin{multline}\label{ukvs1}
p_{0,1,10,11!} ((\cA^i)^t*P_{r_\vs}*_E q_{23}^{-1}\gamma_{r_\vs}*_E q_{34}^{-1} \cQ_{r_\vs f}*_E
q_{45}^{-1}\cP_{f r_\vs}*_E q_{56}^{-1} \zeta_!\zeta^!\gamma_{r_\vs}*_E p_{67}^{-1} \cQ_{r_\vs f}*_E p_{78}^{-1}\cP_{f r_\vs}*_E
q_{89}^{-1} \gamma_{r_\vs}*_E  P_{r_\vs}*_E \cA^j)\\
\to p_{0,1,10,11!} ((\cA^i)^t*_E   P_{r_\vs} *_E (q_{235689}^{-1}\gamma_{r_\vs}|_{\Delta^{(6)}_{\Fr_\vs}}) *_E
 P_{r_\vs}*_E  \cA^j)
\end{multline}

Let $\kappa:\cJ_r\to \Fr_\vs\times \Fr_i\times \Fr_\vs^6\times \Fr_j\times \Fr_\vs$ be the projection.
Denote by $Q$ the image of $\kappa$.   Let us  number the factors from 1 to 10.
denote  by $p_{ij}$ the projection from $Q$ onto the corresponding
pair of factors.  As above, we have projections
$$
q_{34}, q_{56},q_{78}:Q\to \Fr_\vs\times_M \Fr_\vs\to U(N)\times S^1.
$$
Let  $q:Q\to (U(N)\times S^1)^3$ be the product.  Let
$$
\ovQ:=Q\times_{q,(U(N)\times S^1)}\overline{U(N)\times S^1}.
$$

The projection $\kappa$ lifts to a map 
$$
k:\ovJ_r\to \ovQ.
$$

  Because of the geodesic convexity,  $k:\cJ_r\to \ovQ$ is a fibration
with contractible fibers.

Denote 
\begin{multline*}
\cS:=\\
(\cA^i)^t*_E P_{r_\vs}*_E q_{23}^{-1}\gamma_{r_\vs}*_E q_{34}^{-1} \cQ_{r_\vs f}*_E
q_{45}^{-1}\cP_{f r_\vs}*_E q_{56}^{-1} \zeta_!\zeta^!\gamma_{r_\vs}*_E p_{67}^{-1} \cQ_{r_\vs f}*_E p_{78}^{-1}\cP_{f r_\vs}*_E
q_{89}^{-1} \gamma_{r_\vs}*_E  P_{r_\vs}*_E \cA^j)\\
\in \psh(\ovJ_r,\DMetR_0),
\end{multline*}
as in (\ref{ukvs1}).

It follows that 

(1)
$\cS$ is constant along the fibers of $q$, in particular, the map $k^{-1}k_* \cS\to \cS$ is a homotopy equivalence.

(2)  The object $k_*\cS$ is supported on the following subset $\ovQ'\subset \ovQ$ consisting of all points $p$
whose projection onto $Q$,
$
(f_1,f_2,\ldots,f_{10}),
$
satisfies  $B_{\vs}(f_4)\cap B_\vs(f_5)\neq \emptyset$,  $B_\vs(f_6)\cap B_\vs(f_7)\neq \emptyset$.

Let $\cN\subset M^{10}$ consist of all points $(m_1,m_2,\ldots,m_{10})$ satisfying:

$B_{\vs_r}(m_1)\subset B_i(m_2),$

$B_{\vs_r}(m_1)\cap B_{\vs}(m_3) \neq \emptyset,
$
$
m_3=m_4,
$
$
B_{\vs}(m_4)\cap  B_{\vs}(m_5)\neq \emptyset$,
$
m_5=m_6\in U,
$
$
B_{\vs}(m_6)\cap B_\vs(m_7)\neq \emptyset,
$
$
m_7=m_8,
$
$
B_{\vs}(m_8)\cap B_{\vs_r}(m_{10})\neq \emptyset,
$
$B_{i}(m_{9})\supset B_{\vs_r}(m_{10}).
$

According to Assumption \ref{asb}, $U\subset B_{\vs}(m_0)$ for some point $m_0\in M$.
  Let $V:=B_{\vs_{R'_\vs}}(m_0)$ and $W:=B_{\vs_{R''_\vs}}(m_0)$, see Sec. \ref{radius} for definition of $R'_\vs,R''_\vs$.
It follows that

(1)   $\cN\subset V^{10}$;

(2)   if $n\in V$, then $B_{\vs_{R_\vs}}(v)\subset W$.

Choose a section $s$ of $\Fr_\vs|_W$ (which exists because $W$ is a ball).       

Let $\pi:\ovQ'\to \cN$ be the projection.   There exist maps $\lambda_i:\cN\to \ovU(N)$, $i=1,2,3,4$, such that 
for every point $m=(m_1,\ldots,m_{10})\in \cM_M$, the set $\pi^{-1}p$ consists of all points of the form
   \begin{multline*}
(g_1s(m_1),g_2s_i(m_2),g_2\lambda_1(p)s(m_3),c_1\lambda_2(p)s(m_4),c_1s(m_5),c_2s(m_6),c_2\lambda_3(p)s(m_7),\\
g_3\lambda_4(p)s(m_8),g_3 s_j(m_{9}),g_4 s(m_{10})),
\end{multline*}
where $g_i\in \overline{U(N)\times S^1}$, $c_i\in \Re\times \Re=\overline{S^1\times S^1}$.
We thus have an identification
$$
\ovQ'\cong \overline{U(N)\times S^1}^4\times \overline{S^1\times S^1}^2\times \cN.
$$
Let us rewrite (\ref{ukvs1}) under this identification.   Denote by $g_i:\ovQ'\to \overline{U(N)\times S^1}$;
$c_i:\ovQ'\to \overline{S^1\times S^1}$ the corresponding projections.

Next, we have a   graded family of symplectic embeddings 
$$
\overline{ U(N)\times S^1}\times V\to     U(N)\times S^1\times V\cong \Fr_\vs|_{V}\times B_{R_\vs}\to B_{R''_\vs}.
$$
Denote $\ovFr_\vs:=\overline{ U(N)\times S^1}\times V$.   Denote by 
$\cT\in \sh_q(\ovFr_\vs\times E\times E)$ the quantization of this family.

We have projections $\pi_k:\ovQ'\to M\stackrel s\to  \ovFr_\vs$, $k=0,2,3,\ldots,8,9,11$;
$g_i:\ovQ'\to \overline{U(N)\times S^1}$ be the projections.

We can now rewrite (\ref{ukvs1}) as follows (the functions $\lambda_i$ drop out from the expression)

\begin{multline*}
p_{1,2,9,10!}( g_1^{-1}\bbS*_E \pi_1^{-1}\cT^t *_E  \pi_3^{-1}\cT|_E*\gamma_{r_\vs}*\pi_4^{-1}\cT^t*_E
\pi_5^{-1}\cT|_U*_E\gamma_{r_\vs}*_E \pi_6^{-1}\cT^t|_U*_E \pi_7^{-1}\cT*_E \gamma_{r_\vs}*_E \pi_8^{-1}\cT^t*_E \pi_{10}^{-1}\cT*_E g_{4}^{-1}\bbS)\\
\to p_{1,2,9,10!}( g_1^{-1}\bbS*_E \pi_1^{-1}\cT^t *_E \pi_3^{-1}\cT|_U*_E \pi_{34567!}(\ZQ_{\Delta^{(5)}_{\Fr_\vs}\boxtimes \gamma_{r_\vs}}) *_E\pi_8^{-1}\cT^t|_U*_E \pi_{10}^{-1}\cT*_E g_{4}^{-1}\bbS)
\end{multline*}

 Let $p_i:W^5\to W$ be projections.
Let $q_i:W^5\to W\to M\stackrel{s}\to  \Fr_\vs$ be the projections.  The statement reduces to showing that 
the following map is an isomorphism

\begin{multline*}
p_{1,5!}(  q_1^{-1}\cT^t *_E q_2^{-1}\cT*_E\gamma_{r_\vs}*_E q_2^{-1}\cT^t*_E q_3^{-1}\cT|_U*_E\gamma_{r_\vs}*_E q_3^{-1}\cT^t|_U*_E q^{-1}_4\cT*_E \gamma_{r_\vs}*_E q_4^{-1}\cT^t*_E q_5^{-1}\cT)\\
\to p_{1,5!}(q_1^{-1}\cT^t *_E  q_2^{-1}\cT|_U*_E\pi_{234!}(\ZQ_{\Delta^{(3)}_{\Fr_\vs}\boxtimes \gamma_{r_\vs}}) *_E q_4^{-1}\cT^t|_U*_E q_5^{-1}\cT)
\end{multline*}

The statement now follows from the mobility of the family $\cT$, see Condition 3 in Sec  \ref{radius}.

\section{ $\bY(\bbB)$-modules} 

\subsection{Infinitesimal operads}   {\em  An infinitesimal circular operad (ICO)} is 
a colored operad (CO) in the category $\Com_{[0,1]}(\cC)$.    

Given an ICO $\cO$ let $|\cO|
$ denote the CO over $\cC$ obtained
via applying  the totalization functor $\Com_{[0,1]}(\cC)\to \cC$.

Let $\Gr_{[0,1]}(\cC)$ be a category whose every object is a pair $(X^0,X^1)$ of objects of $\cC$.
We set $$
\hom((X^0,X^1);(Y^0,Y^1))=\hom(X^0,Y^0)\oplus \hom(X^1,Y^1). 
$$  

Define a SMC on $\Gr_{[0,1]}(\cC)$ by setting $$
(X^0,X^1)\otimes (Y^0,Y^1):=X^0\otimes Y^0; X^0\otimes Y^1\oplus X^1\otimes Y^0.
$$

We have  functors $\Gr_{[0,1]}(\cC)\stackrel\iota\to \Com_{[0,1]}(\cC)\stackrel\pi\to \Gr_{[0,1]}(\cC)$.

 {\em  A split ICO  } is 
a CO in the category $\Gr_{[0,1]}(\cC)$.   

Given a  split ICO $\cO'$, denote by $X(\cO')$   the set of all maps $D:(\cO')^0\to (\cO')^1[1]$
such that $(\cO',D)$ is an ICO.

\subsection{Categories $\bY(\cO)$, $\bY^\cyc(\cO)$} \label{sco} Let $\cO$ be a CO in $\cC$.  Denote by  $S(\cO)$ the category over $\sets$ whose every 
object is a split ICO $\cO'$  along with an identification of VCO
$(\cO')^0=\cO$.    One can construct a category $\bY(\cO)$ enriched over $\cC$, whose objects are
of the form $(n)^\uncyc$,  $(n)^\cyc$, such that $S(\cO)$ is  equivalent to the category of
functors $\bY(\cO)^\op\to \cC$.  Namely, given $\cU \in S(\cO)$, the corresponding functor
associates $\cU^1(n)^\uncyc$ to $(n)^\uncyc$ and $\cU^1(n)^\cyc$ to $(n)^\cyc$.

 Denote by $\bY^\cyc(\cO)\subset \bY(\cO)$ the full sub-category
consisting of all objects $(n)^\cyc$, $n> 0$.   The category $\bY^\cyc(\cO)$ only depends on the non-cyclic part of $\cO$ .  We therefore can write $\bY^\cyc(\cO^\uncyc)$ instead of $\bY^\cyc(\cO)$.

  It follows that the structure of a circular operad  on a pair $(\cO^\uncyc,\cO^\cyc)$
is equivalent to that of an asymmetric operad on $\cO^\uncyc$ and of a functor $\bY^\cyc(\cO^\uncyc)\to \cC$ on $\cO^\cyc$.

\subsection{Studying $\bY(\gr\bbB)$-modules}

 Let us get back to the sequence of monoidal categories with trace and their 
functors in (\ref{strght}).   Denote by $\cO^{\bbM}$ the full circular operad of $M^\bbM$, and likewise
for $\cO^{\bbH},\cO^{\bbF}$, etc. We therefore have a sequence of circular operads and their maps:
\begin{equation}\label{phiphi}
\gr \bbB\to \cO^{\bbM}\to  \cO^{\bbH}\to \gr \cO^\bbF.
\end{equation}
Hence we have an induced structure of $\bY(\gr\bbB)$-modules and their maps   
\begin{equation}\label{xoduli}
\cO^\bbM\to\cO^\bbH\to\gr\cO^\bbF.
\end{equation}

In this section we will formulate and prove several statements about these modules.
We start with introducing certain endofunctors on the category of $\bY(\gr\bbB)$-modules and their maps.
\subsection{Description of $\bY(\assoc),\bY(\bbB)$}  The full sub-category of $\bY(\assoc)$ consisting 
of objects $(n)^\uncyc$ is isomorphic to $\Delta$; the full sub-category  consisting of all objects $(n)^\cyc$ is 
isomorphic to the cyclic category $\Lambda$.   Denote by $I:\Delta\to \Lambda$ the natural embedding 
$
[n]\mapsto (n+1).
$
 The category $\bY(\assoc)$ is then  equivalent to the cone of $I$.

Let us now describe the category $\bY(\bbB)$.   Let $\alg\in \bR_q$,
$$
\alg:=\bigoplus\limits_{a,b\in \bZ} T_a\unit_{\bR_q}[2b].
$$
We have a unital commutative algebra structure on $\alg$.       We have a functor
$\Hochcyc(\alg):\Lambda^\op\to \bR_q$.  
We have an induced commutative algebra structure on $\Hochcyc(\alg)$.

We have  $$
\hom_{\bY(\bbB)}((m)^\uncyc,(n)^\uncyc)=\Hochcyc_m(\alg)\otimes \hom_{\bY(\assoc)}((m)^\uncyc,(n)^\uncyc);
$$ 
$$
\hom_{\bY(\bbB)}((m)^\uncyc,(n)^\cyc)=\Hochcyc_m(\alg)\otimes \hom_{\bY(\assoc)}((m)^\uncyc,(n)^\cyc);
$$ 
$$
\hom_{\bY(\bbB)}((m)^\cyc,(n)^\cyc)=\Hochcyc_{n-1}(\alg)\otimes \hom_{\bY(\assoc)}((m)^\cyc,(n)^\cyc).
$$ 

Denote by $\bY(\bbB)^{\noncyc}\subset \bY(\bbB)$ the full sub-category  consisting of all objects $(n)^\noncyc$.
Denote by $\bY(\bbB)^{\cyc}\subset \bY(\bbB)$ the full sub-category  consisting of all objects $(m)^\cyc$.
We have a functor $I:\bY(\bbB)^{\noncyc}\to \bY(\bbB)^{\cyc}$, where $I(n)^\noncyc=(n)^\cyc$.     

Let $F:\bY(\bbB)\to \bR_q$ (or  $F:\bY(\gr\bbB)\to \bR_0$) be a functor.  Let $F^{\noncyc}:=F|_{\bY(\bbB)^{\noncyc}}$;
$F^{\cyc}:=F|_{\bY(\bbB)^{\cyc}}$ be restrictions.    The embedding $\assoc\to \bbB$ gives rise
to further restrictions
$$
F^{\noncyc,\Delta},F^{\cyc,\Delta}:\Delta\to \bR_q.
$$

Denote
$$
\Gamma^{\noncyc}(F):=C^*(F^{\noncyc,\Delta});\quad \Gamma^{\cyc}(F):=C^*(F^{\cyc,\Delta}).
$$

\subsection{Condensation} \label{condensc}
  Let $\bbE$ be a colored operad in $\ccA$ with two colors, to be called $\ba$ and $\bm$.
Let $c_1,c_2,\ldots,c_n,c$ be colors. Let $S_\ba,S_\bm\subset\{1,2,\ldots,n\}$ be defined by
$$
S_\ba=\{i|c_i=\ba\};\ S_\bm=\{i|c_i=\bm\}.
$$

if $c=\bm$,  set 
$$
\bbE^\uncyc(c_1,c_2,\ldots,c_n|c)=\ZQ;
$$

if $c=\ba$ and $c_1=c_2=\cdots=c_n=\ba$, set $\bbE^\uncyc(c_1,c_2,\ldots,c_n|c)=\ZQ$. 
Otherwise,  set $\bbE(c_1,c_2,\ldots,c_n|c)=0$.
Finally, set $\bbE^\cyc(c_0,c_1,\ldots,c_n)=\ZQ$  for all $n$.

Let $\pi:\{\ba,\bm\}\to \pt$ be the projection.  We have a natural map of circular  operads 
\begin{equation}\label{bbaa}
\bbE\to\pi^{-1} \bbB.
\end{equation}
Define two functors  $\iota,\epsilon:\bY(\bbE)\to \bY(\bbB)$ as follows.

---$\iota$.   On objects:
 $$
\iota(c_1,c_2,\ldots,c_n)^\uncyc=(n)^\uncyc;\ \iota(c_1,c_2,\ldots,c_n)^\cyc:= (n)^\cyc. 
$$    On maps: induced by (\ref{bbaa});

---$\ve$. On objects $$
\ve(c_1,c_2,\ldots,c_n)^\cyc:=(|S_\bm|)^\cyc,\ \ve(c_1,c_2,\ldots,c_n)^\uncyc:=(|S_\bm|)^\uncyc.
$$
On arrows, cyclic part: we
have
$$
\bY(\bbE)^\cyc((c'_1,c'_2,\ldots,c'_m)^\cyc;(c_1,c_2,\ldots,c_n)^\cyc)=\bigoplus\limits_{f}\bbE(f),
$$
where the direct sum is taken over all $f:\{1,2,\ldots,n\}\to \{1,2,\ldots,m\}$, where $f$ is cyclic and  $f(S_\bm)\subset S'_\bm$.
We thus have an induced map $f':S_\bm\to S'_\bm$.  Next, we have an isomorphism
$\bbE(f)\stackrel\sim\to  \bbB(f')$, whence an induced map 
$$
\bY(\bbE)^\cyc((c_1,c_2,\ldots,c_n)^\cyc;(c'_1,c'_2,\ldots,c'_m)^\cyc)\to \bY(\bbB)^\cyc((|S_\bm|)^\cyc;(|S'_\bm|)^\cyc).
$$

The action of $\ve$ on the non-cyclic arrows is defined in a similar way.

We now have an object
$$
\cK_\con:=\ve\otimes^L_{\bY(\bbE)} \iota \in \swell(\bY(\bbB)^\op\otimes \bY(\bbB)).
$$

Let $\cA,\cB$ be arbitrary categories enriched over $\cC$.   Let $K\in \Doplus(\cA^\op\otimes \cB)$.
Let $F:\cB\to \cC$. 
Define the co-convolution with $K$,  $\hom_\cB(K;F):\cA\to \cC$, where
$$
\hom_\cB(K,F)(a):=\hom_\cB(h_a\circ_{\cA} K; F),\quad a\in A.
$$

Let $F:\bY(\bbB)\to \ccA$ be a functor.
Let  $$\con(F):=\hom_{\bY(\bbB)^\op}(K_\con,F)$$ so that
$\con$ is an endofunctor on the category of functors $\bY(\bbB)\to \ccA$.

We have a  zig-zag natural transformation from $\con$ to $\Id$.  It is defined as follows. 
Denote by $\mu:\bY(\bbB)\to \bY(\bbE)c$,  $\mu((n)^\cyc)=(\bm,\bm,\ldots,\bm)^\cyc$, ;
$\mu((n)^\uncyc)=(\bm,\bm,\ldots,\bm)^\uncyc$, ($n$ occurrences of $\bm$).
We have $\iota\mu=\ve\mu=\Id$, whence the following sequence of maps:
$$
I_{\bY(\bbB)^\op}\otimes^L_{\bY(\bbB)} I_{\bY(\bbB)}\cong \ve\mu\otimes^L_{\bY(\bbB)}  \iota\mu\to
 \ve\otimes^L_{\bY(\bbE)} \iota=K_\con.
$$
where $$I_{\bY(\bbB)^\op}:\bY(\bbB)^\op\to \swell(\bY( \bbB)^\op,\  I_{\bY(\bbB)}:\bY(\bbB)\to \swell \bY(\bbB)$$ are the canonical maps.

Denote $K_I:=I_{\bY(\bbB)^\op}\otimes^L_{\bY(\bbB)} I_{\bY(\bbB)}$ so that we have a  map
\begin{equation}\label{idcon}
K_I\to K_\con
\end{equation}
 We then have 
maps
$$
\con F =\hom_{\bY(\bbB)}(K_\con;F)\to \hom_{\bY(\bbB)}(K_I;F)\stackrel\sim\leftarrow F
$$
which define the required zig-zag.

\begin{Lemma}  Let $V:\bY(\bbB)\to \ccA$. Then the  maps
\begin{equation}\label{kconkiuncyc}
\Gamma^{\uncyc}(\con V)\to \Gamma^{\uncyc}(\hom_{\bY(\bbB)}(K_I;V)),
\end{equation}
\begin{equation}\label{kconkicyc}
\Gamma^{\cyc}(\con V)\to \Gamma^{\cyc}(\hom_{\bY(\bbB)}(K_I;V)),
\end{equation}
induced by (\ref{idcon}),  are
  homotopy equivalences.
\end{Lemma}
{\em Sketch of the proof.} 
Consider the  cyclic case (\ref{kconkicyc}),  the case (\ref{kconkiuncyc}) is treated similarly.
    We have isomorphisms
$$
\Gamma^{\cyc}(\con(V))\cong \hom_{\bY(\bbB)^\cyc}(C_*(K_\con|_{\Lambda^\op\otimes(\bY(\bbB)^\cyc});V);
$$
$$
\Gamma^\cyc(\hom(K_I,V))\cong \hom_{\bY(\bbB)^\cyc}(C_*(K_I|_{\Lambda^\op\otimes \bY(\bbB)^\cyc});V).
$$

The problem now reduces to showing that the following map
$$
C_*(K_\con|_{\Lambda^\op\otimes \bY(\bbB)^\cyc})\to C_*(K_I|_{\Lambda^\op\otimes \bY(\bbB)^\cyc})
$$
induced by (\ref{idcon}),
is a homotopy equivalence.

\subsubsection{Semi-orthogonal decomposition} Let $\ccA$ be a ground category.  
 Call a functor $F:\Delta\to \ccA$ {\em quasi-constant} if every arrow in $\Delta$ acts by a homotopy equivalence.
Equivalently, let $R_\Delta$ be a semi-free resolution of the constant functor $\underline{\ZQ}:\Delta\to \ccA$.
Then $F$ is quasi-constant iff the natural map
$$
\hom_\Delta(R_\Delta;F)\otimes R_\Delta\to F
$$
is a termwise homotopy equivalence.

Denote by  $\cC_{\bY(\bbB)}$ the category of functors
$\bY(\bbB)^\op\to \ccA$, enriched over $\sets$.  

 Let $$
\cC_{\bY(\bbB),\const}\subset \cC_{\bY(\bbB)}
$$
be the full sub-category consisting of all objects $X$, where  $X^{\cyc}|_{\Delta}$ and $X^\uncyc$
are quasi-constant.  

Let $\cC_{\bY(\bbB),\perp}\subset \cC_{\bY(\bbB)}$ be 'the right orthogonal complement' to
 $\cC_{\bY(\bbB),\const}$ --- it consists  of all objects $Y$ such that for every $X\in \cC_{\bY(\bbB),\const}$, we have
$$
\Rhom_{\bY(\bbB)}(X,Y)\sim 0.
$$ 

\begin{Lemma}\label{kritperp}
 We have $X\in \cC_{\bY(\bbB),\perp}$ iff $\Gamma^{\noncyc}(X)\sim 0$ and $\Gamma^{\cyc}(X)\sim 0$.
\end{Lemma}

We have the following semi-orthogonal decomposition.

\begin{Proposition}   Let $X\in \cC_{\bY(\bbB)}$.     Then

1)  $\con X\in \cC_{\bY(\bbB),\const}$;

2) $\Cone(\con X\to \hom_{\bY(\bbB)}(K_I;X))\in \cC_{\bY(\bbB),\perp}.$
\end{Proposition}

\subsubsection{Lemma on the map $\cO^\bbH\to\gr \cO^\bbF$}  Let us get back to the sequence (\ref{xoduli}).
\begin{Proposition}\label{ftoh}
1) The cone of the map
$
\cO^\bbH\to\gr \cO^\bbF
$
belongs to $\cC_{\bY(\bbB),\perp}$;

2)  The induced map $\con \cO^\bbH\to \con \gr\cO^\bbF$ is a term-wise homotopy equivalence.
\end{Proposition}

{\em Sketch of the proof.} Follows from Prop.  \ref{hchsch} and  Lemma \ref{kritperp}.

\subsection{Studying the object $\con(\cO^\bbH)$}

\subsection{The map $I:(\cO^\bbH)^{\uncyc}|_\Delta\to (\cO^\bbH)^{\cyc}|_\Delta$ is a homotopy equivalence}
We have a structure map
\begin{equation}\label{mapIH}
I:(\cO^\bbH)^{\uncyc}|_\Delta\to (\cO^\bbH)^{\cyc}|_\Delta.
\end{equation}

The following Proposition is straightforward

\begin{Proposition}
The map (\ref{mapIH}) is a  term-wise homotopy equivalence.
\end{Proposition}

\subsection{The map $c_1$}\label{diamondc} We have a functor  $\delta: \bY(\bbB)\to \bY(\bbB)\otimes \bY(\assoc)$;
where
$$
\delta(a)=a\otimes a, \ a\in \Ob \bY(\bbB).
$$

For $X:\bY(\bbB)\to \bR_q$;  $V:\bY(\assoc)\to \bR_q$, set  $X\diamond Y:=\delta^{-1}(X\otimes Y)$.

 Let $\cR:\bY(\assoc)\to \ccA$ be a semi-free resolution of the constant functor $\underline{\ZQ}$
We have a termwise homotopy equivalence $$
X\diamond \cR\to X\diamond \underline{\ZQ}= X.
$$

Next, we have an element
$c_1:\cR\to \cR[2]$.

  We have an indudced map
$$
c_1:=c_1^X:X\diamond \cR\to X\diamond \cR[2].
$$
where $\otimes$ is the term-wise product of functors.  We can iterate so as to get maps
$$
(c_1^X)^n:X\diamond \cR\to X\diamond \cR[2n].
$$

Call $X$ {\em $c_1$-nilpotent}, if $(c_1^X)^n$ is homotopy equivalent to 0 for $n$ large enough.

\subsubsection{Formulating the Proposition}

For each  subset $V\in \Re$ and $S\in \SMetR$
we have a   functor $\tau_{V,S}:\bR_0\to \bR_0$, where  $(\gr^a(\tau_{V,S}X)_T=X_S^a$ if $a\in V$ and $T=S$. Set
$\gr^a(\tau_{V,S}X)_T=0$ if either $T\neq S$, or $a\notin V$, or both.

\begin{Proposition}\label{c1}
Let $a>0$ and $S\in \SMetR$.
The object $\tau_{(0,a),S}\cO^\bbH$ is   $c_1$-nilpotent.
\end{Proposition}

The rest of the subsection will be devoted to proving  this proposition.

\subsubsection{Reduction to the case $V$ has only one point}    
It follows that $$
\tau_{(0,a),S} =\bigoplus\limits_{b\in (0,a),b\in \val S.\bZ} \tau_{\{b\},S}.
$$

Denote $\tau_{b,S}:=\tau_{\{b\},S}$.     Therefore, it suffices to show that
for every $a\in \val S.\bZ_{>0}$,  $\tau_{a,S}\cO^\bbH $ is  $c_1$-nilpotent.

\subsubsection{Localizing with respect to a covering of $M$} Choose a triangulation $\cT$ of $M$  as in Sec \ref{triang}.   Let $\cU$ be the open covering of $M$ consisting of all stars of $\cT$..
 We  now have a map
$
\cO^{\bbB_\cU}\to \cO^\bbH.
$
For each  $U\in \cU$, we have  an object $M_U^\bbB\in \cO^{\bbB_\cU}$.   Its full operad is isomorphic to $\bbB$.
We also have its image $M_U^\bbH\in \cO^\bbH$.  Denote by $\cO_U^\bbH$ the full operad of $M_U^\bbH$.
   we have a homotopy equivalence
$$
\cO^\bbH\stackrel\sim\to \holim_{U\in \cU} \cO_U^\bbH.
$$
As the covering $\cU$ is finite,  it suffices to show that
$$
\tau_{a,S} \cO_U^\bbH
$$
is $c_1$-nilpotent for each $U\in \cU$.

\subsubsection{Restriction to $\Delta$}  

  Denote
$$
X^\cyc:=(\tau_{a,S}\cO_U^\bbH)^{\cyc}|_\Lambda.
$$
It follows that $X^{\cyc}|_{\Delta}$ is  quasi-constant.  Let $R_\Delta\to \underline\ZQ$ be  a free resolution in the category 
of functors $\Delta\to \ccA$.   It now follows
that the natural map
$$
R_\Delta\otimes \hom_{\Delta}(R_\Delta;X^{\cyc}|_\Delta)\to X^{\cyc}|_\Delta
$$
is a  termwise homotopy equivalence, that is we have a homotopy equivalence
$$
H\otimes R_\Delta\sim X^{\cyc}|_{\Delta}
$$
for an appropriate $H\in \ccA$.

Let us now pass to $\Lambda$.    Let $\cR$ be a free resolution of the constant functor
$$
\underline{\ZQ}:\Lambda\to \ccA.
$$
 Consider  an object
$$
\hom_\Lambda(\cR;X^{\cyc}|_\Lambda)\otimes \cR.
$$
The endomorphism $u:\cR\to \cR[2]$ acts on both tensor functors. Denote those actions by $u_1$ and $u_2$ respectively.    We can now build a cone
$$
X_\cR:=\Cone(\hom_\Lambda(\cR;X^{\cyc}|_\Lambda)\otimes \cR[-2]\stackrel{u_2-u_1}\longrightarrow\hom_\Lambda(\cR;X^{\cyc}|_\Lambda)\otimes \cR)
$$

We have a natural map $X_\cR\to X$ which is a homotopy equivalence.
This implies:
\begin{Lemma} Suppose $X|_\Delta$  is  quasi-constant.   Then  $c_1^n$-action on $X$ is homotopy equivalent  to 0 
iff such is the induced  $c_1^n$-action  on $R\hom_{\Lambda}(\cR;X)$. 
\end{Lemma}
\subsubsection{Passage to $\ZZ$}   Recall that the category $M^\bbH$ is enriched over the SMC of $\ZZ$-graded objects in $\bR_0$, to be denoted by $\bR_0(\ZZ)$.   

The tensor functor $\pi:\bR_0(\ZZ)\to \bR_0$ gives the induced structure over $\bR_0$.
Denote by $\gr_\ZZ^p$ the corresponding graded component.

\subsubsection{An algebra $\xi$}    

As the set $U$ is contractible,  we can trivialize 
the groupoid 
$$
\Sigma|_U\rightrightarrows \cL_U
$$
so that it becomes isomorphic to
$$
\overline{(S^1\times S^1)^2}\times U\rightrightarrows S^1\times S^1\times U.
$$

Choose a point $a\in S^1\times S^1$.  We then have an embedding of groupoids
$$
i:(\bZ\times \bZ\times U\rightrightarrows U)\to \overline{(S^1\times S^1)^2}\times U\rightrightarrows S^1\times S^1\times U,
$$
where the groupoid structure comes from the group law on $\bZ\times \bZ$.

Set $\xi:=i^!M^\bbH$.    The full operad of the pair $(\xi,M^\bbH)$ is pseudo-contractible
 (in the category $\bR_0(\ZZ)$).  Therefore, via straightening out, we can construct
a zig-zag map  of functors $\Lambda\to \ccA$:
$$
\tau_{a}\cO_\xi^\bbH\to \tau_{a,S}\cO^\bbH,
$$
where $\cO_\xi$ is the full operad of $\xi$.
It follows that the induced map

$$
\hom(\cR;\tau_{a}\cO_\xi^\bbH)\to \hom(\cR; \tau_{a,S}\cO^\bbH))
$$
is  a homotopy equivalence.

Therefore, the problem reduces to showing that the $c_1^n$-action on the LHS is homotopy equivalent to 0.

Next,  it follows, that the full operad of $\xi$ is homotopy equivalent to a full operad
of an algebra $\eta$ in a monoidal category $G$ enriched over $\bR_0(\ZZ)$, to be now defined.

An object of the category $G$ is a $\bZ\times \bZ$-graded object in $\bR$.
Set
$$
\gr^{(a_1,a_2)}_{\ZZ}\hom_G(X,Y):= T_{a_1}\gr^{(a_1,a_2)}\hom(X;Y)[2a_2].
$$
 The tensor product on $G$ is defined as that of graded objects:
$$
\gr^c_G (X\otimes Y)=\bigoplus\limits_{a\in \bZ\times \bZ}\gr^a_GX\otimes \gr^{c-a}_G Y,\quad a,c\in \bZ\times \bZ.
$$

Let us now define the above mentioned algebra $\eta$ in $G$ by setting 
 $\gr^{a,b}_G\eta:=\ZQ$ for all $a,b\in \bZ$.
In other words $\eta=\ZQ[u^{-1},u,v^{-1},v]$, where $u:\ZQ\to T_1\gr^{1,0}\eta$;  $v:\ZQ\to \gr^{0,1}\eta[2]$ are generators. 
The statement now can be proven by direct computation based on HKR, provided that $\ZQ=\bbQ$.
\subsection{Studying the map $\cO^\bbM\to \cO^\bbH$}   
   We have a zig-zag homotopy equivalence between $\cO^\bbM$ and $\cO^\bbB\diamond \cO^{\bbM_0}$.

Let $X:\bY(\assoc)\to \bR_0$;  $N:\bY(\bbB)^\op\to \bR_0$ be  functors.   Denote $$
H(X;N):=R\hom_{\bY(\assoc)}( \cO^\bbB\diamond X;N).
$$
Let  $a\in \bY(\assoc)^\op$.  Let  $h_{a}:\bY(\assoc)\to \bR_0$ be the corresponding free object.   Set
$$
N^{\bZ\times \bZ}(a):=H(h_{a};N).
$$
We have a structure  of a functor $\bY(\assoc)\to \bR_0$ on $N^{\bZ\times \bZ}$.

Let $Z:\bY(\assoc)\to \bR_0$ be a functor.    We have  maps
$$
Z(a)\to R\hom_{\bY(\assoc)}(h_{a},Z)\to R\hom_{\bY(\bbB)}(\cO^\bbB\diamond h_{a};\cO^\bbB\diamond Z),
$$
whence  a map $$
Z\to (\bbB\diamond Z)^{\bZ\times \bZ}.
$$
Therefore, we have a zig-zag map from  $\cO^{\bbM_0}$ to $(\cO^{\bbM})^{\bZ\times \bZ}$.
Denote by $C:\bY(\assoc)\to \bR_0$ the cone of this map.    
\begin{Proposition}  We have $C$ is right orthogonal to every quasi-constant  functor $T:\bY(\assoc)\to \bR_0$,
that is $R\hom(T,C)\sim 0$.
\end{Proposition}

{\em Sketch of the proof}   Let $C_U$ be the cone of 
  the zig-zag map from
$
\cO_U^{\bbM_0}$ to $(\cO^{\bbM}_U)^{\bZ\times \bZ}.
$

It suffices to show a similar statement for $C_U$, which can  be checked directly.

\subsubsection{The structure of the  object $\cO^{\bbM_0}$}  

Call an object $X:\bY(\assoc)\to \bR_0$ {\em strictly constant} if there is an object
$A\in \bR_0$ and a zig-zag homotopy equivalence of functors $\Lambda\to \bR_0$ from $A\otimes \underline{\ZQ}$ to $ X$.

\begin{Proposition} \label{strcon} The object $\cO^{\bbM_0}$ is strictly constant
\end{Proposition}
{\em Sketch of the proof} As above, fix a covering $\cU$ of $M$.
Denote $T:=\cO^{\bbM_0}$;  $S_U:=\cO_U^{\bbM_0}$.
We have  a functor $S:\cU\times \bY(\assoc)\to \bR_0$,  $S(U,F):=S_U(F)$.
 It follows that $S(U,n)\sim \ZQ$ for all $(U,n)\in \cU\times \bY(\assoc)$ which implies that we have a zig-zag termwise homotopy equivalence
between $S$ and $\underline{\ZQ}$.  

We also have a homotopy equivalence of cyclic objects
$$
T\stackrel\sim\to \holim_{U\in \cU} S(U,-).
$$
The statement now follows.
\section{Reformulation}\label{data}
\subsection{$c_1$-Localization} 
\subsubsection{Category $\Fun(\cD)$}    Let $\cD,\cC$ be categories enriched over $\ccA$. Let $\cC$ be tensored over $\ccA$.
Let $\Funct(\cD)$ be the category of  $\ccA$-functors $\cD\to \cC$.   Let $\Fun(\cD):=\Doplus(\cD\otimes \cC)$.
  We have the resolution functor
$\cR:\Funct(\cD)\to \Fun(\cD^\op)$.

2)  Let $$\br:=\Id_{\cD}\otimes^L_\cD \Id_{\cD^\op}\in \Fun(\cD\otimes \cD^\op).
$$
 Let us define an endofunctor $$\lambda:\Fun(\cD)\to \Fun(\cD),$$
$$
\lambda(T):=\br\circ_\cD T.
$$

\subsubsection{$c_1$-Localization in $\Fun((\bY(\cO)^\cyc)^\op)$}\label{c1local}
 Let now $\cD:=\Fun((\bY(\cO)^\cyc)^\op).$

We have a natural transformation $c_1:\br\to \br[2]$ which induces a natural transformation
$c_1:\lambda\to \lambda[2]$.

Let us define an endofunctor
$$
\loc:  \Doplus \Fun((\bY(\cO)^\cyc)^\op)\to  \Doplus\Fun((\bY(\cO)^\cyc)^\op).
$$
$$
\loc(T):=\hocolim (\lambda(T)\stackrel {c_1}\to \lambda(T)[2]\stackrel{c_1}\to \lambda(T) \cR[4]\stackrel {c_1}\to \cdots).
$$

For $F:\bY(\cO)^\cyc\to \cC$,   denote
$$
F_\loc:=\loc(\cR(F)).
$$

It follows that $F_\loc\sim 0$ as long as the $c_1$-action on $F$ is nilpotent.

The functor $\loc$ extends to an endofunctor on
$$
\Doplusprod\Fun((\bY(\cO)^\cyc)^\op).
$$

\subsection{Extending the functor $\con$}

Fix a ground SMC $\cC$ enriched over $\ccA$.
\subsubsection{Functor $\bbD$}  Let $\cG$ be a category enriched over $\cC$. 
Let $\nu\in \cG,V\in \cC$.  Define a functor $V_\nu:\cG\to \cC$,
$$
V_\nu(\mu):=\ihom_\cC(\hom_{\cG}(\mu;\nu);V).
$$
 
Define a functor $$\bbD:\cG\otimes \cC\to \Fun(\cG),
$$ where
$$
\bbD(\nu\otimes V)=V_\nu\otimes^L_{\cG} \Id_{\cG^\op}.
$$

The functor $\bbD$ extends to an endofunctor  on $\Fun(\cG)$.

\subsubsection{Co-convolution with a kernel}  Let $\cD_1,\cD_2$ be categories enriched over $\cC$.
 Let $K\in \Fun (\cD_1^\op\otimes \cD_2).
$
and $F\in \Fun(\cD_2)$.  We then have an object
$$
\hom_{\cD_2}(K,F)\in \Dprod\Fun (\cD_1).
$$

We can now apply the functor $\bbD$ so as to get an object
$$\bbDhom_{\cD_2}(K,F)\in \Dprod\Fun (\cD_1)
$$
Call $\bbDhom$ the co-convolution with $K$.

One extends $\bbDhom(K,-)$ to a functor $$
\Dprod\Fun (\cD_2)\to \Dprod\Fun(\cD_1).
$$

\subsubsection{Composition of co-convolutions} 
Let $$
K_{21}\in   \Fun(\cD_2^\op\otimes \cD_1),\quad K_{32}\in   \Fun(\cD_3^\op\otimes \cD_2).
$$
Consider the composition 
$$
\bbDhom(K_{32},\bbDhom(K_{21};-)):\Dprod\Fun(\cD_3)\to \Dprod\Fun(\cD_1).
$$

and construct  natural transformations  which are termwise homotopy equivalences and whose composition is the identity:
\begin{equation}\label{rettrac}
\bbDhom(K_{32}\circ_{\cD_2} K_{21};-)\to \bbDhom(K_{32},\bbDhom(K_{21};-))\to \bbDhom(K_{32}\circ_{\cD_2} K_{21};-).
\end{equation}

Let us first  do it for
$$K_{32}=\bigoplus\limits_{c\in C} \delta_3^c\otimes d_2^c\otimes U^c,\quad K_{21}=\bigoplus\limits_{b\in B} \delta_2^b\otimes d_1^b\otimes U^b,$$
where $U^c,U^b\in \cC$,  $\delta_3^c\in \cD_3^\op$,  $d_2^c\in \cD_2$, $\delta_2^b\in \cD_2^\op$,
$d_1^b\in \cD_1$.
Let also 
$$
X=\prod\limits_{z\in Z}\bigoplus\limits_{a\in A} d^{za}\otimes U^{za}, \quad d^{za}\in \cD_1,\ U^{za}\in \cC.
$$
We have
$$
\bbDhom(K_{21};X)=\prod\limits_{z\in Z,b\in B}\bigoplus\limits_{a\in A} t\otimes^L_{t\in \cD_2} \hom_\cC\left(U^b\otimes \hom_{\cD_2}(\delta_2^b,t);
U^{za}\otimes \hom_{\cD_1}(d_1^b;d^{za})\right).
$$
  
Next,
\begin{multline*}
\bbDhom(K_{32};\bbDhom(K_{21};X))\\=
\bbD\prod_{z\in Z,b\in B,c\in C}\bigoplus\limits_{a\in A} \delta_3^c\otimes \hom_{\cD_2}(d_2^c,t)\otimes^L_{t\in \cD_2}\hom_\cC\left(U^b\otimes U^c\otimes \hom_{\cD_2}(\delta_2^b,t);
U^{za}\otimes \hom_{\cD_1}(d_1^b;d^{za})\right).
\end{multline*}

In other words we can decompose
$$
\bbDhom(K_{32};\bbDhom(K_{21};X))=GH(K_{32},K_{21},X),
$$
where 
\begin{multline*}
H:(D\bigoplus (\cD_3^\op\otimes \cD_2\otimes \cC))^\op\otimes 
(D\bigoplus(\cD_2^\op\otimes \cD_1\otimes \cC))^\op \otimes D\prod\bigoplus (\cD_1\otimes \cC)\\
\to D\bigoplus(\cD_3^\op\otimes \cD_2\otimes \cC\otimes \cD_2^\op\otimes \cD_1\otimes \cC)^\op\otimes 
D\prod\bigoplus(\cD_1\otimes \cC)\\
\stackrel{\ihom}\longrightarrow
D\prod\bigoplus(\cD_3\otimes  \cD_2^\op\otimes \cC^\op\otimes \cD_2\otimes \cD_1^\op\otimes \cC^\op\otimes \cD_1\otimes \cC).
\end{multline*}
Next, 
$$
G:D\prod\bigoplus(\cD_3\otimes  \cD_2^\op\otimes \cC^\op\otimes \cD_2\otimes \cD_1^\op\otimes \cC^\op\otimes \cD_1\otimes \cC)
\to D\prod\bigoplus(\cD_3\otimes \cC),
$$
is induced by the following functor
$$
G_0:\cD_3\otimes  \cD_2^\op\otimes \cC^\op\otimes \cD_2\otimes \cD_1^\op\otimes \cC^\op\otimes \cD_1\otimes \cC\to D\prod\bigoplus(\cD_3\otimes \cC),
$$
\begin{multline*}
G_0(\delta_3,d_2,U_c,\delta_2,d_1,U_b,d_1^a,U_a)\\
:= \bbD\left(\delta_3\otimes \hom_{\cD_2}(d_2,t)\otimes^L_{t\in \cD_2}\hom_\cC(U_b\otimes U_c\otimes \hom_{\cD_2}(\delta_2,t);
U_a\otimes \hom_{\cD_1}(d_1;d_1^a))\right).
\end{multline*}

Likewise, we can decompose $$\hom_{\cD_1}(K_{32}\circ_{\cD_2} K_{21};X)=EH(K_{32},K_{21},X),
$$
where $E$ is induced by the following functor
$$
E_0:\cD_3\otimes  \cD_2^\op\otimes \cC^\op\otimes \cD_2\otimes \cD_1^\op\otimes \cC^\op\otimes \cD_1\otimes \cC\to D\prod\bigoplus(\cD_3\otimes \cC),
$$
\begin{multline*}
E_0(\delta_3,d_2,U_c,\delta_2,d_1,U_b,d_1^a,U_a)\\:= \bbD\left(\delta_3\otimes \hom_\cC(U_b\otimes U_c\otimes \hom_{\cD_2}(\delta_2,d_2);
U_a\otimes \hom_{\cD_1}(d_1;d_1^a))\right).
\end{multline*}

The retraction (\ref{rettrac}) is induced by the following retraction
\begin{multline*}
\hom_\cC(U_b\otimes U_c\otimes \hom_{\cD_2}(\delta_2,d_2);
U_a\otimes \hom_{\cD_1}(d_1;d_1^a))\\
\to 
\hom_{\cD_2}(d_2,t)\otimes^L_{t\in \cD_2}\hom_\cC(U_b\otimes U_c\otimes \hom_{\cD_2}(\delta_2,t);
U_a\otimes \hom_{\cD_1}(d_1;d_1^a))\\
\to \hom_\cC(U_b\otimes U_c\otimes \hom_{\cD_2}(\delta_2,d_2);
U_a\otimes \hom_{\cD_1}(d_1;d_1^a)) 
\end{multline*}

where the leftmost arrow is induced by letting $t=d_2$ and the rightmost arrow is induced by the composition of hom's.

\subsubsection{The identity kernel $K_I$} Let $$
K_I:=\Id_{\cD^\op}\otimes^L \Id_{\cD}\in \Fun(\cD^\op\otimes \cD) .
$$
For every $F\in \Fun(\cD)$, we have a natural transformation
$$
F\circ_{\cD} K_I\stackrel\sim\to F
$$
which produces an element $$
i_F:\ZQ\to \hom(F\circ_\cD K_I;F)=\bbDhom(F\circ_\cD K_I;F)\to \hom(F;\bbDhom(K_I;F)),
$$
whence a natural transformation
\begin{equation}\label{ki}
F\mapsto \bbDhom(K_I;F).
\end{equation}
This transformation extends onto all $F\in \Dprod\Fun(\cD)$.

We will use the following Lemma
\begin{Lemma}
Let $L\in \Fun(\cD_1^\op\otimes \cD)$ and $F\in \Fun(\cD)$. The following diagram commutes:
$$
\xymatrix{ \bbDhom( L,F)\ar[r]^{(\ref{ki})}\ar[rd]^\iota & \bbDhom(L,\bbDhom(K_I;F))\ar[d]^{(\ref{rettrac})}\\
                                                 &\bbDhom(L\circ_{\cD} K_I;F)},
$$
where the map $\iota$ is induced by the natural map $L\circ_{\cD} K_I\to L$.
\end{Lemma}
\subsubsection{The functor $\con$}  Let $\con:\Dprod\Fun((\bY(\bbB)^\cyc)^\op)\to \Dprod\Fun((\bY(\bbB)^\cyc)^\op)$ be defined
by
$$
\con(F):=\bbDhom(K_\con;F).
$$
The natural transformation $K_\Id\to K_\con$ induces natural transformations
\begin{equation}\label{conid}
\con(F)\to \bbDhom(K_\Id;F)\leftarrow  F.
\end{equation}
\begin{Lemma} \label{conidgamma} The induced maps
$$
\Gamma(\con(F))\to \Gamma(\bbDhom(K_\Id;F))\leftarrow \Gamma( F).
$$
are homotopy equivalences.
\end{Lemma}
\subsubsection{Reformulating the results}
\begin{Claim}\label{c1.2}
We have 
$$
\tau_{>0}\hom\left((\con(\cO^\bbF_\loc);\con(\cO^\bbF_\loc)\right)\sim 0
$$
\end{Claim}
This Claim follows from Prop \ref{ftoh} and \ref{c1}.

Fix a semi-free resolution $\bbBr\to \bbB$.
Let $Z:=\Doplusprod \Fun(\bY(\bbBr)^\op)$. Let $$\cR:\Funct(\bY(\bbBr);\cC)\to \Fun(\bY(\bbBr)^\op)
$$ be the canonical resolution functor.
\begin{Claim}\label{spl}
The following map in $\bR_0$:
$$
\hom_{Z}({\gr\bbBr}; {\gr\cO^\bbF})\to \hom_{Z}({\gr \bbBr};(\gr\cO^\bbF)_\loc),
$$
admits a splitting up-to a homotopy.
\end{Claim}
This Claim follows from Prop \ref{ftoh} and \ref{strcon}.

Fix a homotopy equivalence in $\bR_0$: (from the above Claim):
\begin{equation}\label{hel}
 \hom_{Z}({\gr\bbBr}; {\gr\cO^\bbF})\oplus \cL_0 \to \hom_{Z}({\gr \bbBr};(\gr\cO^\bbF)_\loc).
\end{equation}

\section{Quantization:Introduction}\label{intro2}
So far, we have constructed a CO $\cO^\bbF$ in the category  $\bR_q$ as well as a map
of CO $\red \bbB\to \red \cO^\bbF$ (where $\red$ means that the map is defined over $\bR_0$).
In   the rest of the paper   we deal with the problem of lifting this map onto the level of $\bR_q$, which we also
call the quantization.

More precisely, we start with  defining another CO, to be denoted by $\bbA$, and a map of  CO $\bbA\to \bbB$
(over $\bR_q$), and we lift the through map 
\begin{equation}\label{abf}
\red \bbA\to \red \bbB\to \red \cO^\bbF
\end{equation}
 (over $\bR_0$) to
a map $ \bbA\to \cO^\bbF$ (over $\bR_q$).  So far as $\bbB$ can be interpreted as an operad of $\bZ\times \bZ$-equivariant algebras,   $\bbA$ does then control  $\bZ\times \bZ$-equivariant  $A_\infty$- algebras with 'curvature'.
Putting aside the $\bZ\times \bZ$-equivariance, the structure of an  '$A_\infty$-algebra with curvature' on an object
$X$ of a monoidal category enriched over $\bR_q$ is a Maurer-Cartan element in the Hochschild
complex of  $X$ (viewed as an associative algebra with zero).   

The curvature version of the lifting problem  happens to be more tractable, as it turns out to be controlled by the Hochschild 
complexes   studied in Sec \ref{Hochsch}.  On the other hand,  this weaker version of the lifting problem is still
sufficient for our main task --- building the microlocal category. Indeed, the operad $\cO^\bbF$ being a full operad of
a certain object   in a monoidal cateory $D\bigoplus\bU(\bbF)$, the lifting provides for an $A_\infty$-algebra with curvature
structure on that object over $\bR_q$.   Next, the monoidal category $\Doplus\bU(\bbF)$ acts on the category
$\Doplus \bU(\bbF_R)$, and we define the microlocal category as that of  $A_\infty$ modules  over this algebra in
the category $\Doplus \bU(\bbF_R)$.

Let us outline the main steps of the quantization.   

First, we construct a tensor functor $\pi:\bR_q\to \Com_{\geq 0} \ccA$ (see Sec. \ref{combbb} below).
   The category $\bR_q$ is now enriched over $\Com_{\geq 0}\ccA$ (it is crucial that
the functor $\pi$ induces an isomorphism $$
\hom(\unit_{\bR_q};X)\to \hom(\unit_{\Com_{\geq 0}\ccA};\pi(X)).
$$
Accordingly, we get a functor $\pi_0:\bR_0\to \Gr_{\geq 0}\ccA$, where $\Gr_{\geq 0}\ccA$ is the category
of non-negatively graded objects in $\ccA$.

   We can now solve the resulting quantization problem by induction.  Namely, let $\Com_{[0,k]}\ccA$ be a SMC whose every object is a complex in $\ccA$ concentrated in 
degrees from $0$ to $k$.   The tensor product is defined by
setting
$$
(X\otimes Y)^p:=\bigoplus\limits_{q| 0\leq q\leq p} X^q\otimes Y^{p-q},\quad 0\leq p\leq k.
$$

We have tensor functors $\pr_k:\Com_{\geq 0}\ccA\to \Com_{[0,k]}\ccA$;  $\pr_{lk}:\Com_{[0,l]}\ccA\to \Com_{[0,k]}\ccA$,
$l\geq k$,  where $\pr_k(X)^s=X^s$; $\pr_{lk}(X)^s=X^s$, $0\leq s\leq k$.    One can now proceed similar to the case
of the deformation quantization, exploiting the analogy between  the  category $\Com_{\geq 0}\ccA$ and the category
of (continuous) modules over the ring of formal series in one variable $h$.   The category $\Com_{[0,k]}\ccA$ is then
analogous to the quotient ring $\ZQ[[h]]/(h^{k+1})$.  
Our strategy is thus as follows. First, we have categories $\pr_k \bR_q$ enriched over $\Com_{[0,k]}\ccA$.
Observe that $\pr_0\bR_q$ is isomorphic to the category $\bR_0$ as a category enriched over $\ccA$.
The starting data of the problem (i.e.  a  map of operads $\red \bbA\to \red\bbB\to \red \cO^\bbF$)
can, hence, be formulated within the category $\pr_0 \bR_q$.  This is the base of induction.

 Let the  $k$-th induction step consist
of lifting a map $\pr_k\bbA\to \pr_k\cO^\bbF$ in $\pr_k\bR_q$ to the level of  $\pr_{k+1}\bR_q$, that is to a map $\pr_{k+1}\bbA\to \pr_{k+1}\cO^\bbF$.    One can develop an obstruction theory for this kind of problems
(similar to the deformation quantization theory), the obstructions are controlled by the Hochschild complexes as
in Sec. \ref{Hochsch}),  unfortunately, they {\em do not vanish}.    Let us  describe the obstruction complex for the future 
sake. 
We start with defining the notion of an $A$-module, where $A$ is a CO, and of the object $\cD(M)$ classifying
 the derivations from $A$ to $M$.

Let $A$ be a CO in an SMC $\cC$ enriched over an SMC $\cB$. One first defines the notion of an $A$-module  $M$ as a  collection of objects
$M(n)^\uncyc$, $M(n)^\cyc$,  as well as a CO structure on the direct sum $A\oplus M$
such that the inclusion $A\to A\oplus M$ is a CO map;   and $M$ is an ideal whose square is 0. Equivalently,
$M$ is a $\bY(A)$-module.

One has the  notion of a derivation of $A$ with values in a $A$-module  $M$, denote the set
of such derivations by $\Der(M)$.    $A$ being quasi-free,   the functor $M\mapsto \Der(M)$ is representable and the representing object is denoted by $\Omega_A$ and called
 the  module of Kaehler differential on $A$, provided that $A$.   If $A$ is freely generated
 by a collection $G$, then $\Omega_A$ is freely generated over $A$ by $G$.  The following functor $F:\cB^\op\to \sets$ 
 is therefore
representable $F(T):=\hom(\Omega_A\otimes T; M)$, where $\hom$ is in the category of $A$-modules.
Denote the representing object by $\cD(M)$.  Let $M_\uncyc:=(M^\uncyc,0)$, $M_\cyc:=(0,M^\cyc)$ with the obvious
induced $A$-module structure.   Let $$
\cD(M)^\uncyc:=\cD(M_\uncyc),\  \cD(M)^\cyc:=\cD(M_\cyc).
$$
We have a map $h:\cD(M)^\cyc[-1]\to \cD(M)^\uncyc$ so that $\Cone h\cong \cD(M)$.

We will now define the obstruction complex. 
Let $A,B$ be circular operads in $\bR_q$ and let $f:\red A\to\red B$ be a map of CO.  We then get a $\red A$-module structure
on $\red B$, to be denote by $\red B_f$.
We therefore have an object $\cD(\red B_f)\in \Gr_{\geq 0} \ccA$.   This object is the above mentioned obstruction comlex to lifting $f$.

As was mentioned above,  the obstruction complex for the map $\red \bbA\to \red \cO^\bbF$ as in (\ref{abf}), does not
vanish.
However, one can modify the statement of the problem in order to
achieve the unobstructedness.   As usual, the modification is via replacing the structure of CO with  a richer structure
that we call VCO.  In order to motivate the definition of VCO we will first  perform certain manipulations on $\cO$.

 We start with the following categorical construction.
\subsection{The category $IE$}
 Let $E$ be a (symmetric monoidal) category enriched over $\Com_{\geq 0}\ccA$.
One then defines a (symmetric monoidal) category $IE$ enriched  over $\Com_{\geq 0}\ccA$ as follows.
One starts with defining a category $IE'$ whose every object  $X$ is a collection  of objects $\gr^n X\in E$, $n\in \bZ$.
Let $X,Y\in IE'$. Denote
$$
\cH^k:=\prod\limits_{n\in \bZ} \hom_E(\gr^n X; \gr^{n+k}Y)\in \Com_{\geq 0}\ccA.
$$
$$
\hom_{IE'}(X,Y):=\bigoplus\limits_{k>0} \tau_{\geq k}\cH^{-k} \bigoplus\prod\limits_{k\geq 0} \cH^k,
$$
where $\tau_{\geq k}:\Com_{\geq 0}\ccA\to \Com_{\geq 0}\ccA$ is the stupid truncation:
$(\tau_{\geq k} U)^l=0$ if $l<k$ and $(\tau_{\geq k}U)^l=U^l$ if $l\geq k$.

Observe that we have an isomorphism
$$
\bigoplus\limits_{k>0} \tau_{\geq k}\cH^{-k} \to \prod \limits_{k>0} \tau_{\geq k}\cH^{-k}.
$$
Let now set $IE:=D(IE')$.
 We will also use full sub-categories $I_{\leq 0}E,I_{\geq 0}E\subset IE$ consisting of all objects
$X$ with $\gr^kX=0$ if $k>0$ (resp.  if $k<0$).

\subsection{The operad $\bbB\oplus \cO$ as an object of $I\bR_q$}
 It is convenient to switch from $\cO$ to a direct sum 
\begin{equation}\label{cudef}
\cU:=\bbB\oplus \cO, 
\end{equation}
which is again a CO in $\bR_q$.
Let 
\begin{equation}\label{prop}
p:\cU\to \bbB
\end{equation}
 be the projection.   Let $i:\red \bbB\to \red \cU$ be the diagonal map, that is
$i=\Id_{\bbB}\oplus \iota$, where $\iota$ is as in (\ref{phiphi}).
In particular $pi=\Id_{\bbB}$,  which defines an isomorphism 
\begin{equation}\label{split}
\red \cU(a)\stackrel\sim\to \red \bbB(a)\oplus \red\cO(a)
\end{equation}
for each $a=(n)^\uncyc$ or  $a=(n)^\cyc$,
different from (\ref{cudef})
This isomorphism induces the following one in $\bR_q$:  $\cU(a)\cong (\bbB(a)\oplus \cO(a),D_a)$ for an appropriate differential $D_a$ such that $\red D_a=0$.
We can now define  objects $\cU(a)'\in I\bR_q'$, where $\gr^0\cU(a):=\bbB(a)$;  $\gr^{-1}\cU(a)=\cO(a)$ and
all other components of $\cU(a)$ vanish.    We finally set $\cU(a):=(\cU(a)',D_a)$,  which is well defined since
$\red D_a=0$.    It also follows that $\cU$ is a CO in $I_{\leq 0}\bR_q$.  The projection (\ref{prop})
induces a projection 
$p_{\cU\bbB}:\cU\to \bbB$,  where we view $\bbB$ as an operad in $I\bR_q$ concentrated in the  grading 0.

We are now ready to perform our first manipulation
\subsection{Tensoring with $\bY(\bbB)^\cyc$ over $\bY(\cU)^\cyc$}

 Here is our first manipulation.  Let $(n)^\cyc\in \bY(\bbB)^\cyc$ be an object.  Let $h_n:(\bY(\bbB)^\cyc)^\op\to \bR_q$
be the Yoneda functor:  $$h_n((m)^\cyc):=\hom_{\bY(\bbB)^\cyc}((m)^\cyc,(n)^\cyc).
$$
Via the projection $p_{\cU\bbB}$, one can also view each $h_n$ as a  functor  $\bY(\cU)^\cyc\to I\bR_q$.
Let us define a functor $V:\bY(\bbB)\to I\bR_q$ as follows:
 $$V((n)^\cyc):= h_n\otimes^L_{\bY(\cU)^\cyc} \cU^\cyc.
$$

We have maps
$$V((n)^\cyc)\stackrel{p_{\cU\bbB}}\longrightarrow h_n\otimes \otimes^L_{\bY(\bbB)^\cyc} \bbB^\cyc\to \bbB(n)^\cyc
$$
We have a $\bY(\bbB)$-action on $V$ so that we can define a CO $I\cU$, where $I\cU^\uncyc=\bbB^\uncyc$
and $I\cU^\cyc:=V$. We have a map of CO $\cU\to I\cU$ whose non-cyclic component is the above
defined projection $p_{\cU\bbB}^\uncyc$,  and  the map $\cU(n)^\cyc\to V((n)^\cyc)$
is as follows:
$$
\cU(n)^\cyc\to h_n(n)\otimes \cU^\cyc(n)\to h_n\otimes^L_{\bY(\cU)^\cyc} \cU^\cyc=V((n)^\cyc),
$$
where the first arrow is induced by the identity $\unit\to h_n(n)$.

\subsubsection{Converting $V$ into an object of $\Doplus \Doplus((\bY(\bbB)^\cyc)^\op\otimes \bR_q)$}\label{gcocv}
 Let  us convert $V:\bY(\bbB)^\cyc\to I_{\leq 0}\bR_q$ to a more convenient object.
Let us start with a simpler situation, when we have  a functor
$W:\bY(\bbB)\to I_{\leq 0}\bR_q$ so that we have graded components
$\gr^k W:\bY(\bbB)^\cyc\to \bR_q$.     We can now pass to resolutions
$$
\cR_{\gr^k W}\in \Doplus((\bY(\bbB)^\cyc)^\op\otimes \bR_q).
$$
Denote $$\Fun(\bY(\bbB)^\cyc):=\Doplus((\bY(\bbB)^\cyc)^\op\otimes \bR_q).
$$
 We can now produce an object
$$
\Sigma W:=\bigoplus\limits_k \cR_{\gr^k W}\in \Doplus\Fun(\bY(\bbB)^\cyc).
$$

In our case, when we have a functor $V:\bY(\bbB)^\cyc\to I\bR_q$, we can still define an object
$\Sigma V\in \Doplus\Fun(\bY(\bbB)^\cyc)$ in a similar way.   One can rewrite
$V((n)^\cyc)=(W((n)^\cyc),D_n)$ for an appropriate $W$ as above. It turns out that the differentials $D_n$
induce a differential on $\Sigma W$, to be denoted by $D_\Sigma$, so that one can put
$\Sigma V:=(\Sigma W,D_\Sigma)$.

\subsubsection{Generalizing the definition of CO}
In order to be able to replace $V$ with $\Sigma V$ in the definition of the CO $\cU$,
we need to generalize the notion of CO as follows.  

 Let us define a  generalized CO (GCO) to be the  following structure:   an asymmetric operad $\cO$ and an object
of $\Doplus\Fun(\bY(\cO))$ footnote.  We now have a GCO
$\cV$, where $\cV^\uncyc=\bbB^\uncyc$ and $\cV^\cyc=\Sigma V$.

We now pass to  our next manipulation.
\subsection{
  $c_1$-localization} The goal of this manipulation is   to   simplify the operad $\cV$.  The $c_1$-localization 
was defined in Sec \ref{c1local}.  We repeat it here.

 We generalize
from the case of cyclic objects.   Let $\underline{\ZQ}$ be the constant cyclic object.  Let $\cR\to \underline{\ZQ}$
be a free resolution.  One has an endomorphism $c_1:\cR\to \cR[2]$, 'the first Chern class'.  One can define a similar
endomorphism for any cyclic object $X$.  Indeed,  let $X\diamond \cR$ be a cyclic object, where
$X\diamond \cR(n):=X(n)\otimes \cR(n)$ for every object $n$ of the cyclic category.  
We therefore have a zig-zag
$$
X\stackrel\sim\leftarrow X\otimes \cR\stackrel{c_1}\to X\otimes \cR[2]\stackrel\sim\to X[2].
$$

Let now $\cO$ be a CO.  One generalizes the $c_1$-endomorphism to act on every object $X:\bY(\cO)^\cyc\to \ccA$,  where $\ccA$ is the ground category.  Indeed, if $X:\bY(\cO)^\cyc\to \ccA$ and $T:\Lambda\to \ccA$ (where $\Lambda$ is the cyclic category), one  defines a functor $X\diamond T:\bY(\cO)^\cyc\to \ccA$, where $X\diamond T(n):=X(n)\otimes T(n)$,
whence a $c_1$-map $c_1:X\diamond \cR\to X\diamond \cR[2]$.

Let us now try to invert $c_1$ (this is what  we mean by the term '$c_1$-localization').   The most straightforward approach, where 
one considers a functor $'X_\loc':\bY(\cO)^\cyc\to \ccA$,  
\begin{equation}\label{xlocc}
'X_\loc'(n) =\hocolim (X\otimes \cR\stackrel {c_1}\to X\otimes \cR[2]\stackrel{c_1}\to X\otimes \cR[4]\stackrel{c_1}\to \cdots
\end{equation}
does fail, because such an  $F(n)$ is acyclic for each $n$.   However, one can  define the localization as an ind-object.
One now defines  an object $X_\loc\in \Doplus\Fun(\bY(\cO)^\cyc)$ by the same formula as in  (\ref{xlocc}).    As the hocolim
is taken in the category of formal sums,  $X_\loc$ is now not necessarily acyclic.
Furthermore, one generalizes the operation $X\mapsto X_\loc$ to an endofunctor 
on $\Doplus\Fun(\bY(\cO)^\cyc)$ in the obvious way.

  Let us  now apply the $c_1$-localization to the  GCO $\cV$ from  Sec \ref{gcocv}.

As was mentioned above, we have a map $\cV^\cyc\to \cR_\bbB^\cyc$  in $\Doplus\Fun(\bY(\bbB)^\cyc)$. 
Consider the diagram
$\cV^\cyc_\loc\to (\cR_\bbB^\cyc)_\loc\leftarrow  \cR_\bbB$.   This diagram has a colimit, to be denoted 
by $\cV_l^\cyc$.   We thus has a GCO $\cV_l$, where $\cV_l^\noncyc=\bbB$. 
We have a natural map of GCO $\cV\to \cV_l$, whence a through map of GCO
\begin{equation}\label{tuv}
\bg \cU\to \cV_l.
\end{equation}

  The map $( \cU)^\uncyc\to \bbB^\uncyc$ unduces a functor $t:\bY(\cU)^\cyc\to \bY(\bbB)^\cyc$
The map (\ref{tuv}) induces a map
 \begin{equation}\label{cucvt}
\red t(\bg \cU)^\cyc_l\to \red \cV_l^\cyc
\end{equation}
 in $\red\Doplus\Fun(\bY(\bbB)^\cyc)$.  One can show (Lemma \ref{locsimpl}) that this map is a homotopy equivalence, which is the main reason why we apply the localization.

The object on the LHS of (\ref{cucvt}) is still too complicated, whence our next manipulation. 
\subsection{Condensation}
This operation was defined in Sec \ref{condensc}, we repeat it here.

Call a  co-simplicial object  in $\cC$ {\em quasi-constant} if every arrow in the simplicial category acts by a homotopy
equivalence.   Recall that  the category $\bY(\assoc)^\cyc$ coincides with the cyclic category,  therefore, given 
a functor $X:\bY(\bbB)^\cyc\to \cC$, one can restict it to $\bY(\assoc)^\cyc$ thus getting a co-cyclic object which we 
can further restrict to the simplicial category, thus getting a co-simplicial object, to be denoted by $X|_\Delta$.
Call $X$  quasi-constant if such is $X|_\Delta$.    One defines an endofunctor $\con$ on the category
of functors $\bY(\bbB)^\cyc\to \cC$ (enriched over $\sets$) and a natural transformation $\con \to \Id$,
where $\con X$ is quasi-constant and the induced map $\Rhom(T;\con X)\to R\hom(T;X)$ is a homotopy equivalence
for every quasi-constant $T$.  One extends $\con$ to an endofunctor on $\Doplusprod\Fun(\bY(\bbB)^\cyc)$
endowed with a natural transformation $\con\to \Id$.
We can now consider an object $\con \cV_l^\cyc\in \Doplusprod\Fun(\bY(\bbB)^\cyc).$
The latter object can be identified up-to a homotopy equivalence as follows.  First of all
we have an induced homotopy equivalence
$$
\con \red t(\bg \cU)^\cyc_l\to  \con \cV_l^\cyc.
$$

As follows from Sec \ref{Hochsch} the map of operads $\red\cO^\bbH\to\red \cO^\bbF$ induces  homotopy 
equivalences
$$
\bbB\oplus \red(\cO^\bbH)^\cyc_l\leftarrow\bbB\oplus\con \red(\cO^\bbH)^\loc_l\to  \con \red t(\bg \cU)^\cyc_l\to\red \con \cV_l^\cyc.
$$

This identifies $\red \con \cV_l^\cyc$.    Claim \ref{c1.2}  easily implies that the object
 $\red (\cO^\bbH)^\cyc_l$ is {\em rigid}:   all its liftings to  $\Fun(\bY(\bbB)^\cyc)$ are homotopy equivalent
to each other, which gives us the desired identification of $\cV_l^\cyc$.

We now pass to the last manipulation.

 \subsection{ Splitting the map of obstruction objects $\cD(\cU)\to \cD(\cU_l)$} This map
admits a splitting according to Claim \ref{spl}. We  therefore can write $\cD(\cU_l^\loc)\cong \cD(\cU)\oplus V_0$ for an appropriate object  $V_0$.
We thus have a map $V_0\to \cD(\cU_l^\loc)$.   It follows that this map factors as follows:
$$
V_0\to \cD(\cU_l^\loc)^\cyc
$$
  As was mentioned above, we have a homotopy equivalence
$\cD(\red\con\cV_l^\loc)^\cyc\sim \hom_{\bY(\red \bbB)^\cyc}(\red\bbB^\cyc;\red\con\cV_l^\cyc)$.  One deduces that there
exists a homotopy equivalence $V_0\sim V_1$ and a map $\bbB^\cyc\otimes V_1\to \red \con\cU_l^\cyc$.
According to the above established structure of $\cU_l^\cyc$,  this map lifts to a map
\begin{equation}\label{vv}
\bbB^\cyc\otimes V\to \con\cV_l^\cyc\to \cV_l^\cyc.
\end{equation}
\subsection{Packing the data into a new (rigid) structure}
  Let us now define a structure that  involves the above map as well as the  constructed map $\cU\to   \cV_l$.
We call it GVCO (G=generalized;  V stands for the term  $V$  in (\ref{vv})),  see Sec \ref{VCO},\ref{GVCO} for the definition.
We thus get a GVCO built from the above specified data. 
One redefines the notion of the obstruction space in the setting of GVCO, and it turns out to be acyclic in our case.
This provides for the desired quantization.

We start with the general theory (Sec \ref{gen}-Sec \ref{iter}), where we give a  definition of GVCO, its  Kaehler differential  module, and the obstruction space.  We then  show that vanishing of the obstruction space implies the lifting.

In the next Section \ref{our} we perform the above outlined steps resulting in building a GVCO and producing an
unobstructed quantization problem.

\section{Quantization:general theory}\label{gen}
\subsection{The categories $\bB_q,\bB_0$}\label{combbb}
FIx the notation $\bB_q:=\Com_{\geq 0}\ccA$.    Let also $\bB_0$ be the SMC of $\bZ_{\geq 0}$-graded objects
in $\ccA$.   We have tensor functors $\gr:\bB_q\to \bB_0$, $\bq:\bB_0\to \bB_q$, where the functor $\gr$ annihilates the differential and  $\bq(X)$ is $X$ viewed as a complex with 0 differential.
\subsubsection{Tensor functors $T_q:\bR_q\to \bB_q$,  $T_0:\bR_0\to \bB_0$.}
  We have a ground SMC $\bR_q$.     In order to do the quantization, we are going to build
a lax tensor functor $T_q:\bR_q\to \bB_q$, where  we set
$$
(T_qX)^i=\prod_{S\in \SMetR} \gr^{i\val S} X_S, 
$$
$i\geq 0$.    The differential on $T_qX$ is induced by that on $X$.

  The functor $T_q$ gives $\bR_q$ the structure of a SMC enriched
 over $\bB_q$.

One defines a tensor functor $T_0:\bR_0\to \bB_0$ in a similar way so that $\bR_0$ is enriched over $\bB_0$.

\subsection{VCO}\label{VCO}  {\em A VCO} in a SMC $\cC$, is a collection 
$\cO^\uncyc(n),\cO_1^\cyc(n),\cO_2^\cyc(n),\cO_3^\cyc(n)\in \cC$, $n\geq 0$,  and  $V\in \cC$,
endowed with the following structure:

--- a circular operad structures on $(\cO^\uncyc,\cO_i^\cyc)$, $i=1,2,3$ which induce the same asymmetric
operad structure on $\cO^\uncyc$;

--- maps of circular operads  with their restrictions onto $\cO^\uncyc$ being  the identity:

$$
(\cO^\uncyc,\cO_1^\cyc)\to (\cO^\uncyc,\cO_2^\cyc)\leftarrow (\cO^\uncyc,V\otimes \cO_3^\cyc).
$$

\subsubsection{Definition of GVCO}\label{GVCO}

The VCO structure   can be reformulated as follows:

---  an asymmetric operad $\cO^\uncyc$;

--- $\bY(\cO^\uncyc)^\cyc$-modules $\cO_i^\cyc$; $i=1,2,3$;

---an object $V\in \cC$;

---   maps of $\bY(\cO^\uncyc)^\cyc$-modules $\cO_1\to \cO_2\leftarrow V\otimes \cO_3$.

One therefore can generalize:   let GVCO be the following structure:

---  an asymmetric operad $\cO^\uncyc$;

--- objects $\cO_i^\cyc$ of $\Doplusprod \Fun(\bY^\cyc(\cO^\uncyc)^\op)$; $i=1,2,3$;

--- an object $V\in \cC$;

---  maps  $\cO_1\to \cO_2\leftarrow V\otimes \cO_3$ in $\Doplusprod\Fun (\bY^\cyc(\cO^\uncyc)^\op)$.

GVCO over $\cC$ form a category enriched over $\sets$.

Given a GVCO $\cO$ one gets a VCO operad $[\cO]$, where $[\cO]^\uncyc=\cO^\uncyc$;
$$
[\cO]^\cyc_k((n)^\cyc):=\hom_{\Doplusprod\Fun((\bY(\cO)^\cyc)^\op)}((n)^\cyc;\cO_k^\cyc)
$$
In fact,  $[,]$ is a functor from the category of GVCO to that of VCO. 
We define a map $f$ from a VCO $\cO'$ to a GVCO $\cO$ as a VCO map $f:\cO'\to [\cO]$.

\subsection{Semi-free operads}    Let $\Coll$ be a category (over $\sets$) whose every object is a collection of objects
$X(n)^\uncyc,X(n)_k^\cyc$, $n>0$, $k=1,2,3$, and $X^V\in \cC$. We have a forgetful functor from  the category of VCO
to $\Coll$ which has a left adjoint, to be denoted by $\free$.    A VCO of the form $\free(X)$ is called {\em  free.}

 {\em   A quasi-free }  VCO is that  of the form
$|(\free(X),D)|$, where $D$ is an arbitrary differential.

{\em  A semi-free VCO} is a quasi-free VCO of the form $|(\free(X_\bullet),D)|$,  where  each
$X_\bullet^\uncyc(n),X_\bullet^\cyc(n)_k,X^V\in \Com_{\leq 0}(\cC)$ and $\gr D=0$.

Given a quasi-free  VCO $\cO$,  we get a  GVCO  $\bg \cO$, where 
$$
\bg \cO_k\in \Fun(\bY^\cyc(\cO^\uncyc)^\op).
$$

A  map $\bg\cO_1\to  \cO$,  where $\cO_1$ is a VCO and $\cO$ a GVCO is the same as a map
$\cO_1\to \cO$ as defined above.
\subsection{Infinitesimal operads}   {\em  An infinitesimal VCO (IVCO)} is 
a VCO in the category $\Com_{[0,1]}(\cC)$.    

Given an IVCO $\cO$ let $|\cO|
$ denote the VCO over $\cC$ obtained
via applying  the totalization functor $\Com_{[0,1]}(\cC)\to \cC$.

Let $\Gr_{[0,1]}(\cC)$ be a category whose every object is a pair $(X^0,X^1)$ of objects of $\cC$.
We set $$
\hom((X^0,X^1);(Y^0,Y^1))=\hom(X^0,Y^0)\oplus \hom(X^1,Y^1). 
$$  

Define a SMC on $\Gr_{[0,1]}(\cC)$ by setting $$
(X^0,X^1)\otimes (Y^0,Y^1):=X^0\otimes Y^0; X^0\otimes Y^1\oplus X^1\otimes Y^0.
$$

We have  functors $$\Gr_{[0,1]}(\cC)\stackrel\iota\to \Com_{[0,1]}(\cC)\stackrel\pi\to \Gr_{[0,1]}(\cC).
$$

 {\em  A split IVCO  } is 
a circular operad in the category $\Gr_{[0,1]}(\cC)$.   

Given a  split IVCO $\cO'$, denote by $X(\cO')$   the set of all  elements $D\in \hom^1((\cO')^0, (\cO')^1)$
such that $(\cO',D)$ is an IVCO.

\subsection{The categories  $S^G(\cO),\bY_V(\cO)$} Let $\cO$ be a GVCO in $\cC$.  Denote by  $S^G(\cO)$  be 
 the category over $\sets$ whose every 
object is a  split IGVCO $\cO'$  along with an identification of IVCO
$(\cO')^0=\cO$.

\subsection{Kaehler differentials}
Let $\cO$ be a VCO in $\cC$. 
 We have a functor $\cO'\mapsto X(\cO')$ from  $S(\cO)$ to  $\sets$.   Suppose $X$ is representable.
Call the representing object $\omega_\cO$.   
\begin{Proposition}   Suppose $\cO$ is quasi-free.  Then $X$ is representable so that $\omega_\cO$ is well-defined.
\end{Proposition}

We have a canonical element in $X(\omega_\cO)$ which gives $\omega_\cO$ an IVCO structure.
It also follows that given an IVCO $\cO'$ with $(\cO')^0=\cO$ (as VCO), 
we have a canonical map of split IVCO 
$$
\omega_\cO\to \cO'
$$
which is compatible with the infinitesimal circular operad structures on both sides.

Let us denote $\Omega_\cO:=(\omega_\cO)^1$.   We have a $\bY(\cO)$-module structure on $\Omega_\cO$.
Let $\cO'\in S^G(\cO)$.   We then have  an identification
$$ 
X(\cO')=\hom_{\bY(\cO)}(\Omega_\cO;(\cO')^1).
$$

\subsection{An IVCO $\cX(\cO')$} Let $E$ be a category  over $\sets$( whose every object 
is a pair $(X\in \cC,f:\unit_\cC\to X)$.    An arrow $(X,f)\to (Y,g)$ is a map of complexes $\phi:\Cone f \to \Cone g$
Such that 

--- the composition $$\Cone f\stackrel \phi\to \Cone g\to \unit_\cC[1]$$ 
equals the projection $\Cone f\to \unit_\cC[1]$.

Let $(X,f)\in E$ and let $\cO$ be a VCO  in $\cC$.  Denote
by $\cO^{X,f}$ the following IVCO in $\cC$:
$$
(\cO^{X,f})^0=\cO;\quad (\cO^{X,f})^1=X\otimes \Omega_\cO;
$$
$$
D:\cO\to X\otimes \cO
$$
is defined as the composition
$$
\cO\stackrel {D_\omega}\to \Omega_\cO\stackrel f\to X\otimes\Omega_\cO.
$$

Let $u:(X,f)\to (Y,g)$ be an arrow in $E$ whose components are 
$$
e:\unit_\cC[1]\to Y;\  v:X\to Y;  \Id:\unit_\cC[1]\to \unit_\cC[1].
$$
In order for $u$ to be a map of complexes, one must have:   $dv=0$;  $de+v(f)+g=0$.

Let us define a map
$u_*:|\cO^{X,f}|\to |\cO^{Y,g}|$  via its components:
$$
\Id:\cO\to \cO;\quad D_\omega\otimes e:\cO\to \Omega_\cO\otimes Y;\Id_{\Omega_\cO}\otimes v:\Omega_\cO\otimes X\to \Omega_\cO\otimes Y.
$$
This way, we get a functor $E\to S(\cO)$.     Let $\cO'\in S^G(\cO)$.   We then  get a functor
$H_{\cO'}:E\to \sets$, where 
$$H_{\cO'}((X,f))=\hom_{S(\cO)^G}(\cO^{X,f}:\cO').
$$
\begin{Proposition}  Suppose $\cO$ is quasi-free. Then $H_{\cO'}$ is representable.
\end{Proposition}
Denote the representing object $(\cX(\cO'),D_{\cO'})$.
\section{Lifting}\label{lifting}
Let  $\ccA$ be a ground category.   Let $\cQ_q$ be SMC enriched  over $\Com_{\leq 1} \ccA$.   Let $\cQ_0:= \pi_0 \cQ$. Let $\red:\cQ_q\to \cQ_0$ be the reduction functor (over $\ccA$).

Let $\bbAV,\bbBV$ be  quasi-free  VCO in $\cQ_q$.  Let $\cO$ be a GVCO in $\cQ_q$.   Suppose
we are given a map $\cO\to  \bg\bbBV$ which quasi-splits as well as a map $\alpha:\bbAV\to \bg\bbBV$ (which is the same as a  VCO map $\alpha:\bbAV\to \bbBV$).

Suppose, finally, that we have a map $\beta:\red\bbAV\to \red \cO$ such that the through map
map $\red\bbAV\to \red \cO \to \red\bg\bbBV$ equals $\red \alpha$.

Let us define the following functor $T:S(\bbAV)\to \sets$.   Set $T(\bbAV')$ to consist of all
maps $|\bbAV'|\to \cO$ such that  the composition $|\bbAV'|\to \cO \to \bg\bbBV$ coincides with the composition
$$|\bbAV'|\to \bbAV\stackrel\alpha\to\bg \bbBV
$$
 and
the induced map $\red|\bbAV'|\to \red \cO$ coincides with the composition
$$
\red|\bbAV'|\to \red \bbAV\stackrel\beta\to \red \cO.
$$

Let $T_E:E\to \sets $ be as follows $T_E((X,f)):=T(\bbAV^{X,f}).$
\begin{Proposition}\label{replif} The functor $T_E$ is representable.
\end{Proposition}

 Let $(\cY,g)$ be the representing object. 
\subsubsection{Explicit formula for $\cY$}
Let $G_\ccA$ be a category whose every object is a pair $(X^0,X^1)$ of objects from $\ccA$.  Set
$$
\hom((X^0,X^1);(Y^0,Y^1))=\hom(X^0,Y^0)\oplus \hom(X^1,Y^1).
$$
The category $G_\ccA$ has a symmetric monoidal structure via $$
(X^0,X^1)\otimes (Y^0,Y^1):=(X^0\otimes Y^0,X^0\otimes Y^1\oplus X^1\otimes Y^0).
$$
We have a tensor functor $r:\Com_{\leq 1}(\ccA)\to G_\ccA$, where $r(X^0\to X^1)=(X^0,X^1).$

Let us get back to a sequence of operads 
$$
\red\bbAV\to \red\cO\to \red\bg\bbBV.
$$
We can consider the above sequence as that of maps of $\bY(\red\bbBV)^\op$-modules.    The arrow
$\red\cO\to \red\bg\bbBV$ admits a kernel, to be denoted by $K$, which is a $\bY(\red\bbAV)^\op$-module.

 We now have a well-defined object
$$\cK:=\hom_{\bY(\bbAV)}(\Omega_\bbAV;K)\in G_\ccA.
$$
\begin{Proposition} We have $\cY\cong \cK^1$.
\end{Proposition}

\section{Lifting: Iterations}\label{iter}

\subsection{Preliminaries}

   Let $\bB_q:=\Com_{\geq 0} \ccA$, as above.    Recall that the category $\bR_q$ is enriched over $\bB_q$.

 Let  $F^{\leq i}:\bB_q\to \bB_q$  associate  to 
a complex $(X^\bullet,d)$ the complex $(F^{\leq i}X^\bullet,d)$, where 
$(F^{\leq i}X)^k=X^k$, $k\leq i$,  $(F^{\leq i}X)^k=0$, $k>i$.   We have a lax tensor structure 
on $F^{\leq i}$ so that we have categories $\bR_i:=F^{\leq i}\bR_q$ enriched over $\bB_q$.
We have tensor functors $\red_{NM}:\bR_N\to \bR_M$, $N\geq M$.   We also have quasi-splittings
$
i_{MN}:\hom^k_{\bR_M}(X,Y)\to \hom^k_{\bR_N}(X,Y)
$
of $\red_{NM}$
which are not compatible with the differential.

As above,
let $\bbAV$, $\bbBV$,  be  quasi-free VCO in $\bR_q$ and $\cU$ be a GVCO endowed with the following maps
$$
\cU\to \bg\bbBV;\quad\beta:  \bbAV\to \bg\bbBV\quad \red\bbAV\to \red\cU,
$$
where the composition  $$
\red\bbAV\to \red\cU\to \red\bg\bbBV
$$ coincides with $\red\beta$.

\subsection{Rigidity} We have the following diagram in $\bB_0$:
\begin{equation}\label{condrig}
\hom_{\Dprodoplus\Fun\bY((\red \bbAV)^\op)}(\Omega_{\red\bbAV};\red\cU)\to \hom_{\Dprodoplus(\Fun\bY(\red\bbAV)^\op)}(\Omega_{\red\bbAV};\red\bg\bbBV).
\end{equation}

Call the data from the previous sub-section {\em rigid} if the cone of this diagram is acyclic.

In this sub-section we will construct a VCO $\bbAV'$ in $\bR_q$, a homotopy equivalence
$\bbAV'\to \bbAV$ and a lifting $\bbAV'\to \cU$ such that the induced map
$\red\bbAV'\to \red\cU$ coincides with the composition $\red\bbAV'\to \red\bbAV\to r\cU$ and  the through map
$$
\bbAV'\to \cU\to \bg\bbBV
$$
coincides with the composition $\bbAV'\to \bbAV\to \bg\bbBV$.

\subsection{Categories $\bR_i$}

 Let $\bB_i:=\Com_{[0,i]} \ccA\subset \Com_{\geq 0}\ccA$ be a full sub-category conisting of all complexes
$(X^\bullet,d)$,  where $X^i=0$ for all $i>1$.      Let also $\pi_{\leq i}:\Com_{\geq 0}\ccA\to \Com_{[0,i]}\ccA$ be given by
$$
\pi_{\leq i} X^j=X^j,\  j\leq i;\quad \pi_{\leq i} X^j=0,\ j>i.
$$   
Introduce a tensor structure on $\Com_{[0,i]}\ccA$  by setting 
 $$
(X_1\otimes X_2\otimes \cdots \otimes X_n)_{\Com_{[0,i]}\ccA}:=\pi_{\leq i}(X_1\otimes X_2\otimes \cdots\otimes X_n)_{\Com_{\geq 0}\ccA}.
$$

The functor $\pi_{\leq i}$ has an obvious tensor structure.   Let also $\iota_i:\Com_{[0,i]}\to \Com_{[0,\infty)}$ be the embedding.   We have a lax tensor structure on $\iota_i$, i.e. the natural transormation
$$
\iota_i(X)\otimes \iota_i(Y)\to \iota_i(X\otimes Y)
$$
is not necessarily a homotopy equivalence.

Let $\bR_{\leq i}:=\pi_{\leq i} \bR_q$ so that $\bR_{\leq  i}$ is an SMC enriched over $\Com_{[0,i]}(\bbAV)$, hence over $\bB$.  

Let us define a tensor functor $C:\Com_{[0,i+1]}(\ccA)\to \Com_{[0,1]}(\ccA)$,
where
$$C(X)=(C(X)^0\stackrel D\to C(X)^1)
$$
with $C(X)^0=\pi_{\leq i} X$ and $C(X^1)=X^{i+1}$.  The differential $D$ is chosen so that we have
an isomorphism $|C(X)|\cong |X|$.
This way $\cR_{i+1}$ is enriched over $\Com_{\leq 1}(\ccA)$.

Let now $\rho(\Com_{[0,1]}(\ccA))\to \ccA$ be  a tensor functor defined as follows $\rho(X^0\to X^1):=X^0$. We now have  an equivalence of SMC enriched over $\ccA$: $\rho \cR_{i+1}\cong \cR_i$.

Let us construct:

--- a projective sequence $\cdots\to\bbAV^{(2)}\to \bbAV^{(1)}\to \bbAV^{(0)}=\bbAV$ of  quasi-free $VCO$ over  
$\bR_q$, where  all the arrows are term-wise homotopy equivalences;

---maps  $\alpha_i:\bbAV^{(i)}\to \cU$  over $\bR_i$ compatible with the maps in the above projective sequence,
where the map $\alpha_0$ coincides with the given one;

--- maps $\beta_i:\bbAV^{(i)}\to \bg\bbBV$ over $\bR_q$, compatible with the above projective sequence, where
the map  $\beta_0$ must coincide with the given one.

--- the following diagram must commute in $\bR_i$;

$$\xymatrix{ \cU\ar[r]&\bg\bbBV\\
\bbAV^{(i)}\ar[u]\ar[ur]
}
$$

Let us do this inductively with respect to $i$.  The base $i=0$ is clear.  
The step (from $i$ to $i+1$) is as follows.
By  the assumption, we have the following diagram
 
$$
\xymatrix{\cU\ar[r]^{\bR_q}& \bg\bbBV\\
      \bbAV^{(i)}\ar[u]^{\bR_i}\ar[ur]^{\bR_q}}
$$

where the through map is a term-wise homotopy equivalence.  We can now use the tensor functor $C$ so as to 
apply   Prop \ref{replif}.  We will then get an object $(\cY,g)$ and an operad $\bbAV^{(i+1),\cY,g}$ defined over 
$\bQ_q$ fitting into the following commutative diagram
$$
\xymatrix{&&\cU\ar[rr]^{\bR_q}&& \bg\bbBV\\
\\
                \bbAV^{(i+1),\cY,g}\ar[rr]^{\bR_q}\ar[uurr]^{\bR_{i+1}}\ar[uurrrr]^{\bR_q}&&\bbAV^{(i)}\ar@<-2ex>[uu]_{\bR_i}\ar[uurr]_{\bR_q}&&}
$$
It follows that $\cY$ is acyclic, therefore,the  natural map 

$$
\bbAV^{(i+1),\cY,g}\to \bbAV^{(i)}
$$
is a homotopy equivalence

 We therefore can set $\bbAV^{(i+1)}$ to be a semi-free resolution of $\bbAV^{(i),\cY,g}$. 
This finishes the proof.

\section{Quantization: our case} \label{our}

Let us construct a VCO $\cV$ possessing all the features from the previous section.
\subsection{More on categories: The category  $I \bE$, where $\bE$ is enriched over $\bB_q$}
 \subsubsection{Truncation}
Let $
\tau_{\geq k}:\bB_q\to \bB_q
$
be the truncation ( $\tau_{\geq k}=\Id$ for all $k\leq 0$).
\subsubsection{Definition of $I\bE $}  Let $\bE$ be a category enriched over $\bB_q$.  Define a 
new category $I\bE$, enriched over $\bB_q$.
 
Let us first define a category $I\bE'$, enriched over $\bB_q$, whose every object
is a collection of objects  $X^{n}\in \bE$, $n\in \bZ$.   
Set
$$
 \hom(X,Y):=\prod\limits_{n,m}\tau_{\geq n-m}\hom_\bE(X^{n};Y^{m})\in \bB_q.
$$

Set $I\bE:=D(I\bE')$.  Suppose $\bE$ is an SMC,   we  the have  an SMC structure on $I\bE$,  where 
$$
(X\otimes Y)^{n}=\bigoplus\limits_{m} X^{m}\otimes Y^{(n-m)}.
$$

Let $I_{\leq 0}\bE\subset I\bE$ be the full sub-category of objects $(X^n,D)$, where $X^n=0$
for all $n>0$ and likewise for $I_{\geq 0}\bE\subset I\bE$.

 We have  tensor functors
 $$
\red:\red(I_{\leq 0}\bE)\to \Com_{\leq 0} \bE;\quad \red:\red(I_{\geq 0}\bE)\to \Com_{\geq 0} \bE,
$$
where
$(X,D)\mapsto (\red X,\red D)$.
\subsection{The operad $\bbA$}
\subsubsection{The operad $\MC$}\label{bbabb} View $\ccA$ as a category enriched over $\bB_q$, via
the obvious embedding $i_\bB:\bbA\to \bB_q$, where $i_\bB(X)^0=X$;  $i_\bB(X)^p=0$ if $p\neq 0$.
For $T\in\bbA$, denote by $T^{(k)}\in I\bbA$, the following object: $(T^{(k)})^l=T$ if $l=k$;  $(T^{(k)})^l=0$ otherwise.

Let $\MC$ be a CO in $I_{\geq 0}\ccA$ freely generated by the elements:

--- $m_n:\ZQ^{(0)}[2-n]\to  \MC(n)^\noncyc$, $n\geq 2$;

---$m_n:\ZQ^{(1)}[2-n]\to  \MC(n)^\noncyc$, $n\leq 2$;
 
---$\mu_n:\ZQ\to \MC(n)^\cyc[1-n]$; it is assumed that each  $\mu_n$ is $\bZ/n\bZ$-invariant.

The differential is defined as follows:
$$
dm_n+\sum m_p\{m_{n+1-p}\}=0;\quad d\mu_n+\mu_p\{m_{n+1-p}\}=0.
$$

We have a natural map $\MC\to \assoc$,  which vanishes on all $m_n$, $n\neq 2$, and on all $\mu_n$,  $n\neq 0$.
We have an induced map
$$
\alg\otimes \MC\to \bbB.
$$

Denote by $\bbA$ the canonical semi-free resolution of $\alg\otimes \MC$. We therefore have an induced
map $\bbA\to \bbBr$.


The embedding $\ccA\to \bR_q$,  $X\mapsto X\otimes \unit$, induces an embedding
$$
I_{\geq 0} \ccA\to I_{\geq 0}\bR_q
$$
which is a tensor functor.  This way, $\bbA$ is a CO in $I_{\geq 0}\bR_q$.
\subsection{Formulating the quantization problem}
Based on  the data from Sec. \ref{data}, we will construct a zig-zag map of operads in $I_{\geq 0}\bR_q$:
 $$\bbA\leftarrow \bbA'\to \cO
$$
which lifts the through map
\begin{equation}\label{bbBr2cU}
\bbA\to \bbBr\stackrel\iota\to \gr \cO,
\end{equation}
where $\cO=\cO^\bbF$ and $\iota$ are as in (\ref{phiphi}).

\subsection{A CO  $\cU$  in $I_{\leq 0}\bR_q$} Set $\cU:=\cO\oplus \bbBr$.  
We have a splitting $\vs:\red \bbBr\to\red  \cU\to\red  \bbBr$, where $\vs|_{\red\bbBr}=\iota$ and
$\vs|_{\bbBr}=\Id$.
We therefore can represent
$$
\cU(*)=(\bbBr(*)\oplus K(*),D),
$$
where  $D:K(*)\to \bbBr(*)$ is a differential and $D$ vanishes in $\bR_0$ ,  that is $D\in \tau_{\geq 1}\hom^1(K(*),\bbBr(*))$.
Let us upgrade $\cU$ to a CO in $I_{\leq 0}\bR_q$, where $\cU^0=\bbBr$ and $\cU^{-1}(*)=K(*)$,
all other components of $\cU^\bullet$ vanish.

 Let us also view $\bbBr$ as a circular operad in $I_{\leq 0}\bR_q$ centered in degree 0.
We still  have a map of circular operads 
\begin{equation}\label{alp}
\cU^\bullet\to \bbBr.   
\end{equation}

\subsection{ Tensoring with $\bY^\cyc(\bbBr)$}  Denote $$
\cI^\cyc:=\bY^\cyc(\bbBr)\otimes^L_{ \bY^\cyc(\cU^\bullet)} (\cU^\bullet)^\cyc:\bY^\cyc(\bbBr)\to I_{\leq 0}\bR_q.
$$

We have a through map
$$
\cI^\cyc\to \bY^\cyc(\bbBr)\otimes^L_{ \bY^\cyc(\bbBr)} \bbBr^\cyc\to \bbBr^\cyc.
$$
We have a circular operad structure on the pair $(\bbBr^\noncyc,\cI^\cyc)$ as well as maps of circular operads over $I_{\leq 0}\bR_q$
$$
\cU^\bullet\to (\bbBr^\noncyc,\cI^\cyc)\to \bbBr,
$$
where the composition coincides with the map  (\ref{alp}).

The object $\cI^\cyc$ is quasi-free,  it  therefore defines an object 
\begin{equation}\label{overlinecI}
\overline{\cI}^\cyc\in D\bigoplus ((\bY(\bbBr)^\cyc)^\op\otimes I_{\leq 0}\bR_q).
\end{equation}
Let us define  functor $$
\lambda:\bigoplus(\bY(\bbBr)^\cyc)^\op\otimes I_{\leq 0}\bR_q'\to  \bigoplus\bigoplus(\bY(\bbBr)^\cyc\otimes \bR_q),
$$
where
$$
\lambda \bigoplus\limits_{a\in A}(k_a)\otimes X_a:=\bigoplus\limits_{n\geq 0}\bigoplus\limits_{a\in A} (k_a)\otimes X^{-n}_a.
$$
This functor naturally gives rise to a tensor functor
$$
\lambda:D\bigoplus ((\bY(\bbBr)^\cyc)^\op\otimes I_{\leq 0}\bR_q)\to D\bigoplus\bigoplus((\bY(\bbBr)^\cyc)^\op\otimes \bR_q).
$$

Hence,  $(\bbBr^\uncyc,\lambda(\cI^\cyc))$ is a GCO, to be denoted  below by $I\cU$.

On the classical level, we have the following diagram
\begin{equation}\label{diag||}
\xymatrix{
\red D\bigoplus ((\bY(\bbBr)^\cyc)^\op\otimes I_{\leq 0}\bR_q)\ar[r]^{\red \lambda}\ar[d]^\sim& D\bigoplus\bigoplus((\bY(\bbBr)^\cyc)^\op\otimes \bR_0)\\
D\bigoplus ((\bY(\bbBr)^\cyc)^\op\otimes \Com_{\leq 0}\bR_0)\ar[r]&\Com_{\leq 0} D\bigoplus ((\bY(\bbBr)^\cyc)^\op\otimes \bR_0)\ar[u]^{||}}
\end{equation}
which commutes up-to a natural isomorphism of the functors.

Below, we denote $$\Fun_\bbBr:=D\bigoplus((\bY(\bbBr)^\cyc)^\op\otimes \bR_q).$$

\def\Loc{{\mathbf{Loc}}}
\subsection{
 $c_1$- Localization}

  Consider a diagram in $\Doplus\Fun_\bbBr$ 
\begin{equation}\label{pto}
I\cU^\cyc_\loc \stackrel{p_{\bbB}}\longrightarrow   \bbBr^\cyc_\loc\leftarrow 0.
\end{equation}
This right arrow of this diagram quasi-splits so that it admits a pull-back to be denoted by $I \cO^\cyc_l$.
Let $I\cO_l:=(\bbBr^\noncyc,I\cO^\cyc_l)$.   

As we have a map $\red\bbBr\to \red \cO$, we have a $\bY(\red\bbBr)^\cyc$-module structure on $\cO^\cyc$, hence
an object $\red\cO^\cyc_\loc\in \Doplus\Fun_{\bbBr}.$
  
We now have maps in $\Doplus\Fun_{\red \bbBr}$:
$$
(\red \cO^\cyc)_\loc\to  \red I\cO^\cyc_l.
$$
\begin{Lemma}\label{locsimpl} This map is a homotopy equivalence.
\end{Lemma}

{\em Sketch of the proof}  
We have a functor
$$F:\red D\bigoplus ((\bY(\bbBr)^\cyc)^\op\otimes I_{\leq 0}\bR_q)\to 
\Com_{\leq 0}\Doplus((\bY(\bbBr)^\cyc)^\op\otimes \bR_0)
$$
which coincides with the composition of the left and the bottom arrows in (\ref{diag||}).
The same diagram implies that 
$\red \lambda\overline{ \cI}^\cyc\cong  |F\overline{\cI}^\cyc|$, where $||$ is the totalization functor
from $\Com_{\leq 0}$.

Denote by $$I_k:=(F\overline{\cI}^\cyc)^{-k}\in \Fun((\bY(\bbBr)^\cyc)^\op\otimes \bR_0)
$$
the $-k$-th component of $F\overline{\cI}^\cyc$ as a complex concentrated in non-positive degrees.

As the arrow $p_{\bbB}$ in (\ref{pto}) induces an isomorphism on the 0-th  component, and we have a homotopy equivalence
$$
\cO^\cyc_\loc\stackrel\sim\to (I_1)_\loc,
$$
 it suffices to show
that $(I_k)_\loc\sim 0$ for all $k>1$.   

Next,  we have ($k\geq 2$): \begin{multline*}
I_k
\cong 
\red \bigoplus\limits_{r>0,t_i>0,t_1+t_2+\cdots+t_r=k }
\bY(\bbB)^\cyc\otimes^L_{\bY(\bbB)^\cyc}\bY_{t_1}(\cU)^\cyc\otimes^L_{Y(\bbB)^\cyc} \bY_{t_2}(\cU)^\cyc\otimes^L_{\bY(\bbB)^\cyc}\\ \cdots \otimes^L_{\bY(\bbB)^\cyc} \bY_{t_r}(\cU)^\cyc\otimes^L_{\bY(\bbB)^\cyc} \cO^\cyc\\
\stackrel\sim\longrightarrow \bigoplus\limits_{r>0,t_i>0,t_1+t_2+\cdots+t_r=k } \bY_{t_1}(\cU)^\cyc\otimes^L_{Y(\bbB)^\cyc} \bY_{t_2}(\cU)^\cyc\otimes^L_{\bY(\bbB)^\cyc}\cdots \otimes^L_{\bY(\bbB)^\cyc}\bY_{t_r}(\cU)^\cyc\otimes^L_{\bY(\bbB)^\cyc} \cO^\cyc.
\end{multline*}

It therefore, suffices to show that given any $M:\bY(\bbB)^\cyc\to \bR_0$, the $c_1$-action on
$
\bY_t(\cU)^\cyc\otimes^L_{\bY(\bbB)} M
$
is 0.

Let us give an explicit desctiption of $$\red\bY_t(\cU)^\cyc:(\bY(\bbB)^\cyc)^\op\otimes \bY(\bbB)^\cyc\to \bR_0.$$

Let $(k):=\{0,1,2,\ldots,k\}$ viewed as a cylcically ordered set.  Given a cyclically monotone map $f:(m)\to (n)$,
we have a total order  on every pre-image $f^{-1}i$, $i\in (n)$.

 We have 
$$
\red \bY_t(\cU)^\cyc((n),(m))=\red\bigoplus\limits_{f:(m)\to (n), S\subset (n),|S|=t}\ \ \bigotimes\limits_{s\in S} \cO^\uncyc(f^{-1}s|s)\otimes \bigotimes\limits_{t\notin S} \bbB^\uncyc(f^{-1}t|t),
$$
where the direct sum is taken over all cyclically monotone $f:(m)\to (n)$ and all subsets $S\subset (n)$ of cardinality $t$.
Denote 
$$
\red \bY_t(\cU)^\cyc((n),(m))^S:=\bigotimes\limits_{s\in S} \cO^\uncyc(f^{-1}s|s)\otimes \bigotimes\limits_{t\notin S} \bbB^\uncyc(f^{-1}t|t)
$$

Let $S_t(n)$ be the set of all $t$-element subsets of $(n)$.
Let us now define a functor $\Sigma_t:\bY^\cyc(\assoc)\to \bbQ\text{-mod}$, where 
we set $\Sigma_t((n)):=\bbQ[S_t(n)]$ (the $\bbQ$-span of $S_t(n)$).    Let $f:(n)\to (m)$ be a cyclically monotone map.
We then have an element $f_*\in \bY^\cyc(\assoc)((m),(n))$.  
Let $S\in S_t(n)$,  let $e_S\in \bbQ[S_t(n)]$ be the corresponding element.  Set 
$f_* e_S:=e_{f(S)}$ if $f$ is injective on $S$ and $f_*e_S:=0$ otherwise.
Denote by $$\Sigma'_t:\bY^\cyc(\assoc)^\op\otimes \bY^\cyc(\assoc)\to \bbQ\text{-mod}
$$
the functor $\Sigma'_t((n),(m)):=\Sigma_t((m))$.
Let $\diamond$ be as in Sec \ref{diamondc}.
We now have a retraction of functors $\bY^\cyc(\bbB)^\op\otimes \bY^\cyc(\bbB)\to \bR_0$:
$$
\bY_t(\cU)^\cyc\stackrel\iota\to  \bY_t(\cU)^\cyc\diamond \Sigma'_t\stackrel\pi\to \bY_t(\cU)^\cyc.
$$
This retraction is defined as follows. Let $\ve_S:\bbQ\to \Sigma'_t((m),(n))$  be given by $\ve_S(1)=e_S$.
Set
$$
\iota|_{\bY_t^\cyc(\cU)^S}:=\Id\otimes \ve_S.
$$

Set $\pi|_{\bY_t^\cyc(\cU)^S\otimes e_S}=\Id$;  $\pi|_{\bY_t^\cyc(\cU)^S\otimes e_T}=0$ if $T\neq S$.

The problem now reduces to showing the nilpotence of the $c_1$-action
on $\Sigma_t$,  which can be verified by the direct computation.  Observe that this fact
requires $\bbQ$ as a base ring.

 \begin{Corollary} \label{psycyc} We have a homotopy equivalence in $\Doplus\prod \Fun_{\red \bbBr}$: 
$$
\red \con(\cO^\cyc_\loc)\to  \red \con(I\cO^\cyc_l).
$$
\end{Corollary}

\begin{Corollary}  Let $$A:=\End_{\Doplus\prod\Fun_{\bbBr}}(\con(I\cO^\cyc_l)).$$
Then $\gr^{>0}A$ is acyclic.
\end{Corollary}

\subsubsection{Studying $\hom_{\Doplusprod\Fun_{\bbBr}}({\bbBr^\cyc};\con(I\cO^\cyc_l))$}
Set $$
H:=\hom_{\Doplusprod\Fun_\bbBr}({\bbBr^\cyc};\con(I\cO^\cyc_l)).
$$
$$
A:=\End_{\Doplusprod\Fun_\bbBr}(\con(I\cO^\cyc_I));\quad B:=\End_{\Doplusprod\Fun_\bbBr}({\bbBr^\cyc}).
$$
$H$ is therefore a $B-A$-bimodule.   Let $B=(B_0,d_{B})$,  $A=(A_0,d_{A})$.  We have $d_{B}=0$, $d_{A}=d_{A_0}+\ad x$.
As $\gr^{>0}A_0$ is acyclic,  we have a zig-zag homotopy equivalence between $A$ and $\gr^0 A_0$.
Therefore,  we have a  homotopy equivalence in $\bB_q$ between $H$ and $(H_0,0)$, where $H_0\in  \bB_0$.

Using  Corollary \ref{psycyc}, the maps in (\ref{conid}) and Lemma \ref{conidgamma}, we get the following diagram of homotopy equivalences:
\begin{equation}\label{xym}
\xymatrix{
\hom_{\Doplusprod\Fun_\bbBr}({\bbBr^\cyc};\red \con(\cO^\cyc_\loc))\ar[d]^\sim\ar[r]^\sim&\hom_{\Doplusprod\Fun_\bbBr}({\bbBr^\cyc}; \red \con(I\cO^\cyc_l))=:H_0\\
\hom_{\Doplusprod\Fun_\bbBr}({\bbBr^\cyc};\bbDhom(K_I;\red (\cO^\cyc_\loc)))&
\hom_{\Doplusprod\Fun_\bbBr}({\bbBr^\cyc};\red \cO^\cyc_\loc)\ar[l]^\sim
}
\end{equation}

Let $$G_0:=\hom_{\Doplusprod\Fun_\bbBr}({\bbBr^\cyc};\red \cO^\cyc_\loc).
$$

We have  the following  zig-zag in $\bB_q$
$$
(G_0,0)\stackrel\sim\leftarrow T\to \hom_{\Doplusprod\Fun_\bbBr}({\bbBr^\cyc}; \cO^\cyc_\loc).
$$
for an appropriate $T\in \bB_q$, where the left arrow quasi-splits.

 Let  $\cL_0$ be as in 
 (\ref{hel}) so that we get a map
$$
\cL_0\to \hom_{\Doplusprod\Fun_\bbBr}({\bbBr^\cyc};\red \cO^\cyc_\loc)=:G_0.
$$

Next, we have a map  $\iota:\bbBr\to \red\cO$ as in (\ref{bbBr2cU}) (as well as in (\ref{phiphi})) whence a  commutative diagram

$$\bbBr\to \red \cO\to \red \cO^\cyc_\loc.
$$
Taking pull-back with respect to the arrow $T\stackrel\sim\to (G_0,0)$,  we get diagrams
$$
\xymatrix{
(G_0,0)& T\ar[l]^\sim\\
(\cL_0,0)\ar[u]&V\ar[l]^\sim\ar[u]}
$$
$$
\xymatrix{
(G_0,0)& T\ar[l]^\sim\\
\unit\ar[u]&\cI\ar[l]^\sim\ar[u]}
$$

for appropriate $V,\cI\in \bB_q$.

We therefore have an  induced map
\begin{equation}\label{tridva}
V\otimes  {\bbBr^\cyc}\to I\cO^\cyc_l.
\end{equation}
and a commutative diagram
\begin{equation}\label{ijn}\xymatrix{
                      \red \cO^\cyc_\loc \ar[dr]   &    \\
                            \unit\otimes \bbBr^\cyc\ar[u]&  \red I\cO^\cyc_l\\
                         \red \cI\otimes \bbBr^\cyc\ar[r]\ar[u] &\red T\otimes \bbBr^\cyc\ar[u]}
\end{equation}
\subsection{Constructing VCO's $\cV,\bbAV, \bbBV$}
\subsubsection{$\cV$}
Let us now build a VCO in $I\bR_q$, where 

---$
\cV^\uncyc=\cU^\uncyc;
$

--- $\cV^\cyc_1=\cO^\cyc\oplus \cI\otimes \bbBr^\cyc$;

---$\cV^\cyc_2= I\cO^\cyc_\loc$;

---$\cV^\cyc_3={\bbBr^\cyc}$;

--- $\cV^V=V\otimes \unit_{\bR_q}$.
The map $p_{12}:\cV^\cyc_1\to \cV^\cyc_2$ is as follows:
$$
p_{12}:\cO^\cyc\to I\cO^\cyc_l
$$
is the canonical map and
$$
p_{12}:\cI\otimes \bbBr^\cyc\to I\cO^\cyc_l
$$
coincides with the map in (\ref{ijn}).
The map $\cV^V\otimes \cV^\cyc_3\to \cV^\cyc_2$ is as in  (\ref{tridva}).
\subsubsection{$\bbBV$}
Set $(\bbBV)^\uncyc=\bbBr^\uncyc$;  $(\bbBV)^\cyc_1=\bbBr^\cyc$, $\bbBV^\cyc_2=0$;   $\bbBV^\cyc_3=\bbBr^\cyc$;
$(\bbBV)^V=0$.

\subsubsection{$\bbAV$} Set $(\bbAV)^\uncyc=\bbA^\uncyc$;  $(\bbAV)^\cyc_1=\bbA^\cyc\otimes \cI$,
$\bbAV^\cyc_3=\bbA^\cyc$.
Let $\cA$ be a complex $0\to \cA^{-1}\to \cA^0\to 0$, where $\cA^0=\unit\oplus \unit$ and $\cA^{-1}=\unit$,
the differential $\cA^{-1}\to \cA^0$ is $\Id\oplus(-\Id)$.   Let $i_0,i_1:\unit \to \cA$ be the embeddings
onto the first and the second summand of $\cA^0$. Set
$\bbAV^\cyc_2=\bbA^\cyc\otimes \cA$. Set  $(\bbAV)^V=\unit$.

Set  the map $
\bbAV^\cyc_1\to \bbAV^\cyc_2$ to be  the composition
$$
\bbAV^\cyc_1\to \bbA^\cyc\to \bbAV^\cyc_2,
$$
where the first arrow is induced by the map $\cI\to \unit$ and the second arrow is
$\Id_{\bbA^\cyc}\otimes i_0$;  set the map $\bbAV^\cyc_3\to \bbAV^\cyc_2$
to be $\Id_{\bbA^\cyc}\otimes i_1$.

\subsubsection{The map $p:\cV\to \bbBV$} 

We have $\cV^\uncyc=\cU^\uncyc=(\cO^\bbF)^\uncyc\bigoplus \bbBr$, (see (\ref{data}), whence 
the  projection $p^\uncyc:\cV^\uncyc\to \bbBr^\uncyc$.   
Define the projection $$p_1:\cV_1^\cyc\to \bbBr^\cyc\otimes \cI\to \bbBr^\cyc$$ in a similar way.
Set all the remaining maps to 0.

\subsubsection{The map $\beta:\bbAV\to \bbBV$}
Set 
$$\beta^\uncyc:\bbAV^\uncyc=\bbA^\uncyc\to \bbBr^\uncyc=\bbBV^\uncyc
$$ to be induced by the map  $\bbA\to \bbBr$ as defined in Sec \ref{bbabb}.
Define
$$
\beta_1:\bbAV_1=\bbA^\cyc\to \bbBr^\cyc=\bbBV_1,\ \beta_2:\bbAV_2=\bbA^\cyc\to \bbBr^\cyc=\bbBV_2,\ \beta_3:\bbAV_3=\bbA^\cyc\to \bbBr^\cyc=\bbBV_3,
$$
to be induced by  the same map $\bbA\to \bbBr$ from Sec \ref{bbabb}.

Define $q^V:\unit\to 0$ to be the 0 map.
\subsubsection{The map $\alpha: \red \bbAV\to \red\cV$}
The map $\alpha_1:\red \cI\otimes \bbA^\cyc\to  \red \cV$ is defined as the  direct sum of:
$$
\alpha_{11}:\red \cI\otimes \bbA^\cyc\to \bbA^\cyc\to \cO^\cyc
$$
and
$$
\alpha_{12}=-\beta_{12}:\red \cI\otimes \bbA^\cyc\to \bbA^\cyc\to \bbBr^\cyc,
$$
where $\beta_{12}$ is the tensor product of the above defined arrows 
$\cI\to \unit$ and $\bbA^\cyc\to \bbBr^\cyc$.
As a result, the composition
$$
\red \bbAV_1\to \red \bbBV_1\to \red\bbBV_2
$$
is 0.  We now can set
$\alpha_2=0$;   $\alpha_3:\bbA\to \bbBr$ to be the above defined map, and, finally,
$\alpha^V=0$.

\subsubsection{Checking the rigidity condition}
The rigidity condition (\ref{condrig}) follows from the homotopy equivalence in  (\ref{hel}).

We therefore conclude the existence of a VCO  $\bbA'$, a homotopy equivalence
$\bbAV'\to \bbAV$, and a map $\bbAV'\to  \cV$ such that $\red \bbAV'\to \red \cV$ is homotopy equivalent
to the above defined map $\red\bbAV\to \red \cV$.  Restriction to the $\uncyc$- and $1-$ components
yields the following zig-zag of CO in $\bR_q$:
 $$\bbA\stackrel\sim\leftarrow\bbA'\to \cO,
$$
The reduction of this zig-zag to $\bR_0$ is homotopy equivalent to the arrow 
$$\red\bbA\to \red \bbBr\to \red \bbB\stackrel\phi\to  \red \cO,
$$
where $\phi$ is as in (\ref{phiphi}). This  solves the quantization problem.

\subsection{Building the Microlocal category}
Recall (Sec \ref{strt}) that we have a category $\Doplus\bU(\bbFR)$ enriched over the category $\Dprodoplus\bU(\bbF)$.
Next, we have an object $\cA\in \bU(\bbF)$, where $\cO$ is a quasi-free resolution of the  full operad of $\cA$,
Therefore, $\cA$ carries a $\bbA'$-action.    We have  homotopy equivalences
$\bbA'\to \bbA\to \ba\otimes \MC$,  we can now construct a homotopy equivalent object $\cB\to \cA$ with
a $\ba\otimes \MC$-action on it.   The $\red \bbA\otimes \MC$-action passes through 
the map $\red \bbA\otimes \MC\to\red  \bbB$.

Let $\MC_R$ be a 2-colored operad in  $\ccA$, such that a $\MC_R$-structure on a pair $(A,M)$ is the same as
that of an $A_\infty$-algebra with a curvature on $A$ and that of an $A$-module on $M$.
Let us define an $A_\infty$- category $\cM_M$, enriched over $\bR_q$, whose every object is a $\ba\otimes \MC$-algebra
structure on a pair $(\cB,M)$, $M\in \Doplus \bU(\bbFR)$, whose restriction onto $\cB$ coincides with the existing one.
The  $A_\infty$-structure is defined in the standard way.

\end{document}